\documentclass[12pt]{book}
\usepackage[spanish,english]{babel}
\usepackage[latin1]{inputenc}  
\usepackage{amsmath}
\usepackage{amssymb}
\usepackage{amsthm}
\usepackage{url}
\usepackage{stmaryrd} 
\usepackage{fancyhdr}
\usepackage{tocbibind}

\usepackage{anysize} 
\usepackage{verbatim} 

\usepackage{amscd,tikz,multicol}
\usetikzlibrary{arrows,shapes}

\def\NZQ{\mathbb}               
\def\NN{{\NZQ N}}
\def\TT{{\NZQ T}}
\def\QQ{{\NZQ Q}}
\def\ZZ{{\NZQ Z}}
\def\RR{{\NZQ R}}

\def\PP{{\NZQ P}}
\def\FF{{\NZQ F}}
\def\KK{{\NZQ K}}
\def\MM{{\NZQ M}}
\def\LL{{\NZQ L}}
\def\EE{{\NZQ E}}

\def\Ac{{\mathcal{A}}}
\def\Bc{{\mathcal{B}}}
\def\Uc{{\mathcal{U}}}
\def\B'c{{\mathcal{B'}}}
\def\U'c{{\mathcal{U'}}}
\def\Mc{{\mathcal{M}}}
\def\Nc{{\mathcal{N}}}
\def\Kc{{\mathcal{K}}}
\def\Hc{{\mathcal{H}}}
\def\Fc{{\mathcal{F}}}
\def\Pc{{\mathcal{P}}}
\def\Vc{{\mathcal{V}}}
\def\Cc{{\mathcal{C}}}
\def\Sc{{\mathcal{S}}}
\def\Rc{{\mathcal{R}}}
\def\Tc{{\mathcal{T}}}
\def\Qc{{\mathcal{Q}}}
\def\Ic{{\mathcal{I}}}
\def\Lc{{\mathcal{L}}}
\def\Dc{{\mathcal{D}}}
\def\Xc{{\mathcal{X}}}
\def\Ec{{\mathcal{E}}}

\def\mi{{\mathfrak{i}}}
\def\mj{{\mathfrak{j}}}
\def\mm{{\mathfrak{m}}}
\def\mq{{\mathfrak{q}}}
\def\mp{{\mathfrak{p}}}
\def\mD{{\mathfrak{D}}}
\def\mF{{\mathfrak{F}}}

\def\ab{{\mathbf a}}
\def\bb{{\mathbf b}}

\def\eb{{\mathbf e}}
\def\kb{{\mathbf k}}
\def\ub{{\mathbf u}}

\def\xb{{\mathbf x}}
\def\yb{{\mathbf y}}

\def\Sb{{\mathbf S}}

\def\lex{{\tt lex}\,}
\def\deglex{{\tt deglex}\,}
\def\degrevlex{{\tt degrevlex}\,}
\def\cocoa{{\tt CoCoA}\,}
\def\singular{{\tt Singular}\,}
\def\macaulay{{\tt Macaulay2}\,}
\def\init{{\sf in}}
\def\lt{{\sf lt}}
\def\gin{{\tt gin}}
\def\cocoalib{{\tt CoCoALib}\,}
\def\cpp{{\tt C++}\,}

\def\poly#1#2#3{#1[#2_1,\dots,#2_{#3}]}
\def\pot#1#2{#1[\kern-0.28ex[#2]\kern-0.28ex]}

\def\Exseq#1#2#3{0\rightarrow#1_{#2}\stackrel{#3_{#2}}{\rightarrow}#1_{#2-1}\stackrel{#3_{#2-1}}{\rightarrow}\cdots#1_1\stackrel{#3_1}{\rightarrow}#1_0\rightarrow\Mc\rightarrow 0}
\def\Res#1#2#3{\dots\rightarrow#1_{#2}\stackrel{#3_{#2}}{\rightarrow}#1_{#2-1}\stackrel{#3_{#2-1}}{\rightarrow}\cdots#1_1\stackrel{#3_1}{\rightarrow}#1_0\rightarrow 0}
\def\Implies{\ifmmode\Longrightarrow \else
     \unskip${}\Longrightarrow{}$\ignorespaces\fi}
\def\implies{\ifmmode\Rightarrow \else
     \unskip${}\Rightarrow{}$\ignorespaces\fi}
\def\iff{\ifmmode\Longleftrightarrow \else
     \unskip${}\Longleftrightarrow{}$\ignorespaces\fi}

\let\:=\colon

\newtheorem{Theorem}{Theorem}[section]
\newtheorem{Lemma}[Theorem]{Lemma}
\newtheorem{Corollary}[Theorem]{Corollary}
\newtheorem{Proposition}[Theorem]{Proposition}
\newtheorem{Remark}[Theorem]{Remark}

\newtheorem{Example}[Theorem]{Example}

\newtheorem{Definition}[Theorem]{Definition}

\marginsize{3cm}{2.25cm}{2cm}{2.75cm}
\renewcommand{\headrulewidth}{0.1pt}
\renewcommand{\thesection}{\thechapter.\arabic{section}}
\renewcommand{\thesubsection}{\thesection.\arabic{subsection}}
\renewcommand{\thesubsubsection}{\thesubsection.\arabic{subsubsection}}

\renewcommand{\chaptermark}[1]{\markboth{Chapter \thechapter \quad #1}{}} 
\newcommand{\rrdc}{\mbox{\,\(\Rightarrow\hspace{-9pt}\Rightarrow\)\,}}
\newcommand{\lrdc}{\mbox{\,\(\Leftarrow\hspace{-9pt}\Leftarrow\)\,}}
\newcommand{\lrrdc}{\mbox{\,\(\Leftarrow\hspace{-9pt}\Leftarrow\hspace{-5pt}\Rightarrow\hspace{-9pt}\Rightarrow\)\,}} 

\fancypagestyle{plain}{%
    \fancyhf{} 
    \fancyfoot{}
    \fancyfoot[C]{\it \thepage} 
    \renewcommand{\headrulewidth}{0pt}
    
    }

\makeatletter
\def\cleardoublepage{\clearpage\if@twoside \ifodd\c@page\else
\hbox{} \vspace*{\fill}
\vspace{\fill} \thispagestyle{empty}
\newpage
\if@twocolumn\hbox{}\newpage\fi\fi\fi} \makeatother

\begin{document}

\newcommand\spanishTOCname{\'Indice General}
\makeatletter
\newcommand\spanishTOC{%
    \chapter*{\spanishTOCname}%
      \@mkboth{\MakeUppercase\spanishTOCname}%
              {\MakeUppercase\spanishTOCname}%
    \@starttoc{spt}%
    }
\makeatother

\def\spanishchapter#1{%
  \addcontentsline{spt}{chapter}{\numberline{\thechapter}#1}
}

\def\spanishsection#1{%
  \addcontentsline{spt}{section}{\numberline{\thesection}#1}
}

\def\spanishsubsection#1{%
  \addcontentsline{spt}{subsection}{\numberline{\thesubsection}#1}
}

\def\spanishsubsubsection#1{%
  \addcontentsline{spt}{subsubsection}{\numberline{\thesubsubsection}#1}
}

\def\spanishbibliography#1{%
  \addcontentsline{spt}{chapter}{\numberline{\thebibliography}#1}
}


\renewcommand{\headrulewidth}{0.0pt} 

\begin{titlepage}
  \vspace*{1cm}
  \begin{center}%
    {\Huge {\bf Combinatorial Koszul Homology, \\ Computations and Applications\par}}%
    \vspace*{1.5cm}%
    {\Large \bf Eduardo S\'aenz-de-Cabez\'on Irigaray \par}
    \vspace*{1.5cm}
    {\Large Dissertation submitted for the degree \\ of Doctor of Philosophy\par}
    \vspace*{0.5cm}%
       \par
      \vskip 5em%
\end{center}
\begin{flushright}
   \begin{large}
      \begin{tabular}{l l}
           Supervisors: & Prof. Dr. Luis Javier Hern\'andez Paricio \\
                       & Prof. Dr. Werner M. Seiler
      \end{tabular}
   \end{large}
\end{flushright}
\vspace*{3.5cm}
\begin{center}

{\Large Universidad de La Rioja \\}
{\large Departamento de Matem\'aticas y Computaci\'on \par} 
  
  \vspace*{1cm}

  {\large Logro\~no, December 2007}
\end{center}
\end{titlepage}

\newpage
\thispagestyle{empty} \vspace*{19.5cm}

\noindent This work has been partially supported by project CALCULEMUS (Contract No HPRN-CT-2000-00102) and   GIFT (NEST Contract No 5006) from the European Union. And by \mbox{project ANGI-2005/10} from the Comunidad Aut\'onoma de La Rioja, and \mbox{grants ATUR-04/51}, \mbox{ATUR-05/46}, \mbox{ATUR-06/31} from the Universidad de La Rioja.

\sloppy

\frontmatter 
\fancyhf{}
\fancyfoot[C]{\thepage} 

\chapter*{Acknowledgments}

\begin{flushright}
\mbox{\parbox{0.7\textwidth}{\footnotesize\it A la memoria de mi padre, de quien aprend\'i que s\'olo cuando hacemos caso al coraz\'on por encima de la raz\'on, podemos ser felices.

A mi madre, que me ense\~n\'o que s\'olo el control de la raz\'on puede llevar a buen puerto los impulsos del coraz\'on.

A ambos, de cuya uni\'on he aprendido el necesario equilibrio que hace m\'as sencillo el camino de la vida.

A Elena, mi amor, es decir, mi vida. A Juan, Luc\'ia y M\'ikel, mi vida, es decir, mi amor.}}
\end{flushright}


\vskip1cm

To write here the list of all persons to whom I am grateful for their contribution in a direct or indirect manner to the elaboration of this thesis I would need more time, space and memory than I have available. I will only emphasize the names of Luis Javier Hern\'andez Paricio and Werner M. Seiler, my two supervisors, which have given me their advice and support, and have encouraged me and teached me many things during these years, not only in the mathematical aspects of the thesis. I also want to emphasise Julio Rubio, who acted as a ``third supervisor'' and has lead me in these first steps into the world of mathematics and its circumstances. Special thanks go to Mar\'ia Teresa Rivas for her continuous advice and to Jacques Calmet who hosted me in his Institut in Karlsruhe.

I will just note here the list of a few among all the persons to whom I want to thank for one thing or another. Of course, some of them have been very special and I feel a deeper gratitude toward them. They know who they are, and I hope I am able to show them everyday how much I appreciate everything they give me.

First of all, those who taught me all the mathematics I know and all the mathematics I have forgotten, they gave the support and help I needed many times: Jos\'e Luis Ansorena, Manuel Bello, Manuel Benito, Mar\'ia Pilar Benito, \'Oscar Ciaurri, Luis Espa\~nol, Ignacio Extremiana, Jos\'e Antonio Ezquerro, Juan Carlos Fillat, Jos\'e Javier Guadalupe, Jos\'e Manuel Guti\'errez, Miguel \'Angel Hern\'andez, Jes\'us Antonio Laliena, Laureano Lamb\'an, V\'ictor Lanchares, Eloy Mata, Mar\'ia del Carmen M\'inguez, Jes\'us Mun\'arriz, Ana Isabel Pascual, Jos\'e Mar\'ia P\'erez and Juan Luis Varona.

My dear colleagues, together with whom I learned many things in and outside mathematics: Miriam Andr\'es, Jes\'us Mar\'ia Aransay, C\'esar Dom\'inguez, Elena Fau, Ignacio Garc\'ia Marco, Marcus Hausdorf, Clara Jim\'enez, Fabi\'an Mart\'in, Francisco Javier P\'erez and of course Ana Romero.

The people who made my life easier in Karlruhe in many aspects, most of which have became true friends: Jos\'e Mar\'ia Alonso, Anusch Daemi, Luana Deambrosi, Ignacio Izquierdo, Sabri Caglan Kessim, Tolga Komurcu, Jaime Ponce, Carolina Santoro, Helga Scherer, and Graham Steel. And my friends in Logro\~no who make my life easier every day: Los Caballeros Avutardos, Juan Carlos D\'iez, Raquel Fern\'andez (2), Clotilde L\'opez, Esperanza Madorr\'an, Pacho Moral, Jorge Pad\'in, Raquel San Mart\'in, Elisa Tobalina, Nacho Ugarte and of course, Santiago Urizarna.

Some mathematicians I have encountered in different places who have made direct or indirect contributions to this thesis, in particular John Abbot, Anna Bigatti, Francis Sergeraert and Henry P. Wynn.

My family:  Nieves and Enrique, \'Alvaro, Cristina, Javier and Mar\'ia, Adolfo, Ana, Anuska and Pablo, Cristina, Rub\'en and Carmen, and my beloved Elena, Juan, Luc\'ia and Mikel.

... and many others.
\thispagestyle{empty}

\renewcommand{\headrulewidth}{0.1pt} 

\fancyhf{}
\fancyhead[EL]{\it \thepage} 
\fancyhead[ER]{\itshape Contents} 
\fancyhead[OL]{\itshape Contents}
\fancyhead[OR]{\it \thepage}  
\fancyfoot{}        

\addcontentsline{spt}{chapter}{\'Indice general}
\tableofcontents 


\newpage \thispagestyle{empty}

\mainmatter 


\fancyhf{}
\fancyhead[EL]{\it \thepage} 
\fancyhead[ER]{\itshape Introduction} 
\fancyhead[OL]{\itshape Introduction}
\fancyhead[OR]{\it \thepage}  
\fancyfoot{}        

\include{introduction}


\fancyhf{} 
\fancyhead[EL]{\it \thepage} 
\fancyhead[ER]{\itshape \leftmark} 
\fancyhead[OL]{\rightmark}
\fancyhead[OR]{\it \thepage}  
\fancyfoot{}        

\chapter*{Preface}
\spanishchapter{Prefacio}

\begin{flushright}
\mbox{\parbox{0.5\textwidth}{\footnotesize\it%
 Unfortunately, there seem to be fewer ways to compute Koszul cohomology groups than reasons to compute them.
\rm\cite{G84}}}
\end{flushright}

Monomial ideals are a particular type of ideals in the polynomial ring which have a combinatorial nature. They play an important role in commutative algebra because some problems in this area concerning ideals or modules over the polynomial ring can be reduced to problems about monomial ideals, in particular in the context of Gr\"obner basis techniques. Also, some theoretical properties of certain types of monomial ideals are very relevant in the theory of syzygies and Hilbert functions. Moreover, the combinatorial nature of monomial ideals makes them suitable for applications inside other areas of mathematics and outside mathematics, ranging from graph theory or differential systems to reliability theory. There has been a lot of interest about this type of objects and their applications in the recent years, and they have become a very active area of research.

In this thesis we concern about the homological properties of monomial ideals. Usually, the homological description of a monomial ideal is given by its minimal free resolution, from which one computes the most relevant homological invariants of the ideal. However, it is an open problem to give a closed description of the minimal resolution of a monomial ideal, although some interesting works have dealt with this problem in the past. Here we will use Koszul homology for giving this homological description of monomial ideals. We will see that this homology can give us a good way to describe the homological properties of the ideal as well as some other structural properties of it. Both approaches are in some sense equivalent since they represent two different methods to compute the $Tor$ modules of the ideal.

The combinatorial nature of monomial ideals introduces a combinatorial way to work with Koszul homology, therefore in this context, we speak of {\bf combinatorial Koszul homology} which gives the title to this thesis.

In the first chapter, we introduce the main characters in the play: Koszul homology and monomial ideals. Koszul homology is the homology of a complex first introduced by J-L. Koszul in a geometrical context \cite{K50a,K50b}. It has been an object of interest in commutative algebra for years, and it is also at the merging of important problems in formal theory of differential systems and commutative algebra due to its relation with the Spencer complex \cite{S69}. On the other hand, much work have been done in relation to monomial ideal by algebraists, since they are important objects, in particular in the work in Gr\"obner basis theory, which allows us to reduce many problems related to polynomial ideals to problems concerning monomial ideals, much easier to handle, since they have a combinatorial nature. This chapter will be devoted to present the basic notations and notions, and the main properties of these two objects.

The second chapter is dedicated to describe the homological and structural properties of monomial ideals that can be read from the Koszul homology. First of all, we focus on homological properties and invariants, which are the main goal of this chapter. Then, we treat some algebraic properties of these ideals. These include Stanley decompositions and irredundant irreducible and primary decompositions. We also transfer the results on the homology of monomial ideals to polynomial ideals. Here we need homological perturbation and Gr\"obner basis theory. This makes the methods described in the second chapter applicable in a more general setting, and allows us to follow a program similar to the one used in Gr\"obner basis theory.

The third chapter is devoted to computations. We give an algorithm to compute the Koszul homology of monomial ideals based on different techniques. We use homological and combinatorial techniques, and introduce Mayer-Vietoris trees, which not only allow us to make homological computations on monomial ideals, but are also a new tool to analyze the structure of these ideals. Several types of Mayer-Vietoris trees will be analyzed in this context. Another tool used in this chapter is simplicial homology; for this, some improvements come from the study of discrete Morse functions and dualities and in particular from the application of Stanley-Reisner theory to the Koszul simplicial complexes. A study of the algorithm is provided, together with implementation issues and some experiments and comparisons to other algorithms with a similar purpose. These show that Mayer-Vietoris trees are an efficient alternative to perform homological computations on monomial ideals.

The fourth chapter is devoted to applications. Different applications will be shown: Applications of Mayer-Vietoris trees to several types of monomial ideals, which are themselves applied either inside commutative algebra (Borel-fixed, stable, segment or generic ideals) or in other areas (Valla, Ferrers, quasi-stable ideals). Also, some applications to different fields like the formal theory of differential systems and reliability theory are developed. These applications use the properties of the Koszul homology shown in the second chapter, and the computational tools presented in the third one.

There are concepts from different areas of mathematics that appear in many places of this thesis. Some readers will probably be familiar with some of them and not familiar with others, or vice versa. For this reason, and for readability, we include several appendices in which the relevant definitions are given. They are intended to serve as references for the main concepts, not as introductions or explanations of the different theories involved.  		

\chapter{Koszul Homology and Monomial Ideals}\label{Koszul-homology-monomial-ideals}
\spanishchapter{Homolog\'ia de Koszul e ideales monomiales}

\begin{flushright}
\mbox{\parbox{0.5\textwidth}{\footnotesize\it%
Le complexe de Koszul \'etait, \`a l'origine, une alg\`ebre diff\'erentielle gradu\'ee (ADG) de la forme $B\otimes \wedge(x_1,\dots,x_r)$, $dx_i=b_i\in B$, o\`u $B$ \'etait une sous ADG (dite la \emph{base}).
\rm\cite{H87}}}
\end{flushright}

In this chapter we introduce the main objects we shall deal with: Koszul homology and monomial ideals. The chapter is divided into three sections.

In the first section, we present the definition, origin and main properties of Koszul homology. We pay special attention to its graded and multigraded versions. We also explain in this section the duality between Koszul homology and Spencer cohomology, which is at the origin of the application of the first to differential systems. Main references for this section are \cite{K50a,K50b,S69,S07,S07b}. 

In the second section we define and give the main characterizations and properties of monomial ideals. Based on the combinatorial properties of them, characterizations of their main algebraic and homological invariants are given, and also some algorithms for their computation. We pay special attention to resolutions. Main references for this chapter are \cite{V01,MS04}.

Finally, in the third section we present a non-exhaustive catalog of topological and homological tools applicable to Koszul homology computations on monomial ideals. The main tools introduced are Koszul simplicial complexes, Stanley-Reisner ideals and Alexander duality. Moreover, we apply Mayer-Vietoris sequences and mapping cones to the computation of resolutions and Koszul homology for monomial ideals. Some of these techniques and algorithms are presented here for the first time. Main references for this section are \cite{MS04,B96,S96}.
\section{The Koszul Homology}
\spanishsection{La Homolog\'ia de Koszul}
\subsection{Basic notions}
\spanishsubsection{Nociones b\'asicas}
The Koszul homology of a module over a graded ring is an important object that plays an interesting role in the merging of commutative algebra, algebraic geometry, the formal theory of differential systems and other areas. In commutative algebra and algebraic geometry, the Koszul homology of a module of the polynomial ring is strongly related to its minimal free resolution and all the invariants that can be read from it, like Betti numbers, Hilbert function, Castelnuovo-Mumford regularity, depth, homological dimension etc. There are many applications of Koszul homology in geometry and many interesting problems in which the computation of Koszul homology is an important issue, see for example \cite{G84,G87}. The Koszul complex has also relation with regular sequences and provides a characterization of the modules $Tor(\Mc,\kb)$ for $\Mc$ a graded module and $\kb$ the base field of the polynomial ring, it has been used to characterize regular, Cohen-Macaulay, Gorenstein rings and complete intersections, to name a few examples, see for instance \cite{AB58,B64,BR65,BR64}. In the formal theory of differential systems, the role of Koszul homology shows up because of its duality with respect to Spencer cohomology, which plays a fundamental role in the characterization of involution and formal integrability (see \cite{LS02}). These relations show a parallelism between certain features of complexes in commutative algebra and formal theory, which are more clearly read in a homological algebra context, see section \ref{formal_theory} below.

In the following pages $R$ will be the polynomial ring $\poly \kb x n$ in $n$ variables over a field $\kb$ of characteristic $0$. We will always consider the usual grading and multigrading in $R$. We will be interested in computing the Koszul homology of ideals (considered as $R$-modules) and modules of the form $R/I$ where $I$ is an ideal.

\subsection{The Koszul complex and Koszul homology}\label{koszdef}
\spanishsubsection{El complejo de Koszul y la homolog\'ia de Koszul}
Let $\Vc$ be a $n$-dimensional $\kb$-vector space. Let $S\Vc$ and $\wedge \Vc$ be the symmetric and exterior algebras of $\Vc$ respectively (see Appendix \ref{apal}). We consider the basis of $\Vc$ given by $\{x_1,\dots,x_n\}$; then we can identify $S\Vc$ and $R$ and consider the following complex
$$\KK: 0\rightarrow R\otimes_\kb\wedge^n\Vc\stackrel{\partial}{\rightarrow}R\otimes_\kb\wedge^{n-1}\Vc\stackrel{\partial}{\rightarrow}\cdots R\otimes_\kb\wedge^{1}\Vc\stackrel{\partial}{\rightarrow}R\otimes_\kb\wedge^0\Vc\rightarrow 0$$

Any element of $R\otimes_\kb\wedge^i\Vc$ can be written in two different ways: First, as a $\kb$-linear combination of elements of the form $x_1^{\mu_1}\dots x_n^{\mu_n}\otimes x_1^{j_1}\wedge\cdots\wedge x_n^{j_n}$, where the $\mu_k$ are non-negative integers and the $j_k$ are either $0$ or $1$, and exactly $i$ of them are equal to $1$. In this case, the differentials $\partial$ are given by the rule

$$\partial(x_1^{\mu_1}\dots x_n^{\mu_n}\otimes x_1^{j_1}\wedge\cdots\wedge x_n^{j_n})=\sum_{\{k\vert j_k=1\}} (-1)^{\sigma(k)+1} x_{k}\cdot x_1^{\mu_1}\dots x_n^{\mu_n}\otimes x_1^{j_1}\wedge\cdots\wedge x_k^{j_k-1}\wedge\cdots \wedge x_n^{j_n}$$
where $\sigma(k)$ is the position of $k$ in the ordered set $\{k\vert j_k=1\}$.

Alternatively, elements of $R\otimes \wedge^i\Vc$ can be expressed as $\kb$-linear combinations of elements of the form $ x_1^{\mu_1}\dots x_n^{\mu_n}\otimes x_{j_1}\wedge\cdots\wedge x_{j_i}$ with ${1\leq j_1<\cdots <j_i\leq n}$; then, the differentials have the form

$$\partial(x_1^{\mu_1}\dots x_n^{\mu_n}\otimes x_{j_1}\wedge\cdots\wedge x_{j_i})=\sum_{k=1}^i (-1)^{k+1} x_{j_k}\cdot x_1^{\mu_1}\dots x_n^{\mu_n}\otimes x_{j_1}\wedge\cdots\wedge \widehat{x_{j_k}}\wedge\cdots \wedge x_{j_i}$$

This differential verifies $\partial^2=0$ and makes $\KK$ a complex, which is called the {\bf Koszul complex}. This complex is exact and it is therefore a minimal free resolution of $\kb=R/\mathfrak{m}$, where $\mathfrak{m}=\langle x_1,\dots,x_n\rangle$, the irrelevant ideal in $R$. The exactness of the Koszul complex is a direct consequence of lemma \ref{contracting} later (see remark \ref{rem:contracting}). This exactness is a formulation of the fact that $\kb[x_1,\dots,x_n]$ is a Koszul algebra (see \cite{F99} for the definition and main properties of Koszul algebras).

Given a (multi)graded $R$-module $\Mc$, its {\bf Koszul complex} $(\KK(\Mc), \partial)$ is the tensor product complex $\Mc\otimes_R \KK=\Mc\otimes_R(R\otimes_\kb\wedge \Vc)\simeq\Mc\otimes_\kb\wedge\Vc$:
$$\KK(\Mc): 0\rightarrow \Mc\otimes_\kb\wedge^n\Vc\stackrel{\partial}{\rightarrow}\Mc\otimes_\kb\wedge^{n-1}\Vc\stackrel{\partial}{\rightarrow}\cdots \Mc\otimes\wedge^{1}\Vc\stackrel{\partial}{\rightarrow}\Mc\otimes\wedge^0\Vc\rightarrow 0$$

This complex is no longer acyclic in general, and we define the {\bf Koszul homology} of $\Mc$ as the homology of $\KK(\Mc)$. 

\subsubsection*{Grading, bigrading and multigrading of the Koszul complex}
\spanishsubsubsection{Grado,bigrado y multigrado del complejo de Koszul}
Consider an element of $R\otimes\wedge\Vc$ of the form $x^\mu\otimes x^J$ where $x^\mu=x_1^{\mu_1}\dots x_n^{\mu_n}$ and $x^J=x_1^{j_1}\wedge\dots\wedge x_n^{j_n}$ as before. We say that the {\it total degree} of $x^\mu\otimes x^J$ is $\mu_1+\cdots+\mu_n+j_1+\cdots+j_n$ and that the {\it total multidegree} of $x^\mu\otimes x^J$ is $(\mu_1+j_1,\cdots ,\mu_n+j_n)$. We also say that the \emph{symmetric degree} of $x^\mu\otimes x^J$ is $\mu_1+\cdots+\mu_n$ and that its \emph{exterior degree} is $j_1+\cdots+j_n$. Similarly we obtain the \emph{symmetric multidegree} and \emph{exterior multidegree} of $x^\mu\otimes x^J$.
Equivalently, if $J$ is given in the form $J=1\leq j_1<\cdots<j_i\leq n$ then the total degree of $x^\mu\otimes x^J$ is $\mu_1+\cdots+\mu_n+i$ and the total multidegree is $(\mu_1+[1\in J],\dots,\mu_n+[n\in J])$ where $[i\in J]$ equals $1$ if $i$ is in $J$ and $0$ otherwise. In this case, the \emph{symmetric} and \emph{exterior} degrees are given by $\mu_1+\cdots+\mu_n$ and $i$ respectively; symmetric and exterior \emph{multidegrees} are then found in the obvious way.

It is clear that for these elements, the Koszul differential preserves both the total degree and total multidegree. Thus, we can consider the following (multi)gradings in $\KK$ and $\KK(\Mc)$:
\begin{itemize}
\item {With respect to the total degree,we have

$$\KK=\bigoplus_{d\in \NN} \KK_d \qquad and \qquad \KK(\Mc)=\bigoplus_{d\in \NN} \KK_d(\Mc)$$

where $$\KK_d:0\rightarrow R_{d-n}\otimes\wedge^n\Vc\stackrel{\partial}{\rightarrow}R_{d-n+1}\otimes\wedge^{n-1}\Vc\stackrel{\partial}{\rightarrow}\cdots R_{d-1}\otimes\wedge^{1}\Vc\stackrel{\partial}{\rightarrow}R_{d}\otimes\wedge^0\Vc\rightarrow 0$$

and similarly for $\KK_d(\Mc)$. Note that if $q<d$ then $R_{q-d}\otimes\wedge^q\Vc=0$, so we have

$$\KK_d:0\rightarrow R_{0}\otimes\wedge^d\Vc\stackrel{\partial}{\rightarrow}R_{1}\otimes\wedge^{d-1}\Vc\stackrel{\partial}{\rightarrow}\cdots R_{d-1}\otimes\wedge^{1}\Vc\stackrel{\partial}{\rightarrow}R_{d}\otimes\wedge^0\Vc\rightarrow 0$$

Here, $R_l$ denotes the polynomials of degree $l$ and in $\KK(\Mc)$ we denote by $\Mc_l$ the degree $l$ component of $\Mc$. Because of this grading in $\KK$ and $\KK(\Mc)$, the homologies of them are also graded:
$$H_*(\KK)=\bigoplus_{d\in\NN}H_*(\KK_d)\qquad and \qquad H_*(\KK(\Mc))=\bigoplus_{d\in\NN}H_*(\KK_d(\Mc))$$

For each homological degree $p$ we have $H_p(\KK)=\bigoplus_{d\in\NN}H_p(\KK_d)=\bigoplus_{q+p=d}H_{q,p}(\KK)$; so we have a bigrading and we denote by $H_{q,p}(\KK)$ and $H_{q,p}(\KK(\Mc))$ the respective homology modules at $R_q\otimes\wedge^p\Vc$ and $\Mc_q\otimes\wedge^p\Vc$. We say that $q$ is the {\it symmetric degree} of $H_{q,p}(\KK)$ or $H_{q,p}(\KK(\Mc))$ and $p$ is its {\it exterior degree}.}

\item{With respect to the total multidegree, we have

$$\KK=\bigoplus_{\ab \in \NN^n} \KK_\ab \qquad and \qquad \KK(\Mc)=\bigoplus_{\ab \in \NN^n} \KK_\ab(\Mc)$$
where for every $\ab=(a_1,\dots,a_n)\in\NN^n$ with $a_{j_1},\dots,a_{j_l}\neq 0$ 
$$\KK_\ab:0\rightarrow R_{\ab-(0..1_{j_1}\dots 1_{j_l}..0)}\otimes\wedge^{(0..1_{j_1}\dots 1_{j_l}..0)}\Vc\stackrel{\partial}{\longrightarrow}\bigoplus_{\substack{\nu\subset\{1,\dots,l\}\\ \vert\nu\vert=l-1}}R_{\ab-\nu}\otimes\wedge^{\nu}\Vc\stackrel{\partial}{\longrightarrow}\cdots$$ $$\cdots\stackrel{\partial}{\longrightarrow} \bigoplus_{\substack{\nu\subset\{1,\dots,l\}\\ \vert\nu\vert=1}}R_{\ab-\nu}\otimes\wedge^{\nu}\Vc\stackrel{\partial}{\longrightarrow}R_{\ab}\otimes\wedge^0\Vc\rightarrow 0$$
and similarly for $\KK_\ab(\Mc)$.
Here, $R_\mu$ denotes the set of polynomials of multidegree $\mu$, and $\wedge^\nu\Vc$ denotes the span of $x_1^{\nu_1}\wedge\cdots\wedge x_n^{\nu_n}$. In this case, we have that the homologies of $\KK$ and $\KK(\Mc)$ are also multigraded:
$$H_*(\KK)=\bigoplus_{\ab\in\NN^n} H_{*}(\KK_{\ab})=\bigoplus_{\ab\in\NN^n} H_{\ab}(\KK)\qquad H_p(\KK)=\bigoplus_{\ab\in\NN^n}H_p(\KK_\ab)=H_{p,\ab}(\KK)$$
$$H_*(\KK(\Mc)=\bigoplus_{\ab\in\NN^n} H_{*}(\KK_{\ab}(\Mc))=\bigoplus_{\ab\in\NN^n} H_{\ab}(\KK(\Mc))$$
$$H_p(\KK(\Mc))=\bigoplus_{\ab\in\NN^n}H_p(\KK_\ab(\Mc))=H_{p,\ab}(\KK(\Mc))$$
}
\end{itemize}

\begin{Remark}
In the context of commutative algebra and following the terminology of \cite{AB58,E95,BH98}, the Koszul complex is defined for any commutative ring $R$ and a set of elements of $R$ or even for a set of homogeneous $R$-linear forms. An explicit definition in such context, which appears for example in \cite{Si07} is the following:

For any homogeneous forms $f_1,\dots,f_r\in R$, the Koszul complex $\KK(f_1,\dots,f_r)$ is a complex of free modules $F_i$, $0\leq i\leq r$. Letting $[r]=\{1,\dots,r\}$, we can describe the $i$-th module as
$$F_i=\bigoplus_{\sigma\subseteq[r],\vert\sigma\vert=i}R(-deg(\prod_{j\in\sigma}f_j).$$
The map from $F_i$ to $F_{i-1}$ is given by an ${{r}\choose{i-1}}\times{{r}\choose{i}}$ matrix whose $(\sigma,\tau)$ entry is $0$ if $\sigma$ is not contained in $\tau$. Otherwise $\tau=\sigma\cup\{j\}$ and the $(\sigma,\tau)$ entry is equal to $(-1)^{\vert \tau_{<j}\vert}f_j$ where $\tau_{<j}=\{l \in \tau\vert l<j\}$.

The formulation of the Koszul complex given above and used in this thesis is just $\KK(x_1,\dots,x_n)$.  
\end{Remark}

\subsection{Koszul complex and $Tor$}\label{tor}
\spanishsubsection{El complejo de Koszul y $Tor$}
From the definitions, it is clear that we can identify the Koszul homology modules with $Tor_\bullet^R(\Mc,\kb)$: We have a resolution of $\kb$ (the Koszul Complex) to which we have applied the functor $\Mc\otimes_R-$. The homology of the resulting complex is by definition $Tor_\bullet^R(\Mc,\kb)$ (see Appendix \ref{apal}).

Another way of computing the Koszul homology of $\Mc$ would start with a resolution $\PP: \cdots \rightarrow P_i\stackrel{\delta_i}{\rightarrow} P_{i-1}\rightarrow\cdots\rightarrow P_0$ of $\Mc$ and then compute the homology of $\PP \otimes_R \kb$. This homology is independent of the chosen resolution of $\Mc$ or of $\kb$. If $\PP$ is minimal,the differentials are given by matrices with polynomial coefficients, none of which is a nonzero constant (see Appendix \ref{apal}). Thus, tensoring with $\kb$ yields the zero differential everywhere and then the number of generators of each $P_i$ equals $Tor^R_i(\Mc,\kb)$ and the dimension of the $i$-th Koszul homology of $\Mc$. If $\PP$ is not minimal one could either minimize it with some standard procedure (see for example \cite{CLO98} or the chain complex reduction algorithm below) or compute the homology of the resulting resolution $\PP \otimes_R \kb$. The $Tor$ modules inherit the gradings, bigradings and multigradings we have seen before.

\subsection{Koszul homology and Spencer cohomology}\label{spencer}
\spanishsubsection{Homolog\'ia de Koszul y cohomolog\'ia de Spencer}
Spencer cohomology is an important object in the formal theory of differential systems, in particular in the study of involution from a homological point of view \cite{S07}, and also in the measure of the dimension of the solutions space of PDE systems \cite{KL06}. First introduced by Spencer \cite{S69}, this cohomology is associated to symbol (co)modules of differential systems and a result by Serre relates the computation of this cohomology (or equivalently the {\it Cotor} modules) to the Cartan test for involution. Details can be seen in \cite{LS02} and \cite{S07}. Here we will only make a brief presentation of the Spencer cohomology and the duality between it and Koszul homology.

\subsubsection*{The polynomial de Rham complex and the Koszul complex}
\spanishsubsubsection{El complejo de de Rham polynomial y el complejo de Koszul}
Let again $\Vc$ be a $n$-dimensional $\kb$-vector space, and  $S\Vc$, $\wedge \Vc$ the symmetric and exterior algebras of $\Vc$. The {\it polynomial de Rham complex} at degree $q$, $\Rc_q(S\Vc)$ is given by
$$0\rightarrow S_q\Vc\rightarrow S_{q-1}\Vc\otimes\Vc\stackrel{\delta}{\rightarrow}S_{q-2}\Vc\otimes\wedge^2\Vc\stackrel{\delta}{\rightarrow}\cdots S_{q-n}\Vc\otimes\wedge^n\Vc\rightarrow 0$$

where the differential $\delta$ is given by
$$\delta(x_1^{\mu_1}\dots x_n^{\mu_n}\otimes x_{j_1}\wedge\cdots\wedge x_{j_i})=\sum_{k=1}^n \mu_k\cdot  x_1^{\mu_1}\cdots x_k^{\mu_k-1} \cdots x_n^{\mu_n}\otimes x_{k}\wedge x_{j_1}\wedge\cdots\wedge x_{j_i}$$

This differential satisfies $\delta^2=0$, and it can be seen as the exterior derivative applied to a differential $p$ form with polynomial coefficients, whence the name `polynomial de Rham complex'. The complex is also exact for all values $q>0$, the only nonvanishing cohomology modules are $H^{0,0}(\Rc(S\Vc),\delta)=\kb$. This result is known as the {\it formal Poincar\'e lemma}; the polynomial de Rham complex is then a free coresolution of $\kb$.

A first relation between the polynomial de Rham and the Koszul complexes is given by the fact that their differentials are contracting homotopies for the other complex, i.e.

\begin{Lemma}\label{contracting} We have $(\delta\circ\partial+\partial\circ\delta)(w)=(p+q)w$ for all 
$w\in S_q\Vc\otimes\wedge_p\Vc$
\end{Lemma}

The proof of this lemma is a straightforward verification, it can be found for example in \cite{S07}.

\begin{Remark}\label{rem:contracting}
As a consequence of this lemma we can easily proof the exactness of the Koszul complex: $\delta$ induces a contracting homotopy for $\KK$ connecting the identity and zero maps. The existence of such homotopy implies the exactness of $\KK$.
\end{Remark}

\subsubsection*{Spencer cohomology and Koszul homology}
\spanishsubsubsection{Cohomolog\'ia de Spencer y homolog\'ia de Koszul}
The duality between Spencer cohomology and Koszul homology comes from the duality between the polynomial de Rham and the Koszul complexes. A way to express this duality is given in \cite{S07}, which we follow here: If we apply the functor $Hom_{R}(\cdot,R)$ to some complex of $R$-modules we obtain its dual complex. In our case, we have that there exists a canonical isomorphism between $S_q(\Vc^*)$ and $(S_q\Vc)^*$, and the same holds for exterior products. Thus, we have a canonical isomorphism $S_q(\Vc^*)\otimes \wedge^p(\Vc^*)\simeq (S_q\Vc\otimes \wedge^p\Vc)^* $.
Choosing appropriate dual basis one can show that $\partial$ is the pull back of $\delta$ and one obtains that
 
\begin{Proposition}
 $(\Rc(S\Vc)^*,\delta^*)$ is isomorphic to $(\KK(S(\Vc^*)),\partial)$
\end{Proposition}

To see the duality at the homology level, we need to define Spencer cohomology and Koszul homology in terms of modules and comodules of the symmetric bialgebra \footnote{The definitions of coalgebra, comodule and cotensor product can be found in \cite{EM66,M95}.}:
\begin{Definition}
Let $\Nc$ be a comodule over the symmetric coalgebra $\mathfrak{S}\Vc$. Its Spencer complex $(\Rc(\Nc),\delta)$ is the cotensor product complex $N\boxtimes \Rc(\mathfrak{S}\Vc)$. The Spencer cohomology of $\Nc$ is the corresponding bigraded cohomology; the cohomology group at $\Nc_q\otimes\wedge^p\Vc$ in $(\Rc_{q+p}(\Nc),\delta)$ is denoted by $H^{q,p}(\Nc)$
\end{Definition}
For symmetry with the precedent definition, we recall here the definition of Koszul complex and homology in terms that emphasize the parallelism with the Spencer complex and homology.
\begin{Definition}
Let $\Mc$ be a graded module over the symmetric algebra $S\Vc$. Its Koszul complex $(\KK(\Mc),\partial)$ is the tensor product complex $\Mc\otimes \KK(S\Vc)$. The Koszul homology of $\Mc$ is the corresponding bigraded homology; the homology group at $\Mc_q\otimes\wedge^p\Vc$ in $(\KK_{q+p}(\Mc),\delta)$ is denoted by $H_{q,p}(\Mc)$
\end{Definition}

\begin{Remark}
The reason that we can use the symmetric coalgebra instead of the symmetric algebra in the definition of Spencer cohomology, is that the de Rham differential $\delta$ exploits only the vector space structure of the symmetric algebra $S\Vc$, and thus we may substitute it by $\mathfrak{S}\Vc$ and define $\delta$ on the components of the free $\mathfrak{S}\Vc$-comodule $\mathfrak{S}\Vc\otimes\wedge\Vc$, since both are identical as vector spaces. The duality of the two definitions is then more evident, and this comodule interpretation is more natural in some contexts, see \cite{S07}.
\end{Remark}

A useful application of the duality between Koszul and Spencer (co)homologies is given in section \ref{formal_theory} where the vanishing of certain (co)homology groups is used as a criterion to detect the involution of symbolic systems (see theorem \ref{involution-symbolic-system} and definition \ref{deg_involution}) which is later applied to differential systems. This homological treatment is the basis of the approach in \cite{S07b}.

We give here an algebraic definition of involution for symbolic systems, the application to differential system will be seen in section \ref{formal_theory}:

\begin{Definition}\label{symbolic-system}
Let $N_q\subseteq S_q(\Vc^*)\otimes \Uc$ be a vector subspace ($\Uc$ being a further finite-dimensional vector space). Its prolongation is the subspace
$$N_{q,1}=\{f\in S_{q+1}(\Vc^*)\otimes \Uc\vert \delta(f)\in N_q\otimes \Vc^*\}.$$
A sequence of vector subspaces $(N_q\subseteq S_q(\Vc^*)\otimes \Uc)_{q\in\NN}$ is called a {\it symbolic system} over $\Vc^*$, if $\Nc_{q+1}\subseteq\Nc_{q,1}$ for all $q\in \NN$ note that we can also introduce prolongations as $N_{q,1}=(\Vc\otimes N_q)\cap(S_{q+1}(\Vc^*)\otimes \Uc)$.
\end{Definition}

\begin{Remark}
The new vector space $\Uc$ makes us extend our complexes to tensor product complexes $\Rc(S(\Vc^*)\otimes \Uc)$ and $\KK(S\Vc\otimes \Uc^*)$, but everything remains valid with trivial modifications, since the differentials of the complexes are essentially the same.
\end{Remark}

\begin{Lemma}\label{symbolic-system-comodule}
Let $(\Nc_q)_{q\in\NN}$ be a symbolic system. Then $\Nc=\bigoplus_{q=0}^\infty \Nc_q$ is a graded (right) comodule of the free $\mathfrak{S}(\Vc^*)$-comodule $\mathfrak{S}(\Vc^*)\otimes \Uc$. Conversely, the sequence $(\Nc_q)_{q\in\NN}$ of the components of any graded (right) comodule $\Nc\subseteq\mathfrak{S}(\Vc^*)\otimes \Uc$ defines a symbolic system.
\end{Lemma}

For a symbolic system, we have the following results concerning Spencer cohomology (and hence, Koszul homology):

\begin{Lemma}\label{homology-involutive-system}
Let $\Nc\subseteq \mathfrak{S}(\Vc^*)\otimes \Uc$ be a symbolic system. Then $H^{q,0}(\Nc)=0$ and $dimH^{q-1,1}(\Nc)=dim(\Nc_{q-1,1}/\Nc_q)$ for all $q>0$.
\end{Lemma}

\begin{Theorem}\label{involution-symbolic-system}
Let $N\subseteq\mathfrak{S}(\Vc^*)\otimes \Uc$ be a symbolic system. There exists an integer $q_0\geq0$ such that $H^{q,p}(N)=0$ for all $q\geq q_0$ and $0\leq q\leq n$. Dually, let $\Mc$ be a finitely generated graded polynomial module. There exists an integer $q_0\geq 0$ such that $H_{q,p}(\Mc)=0$ for all $q\geq q_0$ and $0\leq p\leq n$.
\end{Theorem}

The proof of this important theorem is very simple in this context, it is based on the following lemma:
\begin{Lemma}\label{multiplication-zero-map}
Let $\Mc$ be a graded $R$-module. Multiplication by an arbitrary element of $S_+\Vc$ induces the zero map on the Koszul homology $H_*(\KK(\Mc))$
\end{Lemma}

\noindent{\bf Proof: } If $w\in \Mc_q\otimes\wedge_p\Vc$ is a cycle, then for any $v\in \Vc$, the form $vw$ is a boundary. Indeed $\partial(v\wedge w)=-v\wedge(\partial w)+vw=vw$. Since $\partial$ is $S\Vc$-linear, this observation remains true when we take for $v$ an arbitrary element of $S_+\Vc$, i.e. any polynomial without constant term.

\noindent{\bf Proof of Theorem \ref{involution-symbolic-system}:}
The cycles in $\Mc\otimes \wedge^p\Vc$ form a finitely generated $S\Vc$-module. Thus, there exists an integer $q_0$ such that the polynomial degree of all elements in a finite generating set of it is less than $q_0$. All cycles of higher polynomial degree are then linear combinations of these generators with polynomial coefficients without constant terms. By lemma \ref{multiplication-zero-map}, they are therefore boundaries. Hence $H_{q,p}(\KK(\Mc))=0$ for all $q\geq q_0$. $\square$
\begin{Definition}\label{deg_involution}
The {\it degree of involution} of the polynomial comodule $\Nc$ is the smallest value $q_0$ such that $H^{q,p}(\Nc)=0$ for all $q\geq q_0$ and $0\leq p\leq n$. More generally, we say that $\Nc$ is $s$-acyclic at degree $q_0$ if $H^{q,p}(\Nc)=0$ for all $q\geq q_0$ and $0\leq p\leq s$. A comodule that is $n$-acyclic at degree $q_0$ is called involutive at degree $q_0$. Dually, we call a polynomial module $\Mc$ involutive at degree $q_0$ if $H_{q,p}(\Mc)=0$ for all $q\geq q_0$ and $0\leq p\leq n$.
\end{Definition}

\begin{Remark}
The degree of involution of a polynomial module is exactly the Castelnuovo-Mumford regularity, as can be seen from the definition (see appendix \ref{apal}). The equivalence between these two fundamental notions has been almost unnoticed in the literature since the two concepts appear in very distant contexts; however, it was implicitly present in \cite{S02d}, and appears explicitly in \cite{M03}. The subjacent reason for this equivalence, and what made it evident is the duality between Spencer cohomology and Koszul homology, from which the notions of involution and Castelnuovo-Mumford regularity can be defined in their respective contexts.
\end{Remark}

\section{Monomial Ideals}\label{mon}
\spanishsection{Ideales monomiales}
Much attention has been paid in the last decade to monomial ideals and much work has been done around them from different points of view. Their place at the intersection of commutative algebra, algebraic geometry and combinatorics provides with many examples of interaction between algebraic and combinatorial concepts, from which big developments have resulted, see for instance the books by Bruns and Herzog \cite{BH98} and Miller and Sturmfels \cite{MS04}, which together with the references given in the next paragraph can give the interested reader a good view of the interaction between combinatorics and commutative algebra.

At least three different but complementary points of view have been considered when dealing with monomial ideals. Of course the interaction between the three is big and it is a main characteristic of the topic:
\begin{itemize}
\item The first one is the algebraic one. Of course without forgetting the combinatorial nature of monomial ideals, this approach focuses on the properties of monomial ideals and algebras as algebraic objects in relation with the polynomial ring. Here we can consider the book by R. Villarreal \cite{V01} as a basic reference. Rees algebras and their algebraic properties can be considered as one main object in this approach.
\item A second approach, which is somehow considered as the origin of the interest in monomial ideals is the combinatorial one. Again, without forgetting the focus on algebra, this approach studies the combinatorial properties of monomial ideals in relation with simplicial topology, graph theory, etc. Stanley's monograph \cite{S96} can be considered as the main reference, and Stanley-Reisner ideals as the main object in this approach. Stanley-Reisner theory deals mainly with squarefree monomial ideals, a bit of it is presented in section \ref{stanley-alexander}. Another recent and interesting too, used in this context is \emph{discrete Morse theory}, see \cite{B02,OW07}.
\item In the recent years, the interaction of combinatorial, algebraic and computational methods have become more evident, and this interaction has given rise to a third approach to the study of monomial ideals. A representative book of this style of research is the one by Miller and Sturmfels \cite{MS04}, which puts together many of the ideas in their work and those of Bayer, Peeva, etc. (see references later). Resolutions, Betti numbers and Alexander duality are the main objects in this approach, and the computational aspects give it a particular flavour.
\end{itemize}

Here we focus on the Koszul homology of monomial ideals and its properties and computation. Our point of view is probably closest to the third approach described above, since we are interested in Betti numbers, resolutions, and other homological invariants of these ideals. First of all we will make a description of monomial ideals and how their combinatorial nature provides good expressions for many algebraic properties; the approach represented by Villarreal will be therefore used in this part. Second we will see how the Koszul homology of a monomial ideal also interacts very deeply with the combinatorial nature of it and how simplicial complexes are a useful tool to describe this interaction. The approach represented by Miller-Sturmfels' book will be preferred here. Not forgetting Stanley's approach, this will be less used, due to the nature of our goal and methods. Some other considerations will be made on homological methods that can also be applied to get some insight in the nature of the Koszul homology of monomial ideals.

\subsection{Basic terminology}
\spanishsubsection{Terminolog\'ia b\'asica}
After defining  monomial ideals and their basic properties, we see here the identification between monomial ideals of the polynomial ring in $n$ variables and monoid ideals of the monoid $\NN^n$. This identification will be very useful for dealing with monomial ideals. In particular, when $n$ is small, we can draw helpful graphical representations of our ideals, so called \emph{staircase diagrams}. Recall that $R=\kb[x_1,\dots,x_n]$.

\subsubsection{Monomial ideals and multigrading}
\spanishsubsubsection{Ideales monomiales y multigrado}
\begin{Definition}
A \emph{monomial} in $R$ is a product $x^{\ab}=x_1^{a_1}\cdots x_n^{a_n}$ with $a_i\geq 0\quad \forall i$. We say that $\ab=(a_1,\dots,a_n)\in\NN^n$ is the {\rm multidegree} of $x^{\ab}$. An ideal $I\subset R$ is called a {\rm monomial ideal} if it is generated by monomials. If $I$ is a monomial ideal, the quotient ring $R/I$ is called a {\rm monomial ring}.
\end{Definition}

Particularly interesting monomial ideals are {\it squarefree} and {\it face} ideals:

\begin{Definition}
A {\rm squarefree ideal} is a monomial ideal generated by squarefree monomials, i.e. such that their exponent in each variable is either $0$ or $1$. A {\rm face ideal} is an ideal $J$ of $R$ generated by a subset of the set of variables.
\end{Definition}

\begin{Remark}
Some important properties of monomial ideals, that constitute characterizations of them, are the following:
\begin{itemize}
\item If $I$ is a monomial ideal generated by a set of monomials $\{x^\mu\vert \mu\in M\subset\NN^n\}$ then a monomial $x^\nu$ is in $I$ if and only if $x^\nu$ is divisible for some generator $x^\mu$ of $I$.
\item In general, if $f$ is a polynomial in $\poly\kb x n$ then $f$ is in $I$ if and only if every monomial in $f$ lies in $I$, if and only if $f$ is a $\kb$-linear combination of the monomials in $I$.
\item An ideal $I$ in $\poly\kb x n$ is monomial if and only if I is torus fixed, i.e. if $(c_1,\dots,c_n)\in (\kb^*)^n$, then $I$ is fixed under the action $x_i\mapsto c_i x_i$ for all $i$ \cite{HS01}.
\end{itemize}
\end{Remark}

As a consequence of these characterizations, it is easy to see that a monomial ideal is uniquely determined by its monomials, i.e. two monomial ideals are the same if and only if they contain the same monomials. An important result for monomial ideals is the so called Dickson's Lemma, which states that all monomial ideals of $R$ are finitely generated, see for example \cite{CLO96} for a proof. Moreover, this minimal set of monomial generators is unique.

As a vector space, $R=\bigoplus_{\ab\in\NN^n} R_{\ab}$, where $R_\ab$ is the vector subspace generated by $x^\ab$. Since given $\ab,\,\bb\in\NN^n$, $R_\ab\cdot R_\bb \subseteq R_{\ab+\bb}$, we say that $R$ is an $\NN^n$-graded $\kb$-algebra, being the monomial ideals its $\NN^n$-graded ideals. Thus, a monomial ideal $I$ can be expressed as a direct sum of its  $\NN^n$-graded (i.e. multigraded) components.

Another interesting issue is that the family of monomial ideals of $R$, ordered by inclusion, forms a distributive lattice under the operations $I\wedge J=I \cap J$ and $I\vee J=I+J$ (see \cite{V01}).

\subsubsection{Staircase diagrams, lcm-lattice, and the combinatorial nature of monomial ideals}
\spanishsubsubsection{Diagramas en escalera, ret\'iculo de $mcd$ y car\'acter combinatorio de los ideales monomiales}
Let us consider the abelian monoid $(\NN^n,+)$ of $n$-multi indices (or just multi indices when $n$ is clear from the context) with addition defined componentwise. This can, on one side, be identified with the vertices of an $n$-dimensional integer lattice, and on the other side, with the abelian monoid of monomials in $\poly \kb x n$ with the usual product, just by identifying the exponent of a monomial with the correspondent multiindex. Setting $log(x^\mu)=\mu$ gives this correspondence. Given a set $S$ of monomials, we say $log(S)$ is just the set of all $log(x^\mu)$ with $x^\mu\in S$. Thus, the set of multiples of a given monomial $x^\mu$ can be identified with the set $\mu + \NN^n=\{\mu +\nu\vert \nu\in \NN^n\}$. It is a monoid ideal. If we call this the {\it span} of $\mu$ in $(\NN^n,+)$, then the set of monomials in the monomial ideal $I=\langle x^{\mu_1},\dots,x^{\mu_r}\rangle$ can be identified with the union of the spans of its generators, and with the corresponding vertices in the $n$-dimensional integer lattice; we call this the {\it lattice span} of $I$. This union of spans is just $log(I)$.

In the low dimensional cases, where $n$ equals 2 or 3, we can draw the lattice span of a given ideal $I$. These drawings are known as the two or three-dimensional {\it staircase diagrams} of monomial ideals and are useful tools for studying monomial ideals in two or three variables, see \cite{MS04}. Their analogues in higher dimensions are also useful although obviously do not have such nice graphical representations.

Figures \ref{stair1} and \ref{stair2} show some examples of staircase diagrams. In figure \ref{stair1} the colored squares represent the multidegrees of the monomials inside the ideal $I$, whereas in figure \ref{stair2} the colored cubes represent the multidegrees of the monomials not in $I$.

\begin{figure}
\begin{center}
\begin{tikzpicture}[scale=1.5]
\fill[blue!40!white] (0,1.5)rectangle(2.4,2.4);
\fill[blue!40!white] (1,1)rectangle(2.4,2.4);
\fill[blue!40!white] (1.5,0)rectangle(2.4,2.4);
\draw [step=.5cm, very thin](0,0) grid (2.4,2.4);
\fill (0,1.5) circle(1pt) node[left]{$y^3$};
\fill (1.5,0) circle(1pt) node[below]{$x^3$};
\fill (1,1) circle(1pt) node[above left]{$x^2y^2$};
\end{tikzpicture}
\caption{Staircase diagram of $I=\langle x^3,x^2y^2,y^3\rangle$}\label{stair1}
\end{center}
\end{figure}
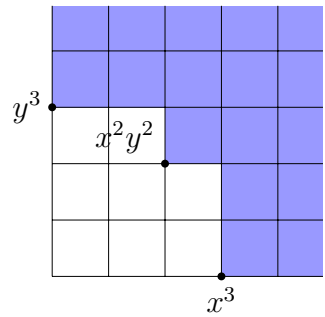

\begin{figure}
\begin{center}
\begin{tikzpicture}[scale=1]
\draw [very thin, ->](0,0)--(-0.7,-0.7);
\draw [very thin, ->](4.5,1.5)--+(1,0);
\draw [very thin, ->](1.5,4.5)--+(0,1);
\foreach \position in {(0,0),(0,1),(0,2),(1.5,0.5),(1.5,1.5),(1.5,2.5),(2.5,0.5),(2.5,1.5),(2.5,2.5)} \filldraw[fill=blue!75,draw=black]\position rectangle +(1,1);
\foreach \position in {(1,0),(1,1),(1,2),(3.5,0.5),(3.5,1.5),(3.5,2.5),(4,1),(4,2),(4,3)}
\filldraw[fill=blue!40!white,draw=black] \position-- ++(0.5,0.5)-- ++(0,1)-- ++(-0.5,-0.5)-- cycle;
\foreach \position in {(0,3),(0.5,3.5),(1,4),(1.5,3.5),(2,4),(2.5,3.5),(3,4)}
\filldraw[fill=blue!15!white,draw=black] \position-- ++(1,0)-- ++(0.5,0.5)-- ++(-1,0)-- cycle;
\fill (0,0) circle(2pt) node[left]{$x^3$};
\fill (1.5,0.5) circle(2pt) node[below right]{$x^2y$};
\fill (4.5,1.5) circle(2pt) node[above right]{$y^3$};
\fill (1.5,4.5) circle(2pt) node[above left]{$z^3$};
\end{tikzpicture}
\caption{Staircase diagram of $I=\langle x^3,x^2y,y^3,z^3\rangle$}\label{stair2}
\end{center}
\end{figure}

In general, the correspondence of monomial ideals and monoid ideals in $(\NN^n,+)$ allows us to interpret some relations and operations such as divisibility and multiplication in terms of addition or substraction of multi indices. This gives us a good tool to easily verify ideal membership of a given monomial or polynomial with respect to a monomial ideal, to compute least common multiples or great common divisors of sets of monomials, or even compute more complicated objects such as intersection of ideals, union, colon ideals... (for detailed descriptions see, for example \cite{B96,MS04,MP01}). More explicitly:
\begin{itemize}
\item A monomial $x^\mu$ divides a monomial $x^\nu$ where $\mu=(\mu_1,\dots,\mu_n)$ and $\nu=(\nu_1,\dots,\nu_n)$ if and only if $\nu_i\geq\mu_i$ for all $1\leq i\leq n$. A monomial $x^\nu$ is in the ideal generated by $\{x^{\mu_1},\dots x^{\mu_r}\}$ if it is divisible by some generator.
\item The {\it least common multiple} of two monomials $x^\mu$ and $x^\nu$ is the monomial $lcm(x^\mu,x^\nu)=x^\rho$ with $\rho=(\max(\mu_1,\nu_1),\dots,\max(\mu_n,\nu_n))$. The \emph{greatest common divisor} of $x^\mu$ and $x^\nu$ is the monomial $gcd(x^\mu,x^\nu)=x^\sigma$ with $\sigma=(\min(\mu_1,\nu_1),\dots,\min(\mu_n,\nu_n))$.
\item If $I=\langle x^{\mu_1},\dots x^{\mu_r}\rangle$ and $J=\langle x^{\nu_1},\dots x^{\nu_s}\rangle$ are monomial ideals, then we have that $I+J=\langle x^{\mu_1},\dots x^{\mu_r},x^{\nu_1},\dots x^{\nu_s}\rangle$ and $I\cap J=\langle\{lcm(x^{\mu_i},x^{\nu_j}) |1\leq i \leq r,1\leq j \leq s\}\rangle$.
\item If $I=\langle x^{\mu_1},\dots x^{\mu_r}\rangle$ is a monomial ideal and $x^\nu$ is a monomial not in $I$, then the quotient ideal $(I:x^\nu)=\{f\in R\vert f\cdot x^\nu\in I\}$ is given by $$(I:x^\nu)=\langle \frac{x^{\mu_1}}{gcd(x^{\mu_1},x^{\nu})},\dots, \frac{x^{\mu_r}}{gcd(x^{\mu_r},x^{\nu})}\rangle=\langle \frac{lcm(x^{\mu_1},x^{\nu})}{x^{\nu}},\dots, \frac{lcm(x^{\mu_r},x^{\nu})}{x^{\nu}}\rangle$$
In particular, $x^\nu$ is said to be {\it regular} for $I$ if it is not a zerodivisor in $R/I$, if and only if $(I:x^\nu)=I$ if and only if $gcd(x^{\mu_i},x^{\nu})=1$ for all $i$, if and only if $lcm(x^{\mu_i},x^{\nu})=x^{\mu_i}\cdot x^\nu$ for all $i$. So, for monomial ideals we have an easy way to characterize {\it regular sequences} (see the definition in appendix \ref{apal}). 
\end{itemize}

Another interesting object related to the combinatorial nature of monomial ideals is the lcm-lattice, see \cite{GPW99}.  Given a monomial ideal $I$ minimally generated by $\{m_1,\dots,m_r\}$, we denote by $L_I$ the lattice with elements labeled by the least common multiples of subsets of $\{m_1,\dots, m_r\}$ ordered by divisibility. We call $L_I$ the {\it lcm-Lattice} of $I$. The minimal element of this lattice is $1$ (lcm of the empty set) and the maximal element is  $lcm(m_1,\dots,m_r)$. We will also denote by $L_{I,i}$ the subset of $L_I$ given by those elements labelled by the $lcm$ of exactly $i$ monomials in $\{m_1,\dots,m_r\}$. Of course, $L_I=\bigcup_{i=0,\dots,r}L_{I,i}$. As we have seen this lattice can be very easily translated into a subset of the integer lattice in $n$ variables. The lcm-lattice $L_I$ will be very useful when discussing Hilbert series, Betti numbers and resolutions of monomial ideals, in particular the Taylor resolution (see later).

\subsection{The algebra of monomial ideals}\label{algebra_monomial}
\spanishsubsection{El \'algebra de los ideales monomiales}
 Many interesting algebraic objects associated with ideals of the polynomial ring have closed form expressions when we deal with monomial ideals. The following is a non-exhaustive list of the most relevant of these objects and/or those which will play some role on the subsequent chapters. Most definitions and properties in this section come from \cite{MS04} and \cite{V01}. For the definition of the several algebraic objects involved, see Appendix \ref{apal}.

\subsubsection{Prime, primary and irreducible monomial ideals.}\label{prime-primary-irreducible}
\spanishsubsubsection{Ideales monomiales primos, primarios e irreducibles}
Prime, primary and irreducible ideals (see appendix \ref{apal}) are particularly important classes of ideals. Their description, characterizations and some important properties can be easily obtained in the case of monomial ideals.

\begin{Proposition}\label{prime}
A monomial ideal is a prime ideal if and only if it is a face ideal i.e. it is generated by some set of variables.
\end{Proposition}

\begin{Proposition}\label{primary}
A monomial ideal $I$ is a primary ideal if and only if, after permutation of the variables, $I$ has the form
$I=\langle x_1^{a_1},\cdots,x_r^{a_r},x^{\bb_1},\dots,x^{\bb_s}\rangle$, where $a_i\geq 1$ and $\cup_{i=1}^s supp(x^{\bb_i})\subset\{x_1,\dots,x_r\}$, $r\leq n$. In particular, the powers of face ideals are primary.
\end{Proposition}

\begin{Proposition}\label{irreducible}
A monomial ideal $I$ is irreducible if and only if up to permutation of the variables, $I$ can be written as $I=\langle x_1^{a_1},\dots,x_r^{a_r}\rangle$ where $a_i>0$ for all $i$, $r\leq n$.
\end{Proposition}

With these descriptions, some interesting properties about these ideals and decompositions of any monomial ideal in terms of special ideals are given. Some properties are more easily stated when speaking about squarefree monomial ideals, which are of great importance due to their relation with simplicial complexes, see \cite{MS04,S78,V01}.
\begin{Proposition}
Let $I$ be a monomial ideal of $\kb[x_1,\dots,x_n]$ with $\kb$ a field, then every associated prime of $I$ is a face ideal.
\end{Proposition}
In the case of squarefree monomial ideals we even have the following relation:
\begin{Corollary}
If $I$ is a squarefree monomial ideal and $J_1,\dots,J_s$ are the associated primes of $I$ then $I=J_1\cap\cdots\cap J_n$. 
\end{Corollary}

\begin{Proposition}
A monomial ideal $I\subset R$ is squarefree if and only if any of the following conditions hold:
\begin{enumerate}
\item $I$ is an intersection of prime ideals.
\item $rad(I)=I$
\item A monomial $m$ is in $I$ if and only if $x_1\cdots x_r\in I$, where $supp(m)=\{x_i\}_{i=1}^r$
\end{enumerate}
\end{Proposition}

\begin{Proposition}
Let $I$ be a monomial ideal, then $I$ has an irredundant primary decomposition $I=J_1\cap\cdots\cap J_r$ where $J_k$ is a primary monomial ideal for all $k$ and $rad(J_k)\neq rad(J_l)$ if $k\neq l$.
\end{Proposition}
This proposition provides with an algorithm for computing minimal primary decompositions of monomial ideals (which need not to be unique), by successive elimination of powers of variables. A detailed description can be found in \cite{V01}.

\begin{Proposition}
If $I$ is a monomial ideal, then there is a unique irredundant decomposition $I=J_1\cap\cdots\cap J_r$ such that every $J_i$ is an irreducible monomial ideal.
\end{Proposition}

\begin{Remark}
Irreducible components of monomial ideals are very closely related to the minimal generators of their Alexander duals. A complete description of the role of this beautiful duality when studying monomial ideals can be found in \cite{MS04} and \cite{M00}, with several examples and other applications. We will meet Alexander duality in section \ref{hom_tools}.
\end{Remark}
\subsubsection{Integral closure of monomial ideals.}
\spanishsubsubsection{Clausura integral de ideales monomiales}
Integral closures (see definition in Appendix \ref{apal}) have become an active area of research within commutative and computer algebra in the recent years (see for example \cite{V05}) in the case of monomial ideals, again their combinatorial properties provide us with characterizations and description that makes this case easier to handle. First of all the following proposition states that the integral closure of a monomial ideal is a monomial ideal, and gives a description (cf. \cite{V01,C02}).
\begin{Proposition}
Let $I\subset R$ be a monomial ideal, then its integral closure $\bar{I}$ is also a monomial ideal. We have that $\bar{I}=\langle m\vert m^l\in I^l \mbox{ for some } l\rangle$
\end{Proposition}
With the help of the staircase diagrams we have seen in the precedent section and the correspondence provided by the $log$ map, we can give a geometric description of the integral closure of a monomial ideal. For this description we need to recall the definition of the (rational) convex hull of a set of elements in the lattice $\NN^n$:
\begin{Definition}
Let $\mu_i=(\mu_{i_1},\dots,\mu_{i_n})\in \NN^n$, then its (rational) convex hull is 
$$conv(\mu_1,\dots,\mu_r)=\left\lbrace \sum_{i=1}^r\lambda_i \mu_i \mbox{ such that } \sum_{i=1}^r\lambda_i=1,\lambda_i\in\QQ_+\right\rbrace $$
\end{Definition}
With this definition, we have 
\begin{Proposition}
Let $I$ be a monomial ideal, then the integral closure of $I$ is given by
$$\bar{I}=(x^\mu\vert \mu\in conv(log(I))\cap\NN^n)$$
\end{Proposition}

\begin{Example}
Consider the ideal in $\RR[x,y]$ given by $I=\langle xy^{15},x^3y^{10},x^6y^5,x^{11}y\rangle$ the ideal itself and its integral closure are depicted in figure \ref{integral2dim}. The integral closure of $I$ is given by
$$\bar I=\langle xy^{15},x^2y^{13},x^3y^{10},x^4y^9,x^5y^7,x^6y^5,x^8y^4,x^9y^3,x^{10}y^2,x^{11}y\rangle$$
In the picture, the new generators are drawn with white dots.
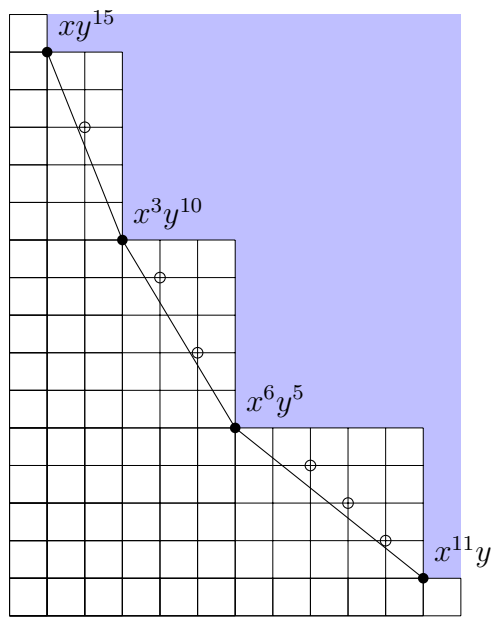
\begin{figure}
\begin{center}
\begin{tikzpicture}[scale=2]
\fill[blue!25!white] (0.25,3.75)rectangle(3,4);
\fill[blue!25!white] (0.75,2.5)rectangle(3,4);
\fill[blue!25!white] (1.5,1.25)rectangle(3,4);
\fill[blue!25!white] (2.75,0.25)rectangle(3,4);
\draw [step=.25cm, very thin](0,0) grid (0.25,4);
\draw [step=.25cm, very thin](0,0) grid (3,0.25);
\draw [step=.25cm, very thin](0,0) grid (0.75,3.75);
\draw [step=.25cm, very thin](0,0) grid (1.5,2.5);
\draw [step=.25cm, very thin](0,0) grid (2.75,1.25);
\fill (0.25,3.75) circle(1pt) node[above right]{$xy^{15}$};
\fill (0.75,2.5) circle(1pt) node[above right]{$x^3y^{10}$};
\fill (1.5,1.25) circle(1pt) node[above right]{$x^6y^5$};
\fill (2.75,0.25) circle(1pt) node[above right]{$x^{11}y$};
\draw (0.5,3.25) circle(1pt);
\draw (1,2.25) circle(1pt);
\draw (1.25,1.75) circle(1pt);
\draw (2,1) circle(1pt);
\draw (2.25,0.75) circle(1pt);
\draw (2.5,0.5) circle(1pt);

\draw [thin](0.25,3.75)-- ++(0.5,-1.25)-- ++(0.75,-1.25)-- ++(1.25,-1);

\end{tikzpicture}
\caption{The ideal $I=\langle xy^{15},x^3y^{10},x^6y^5,x^{11}y\rangle$ and its integral closure}\label{integral2dim}
\end{center}
\end{figure}

\end{Example}

\subsubsection{Hilbert series of monomial ideals.}
\spanishsubsubsection{Series de Hilbert de ideales monomiales}
Hilbert series will be a central topic for us, since its study is closely related to that of resolutions; this relation is object of active research \cite{P07}. Therefore, before going to the monomial case, we stop and give the basic general definitions.

The Hilbert function and Hilbert series of a graded module of the polynomial ring are defined as follows:
\begin{Definition}
Let $\Mc$ be a graded $R$-module, the Hilbert function $H_\Mc:\ZZ\rightarrow \ZZ$ maps each integer $z\in\ZZ$ to the dimension as $\kb$-vector space of the degree-$z$ piece of $\Mc$, i.e.
$$H_\Mc(z):=dim_\kb(\Mc_z)$$
The Hilbert series is defined as
$$HS_\Mc(t):=\sum_{z\in\ZZ}H_\Mc(z)\cdot t^z$$
Through these pages, when it is clear from the context, we will use the notation $H_\Mc$ alternatively for the Hilbert function and series.
\end{Definition}
One main property of the Hilbert function, due to Hilbert himself is that the Hilbert function becomes polynomial for large $z$, thus the information in it can be expressed in a simple way:
\begin{Theorem}\label{hilbertpoly}
If $\Mc$ is a finitely generated graded $R$-module then $H_\Mc(z)$ agrees for large $z$ with a polynomial of degree $\leq n$.
\end{Theorem}
\begin{Definition}
This polynomial is called the {\rm Hilbert polynomial} of $\Mc$ and is denoted by $HP(t)$.
\end{Definition}
Consider now an exact sequence of graded $R$-modules,
$$0\rightarrow\Mc_k\rightarrow\cdots\rightarrow\Mc_0\rightarrow 0$$
using the rank nullity theorem from linear algebra, we have that 
$$HS_{\Mc_k}=\sum_{i=0}^{k-1}(-1)^i HS_{\Mc_i}$$
This fact provides us a good tool for the computation of Hilbert series, and even a way to proof theorem \ref{hilbertpoly}, see \cite{E95} for example.

Another crucial result for computing Hilbert series is the fact that for any degree-preserving monomial ordering $\succeq$ defined in $R$ (see definition of monomial orderings in appendix \ref{apal}) and any ideal $I$ of $R$, we have that the Hilbert series of $I$ and $lt_\succeq(I)$ coincide; thus, we can reduce the computation of these Hilbert series to computations on monomial ideals.

In the case of multigraded modules, we can define a multigraded version of these objects. Multigraded (or $\NN^n$-graded) $R$-modules are those modules $\Mc$ such that $\Mc=\oplus_{\ab\in\NN^n}\Mc_{\ab}$ and $\xb^\ab\Mc_{\bb}\subseteq\Mc_{\ab+\bb}\quad\ab,\bb\in \NN^n$. In this case, one can define the {\it multigraded Hilbert series} as
$$H_\Mc(\xb)=\sum_{\ab\in\NN^n}dim_\kb(\Mc_\ab)\cdot \xb^\ab.$$ Monomial ideals are a particular case of {\it multigraded} $R$-modules.
\begin{Remark}
The multigraded Hilbert series is an element of the formal power series ring $\ZZ[[x_1\dots,x_n]]$ and in this ring we have that $\frac{1}{1-x_i}=1+x_i+x_i^2+\cdots$. Since the multigraded Hilbert series of $R$ is just the sum of all monomials in $R$, we have $H_R(\xb)=\prod_{i=1}^n\frac{1}{1-x_i}$. If we shift the grading of $R$ by $\ab$, i.e. consider the free module generated in multidegree $\ab$, which we will denote $R(-\ab)$, it is isomorphic to $\langle \xb^\ab\rangle$, then $H_{R(-\ab)}(\xb)=\xb^\ab\cdot H_R(\xb)=\frac{\xb^\ab}{1-x_i}$.

If $I$ is a monomial ideal, then the multigraded Hilbert series of the $R$-module $R/I$ is just the sum of all monomials not in $I$.
\end{Remark}

\begin{Remark}
Multigraded Hilbert series of monomial ideals and modules of the form $R/I$ for $I$ a monomial ideal, can be expressed as rational functions of the type $$H_{\Mc}(\xb)=\frac{\Kc_{\Mc}(\xb)}{(1-x_1)\cdots(1-x_n)}$$
the numerator $\Kc_{\Mc}(\xb)$ is known as the $\Kc$-polynomial of $\Mc$, see \cite{MS04}. Most of the time we will be interested in computing the $\Kc$-polynomial of some ideal or some module $R/I$. It is easy to see that if $I$ is a monomial ideal, $\Kc_{R/I}(\xb)=1-\Kc_{I}(\xb)$.
\end{Remark}

When computing the multigraded Hilbert series of a monomial ideal $I$, one can use the $lcm$-lattice in the following way: We need to `count' the monomials in $I$, and we have seen that the multigraded Hilbert series i.e. the sum of all the monomials in an ideal generated by a single monomial $\xb^\ab$ is of the form $H_{\langle \xb^\ab\rangle}(\xb)=\xb^\ab\cdot H_R(\xb)=\frac{\xb^\ab}{\prod_{i=1}^n(1-x_i)}$, so we need to sum these factors for each of the generators of our ideal. But doing so, we would add too many times some of the monomials in $I$, those belonging to the span of more than one generator. To avoid this, we can delete to our sum the corresponding factors of the pairwise intersections of the spans of the generators. These intersections are just ideals generated by the pairwise $lcm$'s. Again, deleting all these, would delete too many times the monomials belonging to more than one of such new ideals, and so on... Thus, we need to perform an inclusion-exclusion procedure, that leads us to the following formula for the multigraded Hilbert function of $I$ in terms of the $lcm$-lattice $L_I$:
$$H_{I}(\xb)=\frac{\sum_{i=1}^n (-1)^{i-1} \sum_{\xb^\ab\in L_{I,i}}\xb^\ab}{(1-x_1)\cdots(1-x_n)} $$
where $L_{I,i}$ is the set of the $lcm$'s of sets of $i$ generators of $I$. If we call $\Kc_{L_I}(\xb)$ to the numerator in this formula, before performing any cancellations among the factors, then we have that 
$$H_{R/I}(\xb)=\frac{1-\Kc_{L_I}(\xb)}{(1-x_1)\cdots(1-x_n)}.$$

\begin{Example}
Consider the ideal $I=\langle x^2,xy,y^2,yz\rangle$; its $lcm$-lattice $L_I$ is formed by: $L_{I,1}=\{x^2,xy,y^2,yz\}$, $L_{I,2}=\{x^2y,x^2y^2,x^2yz,xy^2,xyz,y^2z\}$, $L_{I,3}=\{x^2y^2,x^2yz,x^2y^2z,xy^2z\}$ and $L_{I,4}=\{x^2y^2z\}$ thus, we have 
$$\Kc_{L_I}(\xb)=x^2+xy+y^2+yz-(x^2y+x^2y^2+x^2yz+xy^2+xyz+y^2z)+(x^2y^2+x^2yz+x^2y^2z+xy^2z)-x^2y^2z$$
after cancellation of terms, we obtain 
$$\Kc_{I}(\xb)=x^2+xy+y^2+yz-x^2y-xy^2-xyz-y^2z+xy^2z$$
\end{Example}
The $lcm$-lattice produces as we have seen a closed formula for the $\Kc$-polynomial, but in general it contains much too many terms. Better formulas can be obtained from free resolutions of the ideal, in particular from the minimal free resolution and Betti numbers, as we see in the next paragraph.

\subsection{Homological invariants of monomial ideals}
\spanishsubsection{Invariantes homol\'ogicos de ideales monomiales}
We shall pay special attention to the homology of monomial ideals and its related invariants. The main object here are resolutions. Among them, a special object is the minimal resolution and the ranks of the modules forming it, i.e. the Betti numbers. The fact that monomial ideals are multigraded modules, allows us to obtain multigraded resolutions. This multigrading will be very useful from both the theoretical and computational viewpoints.

\subsubsection{Free resolutions, Betti numbers and Betti multidegrees.}\label{resolutions}
\spanishsubsubsection{Resoluciones libres, n\'umeros de Betti y multigrados de Betti}
Given an $R$-module $\Mc$, a {\it resolution} of it is a complex $\PP$ of free modules that is exact everywhere except at dimension $0$ in which we have $H_0(\PP)\simeq\Mc$ (see Appendix \ref{apal} for the necessary definitions). A well known result by Hilbert states that every finitely generated graded $R$-module has a free resolution of length at most $n$. Among all the possible resolutions of a graded module, a special role is played by the minimal one, which is unique up to isomorphisms. The ranks of the modules in the minimal resolution are called the Betti numbers of $\Mc$ and are the most important homological invariants of an $R$-module. 

As we have seen for the Hilbert series, the fact that monomial ideals are multigraded modules, leads to some extra properties. In the case of free resolutions, the main one is that free resolutions of multigraded modules, are also multigraded.  When taking into account the multigrading of a resolution $\PP: \cdots \Pc_i\stackrel{\delta_i}{\rightarrow}\Pc_{i-1}\rightarrow\cdots\rightarrow\Pc_1\stackrel{\delta_i}{\rightarrow}\Pc_{0}\rightarrow 0$, we will consider the fact that the free modules $\Pc_i$ are $\NN^n$-graded and each homomorphism $\delta_i$ is multidegree preserving. In the case of minimal resolutions of $\NN^n$-graded modules, we have that $\Pc_i=\bigoplus_{\ab \in \NN^n} R^{\beta_{i,\ab}}(-\ab)$ then, we define the $i$-th Betti number of the multigraded module $\Mc$ in multidegree $\ab$ is the invariant $\beta_{i,\ab}=\beta_{i,\ab}(\Mc)$; we will consider only those which are not zero. These are what we will call \emph{multigraded Betti numbers}.
 
Being $I$ a monomial ideal, we have that the $\beta_{i,\ab}(I)$ measure the number of minimal generators required in multidegree $\ab$ for the $i$-th syzygy module of $I$. We note that they can also be characterized by means of {\it Tor} since the $i$-th Betti number of an $\NN^n$-graded module $\Mc$ in degree $\ab$ equals the vector space dimension $dim_\kb Tor^R_{i,\ab}(\Mc,\kb)$. 

It will be useful to collect the multidegrees in which the multigraded Betti numbers are non-zero:
\begin{Definition}
Let $\Mc$ be a multigraded $R$-module. Let us denote by $\Bc(\Mc)$ the set of multidegrees in which the multigraded $Tor$ modules are nonzero: 
$$\Bc(\Mc)=\{\ab\in \NN^n\vert Tor_{i,\ab}(\Mc,\kb)\neq 0 \mbox{ for some }i\}$$
We call $\Bc$ the set of {\rm Betti multidegrees} of $I$. The corresponding sets 
$$\Bc_i(\Mc)=\{\ab\in \NN^n\vert Tor_{i,\ab}(\Mc,\kb)\neq 0\}$$
are called the  {\rm i-th Betti multidegrees} of $\Mc$. Equivalently, the Betti multidegrees are defined to be those multidegrees in which the multigraded Koszul homology of $\Mc$ does not vanish i.e. those in which the multigraded Betti numbers are different from zero. When a given multigraded Betti number is bigger than one, it will sometimes be useful to consider that multidegree as appearing repeated among the Betti multidegrees, then we speak of the \emph{collection} of Betti multidegrees to take these eventual repetitions into account.
\end{Definition}

\begin{Remark}
Given an ideal $I$ considered as an $R$-module, it can be described using a free resolution, and the same can be done with the module $R/I$. Both of them are very related and the information provided for their respective resolutions is equivalent, as is the information given by the Koszul homology of $I$ and $R/I$. The reason is the following: If we have a resolution $\PP$ of $I$ of the form 
$$\PP: \cdots \Pc_i\stackrel{\delta_i}{\rightarrow}\Pc_{i-1}\rightarrow\cdots\rightarrow\Pc_1\stackrel{\delta_i}{\rightarrow}\Pc_{0}\rightarrow 0$$
then we have a resolution of $R/I$
$$\PP: \cdots \Pc_i\stackrel{\delta_i}{\rightarrow}\Pc_{i-1}\rightarrow\cdots\rightarrow\Pc_1\stackrel{\delta_i}{\rightarrow}\Pc_{0}\rightarrow R\rightarrow 0$$
and then for all $i$, we have $\beta_i(I)=\beta_{i+1}(R/I)$ and also for the multigraded Betti numbers: $\beta_{i,\ab}(I)=\beta_{i+1,\ab}(R/I)\,\forall\ab\in\NN^n$. Moreover, we have that
$$H_i(\KK(I))\simeq H_{i+1}(\KK(R/I))$$
and the isomorphism preserves multidegree.

Thus, from a theoretical point of view it is almost equivalent to work with $I$ or with $R/I$. Sometimes it is more convenient to work with resolutions of $R/I$, in particular in the theoretical settings, because some of the concepts are more clear and easy to handle, and it is the most frequent framework in the literature. On the other hand, sometimes is more convenient to work with $I$ and its Koszul complex, because then we can handle the monomials directly, not their equivalence classes modulo $I$, and the generating sets of the Koszul modules can be expressed in a simpler way. In order to recover all the information about $R/I$ from the information over $I$ or vice versa, we need explicit isomorphisms between the Koszul homology modules of $R$ and $R/I$. This is provided by the Spencer differential, in particular observe that if $\mu$ is a generator of the $i$-th Koszul homology group of $I$ then $\delta(\mu)$ is a generator of the $i+1$-st Spencer cohomology group of $R/I$ (see section \ref{spencer}), i.e. the corresponding multidegrees are identical. Therefore, we will speak alternatively of $I$ or $R/I$. Most of the time, we will use $R/I$ in the theoretical part when speaking about resolutions, and we will use $I$ when giving explicit formulas or in Koszul homology computations.
\end{Remark}

Multigraded resolutions are very closely related to multigraded Hilbert functions, as the following proposition (which is an easy consequence of the rank-nullity theorem) shows:

\begin{Proposition}
Let $\PP$ be a multigraded resolution of the monomial ideal $I$ and let $m_{i,\alpha}$ the rank of the multidegree $\alpha$ component of the module $\Pc_i$ in the resolution. Then the $\Kc$-polynomial of $I$ is given by
$$\Kc_I(\xb)=\sum_{i>0,\alpha\in\NN} (-1)^{i-1} m_{i,\alpha} x^\alpha$$
in case $\PP$ is the minimal resolution, we have
$$\Kc_I(\xb)=\sum_{i>0,\alpha\in\NN} (-1)^{i-1} \beta_{i,\alpha} x^\alpha$$
We will denote $\Kc_\PP(\xb)$ the formula $\sum_{i>0,\alpha\in\NN} (-1)^{i-1} m_{i,\alpha} x^\alpha$ before performing any cancellations, and call it the \emph{$\Kc$-polynomial of $I$ in terms of $\PP$}.
\end{Proposition}
In terms of the Betti multidegrees, we have that $\Kc_I$ is the alternating sum of the  $\Bc_i(I)$ in which every multidegree $\ab$ is counted $\beta_{i,\ab}(I)$ times. In order to avoid too many redundant terms and cancellations in the form we obtain for $\Kc_I$ from resolutions, we will prefer smaller resolutions, and the less redundant one is the minimal resolution:

\begin{Proposition}
Let $I$ be a monomial ideal and let $\PP$ be a free resolution of $I$; let $\MM$ be the minimal free resolution of $I$. Any terms that can be cancelled in $\Kc_\MM(\xb)$ can also be cancelled in $\Kc_\PP(\xb)$.
\end{Proposition}

\noindent {\bf Proof: } If a term of multidegree $\ab$ can be cancelled in $\Kc_\MM(\xb)$ then we have that $\ab$ belongs to $\Bc_i(I)$ and $\Bc_j(I)$ with $i\neq j$ each of them having different parity. In other words, the multidegree $\ab$ pieces of $\Mc_i$ and $\Mc_j$, namely $\Mc_{i,\ab}$ and $\Mc_{j,\ab}$ have positive rank. Since by definition of minimal resolution and since the free resolutions of $I$ are multigraded, we have that the rank of $\Pc_{k,\bb}$ is bigger or equal to the rank of $\Mc_{k,\bb}$ for every $k$ and $\bb$ then $\Pc_{i,\ab}$ and $\Pc_{j,\ab}$ have positive rank and hence the same cancellation is possible in $\Kc_\PP(\xb)$ $\square$

\subsubsection{Examples of monomial resolutions}\label{examples-resolutions}
\spanishsubsubsection{Ejemplos de resoluciones monomiales}
The efficient computation of minimal resolutions is a difficult task, even in the monomial case \cite{BPS98,GPW99}, for which no general closed formula has been given. Research in this field has led to define on one hand non-minimal resolutions of general monomial ideals \cite{T60,L98}, and in the other hand, to study special classes of monomial ideals for which their minimal resolutions can be explicitly given or at least easily computed \cite{EK90,MS04}.

Probably the best known resolution of a monomial ideal is Taylor's resolution \cite{T60}, which has a very simple explicit description. In general, Taylor resolution is highly nonminimal. A subresolution of Taylor's is Lyubeznik resolution \cite{L98}, still nonminimal. There are a number of other nonminimal resolutions in the literature cited in the references (for instance, a resolution of minimal length for some types of ideals is described in \cite{S02b}), here we shall describe only cellular resolutions \cite{MS04}, and in particular the Scarf resolution, which is  minimal for a special family of monomial ideals, called generic ideals.
More recently {\it discrete Morse theory} has also been applied to the study of resolutions of monomial ideals \cite{B02,OW07}
\medskip

\noindent \emph{Taylor and Lyubeznik resolutions:}

Let $I$ be a monomial ideal and $\{m_1, \dots, m_r\}$ a generating set of $I$.
For any subset $J=\{j_1,\dots,j_s\}\subseteq \{1,\dots,r\}$, let us denote $m_J=lcm(m_{j_1},\dots m_{j_s})$,and $J^i=\{j_1,\dots,\widehat{j_i},\dots,j_s\}$. We can construct a resolution of $R/I$ in the following way: Let $T_s,\,s\geq 0$ be a free $R$-module generated as a vector space by the set $\{u_J\, s.th.\, |J|=s\}$ and consider the $R$-linear differential
$$d(u_J)=\sum_{i\in J}(-1)^{i-1} \frac{m_J}{m_{J^i}}u_{J^i}$$
it is easy to verify that $d^2=0$. Moreover, this complex is acyclic and it is a resolution of $R/I$. This resolution is due to D.Taylor \cite{T60} and will be denoted by $\TT=(T_i,d_i)$. The length of Taylor's resolution is given by the number of elements in the given generating set of the ideal (normally, we will assume we have a minimal generating set for the ideal), which we denote by $r$. The rank of the $i$-th free module $T_i$ is ${r}\choose{i}$, thus the sum of all these ranks, i.e. $Size(\TT)$ is $2^r$.
\begin{Remark}
Generally, Taylor's resolution is far from minimal, although \cite{B02} gives the following condition for a monomial ideal to be minimally resolved by Taylor's resolution: Let $I$ be a monomial ideal, then the Taylor resolution of $I$ is minimal if and only if for all $J$ in the support of $\TT(I)$ and all $m\in J$ we have $m\nmid lcm(J\smallsetminus\{m\})$.
Other characterizations of monomial ideals for which the Taylor resolution is minimal are given by Fr\"oberg \cite{F78} and Herzog et al. \cite{HHMT06}
\end{Remark}

A subresolution $\LL$ of $\TT$ was given in \cite{L98} and it is known as the {\it Lyubeznik resolution}. It is defined as follows: For a given subset $J\subseteq \{1\dots r\}$ and an integer $1\leq s\leq r$, let $J_{> s}=\{j\in J\vert j>s\}$; then $\LL$ is generated by those basis elements $u_J$ such that for all $1\leq s\leq r$ one has that $m_s$ does not divide $m_{J_{> s}}$. It is clear that, unlike Taylor's, Lyubeznik resolution depends on the ordering in which the generators of the ideal are given.

\begin{Example}\label{taylor-lyu}
Let us consider the following monomial ideal in three variables: $I=\langle x^2y, xy^3, xz, yz\rangle$, the Taylor resolution $\TT(I)$ of $I$ has length $4$, size $16$ and the differentials are given by
\begin{small}
$$
\begin{array}{cc}
{d_1={ \left ( \begin{array}{cccc}
x^2y&xy^3&xz&yz\\
\end{array}\right )}}&
{ d_2={ \left ( \begin{array}{cccccc}
y^2&z&z&0&0&0\\
-x&0&0&z&z&0\\
0&-xy&0&-y^3&0&y\\
0&0&-x^2&0&-xy^2&-x\\
\end{array}\right )}}
\\
{
d_3={ \left ( \begin{array}{cccc}
-z&-z&0&0\\
y^2&0&-1&0\\
0&y^2&1&0\\
-x&0&0&-1\\
0&-x&0&1\\
0&0&-x&y^2\\
\end{array}\right )}
}&
{
d_4={ \left ( \begin{array}{c}
1\\-1\\y^2\\-x\\
\end{array}\right )}
}\\
\end{array}
$$
\end{small}
Lyubeznik's resolution is in this case equal to Taylor's if we keep the ordering $m_1=x^2y$, $m_2=xy^3$, $m_3=xz$, $m_4=yz$. On the contrary if we change the order to $m_1=xz$, $m_2=yz$, $m_3=x^2y$, $m_4=xy^3$, then $\LL$ is generated by $u_1,\,u_2,\,u_3,\,u_4,\,u_{12},\,u_{13},\,u_{14},\,u_{34}$ and $u_{134}$ i.e. the size of $\LL$ is $10$. In this case, $\LL$ is minimal. The differentials are given by
\begin{small}
$$
\begin{array}{ccc}
{\delta_1={ \left ( \begin{array}{cccc}
xz&yz&x^2y&xy^3\\
\end{array}\right )}}&
{ \delta_2={ \left ( \begin{array}{cccc}
y^2&z&0&0\\
-x&0&z&0\\
0&-xy&-y^3&y\\
0&0&0&-x\\
\end{array}\right )}}
&
{
\delta_3={ \left ( \begin{array}{c}
-z\\
y^2\\
-x\\
0\\
\end{array}\right )}
}\\
\end{array}
$$
\end{small}

\end{Example}

\begin{Remark}
In \cite{S02c} Taylor and Lyubeznik resolutions are described in relation with the Buchberger algorithm for the construction of Gr\"obner bases. Taylor resolution is shown to be obtained by repeated application
of the Schreyer theorem from the theory of Gr\"obner bases, and
the Lyubeznik resolution is shown to be a consequence of Buchberger's chain criterion.
\end{Remark}

The Taylor resolution is strongly related to the lcm-lattice of a monomial ideal (see \cite{GPW99}). If we consider Taylor's resolution of $R/I$, we see that all the multidegrees of the generators lie in $L_I$. Then, applying a minimization process to this resolution, to obtain a minimal one (as described for example in \cite{CLO98} or the reduction algorithm described below) we have that no new multidegrees will appear, although some can eventually disappear. Thus, this is a simple way to show that $\beta_{i,\ab}(I)=0$ if $\ab\notin L_I$.

\medskip

\noindent \emph{Cellular resolutions:}

In the thesis \cite{M00} and their book \cite{MS04}, E. Miller and B. Sturmfels define and give some interesting properties and families of what they call {\it cellular resolutions}. These are a type of geometrical resolutions that can be associated to monomial ideals. Here we will only repeat some of the basic definitions, for deeper details we point the interested reader to the original source.
A useful tool to handle monomial resolution are monomial matrices:
\begin{Definition}
A {\rm monomial matrix} is an array of scalar entries $\lambda_{qp}$ whose columns are labeled by {\rm source degrees} $\ab_p$, whose rows are labeled by {\rm target degrees} $\ab_q$, and whose entry $\lambda_{qp}\in\kb$ is zero unless $\ab_p\geq\ab_q$.
\end{Definition}
In the usual unlabeled notation, the entry $\lambda_{qp}$ is replaced by $\lambda_{qp}\cdot\xb^{\ab_p-\ab_q}$.
\begin{Definition}
A {\rm polyhedral cell complex} $X$ is a finite collection of convex polytopes (in a real vector space $\RR^m$) called {\rm faces} of $X$, satisfying two properties:
\begin{itemize}
\item If $\Pc$ is a polytope in $X$ and $F$ is a face of $\Pc$, then $F$ is in $X$
\item If $\Pc$ and $\Qc$ are in $X$, then $\Pc\cap\Qc$ is a face of both $\Pc$ and $\Qc$
\end{itemize}

\end{Definition}

\begin{Definition}
$X$ is a {\rm labeled cell complex} if its $r$ vertices have {\rm labels} that are vectors $\ab_1,\dots,\ab_r$ in $\NN^n$ and the {\rm label} on an arbitrary face $F$ of $X$ is the exponent $\ab_F$ on the least common multiple $lcm(\xb^{\ab_i}\vert i\in F)$ of the monomial labels $\xb^{\ab_i}$ on vertices in $F$.
\end{Definition}

\begin{Definition}
Let $X$ be a labeled cell complex. The {\rm cellular monomial matrix} supported on $X$ uses the reduced chain complex of $X$ for scalar entries, with $\emptyset$ in homological degree $0$. Row and column labels are those on the corresponding faces of $X$. The {\rm cellular free complex $\Fc_X$ supported on $X$} is the complex of $\NN^n$-graded free $R$-modules (with basis) represented by the cellular monomial matrix supported on $X$. The free complex $\Fc_X$ is a {\rm cellular resolution} if it is acyclic (homology only in degree $0$).
\end{Definition}

\begin{Remark}
Observe that Taylor resolution is an example of cellular resolution. In this case, $X$ is the full $(r-1)$-simplex with the $r$ vertices labeled by the exponents of the generators of $I$.
\end{Remark}

Several examples of cellular resolutions can be found in \cite{MS04}, like those for permutohedron ideals, tree ideals, etc.  Of particular importance is the {\it hull resolution}, which is a canonical construction of a resolution of length less than or equal $n$, which is generally nonminimal. An important family of monomial ideals for which the hull resolution is not only minimal but also has a simple description is that of {\it generic} monomial ideals:

\begin{Definition}
A monomial ideal $\langle m_1,\dots,m_r\rangle$ is {\rm generic} if whenever two distinct minimal generators $m_i$ and $m_j$ have the same positive (nonzero) degree in some variable, a third generator $m_k$ strictly divides their least common multiple $lcm(m_i,m_j)$
\end{Definition}
\begin{Definition}
let $I$ be a monomial ideal with minimal generating set $\{m_1,\dots,m_r\}$. The {\rm Scarf complex} $\Delta_I$ is the collection of all subsets of $\{m_1,\dots,m_r\}$ whose least common multiple is unique:
$$\Delta_I=\{\sigma\subseteq\{1,\dots,r\}\vert m_\sigma=m_\tau\Rightarrow\sigma=\tau\}$$
\end{Definition}
The Scarf complex is a simplicial complex of dimension at most $n-1$ and the chain complex supported on it, the {\it algebraic} Scarf complex $\Fc_{\Delta_I}$, is included as a subcomplex in any resolution of $I$. For generic monomial ideals, the Scarf complex provides the minimal resolution, thus, the minimal resolution of these ideals is simplicial:
\begin{Theorem}
If $I$ is a monomial ideal, then its Scarf complex $\Delta_I$ is a subcomplex of the hull complex $hull(I)$. If $I$ is generic then $\Delta_I=hull(I)$ so its algebraic Scarf complex $\Fc_{\Delta_I}$ minimally resolves $I$.
\end{Theorem}

\begin{Example}
The ideal in example \ref{taylor-lyu} is not generic. Its Scarf complex is formed by the sets $\{ \{1\},\{2\},\{3\},\{4\},\{12\},\{34\}\}$ and its algebraic counterpart is not a resolution of $I$.
\end{Example}

\begin{Remark}
In the non-generic case, the Scarf complex can be perturbed to obtain another complex with supports a resolution of $I$ which is non-minimal in general. A simple description of this perturbation procedure can be seen for example in \cite{GW04}.
\end{Remark}
\subsubsection{ Minimizing resolutions}\label{minimalizing}
\spanishsubsubsection{Minimalizaci\'on de resoluciones}
To close this section on free resolutions we refer to the fact that the minimal resolution of an ideal $I$ is contained as a subcomplex in every resolution of $I$. Given some resolution, there are explicit algorithmic ways to reduce it until obtaining the minimal resolution. Inside Theorem $3.15$ in \cite{CLO98} a procedure to simplify a given graded resolution is described. Here we give a simple equivalent description of this process, based on the {\it Chain Complex Reduction Algorithm} (CCR). This algorithm is applied in general to chain complexes and has a geometric interpretation in the case of simplicial, cubical or CW-complexes. A description of it can be found in \cite{KMM04}. We give here first the basic definitions and theorems in which the algorithm is based. After this we apply these results to the case of resolutions of monomial ideals.

\begin{Definition}
Given a free chain complex $\{C,\delta\}$, we say that two chains $c\in C_i$ and $c'\in C_{i+1}$ are a {\rm reduction pair} if $\langle \delta(c),c'\rangle$ has an inverse, where $\langle\cdot,\cdot\rangle$ denotes their scalar product $\langle c_1,c_2\rangle:=\sum_{i=1}^m\alpha_i,\beta_i$ if $c_1=\sum_{i=1}^m \alpha_i\bf{e}_i$ and $c_2=\sum_{i=1}^m \beta_i\bf{e}_i$
\end{Definition}

Assume we have a generator $c$ of the module $C_{i-1}$ and a generator $c'\in C_{i}$ that form a reduction pair in some chain complex $\{C,\delta\}$.  Using this information, we want to define a new chain complex $\{\bar C,\bar\delta\}$ whose basis is a subbasis of $\{C,\delta\}$.
Let $W_k$ denote the set of generators of $C_k$. Say $W_{i}=\{b_1,\dots,b_{d_{i}},b\}$ and $W_{i-1}=\{a_1,\dots,a_{d_{i-1}},a\}$, being $(a,b)$ a reduction pair. Define $\bar{W_{i}}=\{b_1,\dots,b_{d_i}\}$ and $\bar{W_i}=\{a_1,\dots,a_{d_{i-1}}\}$ and $\bar{W_j}=W_j$ for all other $j$. The new differential $\bar\delta$ is defined with the following formulas for every chain in $C$:
$$\bar{\delta_j}(c)=\left\{ \begin{array}{ll}
\delta(c)-\frac{\langle\delta(c),a\rangle}{\langle\delta(b),a\rangle}\delta(b)&\mbox{if }j=i\\
\delta(c)-{\langle\delta(c),b\rangle}b&\mbox{if }j=i+1\\
\delta(c)&\mbox{otherwise}\\
\end{array}
\right.
$$
Observe that in this formula the condition of $\langle \delta(b),a\rangle$ having an inverse is required (which is relevant in case our coefficients do not lie on a field). With this formula, we have that
\begin{Theorem}
$H_*(C)\simeq H_*(\bar{C})$ and the isomorphism is induced in homology by the projection $\bar\pi$, which can be expressed by the following explicit formula:
$$\bar{\pi_j}(c)=\left\{ \begin{array}{ll}
c-\frac{\langle c,a\rangle}{\langle\delta(b),a\rangle}\delta(b)&\mbox{if }j=i-1\\
c-{\langle\delta(c),a\rangle}b&\mbox{if }j=i\\
c&\mbox{otherwise}\\
\end{array}
\right.
$$
\end{Theorem}

This projection is the mathematical expression of the reduction step performed in the original complex via the reduction pair $(a,b)$. The idea of the CCR algorithm is to perform these reduction steps while possible, resulting in each step a smaller complex with the same homology as the original one. The size of the corresponding complex is reduced by two on each reduction step.

In our case, given a monomial ideal $I\subseteq \poly \kb x n$, we describe a free resolution by giving a sequence of modules with their generators, and a sequence of matrices that describe the differentials. Recall that given a free resolution of a monomial ideal, it is minimal if and only if it has no nonzero constant polynomials in the entries of the matrices describing the differentials in the complex. Starting with some non-minimal resolution $\PP$ of $I$, whenever we have a nonzero constant entry $c$ in one of the matrices in the complex, say in the column that expresses the differential applied to generator $e'\in\Pc_{i+1}$ at the row corresponding to generator $e\in\Pc_i$ then clearly $\langle \delta(e),e'\rangle=c$ thus, for every nonzero scalar in a differential, we have a reduction pair and we can apply the standard result described above and we obtain a smaller complex that is still a resolution. It is clear that no new constant terms will be added in this step, thus, after applying the reduction step we have seen for chain complexes, we have a smaller chain complex which the same homology as $\PP$ (i.e. it is still a resolution of $I$) and has both smaller size and less constant terms in the differentials. Thus performing this reduction successively, we will reach a resolution with no constant terms in the differentials, i.e. the minimal resolution.

\section{Topological and homological techniques for monomial ideals}\label{hom_tools}
\spanishsection{T\'ecnicas topol\'ogicas y homol\'ogicas para ideales de monomios}
In this section we give a collection of topological and homological techniques applicable in the study of monomial ideals. Since we are focused on the multigraded homological description of monomial ideals, which can be summarized by the multigraded $Tor$ modules the multigraded minimal resolution and the Betti numbers, the techniques we will explore have this homological flavour. On the other hand, the combinatorial nature of monomial ideals suggests the use of combinatorial techniques. Thus, the first objects that will appear in our collection are simplicial complexes, a distinguished area where homology meets combinatorics. The connection between simplicial complexes and Betti numbers of monomial ideals has been applied among others by Hochster \cite{H77}, Stanley \cite{S96}, Burns and Herzog \cite{BH97}, Eagon and Reiner \cite{ER98}, Bayer \cite{B96} or Miller and Sturmfels \cite{MS04}. The simplicial Koszul complex, introduced by the latter authors, will be our first object, after which we introduce Stanley-Reisner ideals and relate them to Koszul simplicial complexes. Later on we shall move to exact sequences and will explore the Mayer-Vietoris sequence and its application to monomial ideals, as well as the algebraic mapping cone applied to resolutions of monomial ideals. All the objects in this section share their algorithmic capabilities. Our interest in actual computations is clearly fulfilled by the techniques we present here. Some of the techniques in the following pages were presented in \cite{S06b}.

\begin{Remark}
There are a number of simplicial complexes that have been used to compute multigraded Betti numbers of monomial ideals, we focus here on the \emph{simplicial Koszul complex}, other simplicial constructions can be seen in \cite{BH97,BS98,GPW99,P02}.
\end{Remark}

\subsection{The simplicial Koszul complex}\label{simpl-koszul-complex}
\spanishsubsection{El complejo de Koszul simplicial}
Simplicial complexes are well known objects that play an essential role in combinatorial algebraic topology. All the necessary definitions are given in Appendix \ref{aptop}. One way to compute the Betti numbers or Koszul homology of $I$ at a given multidegree $\ab$ is to associate a particular simplicial complex to the ideal and the multidegree and express the Koszul homology of $I$ at $\ab$ in terms of this simplicial complex. This idea is presented in \cite{B96}, \cite{P02} or \cite{MS04} and also in \cite{BH97} and in some sense gives a simplicial meaning to the contribution of each multigraded piece of the homological structure of a monomial ideal.

\subsubsection*{Upper and lower simplicial Koszul complexes}
\spanishsubsubsection{Complejo de Koszul simplicial superior e inferior}
\begin{Definition}
Let $I$ be a monomial ideal, $\ab=(a_1,\dots,a_n)\in \NN^n$ and $x^{\ab}\in I$. Let $l=|\mbox{support}(\ab)|$ i. e.  $l$ is the number of nonzero components of ${\ab}$. We associate to $\ab$ the $l$-simplex $\Delta_{\ab}$, where the vertices are labeled by the variables in the support of $\ab$. Now, we build the subcomplex of $\Delta_{\ab}$ given by 
$$\Delta_{\ab}^I:=\{\tau \in \Delta_{\ab}\vert x^{\ab}/x^{ \tau}\in I\}$$ where $x^\tau$ is a squarefree monomial with exponents given by the variables defining the face $\tau$. This is what is called the {\rm upper Koszul simplicial complex}, in few words, it is described as
$$\Delta_{\ab}^I=\{\mbox{squarefree vectors }\tau \vert \xb^{\ab-\tau}\in I\}$$ 
\end{Definition}

With this definition we have the following result that relates the simplicial homology of the upper Koszul simplicial complex to the multidegree $\ab$ Betti numbers of $I$ and hence with the multidegree $\ab$ piece of the Koszul homology of $I$, see \cite{B96,MS04}.
\begin{Theorem}\label{simp-kos-hom}
$$H_i(\KK(I)_\ab)\simeq \widetilde {H}_{i-1}(\Delta_{\ab}^I)\quad \forall i$$
\end{Theorem}

\noindent {\bf Proof:} The proof is based on the fact that the corresponding chain complexes are equal under the isomorphism that maps an element in $\KK_j(I)$ to a face of $\Delta_{\ab}^I$ . Since this isomorphism between the chain complexes is compatible with the differentials we have the isomorphism in homology. $\square$
\begin {Remark}
The correspondence of the chains of the chain complex associated to $\Delta_\ab^I$ with chains of $\KK(I)_\ab$ that makes explicit this isomorphism  is just the $\kb$-linear extension of the correspondence that associates the face $\tau \in \Delta_{\ab}^I$ with $x^{\ab}/x^\tau\otimes x^{(\tau)}$ where, if $\tau=a_{\tau_1},\dots a_{\tau_k}$  the exterior part (i.e. the right-hand side with respect to the tensor symbol) of the latter has to be understood as $x^{(\tau)}=x^{a_{\tau_1}}\wedge\dots\wedge x^{ a_{\tau_k}}$.
\end{Remark}

Dual to the upper Koszul simplicial complex, we can define the {\it lower Koszul simplicial complex} in the following way, \cite{B96,MS04}:
\begin{Definition}\label{lower-Koszul}
For $\ab\in \NN^n$ define $\ab'$ by subtracting 1 from each nonzero coordinate of $\ab$. The {\rm lower Koszul simplicial complex} of $I$ at $\ab$ is given by
$$\Delta^\ab_{I}=\{\mbox{squarefree vectors }\tau \preceq \ab\vert \xb^{\ab'+\tau}\notin I\}$$ 
\end{Definition}

\subsection{Simplicial computation of Koszul homology}\label{simplicial_Koszul_comp}
\spanishsubsection{C\'alculo simplicial de la homolog\'ia de Koszul}
Koszul simplicial complexes defined at each $\ab\in\NN^n$ together with the fact that $\beta_\ab(I)=0$ if $\ab\notin L_I$ (recall $L_I$ is the lcm-lattice of $I$) provide us with a method of computing the Koszul homology of $I$ in a finite number of steps. Namely, for each multidegree in the lcm-lattice of $I$ (which is a finite set), we construct the corresponding upper or lower Koszul simplicial complex, and compute the corresponding singular homology. Thus, Koszul homology is computable for any monomial ideal. This way of computing Koszul homology can be expressed as an algorithm, it is shown in table \ref{alg_naive_Koszul}.

\begin{table}[!htb]
\centering
\begin{tabular}{|p{13cm}|}
\hline
$$\begin{array}{ll} Algorithm:\, \mbox{Naive Koszul Homology of a Monomial ideal }$I$\\
\hline \hline\\
\ Input: \, \mbox{Minimal generating set of } I=\langle m_1,\dots m_r\rangle\\
\ Output: \mbox{Generators of the Koszul homology of }I\\
\\
\hline\\
\ 1 \quad \mbox{Compute the lcm-lattice of }I\\
\ 2 \quad {\bf for each }\, \ab\in L_I \,{\bf do}\\
\ 3 \quad \quad \mbox{Build the simplicial complex associated to }\ab,\, \Delta_{\ab}^I\\
\ 4 \quad \quad \mbox{Compute the reduced homology of }\Delta_{\ab}^I\\
\ 5 \quad {\bf end for each}\\
\ 6 \quad {\bf return} H(K(I))_*\\
\end{array}
$$
\\
\hline
\end{tabular}
\caption{Algorithm {\bf Naive Koszul Homology}}\label{alg_naive_Koszul}
\end{table}

The problem of this algorithm is its high complexity. The number of multidegrees in $L_I$ is bounded above by $2^r$ and in some cases this bound will be reached, so we are dealing with an exponentially growing set of points in the lattice. On the other hand, at each of these points, we want to explicitly compute the reduced homology of a simplicial complex in $n$ vertices, and this computation depends exponentially on the number of variables $n$. So the algorithm complexity is double exponential in $n$ and $r$. Improvements of this algorithm must then focus in both reducing the number of relevant multidegrees for the computation and stressing the properties of monomial ideals in order to simplify the simplicial computations required at each multidegree. Both ways of improvements will be developed in the next paragraphs.

We have just seen that the problem of computing the Koszul homology of $I$ at a given multidegree is equivalent to the computation of the reduced homology of a subcomplex of an $n$-simplex. Since we want to explicitly compute these reduced homology modules, our problem reduces to compute quotients between the kernels and image subspaces of some matrices with coefficients in $\kb$. The most usual way to do this is via linear algebra algorithms keeping track of the generators of the corresponding modules, in particular Smith normal form algorithms. As our coefficients lie on a field, we can substitute Smith normal form computations with Gaussian elimination, which is more efficient in our case, but we still have high complexity. Being this a hard problem in general, we can take advantage of several facts that will drastically improve our computations: 

\subsubsection*{Computing homology dimension by dimension}
\spanishsubsubsection{C\'alculo de la homolog\'ia dimensi\'on a dimensi\'on}
Instead of computing Koszul homology by exploring the whole lcm-lattice at once, we can proceed dimension by dimension. For the zero dimensional Koszul homology, we know that the generators are just the minimal generators of the ideal. In each subsequent dimension, we start with a candidate set of multidegrees in which we will eventually have non-vanishing Koszul homology. The actual candidate sets used are a subset of $L_I$: we use as candidate set of multidegrees for dimension $i$, the least common multiples of pairs of multidegrees with non-vanishing homology in dimension $i-1$. The reason to do this is that in a syzygy resolution of $I$, if $I$ is monomial, and we have a minimal generating set of the $(i-1)$-th syzygy module, then we can construct a generating set of the next syzygy module such that the multidegree of the new generators are least common multiples of pairs of the former ones. The multidegrees of the minimal generating set of syzygies for the new module will be of course a subset of these. And finally, the multidegrees of the elements in a minimal generating set of the $i$-th syzygy of $I$ module coincide with the multidegrees in which the $i$-th Koszul homology of $I$ does not vanish, for in both cases they are the multigraded Betti numbers. As in each step we obtain the multidegrees of the non-zero Betti numbers, we can use this information for building the set of candidate multidegrees for the next step.
This procedure, on one hand reduces the number of multidegrees to be considered, for we don't need to explore all the lcm-lattice, and on the other hand, it also reduces the complexity of the simplicial complexes for which we have to compute reduced homology:

For each complex we only need to compute homology at a given dimension. For this we will only need a short part of the chain complex $(C_\bullet,\delta)$ associated to $\Delta_\ab^I$, namely, for computing the reduced homology in dimension $i$, we need the generators of $C_{i+1}$, $C_i$ and $C_{i-1}$ and the matrices defining $\delta_{i+1}$ and $\delta_i$.

Now, for our purposes it is enough to consider a subcomplex of $\Delta_\ab^I$: given $i$, let us denote by $\Delta_{\ab,i}^I$ the subcomplex of $\Delta_\ab^I$ whose facets (i.e. the maximal faces) are the $i+1$ and $i$-faces of $\Delta_\ab^I$, note that a simplicial complex can be completely described by giving its facets. Then it is easy to see that $\widetilde {H_i}(\Delta_{\ab}^I)=\widetilde {H_i}(\Delta_{\ab,i}^I)$. Using these smaller complexes, the size of the needed matrices becomes smaller and we can reduce the size of our computations.
\begin{Example}
Consider for example the ideal $I=\langle xyzts,xytu,xyzu,xysu,tsu\rangle$ in the ring $R=\RR[x,y,z,t,s,u]$ . If we take $\ab=(1,1,1,1,1)$ then figure \ref{complexes} shows $\Delta_\ab^I$ and $\Delta_{\ab,2}^I$. We can see in it that $\Delta_{\ab,2}^I$ is contractible, and thus $H_3(\KK(I))_\ab=0$.

\begin{figure}[h]
\begin{center}
$\begin{array}{l@{\hspace{1in}}r}

\begin{tikzpicture}[scale=1]
\filldraw [fill=gray!50!white,draw=black,thick](0,0)-- ++(-2.5,2.5)-- ++(1.5,1.5)--cycle;
\draw [thick](0,0)-- ++(2.5,2.5)-- ++(-1.5,1.5)--cycle;
\fill (0,0) circle(2pt) node[below] {$z$};
\fill (-2.5,2.5) circle(2pt) node[left] {$y$};
\fill (-1,4) circle(2pt) node[above] {$x$};
\fill (2.5,2.5) circle(2pt) node[right] {$s$};
\fill (1,4) circle(2pt) node[above] {$t$};
\fill (2,1) circle(2pt) node[left] {$u$};
\end{tikzpicture} 

&

\begin{tikzpicture}[scale=1]
\filldraw [fill=gray!50!white,draw=black,thick](0,0)-- ++(-2.5,2.5)-- ++(1.5,1.5)--cycle;
\fill (0,0) circle(2pt) node[below] {$z$};
\fill (-2.5,2.5) circle(2pt) node[left] {$y$};
\fill (-1,4) circle(2pt) node[above] {$x$};
\end{tikzpicture}

 \\ [0.4cm]
\Delta_{\ab}^I & \Delta_{\ab,2}^I
\end{array}$
\end{center}
\caption{The simplicial complexes $\Delta_{\ab}^I$ and $\Delta_{\ab,2}^I$ for $\ab=(1,1,1,1,1,1)$. }
\label{complexes}
\end{figure}
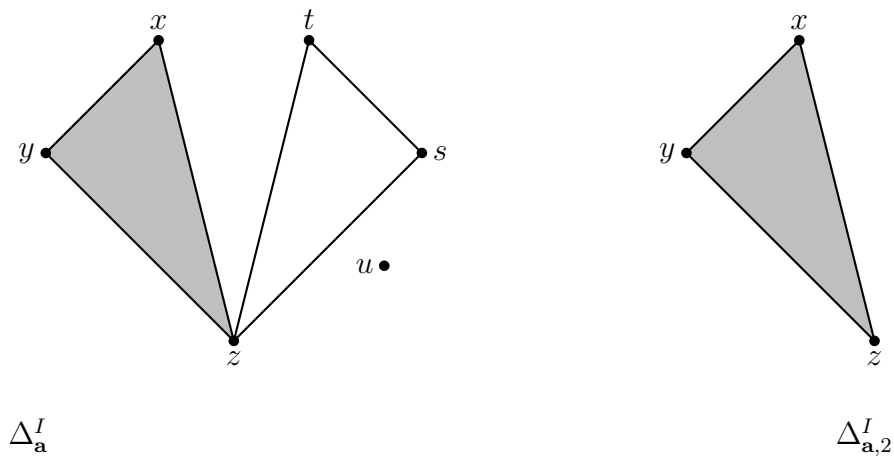

\end{Example}

\subsubsection*{Improvements based on simplicial homology}
\spanishsubsubsection{Mejoras basadas en homolog\'ia simplicial}
An impressive amount of work has been done on the computation of the homology of simplicial complexes. Two main sources of improvement have been developed. On one side, one can improve the techniques for the computations needed on the matrices defining the differentials. Since these matrices have a big degree of sparsity, several techniques for dealing with Smith normal form algorithms and Gaussian elimination in the case of sparse matrices have resulted in fast algorithms for performing these computations \cite{DER97,DSV01}. The second source comes from considerations on the simplicial complexes themselves. Here we will just make two simple considerations on the simplicial complexes we deal with, which can be used to improve our computations.

In first place, sometimes it is easy to detect whether some homology groups vanish or not, without actually compute them. In particular it is easy to detect if $\Delta_{\ab}^I$ is the full $l$-simplex (where $l=\vert supp(\ab)\vert$), and thus contractible: we just check if $x^{\ab-(0,\dots, 1_{i_1},\dots,1_{i_l},\dots,0)}$ is in the ideal $I$, this is performed just by monomial divisions. If it is the case, we know that $H_*(K(I))$ vanishes at $\ab$. On the other hand we can also easily check whether all the $(i+1)$-faces of the $l$-simplex are present in $\Delta_{\ab}^I$ or not. If it is the case, then $\Delta_{\ab,i}^I$ is contractible, and thus the Koszul homology of $I$ in dimesion $i$ vanishes at $\ab$; this is also easy to check.

To illustrate these facts we see in table \ref{contractible_corners} a comparison between the number of elements of the $lcm$-lattice $L_I$, the size of the Taylor resolution, the number of non-contractible Koszul simplicial complexes  and the size of the minimal resolution for some random examples. In the table, $n$ denotes the number of variables, $g$ the number of generators, $md$ the minimal degree of a generator, and $Md$ the maximal degree of a generator. We can see how the simple considerations just exposed produce a severe reduction in the number of complexes to be considered.

\begin {table}[h]
\begin {center}
    \begin{tabular}{|c|c|c|l|l|l|l|}
    \hline
    \small{{\bf n}}& \small{{\bf g}}& \small{{\bf Md/md}} & \small{{\bf lcm-lattice}}& \small{{\bf Taylor}}
    &\small{{\bf non-contractible}}& \small{{\bf Minimal}}\\
    \hline
    3&5&56/26&20&32&12&11\\
    3&12&57/28&158&4096&48&45\\
    6&5&128/45&26&32&19&19\\
    6&12&114/60&467&4096&225&169\\
    9&7&168/73&75&128&43&37\\
    9&8&176/101&157&256&111&103\\
    
    \hline
    \end{tabular}
\caption{Number of non-contractible complexes, compared with $lcm$-lattice and Taylor and Minimal resolutions}
\label{contractible_corners}
\end{center}
\end{table}
\begin{Remark}
In \cite{R07}, the following simplicial complex is defined:
$$\Ac_I:=\{M\subset min(I)\vert \pi(lcm(M))\notin I\}$$
where $min(I)$ is the set of minimal generators of $I$, and $\pi(m)=\frac{lcm(m,x_1\cdots x_n)}{x_1\cdots x_n}$. Observe that the multidegrees of the faces of $\Ac_I$ are exactly those that we have eliminated using the above criterion, i.e. the numbers in the column \emph{non-contractible} in table \ref{contractible_corners} are the number of different multidegrees of the faces of $\Ac_I$. The complex $\Ac_I$ supports a simplicial resolution of $I$, which is a subresolution of Taylor's, the size of which is bounded below by the numbers under the \emph{non-contractible} column in table \ref{contractible_corners}
\end{Remark}

\begin{Remark}
This is yet another easy proof of the fact that the Koszul homology of $I$ vanishes at all multidegrees not in the $lcm$-lattice $L_I$, because if $\ab$ is not in $L_I$ then either $\Delta_{\ab}^I$ is empty (if $\ab \notin I$) or it is the full $l$-simplex (if $\ab\in I$); in both cases, the homology vanishes.
\end{Remark}

\begin{Remark}
We can also use dualities to improve our computations. The simplicial version of Alexander's duality is in particular useful for our purpose, see\cite{B96,MS04}. The {\it Alexander Dual} of a simplicial complex given by a collection $S$ of subsets of the $n$-simplex $\Delta$ is given by
$$S^\vee=\{F|F^c\notin S\}=\{F|F\notin S^c\}$$
where $ ^c$ denotes set complementing. The \emph{Alexander duality theorem} states that $\widetilde{H_i}(S)\simeq \widetilde{H}^{n-i-3}(S^\vee)$ and $\widetilde{H}^i(S)\simeq \widetilde{H}_{n-i-3}(S^\vee)$. Thus, depending on the size of our complex inside the $l$-simplex, it is sometimes better from the computational point of view, to compute in the Alexander dual of $\Delta_\ab^I$. Some other forms of duality can also be of interest for further improvements.
\end{Remark}

With the considerations made above, we can make several improvements to our first \emph{naive} algorithm for compute the Koszul homology of a monomial ideal. Using topological and algebraic arguments, we can reduce the number of corners in which we have to actually perform these computations, moreover, these computations are improved and we only compute the homology of the resulting simplicial complexes at the particular dimensions we need. 

\begin{table}[!htb]
\centering
\begin{tabular}{|p{15cm}|}
\hline
$$\begin{array}{ll} Algorithm:\, \mbox{Dimension-by-dimension Koszul Homology of a Monomial ideal }$I$\\
\hline \hline\\
\ Input: \, \mbox{Minimal generating set of } I:\{m_1,\dots m_r\}\\
\ Output: \mbox{Generators of the Koszul homology of }I\\
\\
\hline\\
\ 1 \quad H_0(K(I))=<m_1,\dots, m_r>\\
\ 2 \quad S\leftarrow \{{\rm multideg}(m_1),\dots,{\rm multideg}(m_r)\}\\
\ 3 \quad i:=1 \\
\ 4 \quad {\bf while }\, i<n \,{\bf and} \,{\vert{S}\vert>1}\, {\bf do}\\
\ 5 \quad \quad lcm(S)\leftarrow \{lcm(s_k,s_j)\vert s_k,s_j\in S\}\\
\ 6 \quad \quad S=\{\},\,H_i(K(I))=\{\}\\
\ 7 \quad \quad {\bf foreach}\, s_{k,j}\in lcm(S)\, {\bf do} \\
\ 8 \quad \quad \quad {\bf if}\, \mbox{Not is-contractible}(s_{k,j})\, {\bf then}\\
\ 9 \quad \quad \quad \mbox{ Build }(\Delta_{s_{k,j}}^I)_i  \\
\ 10 \quad \quad  \quad \mbox{Compute }\hat H_i((\Delta_{s_{k,j}}^I)_i) \\
\ 11 \quad \quad \quad {\bf if}  \hat H_i((S_{s_{k,j}}^I)_i)\neq 0\, {\bf then}\\
\ 12 \quad \quad \quad \quad S\leftarrow S\cup s_{k,j}\\
\ 13 \quad \quad \quad \quad H(K(I))_i\leftarrow \hat H_i(\Delta_{s_{k,j}}^I)\\
\ 14 \quad \quad \quad {\bf endif}\\
\ 15 \quad \quad \quad {\bf endif}\\
\ 16 \quad \quad {\bf endforeach}\\
\ 17 \quad i=i+1\\
\ 18 \quad {\bf return} H_*(K(I))\\
\ 19 \quad {\bf endwhile}\\
\end{array}
$$
\\
\hline
\end{tabular}
\caption{Algorithm {\bf Dimension-by-dimension Koszul Homology}}\label{algorithm_dim_by_dim}

\end{table}
The improved algorithm can be seen in table \ref{algorithm_dim_by_dim}. Observe that in line 9 we use a subroutine {\it Build} that builds the part of the chain complex of $\Delta_{\ab,i}^I$ that we need for our computations, i.e. the generators of the modules of $(i+1)$, $i$ and $(i-1)$-chains, and the matrices of the differentials $\delta_{i+1}$ and $\delta_i$. In line 10, the routine {\it Compute} can make use of the several topological properties and dualities available to improve the computations of the generators of the Koszul homology modules of $I$. 
\subsection{Stanley-Reisner ideals and Alexander duality}\label{stanley-alexander}
\spanishsubsection{Ideales de Stanley-Reiner y dualidad de Alexander}
The theory of Stanley-Reisner rings is one of the most fruitful approaches to monomial ideals \cite{S96}. This theory relates simplicial complexes and squarefree monomial ideals. Although we will not go deep in this direction, we will apply the basics of the theory to the simplicial Koszul complexes we have seen. For this, we will use a simplicial version of the Alexander duality, which has been deeply explored in \cite{MS04}, see also \cite{ER98,RV01}.
\begin{Definition}
Let $\Delta$ be a simplicial complex. The {\rm Stanley-Reisner ideal} of $\Delta$, is the (squarefree) ideal $I_\Delta=\langle \xb^\tau\vert \tau\notin\Delta\rangle$.
\end{Definition}

If we denote by $\mm^\ab$ the prime ideal $\mm^\ab=\langle x_i\vert i\in\ab\rangle$ and $\bar\tau=\{1,\dots,n\}\smallsetminus \tau$ the complement of the set $\tau$, then we have an alternative way to describe $I_\Delta$ in terms of prime components, namely $I_\Delta=\bigcap_{\sigma\in\Delta}\mm^{\bar\sigma}$. Observe how this fact clarifies the natural identification between prime monomial ideals and face ideals we have seen in proposition \ref{prime}.

Alexander duality arises in the homological study of spheres and has a simplicial version \cite{B96} that relates the homology of a simplicial complex with that of its Alexander dual, which is defined as follows:
\begin{Definition}
Let $\Delta$ be a simplicial complex , then the Alexander dual of $\Delta$ is the complex
$$\Delta^*=\{\sigma\vert \bar\sigma\notin\Delta\}$$
\end{Definition}

We can also define the Alexander dual of a squarefree monomial ideal
\begin{Definition}
Let $I=\langle \xb^{\ab_1},\dots,\xb^{\ab_r}\rangle$ be a squarefree monomial ideal. Its Alexander dual is the ideal $I^*=\mm^{\ab_1}\cap\cdots\cap\mm^{\ab_r}$
\end{Definition}

Both Alexander duals for simplicial complexes and ideals are coherent with the operation of taking Stanley-Reisner ideals, i.e. the identity $I_{\Delta^*}=I^*_\Delta$ holds for any simplicial complex $\Delta$.

\begin{Remark}
A generalization of squarefree duality for non-squarefree monomial ideals is given in \cite{MS04} as an expression of the duality between minimal generators and irreducible components of monomial ideals.
\end{Remark}
Our goal here is to relate the homology of a simplicial complex with that of its Stanley Reisner ideal, this requires the following definitions:
\begin{Definition}
The {\rm link} of a face $\sigma$ of the simplicial complex $\Delta$ is given by
$$link_\Delta(\sigma)=\{\tau\in\Delta\vert\tau\cup\sigma\in\Delta\mbox{ and }\tau\cap\sigma=\emptyset\}$$
\end{Definition}
\begin{Definition}
The restriction of $\Delta$ to $\sigma\in\Delta$ is 
$$\Delta\vert_{\sigma}=\{\tau\in\Delta\vert\tau\subseteq\sigma\}$$
\end{Definition}
Using these definitions and Theorem \ref{simp-kos-hom} we have the following dual results:
\begin{Corollary}\label{hochster1}
$\beta_{i-1,\sigma}(I_\Delta)=\beta_{i,\sigma}(R/I_{\Delta})=dim_\kb\tilde H_{\vert\sigma\vert-i-1}(\Delta\vert_\sigma)$
\end{Corollary}
\begin{Corollary}\label{hochster2}
$\beta_{i,\sigma}(I_\Delta)=\beta_{i+1,\sigma}(R/I_{\Delta})=dim_\kb\tilde H_{i-1}(link_{\Delta^*}(\bar\sigma))$
\end{Corollary}

Corollary \ref{hochster1} is known as Hochster's formula \cite{H77,ER98}, corollary \ref{hochster2} is just a reformulation of it.

Take now our Koszul simplicial complexes, and use the Stanley-Reisner machinery together with Hochster's formula, we have the following definition:
\begin{Definition}
Let $I$ be a monomial ideal and $\ab\in\NN^n$, we say that the {\rm upper Koszul ideal of $I$ at $\ab$}, denoted by $KI_\ab$, is the Stanley-Reisner ideal of the upper simplicial Koszul complex of $I$ at $\ab$, $\Delta^I_\ab$. We also define the {\rm lower Koszul ideal of $I$ at $\ab$}, denoted by $KI^\ab$, as the Stanley-Reisner ideal of the lower simplicial Koszul complex of $I$ at $\ab$, $\Delta_I^\ab$. Since upper and lower simplicial complexes are Alexander duals, their corresponding Stanley-Reisner ideals are Alexander dual to each other:
$$KI_\ab=I_{\Delta^I_\ab}=I_{\Delta^\ab_I}^*=(KI^\ab)^*.$$
\end{Definition}
We can use the upper and lower Koszul ideals to make local computations of the Koszul homology of $I$ at a given multidegree:
\begin{Proposition}\label{koszul-stanley-betti}
Let $\ab\in\NN^n$ and let $\sigma=supp(x^\ab)$, then
$$\beta_{i,\ab}(I)=\beta_{i,\sigma}(KI^\ab)=\beta_{\vert\sigma\vert-i-1,\sigma}(KI_\ab)$$
\end{Proposition}
\noindent{\bf Proof: }
Let $\Delta^I_\ab$ be the simplicial Koszul complex of $I$ at $\ab$, it is a complex on the vertices corresponding to the variables in $\sigma$.

From proposition \ref{simp-kos-hom} we have that $\beta_i,\ab(I)=dim_\kb \tilde H_{i-1}(\Delta^I_\ab)$. On the other hand,from corollary \ref{hochster2} we have that $\tilde H_{i-1}(link_{{\Delta^{I}_\ab}}(\emptyset)) =\beta_{i,\sigma}(I_{{\Delta^{I}_\ab}^*})$. Observe that in this case $\bar\sigma=\emptyset$ and since $link_{{\Delta^{I}_\ab}}(\emptyset)=\Delta^{I}_\ab$ we have the result. The last equality comes from the equality $\beta_{\vert\sigma\vert-i,\sigma}(I)=\beta_{i,\sigma}I^*$ which comes immediately from \cite[Theorem 5.48]{MS04} $\square$
\begin{Remark}
As we have seen in the algorithm for computing Koszul homology dimension by dimension, we are usually interested in computing one dimension $k$ of the Koszul homology, and thus we use a subcomplex of $\Delta^I_\ab$, which has the same homology as it in dimension $k$. This $k$-constrained simplicial Koszul complex has Stanley-Reisner ideal, the $k$-constrained Koszul ideal of $I$ at $\ab$. The last proposition can be applied with respect to this smaller ideal.
\end{Remark}

\begin{Example}
Consider the following ideal in $\kb[x,y,z,t,u]$:
$$I=\langle x^{12}y^{2}z^{16}t^{5}u^{4}, xy^{17}z^{3}t^{15}, x^{11}y^{9}z^{9}t^{20}u^{10}, x^{4}y^{13}z^{19}u^{5}, x^{7}y^{14}z^{18}tu^{6}, x^{14}y^{12}z^{9}t^{7}u^{11}\rangle$$ 
The monomial $\xb^\ab$ with $\ab=(12,14,19,20,10)$ is in $L_I$. In order to check whether $H_{*,\ab}(I)$ is zero, we build the simplicial Koszul complex $\Delta^I_\ab$, the facets of which are: $$\{\{1,2,4,5\},\{1,3,4,5\},\{2,3,4,5\},\{1,2,3\}\}$$
and thus its Stanley-Reisner ring is $KI_\ab=\langle xyzt,xyzu\rangle$.
The facets of its Alexander dual are just $\{\{4\},\{5\}\}$
and its Stanley-Reisner ideal is $I_{{\Delta^{I}_\ab}^*}=\langle x,y,z,tu\rangle$.
Applying proposition \ref{koszul-stanley-betti} we have that the dimension of the second reduced homology group of $\Delta^I_\ab$ equals the first Betti number at multidegree $\sigma=xyztu$ of $KI_\ab=I_{{\Delta^I_\ab}}$ and the third Betti number at multidegree $\sigma=xyztu$ of $KI^\ab=I_{{\Delta^I_\ab}^*}$, which are both equal to one, thus, $\beta_{3,\ab}(I)=dim_{\kb}\tilde{H}_2(\Delta^I_\ab)=\beta_{1,xyztu}(IK_\ab)=\beta_{3,xyztu}(IK^\ab)=1$. Observe that this is the only nonzero $\beta_{*,\ab}(I)$.
\end{Example}

\subsection{Mayer-Vietoris sequence in Koszul homology}\label{m-v}
\spanishsubsection{Sucesi\'on de Mayer-Vietoris en homolog\'ia de Koszul}

Sometimes tools coming from topology give some intuitional hints for solving problems in algebra. Algebraic topology is very often used as a way to translate these tools from topology into algebraic tools. In this section, we will describe an analogue of the Mayer-Vietoris sequence coming from topology (see Appendix \ref{aptop}) that can be used for making computations in monomial ideals. The Mayer-Vietoris sequence relates the homology of two spaces, their intersection and their union by means of a short exact sequence of complexes. With the help of Mayer-Vietoris sequences we build in Chapter \ref{computation} a recursive algorithm for homological computations on monomial ideals, which is a main contribution of this thesis. The algorithmic counterpart of Mayer-Vietoris sequences will be called Mayer-Vietoris trees.

\subsubsection{A short exact sequence of Koszul complexes}
\spanishsubsubsection{Una sucesi\'on exacta corta de complejos de Koszul}
In this section we will use the following notation:  $I$ is a monomial ideal in $R=\poly \kb x n$ minimally generated by $\{m_1,\dots,m_r\}$, $I'$ is the ideal generated by $\{m_1,\dots,m_{r-1}\}$ and $\tilde I=I'\cap\langle m_r\rangle$ is the ideal generated by $\{m_{1,r},\dots,m_{r-1,r}\}$, where $m_{i,j}=lcm(m_i,m_j)$. Also, when it is convenient, we will use the following alternative notation: for each $1\leq s\leq r$ denote $I_s:=\langle m_1,\dots m_s\rangle$, $\tilde{I}_s:=I_{s-1}\cap\langle m_s\rangle=\langle m_{1,s},\dots,m_{s-1,s}\rangle$.

For the explicit computations we also need the following notations for elements of Koszul complexes:
\begin{itemize}
\item For each chain $\gamma$ of $\KK(I)$, there is a finite set $\Ac$ of integers such that $\gamma$ can be written as $\gamma=\sum_{i\in\Ac} \alpha_i x^{\mu_i}\otimes x^{J_i}$. 
\item Let $\Bc$ be a finite set of integers such that $\gamma=\sum_{i\in\Bc} \alpha_i x^{\mu_i}\otimes x^{J_i}$ is a chain in $\KK({I'})$. If we denote by $\Bc_{\tilde I}=\{i\in\Bc \,s.th\, x^{\mu_i}\in \tilde I\}$, $\Bc_{I'\smallsetminus\tilde I}=\{i\in\Bc\, s.th \,x^{\mu_i}\in I', x^{\mu_i}\notin\tilde I\}$ then $\Bc=\Bc_{\tilde I}\sqcup \Bc_{I'\smallsetminus \tilde I}$ and $\gamma= \sum_{i\in\Bc_{\tilde I}} \alpha_i x^{\mu_i}\otimes x^{J_i}+\sum_{i \in\Bc_{I'\smallsetminus \tilde I}} \alpha_i x^{\mu_i}\otimes x^{J_i}=\gamma_{\tilde I}+\gamma_{I' \smallsetminus \tilde I}$ being the sum direct.
\item Similarly, if $\gamma\in \KK(\langle m_r \rangle)$ then we have a direct sum decomposition $\gamma=\gamma_{\langle m_r \rangle\smallsetminus \tilde I}+\gamma_{\tilde I}$.
\item Finally, if $\gamma$ is a chain in $\KK(I)$ we have that $\gamma=\gamma_{I'\smallsetminus \langle m_r\rangle}+\gamma_{\tilde I}+\gamma_{\langle m_r \rangle\smallsetminus I'}$.

\end{itemize}With this notation:

\begin{Proposition}\label{mv}
We have the following effective short exact sequence:
\begin{center}
\begin{tikzpicture}
\draw[->] (0,0)node[left]{$0$}--(0.5,0)node[right]{$\KK(\tilde{I})$};
\draw[->] (7.5,0)node[left]{$\KK(I)$}--(8,0) node[right]{$0$};
\draw[->] (1.5,0.1)--(2.5,0.1)node[above,pos=0.5]{$\mi$};
\draw[->] (2.5,-0.1)--(1.5,-0.1)node[below,pos=0.5]{$\rho$};
\draw[->] (5.5,0.1)--(6.5,0.1)node[above,pos=0.5]{$\mj$};
\draw[->] (6.5,-0.1)--(5.5,-0.1)node[below,pos=0.5]{$\sigma$};
\draw (4,0) node{$\KK(I')\oplus\KK(\langle m_r\rangle)$};
\end{tikzpicture}
\end{center}

i.e., we have that $\mi$ and $\mj$ are chain complex morphisms, $\rho$ and $\sigma$ are graded module morphisms, and:
\begin{enumerate}

\item [i)] $\mi$ is injective, $\mj$ is surjective and $ker(\mj)=im(\mi)$

\item [ii)]$\rho\mi=1_{\KK(\tilde I)}$, $\mj\sigma=1_{\KK(I)}$ and $\mi\rho+\sigma\mj=1_{\KK(I')\oplus\KK(\langle m_r\rangle)}$

\end{enumerate}
\end {Proposition}

\noindent {\bf Proof: } The proof of this proposition is a routine verification of the commutativity of some diagrams. The definitions of the involved morphisms are as follows:
Let $\gamma$ a chain in $\KK(\tilde{I})$, $\eta$ a chain in $\KK(I')$, $\eta'$ a chain in $\KK(\langle m_r \rangle)$ and $\zeta$ a chain in $\KK(I)$. Then, we define 

\begin{eqnarray*}
\mi(\gamma)&=&(\gamma, -\gamma)\\
\mj(\eta,\eta')&=&\eta+\eta'\\
\rho(\eta,\eta')&=&\frac{1}{2}\eta_{\tilde I}-\frac{1}{2}\eta'_{\tilde I}\\ \sigma(\zeta)&=&(\zeta_{I'\smallsetminus \langle m_r\rangle}+\frac{1}{2}\zeta_{\tilde I},\zeta_{\langle m_r\rangle\smallsetminus I'}+\frac{1}{2}\zeta_{\tilde I})
\end{eqnarray*} 
Note that for the sequence to be exact we need $\mi$ and $\mj$ to be morphisms of chain complexes, i.e. they must commute with the corresponding Koszul differential. On the other hand, $\rho$ and $\sigma$ do not need to commute with the differential.$\square$

Given a simple element $x^\mu\otimes x^J$ of $\KK(I)$ we say that the multidegree of it is the multidegree of the product $x^\mu\cdot x^J$. As was seen in section \ref{koszdef} it is then clear that the differential $\partial$ in $\KK(I)$ preserves multidegree, and thus
$$\KK(I)=\bigoplus_{\ab\in\NN^n} \KK_{\ab}(I)$$
where $\KK_{\ab}(I)$ is the `piece' of $\KK(I)$ of multidegree $\ab$. In consequence, we have a `multigraded' version of our previous proposition:

\begin{Proposition}\label{mv2}
The following effective sequence of complexes is exact:
\begin{center}
\begin{tikzpicture}
\draw[->] (0,0)node[left]{$0$}--(0.5,0)node[right]{$\KK_{\ab}(\tilde{I})$};
\draw[->] (8.5,0)node[left]{$\KK_{\ab}(I)$}--(9,0) node[right]{$0$};
\draw[->] (1.75,0.1)--(2.75,0.1)node[above,pos=0.5]{$\mi$};
\draw[->] (2.75,-0.1)--(1.75,-0.1)node[below,pos=0.5]{$\rho$};
\draw[->] (6.25,0.1)--(7.25,0.1)node[above,pos=0.5]{$\mj$};
\draw[->] (7.25,-0.1)--(6.25,-0.1)node[below,pos=0.5]{$\sigma$};
\draw (4.5,0) node{$\KK_{\ab}(I')\oplus\KK_{\ab}(\langle m_r\rangle)$};
\end{tikzpicture}
\end{center}
for all $\ab\in\NN^n$
\end {Proposition}

\noindent{\bf Proof: } The proof is exactly the same that the one for proposition \ref{mv} since $\mi$, $\mj$, $\rho$ and $\sigma$ preserve total multidegree.$\square$

\begin{Remark}
The word \emph{effective} in propositions \ref{mv} and \ref{mv2} refers to effective homology (see Appendix \ref{apef}). In this chapter we only make a slight reference to effective homology, which will be of big importance in section \ref{poly} and chapter \ref{computation}. The methods of effective homology will allow us to develop constructive methods from the theory, and therefore implement actual algorithms to make computations with the objects we deal with.
\end{Remark}

\subsubsection{Mayer-Vietoris sequence in Koszul homology}\label{mvs-koszul}
\spanishsubsubsection{Sucesi\'on de Mayer-Vietoris en homolog\'ia de Koszul}
The short exact sequences of propositions \ref{mv} and \ref{mv2} give rise to  long exact sequences in homology, by means of a connecting homomorphism $\Delta$ given by $\rho\partial\sigma$, where $\partial$ is the Koszul differential. In the multidegree $\ab$ case:

\begin{equation}\label{les}
\cdots\longrightarrow H_{i+1}(\KK_{\ab}(I))\stackrel{\Delta}{\longrightarrow}H_{i}(\KK_{\ab}(\tilde I))\longrightarrow$$ $$H_i(\KK_{\ab}(I)\oplus \KK_{\ab}(\langle  m_r\rangle  ))\longrightarrow H_i(\KK_{\ab}(I))\stackrel{\Delta}{\longrightarrow}\cdots
\end{equation}

\begin{Definition}
Given a monomial ideal $I$ minimally generated by $\{m_1,\dots, m_r\}$ we define the (recursive) Mayer-Vietoris exact sequence of $I$ in the following way:

For each $1\leq s\leq r$ we have the following exact sequence of ideals:
$$0\rightarrow \tilde{I}_s\rightarrow I_{s-1}\oplus \langle m_s\rangle\rightarrow I_s\rightarrow 0$$
and the following short exact sequence of Koszul complexes:

$$0\longrightarrow \KK(\tilde {I}_s)\stackrel{\mi}{\longrightarrow} \KK(I_{s-1})\oplus \KK(\langle  m_s\rangle  )\stackrel{\mj}{\longrightarrow}\KK(I_s)\longrightarrow 0, $$
the maps given by
$$\mi(\gamma)=(\gamma, -\gamma),\qquad
\mj(\eta,\eta')=\eta+\eta'
$$
for $\gamma\in\KK(\tilde{I}_s)$, $\eta\in\KK(I_{s-1})$ and $\eta'\in \KK(\langle m_s\rangle)$.

These sequences induce a long exact sequence in Koszul homology for each $s$, the set of all of them is what we call the (recursive) {\rm Mayer-Vietoris Sequence} of $I$.
\end{Definition}

\begin{Remark}\label{pivot}
Strictly speaking, the definition of Mayer-Vietoris sequences of monomial ideals is not fully precise, in the sense that the Mayer-Vietoris sequence associated to a given ideal is not uniquely defined, it depends on how the minimal generators are sorted. The choice of the last generator of the ideal $I$ to be the one which defines the {\it Mayer-Vietoris sequence} is just a matter of convenience in notation. The important fact is that we select {\it some} particular generator to define the sequence. To this one we associate the subindex $s$ which constitutes the {\it breaking point} to generate the sequence. Several selection strategies can be applied to select the distinguished generator, and they can be changed during the process. Thus, strictly speaking, several Mayer-Vietoris sequences are associated to a given monomial ideal. To avoid cumbersome notation, we will give the chosen strategy when necessary; in particular this will be important in computations, see chapter \ref{computation} where this point is treated in more detail.
\end{Remark}

Using the recursive process that naturally comes from these Mayer-Vietoris sequences at every $\ab\in\NN^n$ we could compute the Koszul homology of $I$ from that of ideals with fewer number of generators. For the recursive computations, we need to solve the trivial case, i.e. the computation of the Koszul homology of ideals with just one generator, but the trivial case is of course trivial:
$$H_i(\KK(\langle m_1\rangle))=
\begin{cases}
\kb&i=0,\\
0&otherwise.\\
\end{cases}$$
and the generator of $H_0(\langle m_1\rangle)$ can be identified with $m_1$ itself.

The other ingredient we need for our recursive process is a explicit formula to obtain preimages of the connecting morphism $\Delta$. One way to do this goes as follows:

Let's describe first $\Delta$: $\Delta$ is given by $\rho\partial\sigma$, where $\partial$ is the differential in $\KK(I')\oplus\KK(\langle m_r\rangle)$. More explicitly: given $\gamma\in\KK(I)$ a representative of its homology class in $H_{i+1}(\KK(I))$, we have that
$\sigma(\gamma)=(\gamma_{I'\smallsetminus \langle m_r\rangle}+\frac{1}{2}\gamma_{\tilde I},\gamma_{\langle m_r\rangle\smallsetminus I'}+\frac{1}{2}\gamma_{\tilde I})$. If we denote $\eta=\partial(\gamma_{I'\smallsetminus \langle m_r\rangle}+\frac{1}{2}\gamma_{\tilde I})$ and $\eta'=\partial(\gamma_{\langle m_r\rangle\smallsetminus I'}+\frac{1}{2}\gamma_{\tilde I})$ then
$$\Delta([\gamma])=[\rho\partial\sigma(\gamma)]=[\rho(\eta,\eta')]=[\frac{1}{2}\eta_{\tilde I}-\frac{1}{2}\eta'_{\tilde I}].$$

Now, given $\gamma$ a representative of its homology class in $H_i(\KK(\tilde I))$, and provided $\Delta$ is an epimorphism, we want to find a preimage of $[\gamma]$ by $\Delta$:  First, apply $\mi$ to $\gamma$, i.e. $\mi(\gamma)=(\gamma,-\gamma)$, now we apply the ``inverse'' of $\partial\oplus\partial$ to it (note that each $\partial$ is in a different complex). This can be done using the Spencer differential $\delta$ (see section \ref{spencer}) for which we know that $\partial\delta(\gamma)=c\cdot\gamma$ for some scalar $c$. So, applying the corresponding $\delta$ in each component of $(\gamma,-\gamma)$, we obtain that $(\partial\delta(\gamma),\partial\delta(-\gamma))=(c_1\cdot \gamma,c_2\cdot(-\gamma))$ for some (possibly different) scalars $c_1$, $c_2$. So take $(\eta,\eta')=(\frac{1}{c_1}\delta(\gamma),\frac{1}{c_2}\delta(-\gamma))$ where each $\delta$ lies in the corresponding complex. Finally apply $\mj$ to $(\eta,\eta')$ to obtain $\mj(\eta,\eta')=\eta+\eta'$ and we have the preimage we were looking for.

\subsection{Mapping cones and resolutions}
\spanishsubsection{Conos de aplicaciones y resoluciones}
Many resolutions in commutative algebra arise as iterated {\it mapping cones} (see definitions in Appendix \ref{aptop}). The basic ideas are as follows (see \cite{HT02} or a slightly different approach in \cite{S99}): Let $I$ be a (non necessary monomial) ideal in $R= \poly \kb x n$ generated by $\{f_1,\dots,f_r\}$ and denote $I_j=(f_1,\dots,f_j)$ then, for $j=1\dots,r$ there are exact sequences
$$0\rightarrow R/(I_{j-1}:f_j)\rightarrow R/I_{j-1}\rightarrow R/I_j\rightarrow 0$$
If resolutions of $R/I_{j-1}$ and $R/(I_{j-1}:f_j)$ are known, we can obtain a resolution of $R/I_j$ as a mapping cone (see \cite{M91} for example) of a homomorphism between the known resolutions. This construction provides an inductive procedure to compute a resolution of $R/I$.

More generally, if we have a short exact sequence of $R$-modules 
$$0\rightarrow A\rightarrow B\rightarrow C\rightarrow 0$$
where $A$ and $B$ are simpler than $C$ and we already have resolutions $\PP(A)$ and $\PP(B)$ of $A$ and $B$. The morphism $A\rightarrow B$ induces a morphism between $\PP(A)$ and $\PP(B)$, the cone of which is a resolution of $C$. Note that even if $\PP(A)$ and $\PP(B)$ are minimal, $\PP(C)$ need not to be minimal. If we have that our modules $A$, $B$ and $C$ can be related in a recursive way, we can exploit these ideas to obtain minimal free resolutions from trivial cases.

\begin{Remark}
In the case of Mayer-Vietoris sequences, given a monomial ideal $I=\langle m_1,\dots,m_r\rangle$, and its related ideals $I'=\langle m_1,\dots,m_{r-1}\rangle$ and $\tilde{I}=\langle m_{1,r},\dots,m_{r-1,r}\rangle$, we have the short exact sequence

$$0\rightarrow \tilde{I}\rightarrow I'\oplus\langle m_r\rangle\rightarrow I\rightarrow 0$$

Therefore, the starting point for the recursive process we want to implement is the short exact sequence of ideals associated to the Mayer-Vietoris sequence of an ideal $I$. 
Details on the process of constructing resolutions using Mayer-Vietoris sequences will be given in section \ref{computations}, where the process is explained in an explicit manner.
\end{Remark} 

One of the main difficulties when using the mapping cone techniques to build resolutions is the explicit construction of the comparison morphism, i.e. the chain complex morphism between $A$ and $B$. In the case of resolutions of ideals, we usually have an inclusion between modules or ideals and we want to lift it to a chain complex morphism between the corresponding resolutions. This lifting is theoretically very easy, and in fact very simple existence theorems exist, but on the other hand, the construction of such liftings from given morphisms is not always a trivial task. At this point one can use the constructive techniques from \emph{effective homology}, see \cite{S06} and appendix \ref{apef}. The main tool from effective homology used here is the following theorem:

\begin{Theorem}[$SES_3$ in \cite{S06}] \label{ses3}
Let $(A,\mi,\rho,B,\mj,\sigma,C)$ be an effective short exact sequence of effective chain-complexes:
\begin{center}
\begin{tikzpicture}
\draw[->] (0,0)node[left]{$0$}--(1,0)node[right]{$A$};
\draw[->] (5,0)node[left]{$C$}--(6,0) node[right]{$0$};
\draw[->] (1.5,0.1)--(2.5,0.1)node[above,pos=0.5]{$\mi$};
\draw[->] (2.5,-0.1)--(1.5,-0.1)node[below,pos=0.5]{$\rho$};
\draw[->] (3.5,0.1)--(4.5,0.1)node[above,pos=0.5]{$\mj$};
\draw[->] (4.5,-0.1)--(3.5,-0.1)node[below,pos=0.5]{$\sigma$};
\draw (3,0) node{$B$};
\end{tikzpicture}
\end{center}

where $\mi$ and $\mj$ are chain complex morphisms; $\rho$ and $\sigma$ are graded module morphisms and the following relations hold: $id_A=\rho \mi$, $id_B=\mi\rho+\sigma \mj$ and $id_C=\mj\sigma$,
then, an algorithm constructs a canonical reduction (see \cite{RS02}) between $Cone(\mi)$ and $C$ from the data.
\end{Theorem}

\begin{Remark}
The theorem produces a canonical reduction betwee $Cone(\mi)$ and $C$, a reduction between two chain complexes is essentially an explicit map between a big chain complex and a small one which have the same effective homology. Details about reductions and effective homology can be seen in Appendix \ref{apef}. Here it is enough to say that using this reduction, we can fully describe $C$ as a chain complex and the homology of it from the complex $Cone(\mi)$ and its correspondent homology.
\end{Remark}

To use the precedent theorem in the recursive setting of \emph{mapping cone resolutions} we need {\it effective} resolutions of the leftmost ideals, i.e. we need explicit contracting homotopies for them. Normally, these are easy to give, because we only need an explicit description of them for the trivial cases; the reduction obtained in theorem \ref{ses3} produces effective resolutions in the recursive step. For example, in the case of monomial ideals, the initial step works with ideals generated by one and two monomials, since we know that the Taylor resolution of them is minimal, and moreover, an explicit contracting homotopy is known for it \cite{F78} thus, theorem \ref{ses3} provides us with an {\it effective} resolution of $I_2$ and we can go on with the process.

We have then the following process to build resolutions of monomial ideals:
\begin{itemize}
\item {Consider an \emph{effective} short exact sequence of monomial ideals (also of quotients of $R$ by monomial ideals)
\begin{center}
\begin{tikzpicture}
\draw[->] (0,0)node[left]{$0$}--(1,0)node[right]{$I$};
\draw[->] (5,0)node[left]{$K$}--(6,0) node[right]{$0$};
\draw[->] (1.5,0.1)--(2.5,0.1)node[above,pos=0.5]{$\mi$};
\draw[->] (2.5,-0.1)--(1.5,-0.1)node[below,pos=0.5]{$\rho$};
\draw[->] (3.5,0.1)--(4.5,0.1)node[above,pos=0.5]{$\mj$};
\draw[->] (4.5,-0.1)--(3.5,-0.1)node[below,pos=0.5]{$\sigma$};
\draw (3,0) node{$J$};
\end{tikzpicture}
\end{center}
such that $I$ and $J$ are minimally generated by fewer generators than $K$, so that we can build a recursive process from it. Examples are Mayer-Vietoris sequences, or sequences of the types given in \cite{RS06}, \cite{S99} or \cite{HT99}.}
\item {Start the recursive process with ideals generated by one or two monomials. For them, consider an effective resolution, i.e. a resolution with an explicit contracting homotopy, for example Taylor resolution. Use the effective resolution to lift the map $\mi$ to a chain complex morphism. Then use theorem \ref{ses3} to obtain an effective resolution of the third ideal in the sequence.}
\item {Plug the obtained resolution into the next recursion step and do the same as in the precedent step (lifting the map $\mi$ and applying theorem \ref{ses3}) and proceed forward in the same way. At the end of the process we obtain an \emph{effective} resolution of our ideal.}
\end{itemize}

\begin{Remark}\label{remark-mapping-cones}
The use of effective short exact sequences allows then to overcome the main (algorithmic) difficulty when using this recursive procedure, namely the construction of the comparison maps $\mi$ \cite{HT02}. This procedure gives us also a good recursive description of the differentials involved in the so constructed resolution. We know that $Cone(\mi)_q=B_q\oplus A_{q-1}$ and the differential is given by $d_q^{Cone(\mi)}=\left ( \begin{array}{cc}
d_q^B&\mi_{q}\\0&-d_{q-1}^A
\end{array}\right )$. Thus, if we keep minimality at each step, we know that the only possible part of the matrix which can be reduced is that corresponding to $\mi_q$, and the minimization process is improved. Moreover, in the case of multigraded resolutions (e.g. those of monomial ideals)  as we keep track of the multidegrees involved, only when the same multidegree appears in the resolutions of both $\tilde I_s$ and $I_{s-1}\oplus\langle m_s\rangle$ at the same homological dimension we can have some non-minimality on the resolution of $I_s$, i.e. these are the only multidegrees in which we can eventually find a reduction pair in the resolution.
\end{Remark}
\fancyhead[ER]{\itshape Chapter 2 \ \ Koszul Homology and the Structure of Monomial Ideals} 

\chapter{Koszul Homology and the Structure of Monomial Ideals}\label{structure}
\spanishchapter{Homolog\'ia de Koszul y estructura de ideales monomiales}

In this chapter we investigate to what extent the knowledge of the Koszul homology of a monomial ideal gives us a knowledge of the ideal itself. Of course, there are different levels of knowledge of the Koszul homology of a monomial ideal $I$, ranging from the ranks of the homology modules to the explicit knowledge of a generating set of each of the modules or to the knowledge of the \emph{effective Koszul homology} of the ideal, in the sense of appendix \ref{apef}.
The first natural aspect to investigate is the relation between the Koszul homology of $I$ and its other homological invariants, in particular the minimal free resolution, and invariants related to it, like the Betti numbers, projective dimension, Castelnuovo-Mumford regularity, Hilbert function, etc. This is explored in the first section.
Not so obvious appears to be the relation between the Koszul homology of $I$ and the algebraic structure of it, given by combinatorial and irreducible decompositions, as well as primary decompositions, associated primes or height. We will see that the combinatorial nature of monomial ideals allows us to discover direct relations between these different aspects of them.

In the first section we describe how to obtain Betti numbers and resolutions from Koszul homology. The work is based on two bicomplexes, namely the $Tor$ bicomplex and the \emph{Aramova-Herzog} bicomplex \cite{AH95}. Using homological perturbation techniques explicit methods are obtained to compute resolutions from Koszul generators. Main references here are \cite{AH95,AH96} and \cite{RS06,RSS06}.

In the second section, given a monomial ideal $I\subseteq R=\kb[x_1\dots,x_n]$ and its Koszul homology, we obtain a Stanley decomposition of $R/I$. Here we treat the zero-dimensional case, and based on it, the general one. In this type of decompositions, the $(n-1)$-st Koszul homology plays an important role, due to its relation to the boundary of the ideal.

In the next section we face the problem of obtaining an irredundant irreducible decomposition of a monomial ideal from its Koszul homology. Again, the artinian case its treated first, and for the non-artinian case we give two different procedures, one based on the artinian case and an independent one. These methods are an alternative to other ways of obtaining irreducible decompositions existing in the literature \cite{MS04,M04,R07}.

In the fourth section, a simple procedure for obtaining a minimal primary decomposition of a monomial ideal from its Koszul homology is given. Having this decomposition we can read the height of the ideal and its associated primes. Again, the $(n-1)$-st homology plays a crucial role.

Finally, in the last section we face polynomial ideals. The initial ideal of a polynomial ideal $I$ is a monomial ideal to which we can apply the known techniques to compute its Koszul homology and/or minimal resolution. From this, using homological perturbation and Gr\"obner bases techniques, we obtain the Koszul homology and/or resolutions of $I$. This technique has been used in two forms in the literature, one by L. Lambe and others \cite{JLS02,LS02} perturbing Taylor and Lyubeznik resolutions. The other by F. Sergeraert \cite{RS06,S06,RSS06} using effective mapping cones.

\section{Minimal resolutions and Betti numbers}\label{structure_minres}
\spanishsection{Resoluciones m\'inimas y n\'umeros de Betti}
Even if, as we saw in section \ref{tor} the equivalence between the Koszul homology of an $R$-module and its minimal free resolution is quite evident, it is not trivial to obtain an explicit expression of one of them given an explicit expression of the other or vice versa. The two directions of this equivalence give raise to two different problems:

In the first one, given an explicit minimal free resolution of an $R$-module, we want to obtain an explicit set of generators of its Koszul homology modules. This problem was treated by J. Herzog in \cite{H92}. We will not be concerned with this problem, although in chapter \ref{computation} we present a structure that will allow us to obtain both minimal resolutions and generating sets of the Koszul homology for monomial ideals.

The second problem is the opposite one, namely given an explicit set of generators of the Koszul homology, obtain an explicit free resolution. We present here two different approaches to this problem. The first one is contained in the work of A. Aramova and J. Herzog \cite{AH95,AH96}, and makes use of homological tools such as bicomplexes and spectral sequences. The second, more recent approach is more computationally oriented, and using the work of Aramova and Herzog as starting point, develops actual computations using effective homology. This approach is due to F. Sergeraert, and is presented in \cite{RS06,RSS06}.

We start the section with a short clarification of what one can understand by ``knowing'' the Koszul homology of an ideal and what can be extracted from this knowledge, and with an explicit example that illustrates the basic facts of the equivalence between Koszul homology and minimal resolutions:
\begin{itemize}
\item Assume we know the dimensions of the different $H_i(\KK(I))$, then due to the $Tor^R_\bullet(I,\kb)$ equivalence, we have the Betti numbers of $I$, i.e. the ranks of the modules in the minimal resolution, the same holds for (multi)graded Betti numbers. If we know $H_{i,\ab}(\KK(I))$ for every $\ab\in\NN^n$ and $i\in \NN$ then we have the multigraded Betti numbers of $I$, i.e. the multigraded components of the modules in the minimal resolution. Some other invariants can be read automatically from the (graded, multigraded) Betti numbers, such as the projective dimension, Castelnuovo-Mumford regularity, depth, dimension, etc.
\item If we have an explicit generating set for every Koszul homology module, we of course know the graded and multigraded Betti numbers of $I$, but the situation is more favorable. In this case, from the multidegree of every generator of the $H_i(\KK(I))$ we can build the multigraded components of the modules in the minimal resolution of $I$. We still lack the differentials. On the other hand, we have an explicit knowledge of the generators of the homology modules, i.e. the explicit structure of them. From this information we still need to work on the \emph{Tor bicomplex} to obtain explicit differentials of the minimal resolution. We know the inner structure of the modules in the Koszul homology, but we still need to describe somehow the relations among the generators of the different modules. Here comes the different versions of the work on the $Tor$ bicomplex explained in \cite{AH95,AH96,RS06,RSS06}.
\end{itemize}
Let us start with the following example:
\begin{Example}\label{example-complete}
Consider the following ideal in $R:=\kb[x_1,\dots,x_5]$:
\[ I=\langle x_1^{19}x_2^{13}x_3^{5}x_4^{3}x_5^{6}, x_1^{20}x_2^{14}x_3^{4}x_4^{4}x_5^{5}, x_1^{12}x_3x_4^{7}x_5^{3}, x_1^{2}x_2^{14}x_3^{17}x_5^{3}, x_1^{13}x_2^{11}x_3^{12}x_4^{4}x_5^{8} \rangle \]
Bases of its Koszul homology modules are given by the following generators:
\begin{itemize}
\item Generators in dimension $0$:
$$\begin{array}{l}
g^0_1=x_1^{19}x_2^{13}x_3^{5}x_4^{3}x_5^{6}\\
g^0_2=x_1^{20}x_2^{14}x_3^{4}x_4^{4}x_5^{5}\\
g^0_3=x_1^{12}x_3x_4^{7}x_5^{3}\\
g^0_4=x_1^{2}x_2^{14}x_3^{17}x_5^{3}\\
g^0_5=x_1^{13}x_2^{11}x_3^{12}x_4^{4}x_5^{8}
\end{array}
$$
\item Generators in dimension $1$:
\begin{equation}\begin{split}
g^1_1&=x_1^{12}x_2^{13}x_3^{17}x_4^{7}x_5^{3}\otimes x_2-x_1^{12}x_2^{14}x_3^{17}x_4^{6}x_5^{3}\otimes x_4\\
g^1_2&=x_1^{19}x_2^{12}x_3^{5}x_4^{7}x_5^{6}\otimes x_2-x_1^{19}x_2^{13}x_3^{5}x_4^{6}x_5^{6}\otimes x_4\\
g^1_3&=x_1^{20}x_2^{14}x_3^{4}x_4^{6}x_5^{5}\otimes x_4-x_1^{20}x_2^{13}x_3^{4}x_4^{7}x_5^{5}\otimes x_2  \\
g^1_4&=x_1^{13}x_2^{11}x_3^{12}x_4^{6}x_5^{8}\otimes x_4-x_1^{13}x_2^{10}x_3^{12}x_4^{7}x_5^{8}\otimes x_2\\
g^1_5&=x_1^{19}x_2^{13}x_3^{17}x_4^{3}x_5^{6}\otimes x_2-x_1^{19}x_2^{14}x_3^{17}x_4^{2}x_5^{6}\otimes x_4  \\
g^1_6&=x_1^{20}x_2^{14}x_3^{16}x_4^{4}x_5^{5}\otimes x_3-x_1^{20}x_2^{14}x_3^{17}x_4^{3}x_5^{5}\otimes x_4\\
g^1_7&=x_1^{13}x_2^{13}x_3^{17}x_4^{4}x_5^{8}\otimes x_2-x_1^{13}x_2^{14}x_3^{17}x_4^{3}x_5^{8}\otimes x_4  \\
g^1_8&=x_1^{20}x_2^{14}x_3^{4}x_4^{4}x_5^{6}\otimes x_3-x_1^{20}x_2^{13}x_3^{5}x_4^{4}x_5^{6}\otimes x_2\\
g^1_9&=x_1^{19}x_2^{13}x_3^{12}x_4^{4}x_5^{7}\otimes x_5-x_1^{19}x_2^{12}x_3^{12}x_4^{4}x_5^{8}\otimes x_2 
\nonumber
\end{split}\end{equation}

\item Generators in dimension $2$:
\begin{equation}\begin{split}
g^2_1&=x_1^{19}x_2^{14}x_3^{17}x_4^{6}x_5^{5}\otimes x_4x_5-x_1^{19}x_2^{13}x_3^{17}x_4^{7}x_5^{5}\otimes x_2x_5 +x_1^{19}x_2^{13}x_3^{17}x_4^{6}x_5^{6}\otimes x_2x_4\\
g^2_2&=x_1^{20}x_2^{13}x_3^{16}x_4^{7}x_5^{5}\otimes x_2x_3-x_1^{20}x_2^{13}x_3^{17}x_4^{7}x_5^{4}\otimes x_2x_5 +x_1^{20}x_2^{14}x_3^{16}x_4^{6}x_5^{5}\otimes x_3x_4\\
g^2_3&=x_1^{13}x_2^{14}x_3^{17}x_4^{6}x_5^{7}\otimes x_4x_5-x_1^{13}x_2^{13}x_3^{17}x_4^{7}x_5^{7}\otimes x_2x_5 +x_1^{13}x_2^{13}x_3^{17}x_4^{6}x_5^{8}\otimes x_2x_4\\
g^2_4&=x_1^{20}x_2^{13}x_3^{4}x_4^{7}x_5^{6}\otimes x_2x_3-x_1^{20}x_2^{13}x_3^{5}x_4^{6}x_5^{6}\otimes x_2x_4 +x_1^{20}x_2^{14}x_3^{4}x_4^{6}x_5^{6}\otimes x_3x_4\\
g^2_5&=x_1^{19}x_2^{12}x_3^{11}x_4^{7}x_5^{8}\otimes x_2x_3-x_1^{19}x_2^{12}x_3^{12}x_4^{6}x_5^{8}\otimes x_2x_4 +x_1^{19}x_2^{13}x_3^{11}x_4^{6}x_5^{8}\otimes x_3x_4\\
g^2_6&=x_1^{20}x_2^{13}x_3^{16}x_4^{4}x_5^{6}\otimes x_2x_3-x_1^{20}x_2^{13}x_3^{17}x_4^{3}x_5^{6}\otimes x_2x_4 +x_1^{20}x_2^{14}x_3^{17}x_4^{3}x_5^{5}\otimes x_4x_5\\
&\quad+x_1^{20}x_2^{14}x_3^{16}x_4^{4}x_5^{5}\otimes x_3x_5\\
g^2_7&=x_1^{18}x_2^{14}x_3^{16}x_4^{4}x_5^{8}\otimes x_1x_3-x_1^{19}x_2^{13}x_3^{16}x_4^{4}x_5^{8}\otimes x_2x_3 +x_1^{18}x_2^{14}x_3^{17}x_4^{4}x_5^{7}\otimes x_1x_5\\
&\quad+x_1^{19}x_2^{13}x_3^{17}x_4^{4}x_5^{7}\otimes x_2x_5
\nonumber
\end{split}\end{equation}
\item Generators in dimension $3$:
\begin{equation}\begin{split}
g^3_1&=x_1^{20}x_2^{14}x_3^{16}x_4^{6}x_5^{5}\otimes x_3x_4x_5-x_1^{19}x_2^{14}x_3^{17}x_4^{6}x_5^{5}\otimes x_1x_4x_5 +x_1^{19}x_2^{14}x_3^{16}x_4^{7}x_5^{5}\otimes x_1x_3x_5\\
&\quad-x_1^{19}x_2^{14}x_3^{16}x_4^{6}x_5^{6}\otimes x_1x_3x_4\\
g^3_2&=x_1^{18}x_2^{14}x_3^{16}x_4^{6}x_5^{8}\otimes x_1x_3x_4-x_1^{18}x_2^{14}x_3^{16}x_4^{7}x_5^{7}\otimes x_1x_3x_5 -x_1^{19}x_2^{14}x_3^{16}x_4^{6}x_5^{7}\otimes x_3x_4x_5\\
&\quad+x_1^{18}x_2^{14}x_3^{17}x_4^{6}x_5^{7}\otimes x_1x_4x_5
\nonumber
\end{split}\end{equation}
\end{itemize}
So, we have that $$\beta_0(I)=5,\quad \beta_1(I)=9,\quad \beta_2(I)=7 \mbox{ and }\beta_3(I)=2.$$
We can also read the graded and multigraded Betti numbers of $I$, considering the total degree and multidegree of each generator. We display the graded Betti numbers in the usual Betti diagram (see appendix \ref{apal}):
\begin{center}
\begin{tabular}{lrrrr}
&$0$&$1$&$2$&$3$\\ \hline
$23:$&$1$&$-$&$-$&$-$\\
$36:$&$1$&$-$&$-$&$-$\\
$46:$&$1$&$-$&$-$&$-$\\
$47:$&$1$&$-$&$-$&$-$\\
$48:$&$1$&$1$&$-$&$-$\\
$49:$&$-$&$2$&$-$&$-$\\
$50:$&$-$&$1$&$1$&$-$\\
$52:$&$-$&$1$&$-$&$-$\\
$55:$&$-$&$2$&$-$&$-$\\
$57:$&$-$&$-$&$2$&$-$\\
$58:$&$-$&$1$&$-$&$-$\\
$59:$&$-$&$1$&$1$&$-$\\
$60:$&$-$&$-$&$1$&$-$\\
$61:$&$-$&$-$&$2$&$1$\\
$62:$&$-$&$-$&$-$&$1$\\ \hline
Tot:&$5$&$9$&$7$&$2$\\
\end{tabular}
\end{center}
The multidegrees of the nonzero multigraded Betti numbers are:
\begin{center}
\begin{tabular}{|r|l|}
\hline
$i$&$\ab\in\NN^n\mbox{ such that } H_{i,\ab}(\KK(I))\neq 0$\\ \hline
$0$&$(19,13,5,3,6),(20,14,4,4,5),(12,0,1,7,3),(2,14,17,0,3),(13,11,12,4,8)$\\
$1$&$(12,14,17,7,3),(19,13,5,7,6),(20,14,4,7,5),(13,11,12,7,8),(19,14,17,3,6)$\\
&$(20,14,17,4,5),(13,14,17,4,8),(20,14,5,4,6),(19,13,12,4,8)$\\
$2$&$(19,14,17,7,6),(20,14,17,7,5),(13,14,17,7,8),(20,14,5,7,6),(19,13,12,7,8)$\\
&$(20,14,17,4,6),(19,14,17,4,8)$\\
$3$&$(20,14,17,7,6),(19,14,17,7,8)$\\
\hline
\end{tabular}
\end{center}
The \emph{projective dimension} of $I$ is then $3$ and using the graded Auslander-Buchsbaum formula $n=depth(I)+projdim(I)$ we have that the depth is $2$. Finally, the Castelnuovo-Mumford regularity of $I$ is $62$. The size of the minimal resolution is $5+9+7+2=23$.

Observe that we can also write the multigraded Hilbert series of $I$ out of the multidegrees of these generators, since this is the $\Kc$-polynomial associated to the minimal free resolution of $I$, $\sum_i(-1)^i \sum_{\ab}\beta_{i,\ab}x^\ab$. Other invariants such as the Krull dimension or multiplicity can also be read from the Koszul generators, as will be seen in section \ref{structure-combinatorial}.
\end{Example}
Thus, Koszul homology gives us an alternative way to obtain graded and multigraded Betti numbers, depth, dimension, Castelnuovo-Mumford regularity, and other invariants, the usual way being computing a minimal resolution. Since the actual computation of a full minimal resolution is in general a complicated task, if we are able to avoid it and produce algorithms to compute the relevant data from the Koszul homology, the latter would be quite useful. This issue is discussed in chapter \ref{computation}.

The rest of this section will be devoted to the \emph{constructive} transition between the (effective) Koszul homology of $I$ to its minimal resolution, \emph{constructive} can be understood as meaning as explicit as possible. 

\subsection{The $Tor$ bicomplex}
\spanishsubsection{El bicomplejo $Tor$}
The natural equivalence between the minimal resolution of an $R$-module $M$ and its Koszul homology is described by the two equivalent ways to compute $Tor^R_\bullet(\kb,M)$, as we have already seen. This equivalence is the basis of most of the computations and considerations we make. A good way to describe the situation is given by the $Tor$ bicomplex (recall the definition of bicomplex from Appendix \ref{apal}). In this case we build a bicomplex form the Koszul complex $(\KK(R),\partial)$ and the minimal resolution $\FF:(\Fc,d)$ of $I$, see figure \ref{Tor-bicomplex}.

\begin{figure}
\begin{center}
\begin{tikzpicture}[scale=1.1]

\draw (4,0) node[right]{$R\otimes\wedge^{p-1} \Vc\otimes \Fc_{q-1}$};
\draw (9,0) node[right]{$R\otimes\wedge^p \Vc\otimes \Fc_{q-1}$};
\draw (4,2) node[right]{$R\otimes\wedge^{p-1} \Vc\otimes \Fc_q$};
\draw (9,2) node[right]{$R\otimes\wedge^p \Vc\otimes \Fc_q$};
\draw[->] (4,0)--(3,0)node[above,pos=0.5]{$\partial$}node[left]{$\cdots$};
\draw[->] (9,0)--(7.3,0)node[above,pos=0.5]{$\partial$};
\draw[->] (9,2)--(7,2)node[above,pos=0.5]{$\partial$};
\draw[->] (4,2)--(3,2)node[above,pos=0.5]{$\partial$}node[left]{$\cdots$};
\draw[->] (13,0)node[right]{$\cdots$}--(12,0)node[above,pos=0.5]{$\partial$};
\draw[->] (13,2)node[right]{$\cdots$}--(12,2)node[above,pos=0.5]{$\partial$};
\draw[->] (5.5,-0.2)--(5.5,-1.2)node[right,pos=0.5]{$d$}node[below]{$\vdots$};
\draw[->] (10.5,-0.2)--(10.5,-1.2)node[right,pos=0.5]{$d$}node[below]{$\vdots$};
\draw[->] (5.5,1.8)--(5.5,0.3)node[right,pos=0.5]{$d$};
\draw[->] (10.5,1.8)--(10.5,0.3)node[right,pos=0.5]{$d$};
\draw[->] (5.5,3.3)node[above]{$\vdots$}--(5.5,2.3)node[right,pos=0.5]{$d$};
\draw[->] (10.5,3.3)node[above]{$\vdots$}--(10.5,2.3)node[right,pos=0.5]{$d$};
\end{tikzpicture}\caption{The $Tor^R(I,\kb)$ bicomplex.}\label{Tor-bicomplex}
\end{center}
\end{figure}

The vertical arrows are exact in every column $p$ except at positions $(p,0)$, in which we obtain a homology group $H^d_{(p,0)}$, which, since $\FF$ is a resolution of $I$, is isomorphic to $R\otimes_\kb\wedge^p\Vc\otimes_R I$. The horizontal differential $\partial$ induces here a differential $\partial'$, which makes $(H^d_{(p,0)},\partial'_p)$ a chain complex. Simmetrically, using rows instead of columns, we have a second chain complex $(H^\partial_{(0,q)},d'_q)$. The homology groups of both chain complexes are isomorphic, which can be proved using the spectral sequences of the bicomplex (see for instance Appendix \ref{apal} and \cite{E95}). This is, in a more general setting, the proof of the symmetry of $Tor$ modules, i.e. $Tor^R(M,N)\simeq Tor^R(N,M)$. 

If we have an explicit generating set of the Koszul homology modules, then the data we have available about this bicomplex are the following:
\begin{itemize}
\item[-] The Koszul complex, i.e. the modules $R\otimes_\kb\wedge^p \Vc$ and the differential $\partial$.
\item[-] The modules $\Fc_q$ and their decomposition in multigraded components. We know these because having the explicit generating set of the Koszul homology, we have the multigraded Betti numbers, which describe  the $\Fc_i=\oplus_{\mu\in\NN^n}\Fc_{i,\mu}$.
\item[-] Explicit generating sets for the homology of $\KK(I)$, i.e. we have an explicit description of $H^d_{(p,0)}$ for all $p$.
\end{itemize}

What we lack is the differential $d$ of the minimal resolution. In this case, the $Tor$ bicomplex, although useful to describe the situation, is not of enough use to actually compute the lacking differential, because the induced differential in the vertical spectral sequence, i.e. the differential induced in $(H^\partial_{(0,q)},d'_q)$, is the zero differential, and thus most of the information is lost. A solution was proposed in \cite{AH95} and \cite{AH96}. The first of these two papers treats the case $I$ has a \emph{pure} resolution, and the second treats the general case, giving explicit resolutions in some cases, in particular for stable ideals. In this paper we can see that even if the full Koszul homology is explicitly known, it is not easy to make explicit the needed isomorphisms and the identification of the differentials in the resolution. However, this construction allows to generalize the resolution given by Eliahou and Kervaire \cite{EK90}. We summarize here the construction in \cite{AH96}.
First of all, given an $R$-module $\Mc$, from the $Tor^R(\Mc,\kb)$ bicomplex, we know that the minimal free resolution of $\Mc$ can be written 

$$\cdots\stackrel{d_i}{\longrightarrow}R\otimes_\kb H_{i-1}(\KK(\Mc))\stackrel{d_{i-1}}{\longrightarrow}\cdots\rightarrow R\otimes_\kb H_1(\KK(\Mc))\stackrel{d_1}{\longrightarrow}R\otimes_\kb H_0(\KK(\Mc))\stackrel{d_0}{\longrightarrow} \Mc\rightarrow 0$$

In order to construct the differentials, we can proceed inductively, assuming the $d_j$ are known for $j<i$, and starting with $d_0([g^0_k])=g^0_k$ for all $k$, here the $g^i_k$ are the generators of $H_i(\KK(\Mc))$. Next, consider the bicomplex $\KK\otimes_R \Fc_{<i}$, where the $\Fc_j$ are the modules in the minimal resolution described above, i.e. $\Fc_j=R\otimes_\kb H_i(\KK(\Mc))$, it is just the part of the $Tor$ bicomplex, in which all the differentials are already known; see figure \ref{Tor-bicomplex}.

Take now a generator $g^i_k$ of $H_i(\KK(\Mc))$. Since $\Fc_{<i}$ is exact, we have that there exists an element $f_0\in\KK_i\otimes_\kb\Fc_0$ such that $g^i_k=(id\otimes d)(f_0)$. Moreover, there exist elements $f_j\in \KK_{i-j}\otimes \Fc_j$, for $j=1,\dots,i-1$ such that $(\partial_{i-j}\otimes id)(f_j)=(id\otimes d_{j+1})(f_{j+1})$ for all $j\in\{0,\dots,i-2\}$. Then, we define $(id \otimes d_i)(1\otimes[g^i_k])=(\partial_1\otimes id)(f_{i-1})$. See the diagram in figure \ref{AH2-bicomplex}:

\begin{figure}
\begin{center}
\begin{tikzpicture}[black!40,scale=1.1]
\draw (4,0) node[text=black,right]{$R\otimes\wedge^{i-1} \Vc\otimes \Fc_0$};
\draw (9,0) node[text= black,right]{$\mathbf{f_0}\in R\otimes\wedge^{i} \Vc\otimes \Fc_{0}$};
\draw (4,2) node[text=black,right]{$\mathbf{f_1}\in R\otimes\wedge^{i-1} \Vc\otimes \Fc_{1}$};
\draw (9,2) node[right]{$R\otimes\wedge^i \Vc\otimes \Fc_1$};
\draw (0,2) node[text=black,right]{$R\otimes\wedge^{i-2} \Vc\otimes \Fc_1$};
\draw (0,4) node[text=black,right]{$\mathbf{f_2}\in R\otimes\wedge^{i-2} \Vc\otimes \Fc_{2}$};

\draw[->,black] (0,4)--(-1,4)node[above,pos=0.5]{$\partial_{i-2}$}node[left]{$\cdots$};
\draw[->] (0,2)--(-1,2)node[above,pos=0.5]{$\partial_{i-2}$}node[left]{$\cdots$};
\draw[->] (4,0)--(3,0)node[above,pos=0.5]{$\partial_{i-1}$}node[left]{$\cdots$};
\draw[->,black] (9,0)--(7.3,0)node[above,pos=0.5]{$\partial_i$};
\draw[->] (9,2)--(7.6,2)node[above,pos=0.5]{$\partial_i$};
\draw[->,black] (4,2)--(3,2)node[above,pos=0.5]{$\partial_{i-1}$};
\draw[->] (13.3,0)node[right]{$\cdots$}--(12.3,0)node[above,pos=0.5]{$\partial_{i+1}$};
\draw[->] (13.3,2)node[right]{$\cdots$}--(11.7,2)node[above,pos=0.5]{$\partial_{i+1}$};

\draw[->] (5.5,-0.2)--(5.5,-1.2)node[right,pos=0.5]{$d_0$}node[below]{$\vdots$};
\draw[->,black] (10.5,-0.2)--(10.5,-1.2)node[right,pos=0.5]{$d_0$}node[text=black,below]{$g^i_k\in \KK_i(\Mc)$};
\draw[->] (1.5,1.8)--(1.5,0.3)node[right,pos=0.5]{$d_1$};
\draw[->,black] (5.5,1.8)--(5.5,0.3)node[right,pos=0.5]{$d_1$};
\draw[->] (10.5,1.8)--(10.5,0.3)node[right,pos=0.5]{$d_1$};
\draw[->,black] (1.5,3.3)--(1.5,2.3)node[right,pos=0.5]{$d_2$};
\draw[->] (5.5,3.3)node[above]{$\vdots$}--(5.5,2.3)node[right,pos=0.5]{$d_2$};
\draw[->] (10.5,3.3)node[above]{$\vdots$}--(10.5,2.3)node[right,pos=0.5]{$d_2$};
\draw[->] (1.5,5.3)node[above]{$\vdots$}--(1.5,4.3)node[right,pos=0.5]{$d_3$};
\end{tikzpicture}\caption{Construction of $d_i([g^i_k])$.}\label{AH2-bicomplex}
\end{center}
\end{figure}

Then, to construct the minimal free resolution of $I$ once we know a set of cycles $g^i_k$ the homology classes of which generate $H_i(\KK(I))$ for all $i$, the procedure proposed in \cite{AH96} to determine the differentials $d_i$ works inductively by finding for each $g^i_k$ elements $f_0,\dots,f_{i-1}$ in $\KK\otimes G$ satisfying the conditions exposed above.

\begin{Remark}
Note that this procedure relies on the existence of the elements $f_j$, but in order to explicitly find them, ad-hoc methods should be developed for different situations. One such method is given in \cite{AH96} for stable monomial ideals and is also used for more general notions of stability.
\end{Remark}

We reproduce here the proof of this construction given in \cite{AH96}, which relies in the natural isomorphism $\phi_i:H_i(I)\rightarrow Tor(I,\kb)$. This is defined as follows: If $(\Fc,d)$ is the minimal free resolution of $I$, then given an element $[z]$ of $H_i(\KK(I))$ there is some $a_0\in \KK_i\otimes\Fc_0$ with $d_0(a_0)=z$ and elements $a_j\in \KK_{i-j}(\Fc_j),\, j=1,\dots,i-1$ with
$$\partial_{i-j}(a_j)=(-1)^{j+1}d_{j+1}(a_{j+1})$$ 
then we have $\phi_i([z])=\bar{a}_i$ where $\bar{a}_i$ is the image of $a_i\in\KK_0\otimes\Fc_i=\Fc_i$ in $\Fc_i/\mm\Fc_i$. We have that $(a_0,\dots,a_i)\in\bigoplus_{j=0}^i\KK_{i-j}(\Fc_j)$ is a cycle in $\KK\otimes\Fc$; it is called a lifting of $z$. Identifying $\Fc_{<i}$ and $G_{<i}$ then, for some $k$ we may take $a_{jk}=(-1)^{j(j+1)/2}f_j$ for $j=0,\dots,i-1$ to obtain a lifting of $g^i_k$. Since $\phi([g^i_k])=\bar{a}_{ik}$, the $a_{ik}$ form a basis of $\Fc_i$. Then, the elements $\partial_1(a_{(i-1)k}$ generate $Im(d_i)=Ker(d_{i-1})$. Hence, the elements $(id\otimes d_i)(1\otimes[g^i_k])=\partial_1g_{(i-1)k}=(-1)^{i(i-1)/2}\partial_1(a_{(i-1)k})$ generate $Ker(d_{i-1})$. This proofs that $Im(d_i)=Ker(d_{i-1})$ which suffices to see that $G_{<i+1}$ is exact.

\subsection{The Aramova-Herzog bicomplex}
\spanishsubsection{El bicomplejo de Aramova-Herzog}
To present the second approach to the problem, we start by describing the Aramova-Herzog bicomplex introduced in \cite{AH95}. In this paper, A. Aramova and J. Herzog show that the minimal free resolution of a finitely generated $R$-module $M$ ``can be partially recovered from the cycles generating the Koszul homology'' (\cite{AH95}, p.1). The main result of the paper is an isomorphism between two spectral sequences that arise in the study of a certain bicomplex. This spectral sequence isomorphism is used to prove that the Koszul cycles give the \emph{pure part} of the minimal resolution of $\Mc$. The bicomplex on which they base, will be denoted here as \emph{Aramova-Herzog bicomplex}, $AH(\Mc)$. Based on this bicomplex, an approach from the \emph{effective homology} point of view can give us the explicit construction of an effective resolution of $I$ form its \emph{effective} Koszul homology. This was presented in \cite{RS06,RSS06}. In this subsection some notions from effective homology are very important, in particular those of \emph{reduction}, \emph{equivalence} and the \emph{basic perturbation lemma}; a description of the basics of effective homology can be seen in appendix \ref{apef} if the reader is unfamiliar with the subject.

\begin{figure}
\begin{center}
\begin{tikzpicture}[scale=0.65]
\draw (4,0) node[right]{$\Mc\otimes\wedge^{0} \Vc\otimes R_0$};
\draw (4,2) node[right]{$\Mc\otimes\wedge^{1} \Vc\otimes R_0$};
\draw (4,4) node[right]{$\Mc\otimes\wedge^{2} \Vc\otimes R_0$};
\draw (4,6) node[right]{$\Mc\otimes\wedge^{3} \Vc\otimes R_0$};

\draw (10,2) node[right]{$\Mc\otimes\wedge^{0} \Vc\otimes R_1$};
\draw (10,4) node[right]{$\Mc\otimes\wedge^{1} \Vc\otimes R_1$};
\draw (10,6) node[right]{$\Mc\otimes\wedge^{2} \Vc\otimes R_1$};

\draw (16,4) node[right]{$\Mc\otimes\wedge^{0} \Vc\otimes R_2$};
\draw (16,6) node[right]{$\Mc\otimes\wedge^{1} \Vc\otimes R_2$};

\draw (22,6) node[right]{$\Mc\otimes\wedge^{0} \Vc\otimes R_3$};

\draw[->] (5.5,-0.3)--(5.5,-1.2)node[below]{$0$};
\draw[->] (5.5,1.7)--(5.5,0.3)node[right,pos=0.5]{$\partial ''$};
\draw[->] (5.5,3.7)--(5.5,2.3)node[right,pos=0.5]{$\partial ''$};
\draw[->] (5.5,5.7)--(5.5,4.3)node[right,pos=0.5]{$\partial ''$};
\draw[->] (5.5,7.7)node[above]{$\vdots$}--(5.5,6.3)node[right,pos=0.5]{$\partial ''$};

\draw[->] (11.5,1.7)--(11.5,0.3)node[below]{$0$};
\draw[->] (11.5,3.7)--(11.5,2.3)node[right,pos=0.5]{$\partial ''$};
\draw[->] (11.5,5.7)--(11.5,4.3)node[right,pos=0.5]{$\partial ''$};
\draw[->] (11.5,7.7)node[above]{$\vdots$}--(11.5,6.3)node[right,pos=0.5]{$\partial ''$};

\draw[->] (17.5,3.7)--(17.5,2.3)node[below]{$0$};
\draw[->] (17.5,5.7)--(17.5,4.3)node[right,pos=0.5]{$\partial ''$};
\draw[->] (17.5,7.7)node[above]{$\vdots$}--(17.5,6.3)node[right,pos=0.5]{$\partial ''$};

\draw[->] (23.5,5.7)--(23.5,4.3)node[below]{$0$};
\draw[->] (23.5,7.7)node[above]{$\vdots$}--(23.5,6.3)node[right,pos=0.5]{$\partial ''$};

\draw[->] (3,0)node[left]{$0$}--(4,0);
\draw[->] (3,2)node[left]{$0$}--(4,2);
\draw[->] (3,4)node[left]{$0$}--(4,4);
\draw[->] (3,6)node[left]{$0$}--(4,6);

\draw[->] (8.5,0)--(11,0);
\draw[->] (8.5,2)--(10,2)node[above,pos=0.5]{$\partial '$};
\draw[->] (8.5,4)--(10,4)node[above,pos=0.5]{$\partial '$};
\draw[->] (8.5,6)--(10,6)node[above,pos=0.5]{$\partial '$};

\draw[->] (14.5,2)--(17,2);
\draw[->] (14.5,4)--(16,4)node[above,pos=0.5]{$\partial '$};
\draw[->] (14.5,6)--(16,6)node[above,pos=0.5]{$\partial '$};

\draw[->] (20.5,4)--(23,4);
\draw[->] (20.5,6)--(22,6)node[above,pos=0.5]{$\partial '$};

\draw[->] (26.5,6)--(27.5,6)node[right]{$0$};
\end{tikzpicture}\caption{The Aramova-Herzog bicomplex.}\label{AH-bicomplex}
\end{center}
\end{figure}

Let $\Mc$ be an $R$-module, consider the bicomplex $AH(\Mc)$ in figure \ref{AH-bicomplex}, in which the $R_p$ denote the homogeneous components of degree $p$ of $R$ and $\Vc$ is a $\kb$-vector space with a basis in correspondence with the variables $x_1,\dots,x_n$ of $R=\kb[x_1\dots,x_n]$. Each $\wedge^q\Vc$ is a homogeneous component of $\wedge \Vc$. Note that this bicomplex can be seen as $\KK(\Mc)\otimes R$, being the vertical differential essentially the Koszul differential: $\partial''=\partial_{\KK(M)}\otimes id_R$. On the other hand, $AH(\Mc)$ can be seen as $\Mc\otimes\KK(R)$, being the horizontal differential $\partial'=id_{\Mc}\otimes\partial_{\KK(R)}$.

Since $\KK(R)$ is acyclic, we can construct a reduction $AH(\Mc)\rrdc\Mc$ in the following way: Every row of $AH(M)$ is acyclic, except the bottom one $0\rightarrow\Mc\otimes\wedge^0\otimes R_0\rightarrow 0$. Applying the functor $\Mc\otimes -$ gives a \emph{horizontal} reduction $(AH(\Mc),\partial')\rrdc\Mc$, where the vertical differential has been removed. Using the basic perturbation lemma, we reinstall the vertical differential $\partial''$. In order to apply the basic perturbation lemma, the so called nilpotency condition (see appendix \ref{apef}) must be satisfied; it is indeed the case here, as can be seen in \cite{RS06}.

Similarly, we perform a \emph{vertical} reduction: We begin with an equivalence $\KK(I)\lrrdc \Hc$, where $\Hc$ is a chain complex describing the homology of $\KK(\Mc)$. Note that this is our starting point, i.e. the knowledge of the \emph{effective} homology of $\Mc$ (in our case $\Mc=I$ is a monomial ideal, but the construction admits more generality). This equivalence can be applied to every column of $AH(\Mc)$ and produces a new equivalence $(AH(\Mc),\partial'')\lrrdc \Hc\otimes R$ where the horizontal differential is set to zero. Again, the application of the basic perturbation lemma, allows us to reinstall the horizontal differential in $AH(\Mc)$. The resulting equivalence inserts a differential in $\Hc$, which makes it a chain complex $(\Hc,d)$; this is the resolution we want. The fact that this resolution is minimal is given by the minimality of the ranks of the modules in the resolution. It might seem that these considerations leading to the construction of the minimal resolution are far from being constructive, but one must note that the application of effective homology, and in particular of the basic perturbation lemma, is essentially algorithmic. In particular, an implementation of these techniques is available in the computer algebra program Kenzo \cite{Kenzo}. In our case, the actual description obtained for the differential in the minimal free resolution, might be complicated but it is automatically produced by the formulas of the basic perturbation lemma, therefore algorithms to compute it are available.

\section{Combinatorial decompositions}\label{structure-combinatorial}
\spanishsection{Descomposiciones combinatorias}
Given an ideal $I$ in the polynomial ring $R$, the factor algebra corresponding to $I$ is a $\kb$-algebra $S=R/I$. If $I$ is a monomial ideal, then $S$ can be seen as the complement in the set of all monomials of the monomials contained in $I$, i.e. the standard monomials of $I$. As an abuse of notation, and since we are dealing with monomial ideals, we will identify standard monomials with their equivalence classes in $S=R/I$, therefore, we will use the notation $x^\mu\in R/I$ to say that the equivalence class of $x^\mu$ in $R/I$ is different from zero. We call a combinatorial decomposition of $S$ to a representation of it as a finite sum of $\kb$-vector spaces of the form $\xb^\mu\cdot \kb[\xb_\mu]$, where $\xb_\mu$ is a subset of the variables defining $R$. If the sum is direct, this is a {\it Stanley decomposition} \cite{S78}: \begin{equation}\label{stanley}
S=\bigoplus_{\mu\in \Fc}\xb^\mu\kb[\xb_\mu]
\end{equation}
where $\Fc$ is a finite subset of $\NN^n$. Note that Stanley decompositions need not to be unique.

The existence of such decompositions is guaranteed by the following result, see \cite{CLO96}:
\begin{Proposition}\label{stanley-monoid}
Let $I\subseteq \NN^n$ be a monoid ideal and $\bar{I}=\NN^n\smallsetminus I$ its complementary set. Then, there exists a finite set $\bar{N}\subseteq \bar{I}$ and for each $\nu\in\bar{N}$ a set of indices $N_\nu\subseteq\{1,\dots,n\}$ such that
$$\bar{I}=\bigcup_{\nu\in\bar{N}}(\nu+\NN^n_{N_\nu})$$
and $(\nu+\NN^n_{N_\nu})\cap (\mu+\NN^n_{N_\mu})=\emptyset$ for all $\mu,\nu\in\bar{N}$.
\end{Proposition}

This decomposition of the complementary set of a monoid ideal automatically guarantees the existence of Stanley decompositions of complementary sets of monomial ideals, since the set of monomials in a monomial ideal is a monoid ideal in $\NN^n$. Note that this applies also to general polynomial ideals, since given a degree compatible term order $<$ and a polynomial ideal $I$, the $k$-vector dimension of $R/I$, and $R/lt_{<}I$ is the same, therefore, a Stanley decomposition of the complement of a polynomial ideal is defined as a Stanley decomposition of its leading ideal. Stanley decompositions provide easy ways to obtain some important invariants of graded $\kb$-algebras, such as the Krull dimension or Hilbert series, in particular:
\begin{Proposition}[\cite{S78,SW91}]
Let $S=R/I$ have a Stanley decomposition as in \emph{(\ref{stanley})}, and let $d$ be the maximum of the numbers $\vert \xb_\mu\vert$, $\mu\in \Fc$. Then
\begin{enumerate}
\item $d$ is the Krull dimension of $S$
\item The Hilbert series of $S$ is given by
$$H(S;t)=\sum_{\mu\in\Fc}\frac{t^{\vert\mu\vert}}{(1-t)^{\vert x_\mu\vert}}$$
\end{enumerate}

\end{Proposition}
In this section we give some procedures to obtain Stanley decompositions of $S=R/I$ from the Koszul homology of $I$ in the case $I$ is a monomial ideal.

We need some previous definitions and notations that will be used later:
\begin{Definition}\label{boundary}
Let $x^\mu\in I$, and let $supp(x^\mu)=\{i_1,\dots,i_s\}$. We say that $x^\mu$ is {\rm on the boundary or wall} of $I$ if $x^{\mu'}=x^\mu/x_{i_1}\cdots x_{i_s}\notin I$. We denote $bound(I)$ to this boundary. If on the contrary $x^{\mu'}=x^\mu/x_{i_1}\cdots x_{i_s}\in I$ we say that $x^\mu$ is {\rm inside} the ideal.

Let $x^\mu\in bound(I)$, for each minimal generator $x^\nu$ of $I$ that divides $x^\mu$ consider the variables $0_{\mu,\nu}=\{x_i\vert i\in supp(x^\mu)\smallsetminus supp(x^\mu/x^\nu)\}$ that become zero when dividing $x^\mu$ by $x^\nu$, i.e. all $x_i$ such that $\mu_i=\nu_i$. We say that $x^\mu$ is on a $k$-wall of $I$ if $\vert \bigcup_{x^\nu\vert x^\mu}0_{\mu,\nu}\vert=s-k$
\end{Definition}

\begin{Remark}\label{boundary-roune}
An alternative and elegant definition of the boundary of a monomial ideal is given in\cite{R07}. We reproduce it here:
$$bound(I):=\{m\in I\vert m \mbox{ is a monomial and }\pi(m)\notin I\}$$
where the function $\pi$ is defined by
$$\pi(m):=\frac{lcm(m,x_1\cdots x_n)}{x_1\cdots x_n}$$
The close relation of this definition and the irreducible decomposition of $I$ is showed in \cite{R07} and will be recalled later in section \ref{irreducible_decomposition}.
\end{Remark}

\begin{Remark}
Note that $x^\mu$ is in the boundary or wall of $I$ if and only if it is in some $k$-wall. Note also that the generators of $I$ are always in $0$-walls. The $0$-walls correspond to the \emph{corners} defined in \cite{B96}. Corners are defined there as the points $\mu$ in $\NN^n$ such that the appearance of the staircase diagram of $I$ near $\mu$ changes in each coordinate direction. It is also stated that there exists some $i$ such that $H_{i,\mu}(\KK(I))\neq 0$ only if $\mu$ is a corner (a $0$-wall in our language). The converse is not true, since there are corners ($0$-walls) with null homology in every dimension, take for example the monomial $xyzt$ in the ideal $xy,yt,zt$, which is a corner with no Koszul homology. However, every corner of dimension $3$ or less has some not-null Koszul homology group.
\end{Remark}
\begin{Definition}

Let $x^\mu\in bound(I)$, and let $supp(x^\mu)=\{i_1,\dots,i_s\}$. We say that $x^\mu$ is {\rm a closed corner} of $I$ if $x^\mu/x_{i_1}\cdots \hat x_{i_l}\cdots x_{i_s}\in I$ for all $l\in\{1\dots s\}$. If the support of $x^\mu$ is $\{1,\cdots,n\}$ (i.e. $x^\mu$ has {\rm full support}) and $x^\mu$ is a closed corner, we say that it is a {\rm maximal corner} of $I$.
\end{Definition}

\begin{Remark}
The number of maximal corners of a monomial ideal $I$ is calculated in \cite{A97} by computing the number of maximal standard monomials w.r.t divisibility, modulo $I$. It is shown there that if we call $c_n(p)$ to the number of maximal standard monomials module the ideal $I$ with $p$ minimal generators in $n$ variables, then $c_n(p)$ is not a polynomial in $p$.
\end{Remark}

\begin{Example}
Consider the monomial ideal $I=\langle x^3,x^2y,xz,y^3,z^3\rangle$, and the monomials $xy^2z, xyz^2$ and $x^3yz$, which are drawn in figure \ref{wall-picture}. We have that
\begin{itemize}
\item The only generator that divides $xy^2z$ is $xz$, and $0_{xy^2z,xz}=\{x,z\}$, thus, $xy^2z$ is in a $1$-wall of $I$
\item The only generator that divides $xyz^2$ is $xz$, and $0_{xyz^2,xz}=\{x\}$, thus, $xyz^2$ is in a $2$-wall of $I$
\item $x^3yz$ is divisible by the generators $x^3,x^2y$ and $xz$, and $0_{x^3yz,x^3}=\{x\}$, $0_{x^3yz,x^2y}=\{y\}$ and $0_{x^3yz,xz}=\{z\}$thus, $x^3yz$ is in a $0$-wall of $I$
\end{itemize}

Observe that all the corners drawn in the staircase diagram of $I$ form the boundary or wall of the ideal (strictly speaking, those parts of it that are not contained in any coordinate plane $x_i=0$). In the staircase diagram one can also see the closed and maximal corners. Here, the closed corners are $x^3z,x^3y,x^2y^3,xz^3$ and $y^3z^3$, and the maximal corners are $xy^3z^3,x^2y^3z$ and $x^3yz$.

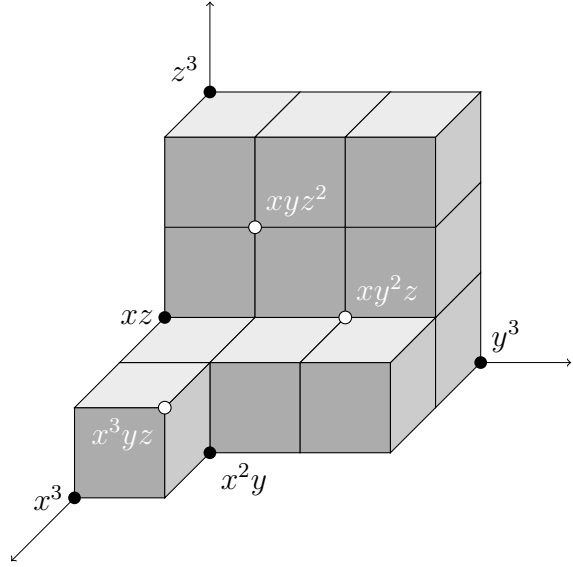
\begin{figure}
\begin{center}
\begin{tikzpicture}[scale=1.2]
\draw [very thin, ->](0,0)--(-0.7,-0.7);
\draw [very thin, ->](4.5,1.5)--+(1,0);
\draw [very thin, ->](1.5,4.5)--+(0,1);
\foreach \position in {(0,0),(1.5,0.5),(2.5,0.5),(1,2),(1,3),(2,2),(2,3),(3,2),(3,3)} 
\filldraw[fill=gray!65,draw=black]\position rectangle +(1,1);
\foreach \position in {(1,0),(3.5,0.5),(4,1),(4,2),(4,3)}
\filldraw[fill=gray!40!white,draw=black] \position-- ++(0.5,0.5)-- ++(0,1)-- ++(-0.5,-0.5)-- cycle;
\foreach \position in {(0,1),(0.5,1.5),(1.5,1.5),(2.5,1.5),(1,4),(2,4),(3,4)}
\filldraw[fill=gray!15!white,draw=black] \position-- ++(1,0)-- ++(0.5,0.5)-- ++(-1,0)-- cycle;
\fill (0,0) circle(2pt) node[left]{$x^3$};
\fill (1.5,0.5) circle(2pt) node[below right]{$x^2y$};
\fill (1,2) circle(2pt) node[left]{$xz$};
\fill (4.5,1.5) circle(2pt) node[above right]{$y^3$};
\fill (1.5,4.5) circle(2pt) node[above left]{$z^3$};

\fill[color=white, draw=black](3,2) circle(2pt) node[above right]{$xy^2z$};
\fill[color=white,draw=black](2,3) circle(2pt) node[above right]{$xyz^2$};
\fill[color=white,draw=black](1,1) circle(2pt) node[below left]{$x^3yz$};

\end{tikzpicture}
\caption{Some elements on the wall of $I=\langle x^3,x^2y,xz,y^3,z^3\rangle$}\label{wall-picture}
\end{center}
\end{figure}
\end{Example}

\begin{Proposition}\label{closed-maximal}
If $x^\mu$ is a closed corner of $I$ with $supp(x^\mu)=\{i_1,\dots,i_s\}$, then $H_{\mu,s-1}(\KK(I))\neq0$. In particular, $H_{\mu,n-1}(\KK(I))\neq0$ if and only if $x^\mu$ is a maximal corner of $I$.
\end{Proposition}

\noindent {\bf Proof:} For the first statement, consider the upper simplicial Koszul complex of $I$ at $\mu$. It is a simplicial complex on $s$ vertices, which has as facets all possible $s-1$ faces, but does not contain the (unique) $s$-face. Then, it is clear that $H_{\mu,s-1}(\KK(I))\simeq\tilde H_s(\Delta^I_\mu)\simeq \kb$.
For the second statement, the {\it if} implication is just a particular case of the first statement with $s=n$. The {\it only if} part comes form the easy observation that the only possibility for a simplicial Koszul complex to have homology in dimension $n-1$ is if all $n-1$ faces are present in the complex, and the $n$-face is not. Observe that this is not the case if $s<n$ $\square$

\begin{Remark}
The last lemma gives us an alternative definition of a maximal corner of $I$, namely those $x^\mu$ such that $H(\KK)_{\mu,n-1}\neq0$.

Note also that the closed corners are those multidegrees at which the lower Koszul simplicial complex is just the empty set: $\Delta^{\mu}_I=\{\emptyset\}$. The Koszul ideals of a maximal corner $\mu$ are given by $KI_\mu=\langle x_1\cdots x_n\rangle$ and $KI^\mu=\langle x_1,\dots,x_n\rangle$. In the case of a closed corner $\mu$ with support given by $x_{i_1},\dots,x_{i_k}$, $KI_\mu=\langle x_{i_1}\cdots x_{i_k}\rangle$ and $KI^\mu=\langle x_{i_1},\dots,x_{i_k}\rangle$.
\end{Remark}

\begin{Definition}
Let $I$ be a monomial ideal, and $\mu\in \NN^n$ such that $x^\mu\in I$. Recall that the lower Koszul simplicial complex of $I$ at $\mu$, $\Delta^\mu_I$ is given by the set of squarefree vectors $\tau$ such that $x^{\mu'+\tau}$ is not in $I$, where $\mu'$ is defined by subtracting $1$ from each nonzero coordinate of $a$ (see definition \ref{lower-Koszul}). We say that a set of variables $\{x_{j_1},\dots,x_{j_k}\}$ is a \emph{cone of locally free directions} of $x^\mu$ (with respect to $I$) if $\tau=\{j_1\cdots j_k\}\in \Delta^\mu_I$, i.e. if $x^{\mu'}\cdot x^\sigma\notin I$ for every subset $\sigma\subseteq \tau$, this cone will be denoted by $[x_{j_1},\dots,x_{j_k}]$. The set of cones of locally free directions of $\xb^\mu$ will be denoted $LF(\xb^\mu)$ and is given by the maximal locally free cones of $x^\mu$ with respect to inclusion (i.e. the facets of $\Delta^\mu_I$). If a variable $x_i$ is such that $\{i\}\notin\Delta^\mu_I$ then $x_i$ is a \emph{non-free direction} of $x^\mu$ (with respect to $I$).

We say that a set of variables $\{x_{j_1},\dots,x_{j_k}\}$ is a \emph{cone of true free directions}, denoted $[x_{j_1},\dots,x_{j_k}]$, of $x^\mu$ (with respect to $I$) if $x^{\mu'}\cdot x^\sigma\notin I$ for every monomial $x^\sigma\in \kb[x_{j_1},\dots,x_{j_k}]$. The set of cones of true free directions of $\xb^\mu$ will be denoted $TF(\xb^\mu)$. If there exists any integer $p>0$ such that $x^{\mu'}\cdot x_i^p\in I$ then $x_i$ is not a true free direction of $x^\mu$ (with respect to $I$) and it is not contained in any cone of free directions of $x^\mu$.

For a multidegree $\mu\in \NN^n$ such that $x^\mu\notin I$ we say that $\{x_{j_1},\dots,x_{j_k}\}$ is a \emph{cone of true free directions}, of $x^\mu$ (with respect to $I$) if $x^{\mu}\cdot x^\sigma\notin I$ for every monomial $x^\sigma\in \kb[x_{j_1},\dots,x_{j_k}]$. Note that if $x^\mu$ is in the boundary of $I$ then the cones of free directions of $x^\mu\in I$ and $x^{\mu'}\notin I$ coincide by definition.
\end{Definition}

Observe that even if the definition of a true free direction involves an infinite condition, it can be detected with a finite criterion. Let $x^\lambda$ be the least common multiple of all generators of $I$. If $x^{\mu'}\cdot x_i^p$ is not in $I$ for $\mu'_i+p=\lambda_i$ then $i$ is a true free direction for $x^{\mu}$.
 
\begin{Example}
Consider the ideal in the previous example, $I=\langle x^3,x^2y,xz,y^3,z^3\rangle$. We have that $\Delta^{xyz}_I=\{\{x\},\{y\},\{z\},\{x,y\},\{y,z\}\}$ thus,  $LF(xyz)=\{[x,y],[y,z]\}$. However, $TF(xyz)=\emptyset$ since $(xyz)'=1$ and $1\cdot x^3\in I$, $1\cdot y^3\in I$ and $1\cdot z^3\in I$.

Note that if we consider instead $I=\langle x^2y,xz,y^3,z^3\rangle$ then $LF(xyz)$ remains unchanged, but $TF(xyz)=\{[x]\}$ since $1\cdot x^p\notin I$ for all $p\geq 0$.
\end{Example}

With these definitions, we are now in a position to give procedures to obtain a Stanley decomposition of $R/I$ given the Koszul homology of $I$. We treat the artinian and non-artinian case in a different way:

\subsection{Artinian case}
\spanishsubsection{Caso artiniano}
\begin{Definition}
We say that a monomial ideal $I$ is {\rm artinian} if the ring $R/I$ is zero-dimensional. Artinian monomial ideals can be characterized as those having among their minimal generators, elements of the form $x_i^{a_i}$ with $a_i>0$ for all $i\in\{1\dots n\}$.
\end{Definition}

\begin{Lemma}\label{artinian}
Let $I$ be any proper artinian monomial ideal, then $H_{n-1}(\KK(I))\neq 0$
\end{Lemma}

\noindent {\bf Proof:} An easy proof of this lemma can be given using quasi stable ideals and Pommaret bases, see section \ref{quasi-stable} and \cite{S02b} for the necessary definitions and properties). Artinian ideals are in particular quasi-stable ideals, which implies that their minimal generating set can be completed to a Pommaret basis. In particular, $x_1^{a_1}$ is an element of the Pommaret basis of the ideal $I$, and this implies (see \cite{S02b}) that $depth(I)=1$. Using now the Auslander-Buchsbaum formula, we have that $projdim(I)=n-depth(I)=n-1$. Since $dim(H_{n-1}(\KK(I)))$ equals the rank of the $(n-1)$-st module in any minimal resolution of $I$, we know that it is greater than zero $\square$

Artinian monomial ideals have finite complementary sets, thus their Krull dimension is zero. In this case it is easy to derive a Stanley decomposition of $R/I$ from the knowledge of $H_*(\KK(I))$:
\begin{Proposition}\label{Stanley-artinian}
Let $I$ be an artinian monomial ideal.  Let $\Bc_{n-1}(I)$ be the set of multidegrees $\nu$ such that $H_{\nu,n-1}(\KK(I))\neq 0$. A Stanley decomposition of $R/I$ is given by
$$R/I\simeq \bigoplus_{\substack{\text{$x^\mu$ divides $x^{\bar{\nu}}$} \\\text{for some $\nu \in \Bc_{n-1}(I)$}}} \kb\cdot x^\mu$$
where $\bar{\nu}_i=\nu_i-1$ if $i\in supp(x^\nu)$ and $\bar{\nu}_i=\nu_i$ otherwise.
\end{Proposition}

\noindent {\bf Proof: } Take a multidegree $\mu$ such that $x^\mu\in R/I$, we have to see that there exists at least one maximal corner $x^\nu$ such that $x^\mu\vert x^\nu$: choose $\rho$ such that $x^\mu\cdot x^\rho$ is on the wall of the ideal, but $x^{\mu+\rho}\cdot x_i$ is inside the ideal for all $i$ note that such a $\rho$ always exists because the ideal is zero dimensional and thus the number of multidegrees in the wall of $R/I$ is obviously finite. Clearly $x^{\mu+\rho}$ is a maximal corner of $I$, use now Lemma \ref{closed-maximal} and we are done $\square$

\begin{Remark}
An algorithmical enumeration of the elements in the index of the sum in this proposition would proceed by ascending degree, starting from degree zero and checking whether each monomial considered is in the ideal or not, until we reach the maximum of the degrees of the maximal corners of $I$.
\end{Remark}

\begin{Remark}
Observe that in proposition \ref{Stanley-artinian} we only use the $(n-1)$-st Koszul homology of $I$. This homology is very easy to compute at each multidegree. In fact, one needs only to check whether the corresponding simplicial Koszul complex has the unique $n$ face; if so, then we know there's no $(n-1)$-st homology at this multidegree. If this maximal face is not present, we check whether any on the $n$ $(n-1)$-faces fails to be in the complex. If any of them fails, we have no $(n-1)$-st homology. Only if all of them are present we have nonzero homology at this multidegree. The vertices to be checked are those in $L_{I,n-1}$, the number of which is bounded above by ${r}\choose{n-1}$, being $r$ the number of minimal generators of $I$.
\end{Remark}

\subsection{Non-artinian case}
\spanishsubsection{Caso no artiniano}

If our ideal $I$ is not artinian, one procedure to compute a Stanley decomposition of $I$ consists in associating to it an artinian ideal $\hat I$, compute the $(n-1)-st$ Koszul homology of $\hat I$, and then use this result to compute a Stanley decomposition of $I$.

We introduce the following ideal: Let $\xb^\lambda$ be the least common multiple of all the generators of $I$. Let $\lambda+1=(\lambda_1+1,\dots,\lambda_n+1)$ and $\mm^{\lambda+1}=\langle x_1^{\lambda_1+1},\dots,x_n^{\lambda_n+1}\rangle$. We call $\hat{I}=I+\mm^{\lambda+1}$. Now we have that $\hat{I}$ is artinian. Lets call $\hat I$ the {\it artinian closure of $I$}, since it is the largest artinian monomial ideal such that all the minimal generators of $I$ are minimal generators of $\hat I$. Compute now the complementary decomposition of $\hat{I}$. It is clear that every $x^\mu\in R/\hat I$ is also in $R/I$ thus, we call the decomposition of $R/\hat I$ the {\it inner part of the decomposition of $I$}.

Take all the maximal corners $x^\nu$ of $\hat{I}$, such that there is some $i$ for which $\nu_i=\lambda_i+1$. Consider here the lower Koszul simplicial complex (see section \ref{complexes} above), lets call it $\Delta^I_{\nu}$. 

Consider now as {\it points of the skeleton} the following set of monomials: every monomial $x^\tau\cdot x^\nu/supp(x^{\nu})$ such that $\tau\in\Delta^I_\nu$ and the {\it backward cone} in their nonfree directions, i.e. all the divisors of $x^\tau\cdot x^\nu/supp(x^\nu)$ in its nonfree directions $\{1,\dots,x_n\}\smallsetminus \tau$. Coming from these points, we add to the Stanley decomposition of $R/I$ the cones in the free directions of $x^\nu$ i.e. all multiples of them in the variables given by the corresponding $\tau$, $x^\tau\cdot x^\nu/supp(x^\nu)\cdot\kb[\xb_\tau]$.

Observe that if a point is a {\it nonfree divisor} of several points in the skeleton, then all these point have the same free directions. To prove this, assume $x^\rho$ divides both $x^\mu$ and $x^\nu$ and let $M$ be the set of free directions of $x^\mu$ and $N$ the set of free directions of $x^\nu$. We know that $M\cap supp(\mu-\rho)=\emptyset$ and $N\cap supp(\nu-\rho)=\emptyset$. Assume now that there exists $i\in M-N$, then we must have $\rho_i=\mu_i=\lambda_i+1$ and $\nu_i<\lambda_i+1$ and we have a contradiction, because in that case $x^\rho$ would not be a divisor of $x^\nu$.

Now, let $x^\mu\in R/I$, then we can have two cases, first, if $x^\mu$ divides $x^{\lambda+1}$ then it is in the inner part of $\hat{I}$, and from the first considerations above we know that we have collected all the inner part of the Stanley decomposition of $R/I$. Second, if $x^\mu$ does not divide $x^{\lambda+1}$ then exists $i$ such that $\mu_i>\lambda_i+1$, for each such $i$ do $\bar{\mu}=\mu-(\mu_i-\lambda_{i}+1))_i$, and then we are back in the first case at a point in which all $i$ are free directions, and we know it is in some of the considered cones.

All these considerations prove the following proposition:
\begin{Proposition}\label{Stanley-nonartinian}
Let $I$ be a monomial ideal, and let  $\hat I$ be its artinian closure. Then, provided we know $\Bc_{n-1}(\hat I)$ there is a procedure to construct a Stanley decomposition of $S=R/I$.
\end{Proposition}

\begin{Example}
Consider the ideal $I= \langle x^2y,xz,y^3,z^3\rangle$; $\lambda+1=(3,4,4)$, its integral closure is $\hat I=I+\langle x^3,y^4,z^4\rangle=\langle x^3, x^2y,xz,y^3,z^3\rangle$. 

We start with the combinatorial decomposition of $R/\hat I$. For this, we look at the multidegrees of the generators of $H_2(\KK(\hat I))$. These are $(3,1,1),(2,3,1)$ and $(1,3,3)$. Using them as in the artinian case, we obtain the \emph{inner part} of the decomposition of $R/I$. 

Now we select $(3,1,1)$, since it is the only maximal corner with some of the exponents of $\lambda+1$. The lower simplicial Koszul complex at $(3,1,1)$ is $\Delta^I_{(3,1,1)}=\{\emptyset,\{x\}\}$. The point in the skeleton are then $x^2$ and $x^3$. The element $x^{2}$ and its divisors $1,x$ have already been obtained when computing the inner part of the decomposition of $R/I$ thus, we only need to add $x^3\cdot \kb[x]$.

Therefore we have the following Stanley decomposition
$$R/I=1\oplus x\oplus x^2\oplus xy\oplus xy^2\oplus y\oplus y^2\oplus yz \oplus yz^2\oplus y^2z\oplus y^2z^2\oplus z\oplus z^2\oplus x^3\cdot\kb[x]$$

from which the Krull dimension of $R/I$ is $1$, and the Hilbert function is given by
$$H(R/I,t)=1+3t+5t^2+3t^3+t^4+\frac{t^3}{(1-t)}$$
\end{Example}

\section{Irreducible decomposition}\label{irreducible_decomposition}
\spanishsection{Descomposici\'on irreducible}
An irreducible ideal of a ring is an ideal such that is not the intersection of two ideals which properly contain it. We have seen in Proposition \ref{irreducible} a characterization of irreducible monomial ideals as those monomial ideals generated by powers of variables. An \emph{irreducible decomposition} of a monomial ideal $I$ is an expression of it as intersection of irreducible monomial ideals:
$$I=\mathfrak m^{\ab_1}\cap\cdots\cap\mathfrak m^{\ab_s}$$
where $\mm^{\ab}=\langle x_i^{a_i}\vert a_i>0\rangle$, being $\ab=(a_1,\dots,a_n)\in\NN^n$.
Every monomial ideal has an irreducible decomposition. Such decompositions are called \emph{irredundant} if no intersectands can be omitted. The irredundant irreducible decomposition of a monomial ideal is unique \cite{V01}. We say that $\mm^\ab$ is the irreducible ideal associated with the multiindex or multidegree $\ab$.

In this section we address the problem of obtaining an irredundant irreducible decomposition of a monomial ideal provided we know its Koszul homology. First, we give a procedure to obtain a irreducible decomposition of an artinian monomial ideal $I$ from its Koszul homology. Second, we obtain an irredundant decomposition of $I$ from the Koszul homology of its artinian closure $\hat I$. Finally, we give an algorithm to read the irreducible decomposition of $I$ from the Koszul homology of $I$ itself. 

Before the description of the different procedures relating Koszul homology and irreducible decompositions, we give a collection of auxiliary lemmas. These lemmas relate the monomials dividing a certain multidegree with the irreducible ideal associated to it. They will be used along the section.

\begin{Lemma}
Let $x^\nu\in \NN^n$ and $C=\{x_{i_1},\dots,x_{i_k}\}\subseteq \{1\dots,n\}$ such that $C\cap supp(x^\nu)=\emptyset$ i.e. $C$ is the complement of the support of $x^\nu$ in the set of all variables. Then
$$x^\mu\notin \mm^\nu\iff x^\mu\in \bigcup_{x^\rho\vert x^{\bar \nu}}x^\rho\cdot\kb[x_{i_1},\dots,x_{i_k}]$$
where $\bar{\nu}_i=\nu_i-1$ if $i\in supp(x^\nu)$ and $\bar{\nu}_i=\nu_i$ otherwise.
\end{Lemma}

\noindent{\bf Proof: } For the direct statement, assume $x^{\mu}\notin \mm^\nu$. This means $\mu_i<\nu_i$ for all $\nu_i\neq0$ i.e. for all $\nu_i$ s.th. $i\notin C$. Then, clearly $x^\mu\in \bigcup_{x^\rho\vert x^{\bar \nu}}x^\rho\cdot\kb[x_{i_1},\dots,x_{i_k}]$. On the other hand, for the reverse statement, take $x^\mu\in\mm^\nu$, then $\mu_i\geq\nu_i$ for all $i\notin C$, which means $x^\mu\notin \bigcup_{x^\rho\vert x^{\bar \nu}}x^\rho\cdot\kb[x_{i_1},\dots,x_{i_k}]$ $\square$

\begin{Corollary}\label{divisib-irreduc}
Let $x^\mu,x^\nu$ be two different monomials with full support, then
$$x^\mu\vert x^{\bar\nu}\iff x^\mu\notin \mm^\nu$$
\end{Corollary}

For the non-artinian case, we need the following lemmas:

\begin{Lemma}\label{support-irred}
Let $x^\mu,x^\nu\in\NN^n$ such that $supp(x^\mu)=supp(x^\nu)$ and $x^\mu\vert x^\nu$. Then $\mm^\mu\cap\mm^\nu=\mm^\nu$
\end{Lemma}
\noindent {Proof: } Note that since the supports coincide, $\mm^\mu\cap\mm^\nu=\langle x_{i_k}^{max(\mu_k,\nu_k)}\vert k\in supp(x^\mu)=sup(x^\nu)\rangle$; but since $x^\mu\vert x^\nu$, we have $max(\mu_k,\nu_k)=\nu_k$ for all $k$ $\square$

\begin{Lemma}\label{freecone-irreduc}
Let $x^\nu\in \NN^n$, let $C=\{x_{c_1},\dots,x_{c_k}\}\in\{1,\dots,n\}$ a set of variables, let $x^{{\nu'}}=x^\nu/\prod_{i\notin C}x_i^{\alpha_i}$ for some $\alpha_i>0\,\forall i$. Take $x^\mu\in\NN^n$ s.th $x^\mu\nmid x^\nu$; then
$$x^\mu\in x^{{\nu'}}\cdot\kb[x_{c_1},\dots,x_{c_k}]\iff x^\mu\notin\mm^{\underline{\nu}}$$
where $\underline{\nu}_i=0$ if $i\in C$ and $\underline{\nu}_i=\nu_i$ otherwise.
\end{Lemma}

\noindent{\bf Proof: } Because of lemma \ref{support-irred}, we can take, without loss of generality, $\alpha_i=1\forall i$. $\square$

\subsection{Irreducible decomposition of $I$ from $H_*(\KK(\hat{I}))$}
\spanishsubsection{Descomposici\'on irreducible de $I$ a partir de  $H_*(\KK(\hat{I}))$}
Just as we have seen for combinatorial decompositions, a key role in the relation between Koszul homology and irreducible decompositions is played by the $(n-1)$-st Koszul homology of the ideal. This relation comes from the easy but decisive proposition \ref{closed-maximal}.

\begin{Proposition}\label{irreducible-artinian}
Let $I$ be an artinian monomial ideal, and let $\Bc_{n-1}(I)$ be the set of multidegrees $\nu$ such that $H_{\nu,n-1}(\KK(I))\neq 0$. The irredundant irreducible decomposition of $I$ is given by
$$I=\bigcap_{\nu\in\Bc_{n-1}(I)}\mm^{\nu}$$
\end{Proposition}
\noindent {\bf Proof: } From proposition \ref{Stanley-artinian} we know  that a monomial $x^\mu$ is not in $I$ if and only if $x^\mu$ divides $x^{\bar{\nu}}$ for some $\nu\in\Bc_{n-1}(I)$. On the other hand, by corollary \ref{divisib-irreduc}, $x^\mu$ divides $x^{\bar \nu}$ if and only if $x^\mu\notin\mm^\nu$. Therefore $x^\mu\notin I\iff x^\mu\notin\bigcup_{\nu\in\Bc_{n-1,I}}\mm^\nu$, thus
$$x^\mu\in I\iff x^\mu\in\bigcap_{\nu\in\Bc_{n-1}(I)}\mm^\nu$$
Thus, we have an irreducible decomposition of $I$, and since $x^\nu\nmid x^{\nu'}\,\forall x^\nu,x^{\nu'}\in\Bc_{n-1}(I)$, it follows that the decomposition is irredundant $\square$

\begin{Proposition}\label{irreducible-nonartinian}
Let $I$ be a monomial ideal, and $\hat{I}$ its artinian closure; let $\Bc_{n-1}(\hat{I})$ be the set of multidegrees $\nu$ such that $H_{\nu,n-1}(\KK(\hat{I}))\neq 0$. The irredundant irreducible decomposition of $I$ is given by
$$I=\bigcap_{\nu\in\Bc_{n-1}(\hat{I})}\mm^{\underline{\nu}}$$
where $\underline{\nu}_i=0$ if $i$ is in some $C\in TF(x^\nu)$ and $\underline{\nu}_i=\nu_i$ otherwise.
\end{Proposition}

\noindent {\bf Proof:}
On one hand, we have that if the monomial $x^\mu$ is not in $\hat{I}$, then $x^\mu$ is not in $I$. But if $I$ is not artinian, then it is strictly smaller than $\hat{I}$. All the elements that divide the maximal corners $x^\nu$ that are in $\Bc_{n-1}(I)\cap \Bc_{n-1}(\hat{I})$ are not in $\hat{I}$ therefore, take $x^\nu\in\Bc_{n-1}(I)$ such that $x^\nu\notin\Bc_{n-1}(\hat{I})$. We have seen in the procedure that lead to proposition \ref{Stanley-nonartinian} that $TF(x^\nu)\neq\emptyset$ and then for each $C\in TF(x^\nu)$, $x^\mu\in x^\nu\cdot\kb[C]$ if and only if $x^\mu\notin
\mm^{\bar{\nu}}$. Observe now that $\mm^\nu\cap\mm^{\underline{\nu}}=\mm^{\underline{\nu}}$, and thus
$$x^\mu\in I\iff x^\mu\in\bigcap_{\nu\in\Bc_{n-1}(\hat{I})}\mm^{\underline{\nu}}$$
and we obtain the desired decomposition.$\square$

\begin{Example}
Take the artinian ideal $I=\langle x^3,x^2y,xz,y^3,z^3\rangle$, we have that $\Bc_{n-1}(I)=\{(3,1,1),(2,3,1),(1,3,3)\}$, then
$$I=\langle x^3,y,z\rangle\cap\langle x^2,y^3,z\rangle\cap\langle x,y,z^3\rangle$$
Take now $I'=\langle x^2y,xz,y^3,z^3\rangle$, which is not artinian, being $\hat{I'}=I$. The only element of $\Bc_{n-1}(I)$ that is not in $\Bc_{n-1}(I')$ is $(3,1,1)$, and we have that $TF(x^3yz)=\{[x]\}$ therefore 
$$I=\langle y,z\rangle\cap\langle x^2,y^3,z\rangle\cap\langle x,y,z^3\rangle$$

\end{Example}

\subsection{Irreducible decomposition of $I$ from $H_*(\KK(I))$}
\spanishsubsection{Descomposici\'on irreducible de $I$ a partir de  $H_*(\KK(I))$}
Even if it is easy to obtain an irreducible decomposition of $I$ from its artinian closure, it is convenient to have a procedure that computes such decompositions from the Koszul homology of $I$ itself. Algorithm \ref{Koszul->Irreducible} below solves the problem of reading the irredundant irreducible decomposition of $I$ form its Koszul homology. It receives as input the multidegrees of the generators of the Koszul homology of $I$ and returns the irredundant irreducible decomposition of $I$. The algorithm consists basically on four loops:
\begin{itemize}
\item The first loop is described in lines $1-8$ and runs on the generators of the Koszul homology. For each multidegree $\mu$ of them, we compute the corresponding $\overline{\mu}$; these can be equal for several different $\mu$. Now, for each $\overline{\mu}$ we first look whether it has no locally free direction. If it is the case, we are at a closed (or maximal) corner, and these automatically correspond to irreducible components of $I$; thus, we store the corresponding pairs $(\mu,\emptyset)$ in a list of components; note that in this case there is only one $\mu$ that corresponds to $\overline{\mu}$. For the other multidegrees, we find which of the locally free directions are true free directions, and in case there is any, we store the corresponding pairs $(\nu,TF(\overline{\mu}))$ as candidates. Here, $\nu$ is the $lcm$ of all the $\mu$ that correspond to the same $\overline{\mu}$ that we are analyzing. When there are several $\mu$ that correspond to the same $\overline{\mu}$ we have to treat them as a single item, in which the element `out of the ideal' is $\overline{\mu}$ and corresponds to only one element `in the ideal', in this case $\nu$, observe that $\overline{\nu}=\overline{\mu}$.
\item The second loop, which is described in lines $9-11$ is just an auto-reduction of the set of candidates, so that among all the candidates $(\nu,C)$ for each $C\in TF(\overline{\mu})$, we keep only those maximal elements with respect to divisibility of the multidegrees. The so reduced set of candidates is then added to the components in line $12$.
\item The third loop runs in lines $13-15$ and transforms the pairs $(\nu,C)$ in a multidegree $\underline{\nu}$ that will represent the corresponding irreducible component: For each element $C=\{x_{i_1}\dots x_{i_k}\}$ in $TF(\overline{\mu})$ we have an irreducible component $\mm^{\underline{\nu}}$ with $\underline{\nu_{i}}=0$ if $x_i\in C$ and $\underline{\nu_{i}}=\nu_i$ otherwise. These irreducible components form an irreducible decomposition of $I$. If one such $\underline{\nu}$ is equal to $1$, we do not store it since $\langle 1\rangle\cap I=I$ for every ideal $I$, and it would be redundant.
\item Finally, in lines $16-18$ we reduce the set of irreducible components to obtain an irredundant irreducible decomposition of $I$, which is the output of the algorithm (line $19$).
\end{itemize}

\begin{table}[!htb]
\centering
\begin{tabular}{|p{14cm}|}
\hline
$$\begin{array}{ll} Algorithm:\, \mbox{Irreducible decomposition of a monomial ideal }$I$\\
\hline \hline\\
\ Input: \, H_*(\KK(I))\\
\ Output: \mbox{ Irredundant irreducible components of }$I$ \\
\\
\hline\\
\ 1 \quad {\bf for each }\, \overline{\mu}\in \NN^n {\bf s.th. }\, H_{\mu,i}(\KK(I))\neq 0\,\mbox{ for some } i>0\,{\bf do}\\
\ 2 \quad \quad \mbox {compute}\, LF(\overline{\mu}) \\
\ 3 \quad \quad {\bf if}\, LF(\overline{\mu})=\emptyset\,{\bf then}\mbox{ store }(\mu,\emptyset)\mbox{ as a {\bf component}}\\
\ 4 \quad \quad {\bf else}\\
\ 5 \quad \quad\quad \mbox{ compute } TF(\overline{\mu})\\
\ 6 \quad \quad\quad \mbox {store as \bf{candidate} } (\nu,TF(\overline{\mu})) \mbox { if } TF(\overline{\mu})\neq\emptyset\\
\ 7 \quad\quad{\bf endif}\\
\ 8 \quad {\bf end foreach}\\
\ 9 \quad  {\bf while}\, \exists\,(\nu,C)\mbox{ and } (\nu',C) \mbox {in \bf{candidates} s. th. } \nu\vert\nu'\\
\ 10 \quad  \quad \mbox{ delete } (\nu,C)\mbox{from {\bf candidates}}\\
\ 11 \quad {\bf end while}\\
\ 12 \quad  {\bf components}={\bf components}\cup{\bf candidates}\\
\ 13 \quad {\bf for each}\,(\nu,C)\mbox{ in {\bf components}}\\
\ 14 \quad \quad {\bf irred\_components}\leftarrow \underline \nu\quad {\bf if } \underline{\nu}\neq 1\\
\ 15 \quad {\bf end foreach}\\
\ 16 \quad  {\bf while}\, \exists\, \underline{\nu}, \underline{\nu'} \mbox {in \bf{irred\_components} s. th. } supp(\underline{\nu})=supp(\underline{\nu'})\mbox{ and }\underline{\nu}\vert\underline{\nu'}\\
\ 17 \quad  \quad \mbox{ delete } \underline{\nu}\mbox{ from {\bf irred\_components}}\\
\ 18 \quad {\bf end while}\\
\ 19 \quad {\bf return(irred\_components)}\\
\end{array}
$$
\\
\hline
\end{tabular}\caption{Algorithm Koszul$\rightarrow$Irreducible Decomposition}\label{Koszul->Irreducible}
\end{table}

\begin{Remark}
The computation of $LF(\mu)$ and $TF(\mu)$ can always be performed in any case. The first one is based on the computation of the lower Koszul simplicial complex of $I$ at $\mu$ and check whether some multidegrees are in the ideal or not, which always finishes, because our ideals are finitely generated. For the computation of $TF(\mu)$ we just have to check whether multiplication by some given variable lies in the artinian closure of $I$ or not, which is also made in finitely many steps.
\end{Remark}

\begin{Example}\label{example-irreducible}
let us follow the flow of the algorithm with a particular example.
Take the ideal $I=\langle x^2y,xz,y^3,z^3\rangle$. The multidegrees of its Koszul homology are given in the following table:
\begin{center}
\begin{tabular}{|r|l|}
\hline
$i$&$\mu\in \NN^n\mbox{ s.th. }H_{\mu,i}(\KK(I))\neq0$\\
\hline
$0$&$(2,1,0),(1,0,1),(0,3,0),(0,0,3)$\\
$1$&$(2,1,1),(2,3,0),(1,0,3),(1,3,1),(0,3,3)$\\
$2$&$(2,3,1),(1,3,3)$\\
\hline
\end{tabular}
\end{center}
First of all, the first loop, running on the generators of the Koszul homology of the ideal gives us the following lists of candidates and components

$${\bf candidates}=\{((2,1,0),[x,y]),((1,0,1),[x,z]),((0,3,0),\{[x,y],[y,z]\}),((0,0,3),\{[x,z],[y,z]\}),$$
$$((2,1,1),[x]),((0,3,3),[y,z]),((2,3,0),[x,y]),((1,0,3),[x,z])\}$$
$${\bf components}=\{((2,3,1),\emptyset),((1,3,3),\emptyset)\}.$$
The second loop performs some reductions in the {\bf candidates} list:
$${\bf candidates}=\{((2,1,1),[x]),((0,3,3),[y,z]),((2,3,0),[x,y]),((1,0,3),[x,z])\}$$
This final list of candidates is joined to the list of components.

The third loop transforms each pair $(\mu,C),\,C\in TF(\mu)$ to the corresponding $\overline{\mu}$, and stores the result in {\bf irred\_components}:
$${\bf irred\_components}=\{yz,x^2y^3z,xy^3z^3\}$$
In this case, there is no possible reduction in {\bf irred\_components} and thus, we have the following irredundant irreducible decomposition of $I$:
$$I=\mm^{(0,1,1)}\cap\mm^{(2,3,1)}\cap\mm^{(1,3,3)}=\langle y,z\rangle\cap\langle x^2,y^3,z\rangle\cap\langle x,y^3,z^3\rangle$$
\end{Example}

\noindent {\bf Proof of algorithm \ref{Koszul->Irreducible}:} Let us start with the case $I$ is artinian. In this case, the algorithm stores in the {\bf  components} list only the closed and maximal corners, since they have no true free directions. This is done in the second loop; the third loop stores nothing, because there is no element of $I$ with true free directions. Now, the reduction of the {\bf components} list, leaves only the maximal corners, since every other closed corner divides at least one of them, and the cone of free directions is empty for every component. The fourth loop stores as irreducible components just the multidegrees of the list {\bf components}. No further reduction is performed, since all the maximal corners have full support, and they cannot be divisible by each other. So we have that the output of the algorithm is the list of multidegrees of the maximal corners, which by proposition \ref{irreducible-artinian} gives the irredundant irreducible decomposition of $I$.

In the non-artinian case, from proposition \ref{irreducible-nonartinian}, we just have to proof that for each maximal corner $x^\nu$ of $\hat I$ that is not a maximal corner of $I$, and for each $C\in TF(x^\nu)$, we store in the second loop of the algorithm an element $x^{\nu'}$ such that $x^{\nu'}$ divides $x^{\nu}$, there is some $i$ such that $H_{i,\nu'}(\KK(I))\neq0$ and there is $C'\in TF(x^{\nu'})$ such that $C'\subseteq C$ as sets. If such $x^{\nu'}$ exists, then the intersection of the ideals given by the output of the algorithm is the same as the one given in the right-hand side of the equality in proposition \ref{irreducible-nonartinian}.

To see this, let $x^\nu$ be a maximal corner of $\hat I$, and $C\in TF(x^\nu)$, $C=\{x_{C_1},\dots,x_{C_k}\}$. Then, we know that $\mm^{\underline{\nu}}$ is an irreducible component of $I$. Therefore, $x^{\lambda+1\smallsetminus \underline{\nu}}$ is a minimal generator of the Alexander dual of $I$ with respect to $\lambda+1$, $I^{[\lambda+1]}$. We have that there is homology in $x^{\nu}-\kb[C]\cap I$ if and only if there is nonzero homology in $x^{\lambda+1\smallsetminus\underline{\nu}}\cdot\kb[C]$, note that $supp(x^{\lambda+1\smallsetminus\underline{\nu}})=[n]\smallsetminus C$. In order to look for homology in the plane $x^{\lambda+1\smallsetminus\underline{\nu}}\cdot\kb[C]$, we first need some $lcm$ of generators of $I^{[\lambda+1]}$ in this plane. It is easy to see that there is no such $lcm$ only if we are in one of the following situations:
\begin{itemize}
\item The support of every other generator $m$ of $I^{[\lambda+1]}$ is in $[n]\smallsetminus C$. But then $I^{[\lambda+1]}\subseteq \kb[[n]\smallsetminus C]$, if and only if $I\subseteq \kb[[n]\smallsetminus C]$, but this contradicts the fact that $x^\nu$ is a maximal corner.
\item For every generator $m$  of $I^{[\lambda+1]}$ such that $supp(m)\cap C\neq \emptyset$ there is some $i\in C$ such that $m_i>(\lambda+1-\underline{\nu})_i$. This means that all the maximal corners and irreducible components of $I^{[\lambda+1]}$ satisfy this same property. This implies that all the generators $\underline{m}$ of $I$ (which correspond to these components), satisfy that there exists $i\in C$ such that $\underline{m}_i<\nu_i$ for every generator $\underline{m}$ such that $\underline{m}_j<\nu_j$ $\forall j\notin C$. Thus, we have that $x^\nu$ has further true free directions, which is a contradiction.
\end{itemize}
Then we know that there is some $lcm$ located in the plane  $x^{\lambda+1\smallsetminus\underline{\nu}}\cdot\kb[C]$, i.e. there is some syzygy in it. Consider the syzygies in this plane. Since we are dealing with monomial ideals, every syzygy lies at a multidegree (in particular, lies at the lcm-lattice of $I$). Take a syzygy in our plane with minimal multidegree with respect to divisibility. Assume it is not  a minimal syzygy, then it can be expressed as a linear combination of other syzygies, not lying in our plane, i.e. $S=\alpha_1S_1+\cdots+\alpha_kS_k$, $deg(\alpha_i)>0$. Then the multidegree of $S_1$, for instance, lies at the $lcm$ of some generators, being $m'$ one of them. Now, the multidegree of $S$ is a multiple of $m'$, but then, there is a syzygy the multidegree of which is the $lcm$ of $m'$ and $x^{\lambda+1-\underline{\nu}}$ is at the $lcm$ of them, and it is in our plane. Since $S$ is a multiple of both, it is divisible by the multidegree of this new syzygy, which is a contradiction.
Therefore $S$ is a minimal syzygy, and thus there is homology at this multidegree.$\square$

The close relation between the irreducible decomposition of $I$ and an artinian ideal $I'$ derived from $I$ is already present in \cite{BPS98}, and is also used in \cite{R07}. We have included the computation of the irreducible decomposition of $I$ directly from the Koszul homology of $I$ itself because of completeness from a theoretical point of view. Algorithmically, it looks more efficient to use artinian ideals, because then we only need the $(n-1)$-st homology groups, which are easy to compute, and furthermore, their number can be computed using the calculations in \cite{A97}. Also, an alternative for the computation of irreducible decomposition is the efficient algorithm recently presented in \cite{R07}, which is based on a simplicial complex that stores the structure of the boundary of a monomial ideal. This algorithm is known as the \emph{label algorithm}. However, going to the artinian closure first is not always the more efficient solution, in particular when the number of variables is high and we have to introduce a significant number of new generators.

\section{Primary decompositions, associated primes, height}\label{primary-decomposition}
\spanishsection{Descomposiciones primarias, primos asociados, altura}
Every monomial ideal has a primary decomposition in which all components are monomial ideals. Recall that a monomial ideal is primary if it is of the form $\langle x_{i_1}^{a_1},\dots,x_{i_j}^{a_j},x^{b_1},\dots,x^{b_k}\rangle$ where $supp(x^{b_l})\subseteq \{i_1,\dots,i_j\}$ for all $l$. If for any two components $\mq_i,\mq_j$ in a (monomial) primary decomposition of a monomial ideal $I$ we have that $rad(\mq_i)\neq rad(\mq_j)$ if $i\neq j$ then we say that the decomposition is irredundant. Even minimal irredundant primary decompositions are not unique for monomial ideals. What is unique in such a decomposition is the number of its terms and the primary components that correspond to minimal associated primes \cite{V01}. Algorithms for computing primary decompositions are given in \cite{V01}, using elimination of powers of variables; and in \cite{HS01}, using Alexander duals; the latter is implemented in \macaulay \cite{M2}.

In any case, the relation between irredundant irreducible and primary decompositions is very close. Given an irredundant irreducible decomposition of $I$, it is easy to obtain an irredundant primary decomposition: Take all irreducible components that share the same support and intersect them. Each of these intersections will be a primary ideal and the intersection of all of them equals $I$, therefore we have a primary decomposition of $I$ which is clearly irredundant.

To obtain the associated primes of $I$ we can either start from the irreducible or from the primary decomposition. If starting from the irreducible decomposition, just take the different supports of the irreducible components. If starting from the primary decomposition, take the radical of each primary component, which amounts to consider the support of it.

The height of a prime ideal is defined as the supremum of the lengths of all chain of prime ideals contained in $P$. In the monomial case it amounts to the number of variables of the support of $P$. Then the height of an ideal $I$ is
$$ht(I)=min\{ht(P)\vert I\subset P \mbox{ and } P \in Spec(R)\}$$
In the case of monomial ideals it is just the minimum of the number of variables in the support of the associated primes of $I$ (i.e. the minimum of the number of variables in the support of any component in an irredundant irreducible or primary decomposition of $I$).

We have seen that we can obtain an irredundant irreducible decomposition of $I$ from its Koszul homology, and it is easy to compute another primary decomposition of $I$ from which we can directly read the associated primes and height of the ideal. Thus, this information is also very easily accessible from the Koszul homology of the ideal.

\begin{Example}\label{example-primary}
In example \ref{example-irreducible} we obtained the irreducible decomposition
$$I=\langle x^2y,xz,y^3,z^3\rangle=\langle y,z\rangle\cap\langle x^2,y^3,z\rangle\cap\langle x, y^3,z^3\rangle$$
The first irreducible component is the only one having only $\{y,z\}$ as support and the other two have full support, so taking their intersection $\langle x^2,y^3,z\rangle\cap\langle x, y^3,z^3\rangle=\langle x^2,y^3,z^3,zx \rangle$ we have a primary decomposition $I=\langle y,z\rangle\cap\langle x^2,y^3,z^3,zx \rangle$. The corresponding associated primes are then $\langle y,z\rangle$ and $\langle x,y,z\rangle$ the first being minimal and the second embedded. The height of this ideal is then clearly $2$.
\end{Example}

This process can be summarized in a tree the nodes of which contain sets of multidegrees and are labeled in the following form:
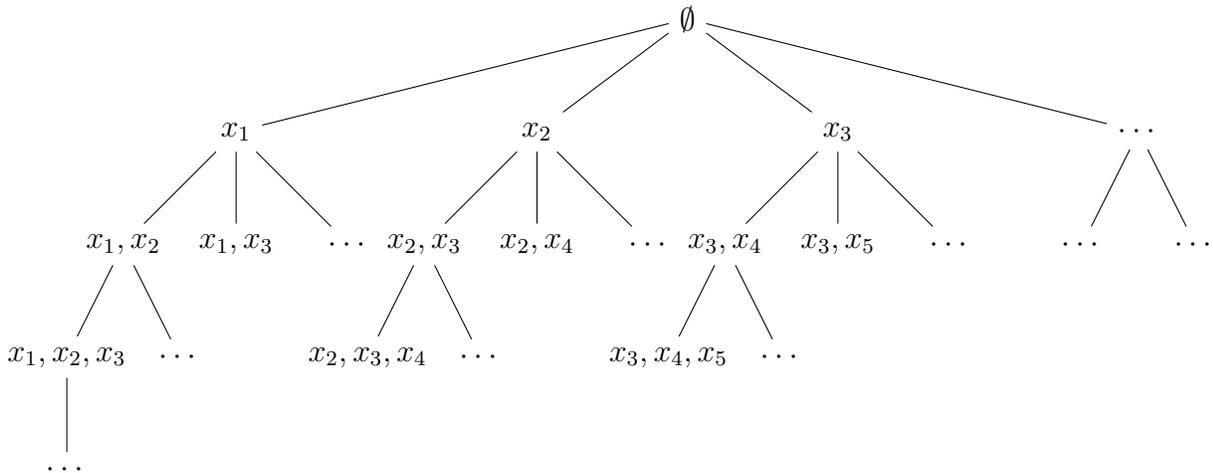
\begin{center}
\begin{tikzpicture}[scale=1]
\tikzstyle{level 1}=[sibling distance=4cm]
\tikzstyle{level 2}=[sibling distance=1.5cm]
\node{$\emptyset$}
child{ node {$x_1$}
	child{node{$x_1,x_2$}
		child{node{$x_1,x_2,x_3$}
			child{node{$\cdots$}}}
		child{node{$\cdots$}}}
	child{node{$x_1,x_3$}}
	child{node{$\cdots$}}}
child{ node {$x_2$}
	child{node{$x_2,x_3$}
		child{node{$x_2,x_3,x_4$}}
		child{node{$\cdots$}}}
	child{node{$x_2,x_4$}}
	child{node{$\cdots$}}}
child{ node {$x_3$}
	child{node{$x_3,x_4$}
		child{node{$x_3,x_4,x_5$}}
		child{node{$\cdots$}}}
	child{node{$x_3,x_5$}}
	child{node{$\cdots$}}}
child{ node {$\cdots$}
	child{node{$\cdots$}}
	child{node{$\cdots$}}};
\end{tikzpicture} 
\end{center}
In this tree, each node is labelled with a set of variables $x_{i_1},\dots,x_{i_k}$, and contains all the multidegrees of the $(n-1)$-st Betti multidegrees of $\hat I$ (the artinian closure of $I$) such that the exponents of the variables $x_{i_1},\dots,x_{i_k}$ are bigger than $\lambda_i$ where $\lambda$ is the least common multiple of the generators of $I$. If we make zero such variables in the multidegrees $\mu$ in the node, obtaining $\underline {\mu}$, then the $\mm^{\underline{\mu}}$ are the irreducible components of $I$ with support $[n]\smallsetminus\{x_{i_1},\dots,x_{i_k}\}$. The intersection of the irreducible components in each node gives us one of the primary components in the irredundant primary decomposition of $I$. The nonempty labels correspond to the complements of the associated primes, and those which are leaves give the minimal associated primes. The height of $I$ is given by subtracting the maximum height of a leaf in the tree from $n$.

\begin{Example}
The tree of the ideal $I=\langle x^3y^5z,y^5z^4,y^3z^5,xyz^5,x^2z^5,x^4z^3,x^4y^2z^2,x^4y^4z\rangle$ \cite[Example 5.22]{MS04} is
\begin{center}
\begin{tikzpicture}[scale=1]
\tikzstyle{level 1}=[sibling distance=4cm]
\tikzstyle{level 2}=[sibling distance=1.5cm]
\node{$\emptyset: (4,5,5)$}
child{ node {$x:(5,2,3),(5,4,2)$}
	child{node{$x,y:(5,6,1)$}}}
child{ node {$y:(3,6,4)$}}
child{ node {$z:(2,1,6),(1,3,6)$}};
\end{tikzpicture} 
\end{center}
Which gives the irredundant primary decomposition
$$I=\langle x^4,y^5,z^5\rangle\cap\langle y^4,z^3,y^2z^2\rangle\cap \langle z\rangle\cap\langle x^3,z^4\rangle\cap\langle x^2,y^3,xy\rangle$$
and $ht(I)=1$.
\end{Example}
A procedure similar to that exposed in section \ref{irreducible_decomposition} can be followed to obtain an irredundant primary decomposition directly from the Koszul homology of $I$, using the sets of free directions of the Koszul homology generators of $I$; this procedure will not be exposed here.


\section{Koszul Homology for Polynomial Ideals}\label{poly}
\spanishsection{Homolog\'ia de Koszul para ideales polinomiales}
In this section we move for a moment from the study of monomial ideals to that of polynomial ideals. The question we face is whether the techniques we are exploring for computing the Koszul homology and free resolutions of monomial ideals are of any use when computing resolutions of general polynomial ideals. 

Minimal free resolutions of polynomial ideals have already been studied by many authors, see for example \cite{S80,LS98,S99,KR05,GP02}. There is a strong relation between the minimal free resolution of a polynomial ideal and that of its leading ideal, see \cite{H04} for example, and the references therein. Algorithms for the computation of minimal free resolutions have been implemented in the main computer algebra systems, like CoCoA \cite{cocoa}, Macaulay \cite{M2} or Singular \cite{Singular}. These algorithms make use of different strategies to perform the computations; what we explore in this section is a different strategy which is not implemented in the mentioned computer algebra systems, namely transfer the problem of computing the minimal resolution of a polynomial ideal to making the computation in a monomial setting then perform the computations using the techniques we have available, and finally transfer back the results to the general polynomial setting using homological perturbation. 

This strategy has been presented so far in two different approaches. One is represented by L. Lambe et al. \cite{JLS02,LS02}, the other one by F. Sergeraert et al. \cite{S06,RS06,RSS06}. Both approaches have in common the crucial use of the \emph{basic perturbation lemma}, which is the central result in \emph{effective homology} to perform the transfer between the monomial and the polynomial settings. The differences are given by the different complexes and techniques used to perform this transfer: the \emph{Bar} resolution is the main tool in one case, \emph{mapping cone} resolutions are preferred in the other. 

In this section, \emph{effective homology} plays a crucial role. The basics of it should be read in Appendix \ref{apef} and the references therein if the reader is not familiar with them. Another main ingredient here is Gr\"obner basis theory, although we will only use the basics of it. Again, the reader unfamiliar with this subject can read the simplest concepts in Appendix \ref{apal}.

Both approaches considered in this chapter share the goal of producing actual algorithms and computations. The algorithms produced in the second one are conceptually clear and computationally efficient, being the implementations available in the \emph{Kenzo} system \cite{Kenzo}.

Making use of the Basic Perturbation Lemma, we can design a strategy for computing the Koszul homology or minimal resolutions of polynomial ideals. The basic procedure would work as follows: 

\begin{itemize}
\item Let $I\subseteq \poly \kb x n $ be an ideal of the polynomial ring. We compute a Gr\"obner basis for $I$ with respect to some monomial ordering $>$ and obtain the initial ideal $in_>(I)=J$. We have that $\KK(R/I)$ and $\KK(R/J)$ are canonically isomorphic as graded $\kb$-vector spaces, but they have different differentials.

\item Since $J$ is a monomial ideal and using the techniques in sections \ref{mon} and \ref{hom_tools} or in chapter \ref{computation}, we compute the Koszul homology and/or minimal resolution of $J$. We can also use any monomial resolution of $J$. As we will see later, we need an explicit contracting homotopy of the used resolution, in order to apply the basic perturbation lemma. This lemma needs two complexes (resolutions in our case) that are in some sense homologically equivalent. One of them plays the role of the \emph{big} complex and the other is the \emph{small} complex. If we perturb the differential in the big complex, we can under some conditions, obtain a perturbation of the differential in the small complex such that the perturbed complexes are still homologically equivalent. In our case this \emph{perturbation} will perform the transfer between a resolution of the initial ideal of a polynomial ideal and a resolution of the latter.

\item Use the basic perturbation lemma to deduce the Koszul homology/minimal resolution of $R/I$. We will in general obtain non-minimal resolutions of $R/I$ but there exist standard procedures to minimize resolutions, as we have seen in section \ref{resolutions}.

\end{itemize}

\subsection{Perturbing the Lyubeznik resolution}
\spanishsubsection{Perturbaci\'on de la resoluci\'on de Lyubeznik}

The final step of this strategy, i.e. the actual use of the basic perturbation lemma can be done in different ways, using different reductions. One approach uses the $Bar$ resolution as the \emph{big} complex together with Lyubeznik resolution, which plays the role of the \emph{small} complex. This procedure is developed in detail in \cite{JLS02} and \cite{LS02}. The main theorem in the first of these papers is a formulation of this strategy:

\begin{Theorem}[\cite{JLS02},Theorem 5.1]\label{barlyu}
Let $I=(f_1,\dots,f_r)$ be  an ideal in $R$ and choose a monomial order. Let $G=\{g_1,\dots,g_s\}$ be a Gr\"obner basis for $I$ with respect to the given order. Let $M=lt(G)=(m_1,\dots,m_s)$ be the ideal generated by the leading terms of $G$. Let $L$ be the Lyubeznik resolution of $R/M$ over $R$. There is an explicit perturbation $p$ of the differential $d$ in $L$ so that $(L,d+p)$ is a resolution of $R/I$ over $R$.
\end{Theorem}

The heart of the proof is conceptually simple, we just need to apply the basic perturbation lemma to a reduction between two resolutions of $R/M$ that satisfy certain conditions. These conditions are satisfied by the Lyubeznik and Bar resolutions. Let us briefly explain the details:

\subsubsection*{The Bar resolution.}
\spanishsubsubsection{La resoluci\'on Bar}

Let us consider the Bar resolution for algebras \cite{M95}. In our case, we have the polynomial ring $R$, over $\kb$ which has a $\kb$-algebra structure. Let us consider an $R$-module $N$ and the relatively free $R$-module
$$B_n(R,\Mc)=R\otimes(R/\kb)\otimes\stackrel{n\,times}{\cdots}\otimes(R/\kb)\otimes N$$
where $\otimes$ is $\otimes_{\kb}$ and $R/\kb$ is the cokernel of the unit map $\kb\rightarrow R$. The $B_n(R,\Mc)$ are spanned by elements of the form
$f[f_1|\dots|f_n]m\equiv f\otimes[(f_1+\kb)\otimes\dots\otimes(f_n+\kb)]\otimes m$.

If we write $\bar{B}(R)=\sum_{n=0}^\infty\bar{B}_n(R)$ with $\bar{B}_0(R)=\kb$ and $\bar{B}_n(R)=\otimes^n(R/\kb)$ then we have the following $R$-module
$$B(R,\Mc)=R\otimes\bar{B}(R)\otimes \Mc$$. To have an $R$-complex, we give $B(R,\Mc)$ the following differential: 

\begin{eqnarray*}
\partial_n(f[f_1|\dots|f_n]m)&=&ff_1[f_2|\dots|f_n]m+\sum_{i=1}^{n-1}(-1)^if[f_1|\dots|f_if_{i+1}|\dots|f_n]m\\
& & +(-1)^nf[f_1|\dots|f_{n-1}]f_nm
\end{eqnarray*}

Consider now $B(R,\Mc)\stackrel{s_{\Mc}}{\rightarrow}B(R\Mc)$ given by
$$s_{\Mc}(f[f_1|\dots|f_n]m)=[f|f_1|\dots|f_n]m,$$
and in the case $\Mc=R/I$ for some ideal $I$ of $R$ , let $R/I\stackrel{nf_I}{\rightarrow}R$ given by normal form with respect to a chosen term ordering $<$. Then, we have a reduction \cite{LS02},\cite{JLS02}
\begin{center}
\begin{tikzpicture}
\draw[->](0,0) node(x){$B(R,R/I)$}	
	(x) ..controls +(2,0.5)and +(2,-0.5) ..(x) node[pos=0.5,right]{$s_{R/I}$}; 
\draw[->](-0.1,-0.3)--(-0.1,-1.5) node[pos=0.5,left]{$nf_I$} node[below]{$R/I$};
\draw[->](0.1,-1.5)--(0.1,-0.3) node[pos=0.5,right]{$\sigma_{R/I}$};
\end{tikzpicture}
\end{center}
we call $B(R,\Mc)$ the Bar resolution of $R/I$ over $R$.

\subsubsection*{Lyubeznik resolution}
\spanishsubsubsection{La resoluci\'on de Lyubeznik}

We have already met Lyubeznik's resolution \cite{L98} in section \ref{examples-resolutions}. Recall that it is a subresolution of Taylor's defined as follows: For a given subset $J\subseteq \{1\dots r\}$ and an integer $1\leq s\leq r$, let $J_{> s}=\{j\in J\vert j>s\}$; then the Lyubeznik resolution, $\LL$ is generated by those basis elements $u_J$ such that for all $1\leq s\leq r$ one has that $m_s$ does not divide $m_{J_{> s}}$. Note that this resolution depends on the order in which the generators of the ideal are given.

Lyubeznik resolution is an example of {\it relatively free} complex, i.e. it is of the form $X=R\otimes \bar{X}$ where $\bar{X}$ is free over $\kb$ \cite{M95}. In \cite{F78}, an explicit contracting homotopy was given for Taylor resolution, which can be restricted to Lyubeznik resolution \cite{JLS02}. This homotopy works as follows: for a monomial $x^\ab$ and a basis element $u_J$, let $\iota(x^\ab u_J)=\min\{i|m_i\leq x^\ab m_J\}$. Define now a $\kb$-linear endomorphism in $\LL$ by
$$h(x^\ab u_J)=[\iota < j_1]\frac{x^\ab m_J}{m_{\{\iota\}\cup J}}u_{\{\iota\}\cup J}$$ where $\iota=\iota(x^\ab u_J)$ and $[p]$ equals $0$ if $p$ is false and $1$ otherwise. An important property of $\LL$ is that $h(\bar{\LL})=0$ \cite{JLS02}, with this property, we can use the following lemma
\begin{Lemma}
Suppose that $(X,d_X)$ and $(Y,d_Y)$ are relatively free resolutions and that the contracting homotopy $h_X$ satisfies $h(\bar{X})=0$ and $d_X(\bar{X})\cap \bar{X}=0$ and the homotopy $h_Y$ satisfies $h_Y(Y)\subseteq(\bar{Y})$, then the constructions above give a reduction
\begin{center}
\begin{tikzpicture}
\draw[->](0,0) node(x){$Y$}	
	(x) ..controls +(1.5,0.5)and +(1.5,-0.5) ..(x) node[pos=0.5,right]{$h$}; 
\draw[->](-0.1,-0.3)--(-0.1,-1.5) node[pos=0.5,left]{$f$} node[below]{$X$};
\draw[->](0.1,-1.5)--(0.1,-0.3) node[pos=0.5,right]{$\nabla$};
\end{tikzpicture}
\end{center}
\end{Lemma}
We have just seen that the Lyubeznik resolution, the Bar resolution and their corresponding homotopies, satisfy the hypotheses in this lemma. Therefore we can now use theorem \ref{barlyu} to build a resolution of a polynomial ideal $I$ from the Lyubeznik resolution of its initial ideal $lt(I)$. The proof is given in \cite{JLS02}; the main ideas are the following: There is a vector space isomorphism $R/I\cong R/lt(I)$ and from it we have an isomorphism of $B(R,R/I)$ and $B(R,R/lt(I))$ {\it as vector spaces}. Then we have that $B(R,R/lt(I))$ supports two differentials, the one for $R/lt(I)$ and the one for $R/I$. Thus the last differential can be seen as a perturbation of the first one. On the other hand, from the precedent lemma, we have a reduction between $B(R,R/lt(I))$ and $\LL$. Now, using the perturbation lemma, we obtain a reduction between $B(R,R/I)$ and the perturbed $\LL$ which becomes then a resolution of $R/I$. The only point to check is the nilpotency condition, which, as seen in \cite{JLS02}, is satisfied. 

\begin{Remark}
Lyubeznik resolutions of monomial ideals are not minimal in general. And even if it is the case, the minimality of the Lyubeznik resolution $\LL$ of $R/M$ does not imply minimality of $\LL$ with the new differential, thus, a further step should be performed to obtain the minimal resolution of $R/I$.
\end{Remark}

\begin{Remark}
The resolution provided by the basic perturbation lemma is \emph{effective}, in the sense that an explicit contracting homotopy is provided for it as a byproduct. One does not obtain explicit closed form formulas for the differential in the new resolution, but the basic perturbation lemma results in algorithmic formulas that can be easily implemented.
\end{Remark}

\subsection{Using mapping cone resolutions}
\spanishsubsection{Uso de resoluciones de conos de aplicaciones}

The second procedure we will examine is given by F. Sergeraert et al. in \cite{S06,RS06,RSS06}. It uses a recursive procedure to compute the minimal resolution of $R/lt(I)$ based on mapping cones. At each step of the recursion, a reduction is produced between  a free $\kb$-chain complex of finite type (the {\it small} complex in this case) and a certain complex, $\widehat{C_*}$ produced by a recursive procedure (the {\it big} complex). This big complex is produced using iterated mapping cones in an effective homology framework. The basic idea is to use theorem \ref{ses3} which produces effective resolutions. Another reduction is produced between this big complex and the Koszul complex of $R/lt(I)$. The differential of this Koszul complex can be perturbed to obtain the Koszul complex of $R/I$; using the {\it easy basic perturbation lemma}, this new differential is transferred to $\widehat{C_*}$ and then, the basic perturbation lemma is applied again, this time between $\widehat{C_*}$ with the new differential, and the small $\kb$-complex, which now describes the Koszul homology of $R/I$. Observe that the context of this approach is Koszul homology computation, therefore, the minimality of the used resolutions is important. Let us go a bit into the details (full details and completely described examples are given in the references).

The first step of this approach is the same one we have just seen in the precedent section, namely, we start with our polynomial ideal $I$ and compute a Gr\"obner basis of it with respect to some selected term order. We obtain then a monomial ideal $lt(I)=M$ and we work in a monomial setting. At this point things change; instead of choosing a given resolution with an explicit contracting homotopy, we actually build the \emph{effective} minimal resolution of $M$ in a recursive way, in fact making extensive use of the basic perturbation lemma. This is done in a similar way as in sections \ref{resolutions} and \ref{computations}: The process is recursive on the number of minimal generators of the ideal. Given a monomial ideal $M=\langle m_1,\dots,m_r\rangle$, the short exact sequence of $R$-modules
$$0\rightarrow R/\langle M:m_r\rangle\stackrel{\times m_{r-1}}{\rightarrow} R/\langle m_1,\dots,m_{r-1}\rangle\stackrel{pr}{\rightarrow}R/M\rightarrow 0$$
induces a short exact sequence involving the corresponding Koszul complexes. Using theorem \ref{ses3} we have that the effective homology of $R/M$ is given by the effective homologies of $R/\langle M:m_r\rangle$ and $R/\langle m_1,\dots,m_{r-1}\rangle$. The first recursion step is given by the short exact sequence $0\rightarrow R\stackrel{\times m_1}{\rightarrow}R\rightarrow R/\langle m_1\rangle\rightarrow 0$; therefore, we only need the homology of $\KK(R)$ to start the process. What we obtain using this recursive procedure is an \emph{effective} resolution, and therefore we have an explicit homotopy, so that homological perturbation can be applied. In fact what we obtain is an equivalence

$$(\KK(R/M),\partial_M)\lrdc (C_*,d_C)\rrdc (\kb_\bullet,\delta)$$

where the left-hand complex is the Koszul complex of $M=lt(I)$, the central complex is a complex we have constructed in the recursive process, and the right-hand complex is a $\kb$-chain complex of finite type.

We now have that since $R/I$ and $R/M$ are isomorphic as $\kb$-vector spaces, their Koszul complexes  $(\KK(M),\partial_M)$ and $(\KK(I),\partial_I)$ are also isomorphic as graded $\kb$-vector spaces, but they have different differentials. Therefore, we can apply the basic perturbation lemma, provided the nilpotency condition is satisfied. But this is indeed the case, the details can be seen in \cite{RS06} and are based in the fact that the different between the two differentials strictly decreases the multigrading of the term it is applied to. Our goal is then to obtain an equivalence

$$(\KK(R/I),\partial_I)\lrdc (C'_*,d_{C'})\rrdc (\kb_{\bullet}',\delta')$$

This is done in two steps:
\begin{itemize}
\item First we apply the \emph{easy} basic perturbation lemma (see Appendix \ref{apef}) between $\KK(M)$ and $\KK(I)$ obtaining $(C'_*,d_{C'})$ all the hypothesis here are trivially satisfied, so this step is straightforward. Note that $(C_*,d_{C})$ and $(C'_*,d_{C'})$ are just the same graded $\kb$-vector spaces, with different differentials.
\item Second, apply now the basic perturbation lemma between $(C_*,d_{C})$ and $ (C'_*,d_{C'})$ to obtain a perturbation in $(\kb_\bullet,\delta))$ that gives the looked for equivalence. Again, the nilpotency condition must be satisfied here, and again a multidegree argument shows that it is the case \cite{S06,RS06}.
\end{itemize}

With this consideration, we have the following

\begin{Theorem}[\cite{RS06}, Theorem 95]
The homological problem \footnote{The meaning of \emph{homological problem} is described in Appendix \ref{apef}.}of $\KK(R/I)$ is solved. 
\end{Theorem}

\begin{Remark}
Observe that what we obtain in this approach is an equivalence between the Koszul complex of $I$ and a finite effective complex, i.e. we have computed the Koszul homology of $I$. If what we are looking for is the minimal resolution of $I$, an explicit description of it can be obtained by similar methods from the \emph{effective} Koszul homology of $I$. This have been already seen in section \ref{structure_minres}, since the methods exposed there are also valid for polynomial ideals.
\end{Remark}
\fancyhead[ER]{\itshape Chapter 3 \ \ Computation of Koszul homology}
\chapter{Computation of Koszul Homology}\label{computation}
\spanishchapter{C\'alculo de la homolog\'ia de Koszul}
In this chapter we introduce a new technique for the computation of Koszul homology, (multigraded) Betti numbers and resolutions of monomial ideals. This technique is based on Mayer-Vietoris sequences (see section \ref{m-v}) and Mayer-Vietoris trees associated to monomial ideals.

In the first section of the chapter we introduce Mayer-Vietoris trees, from the Mayer-Vietoris sequences we saw in chapter one. Mayer-Vietoris trees can be considered, on one hand as a description of the ideal together with the relevant part of its $lcm$-lattice, and on the other hand as an algorithm that allows us to compute the homological invariants of the ideal. In particular, it gives us bounds for the multigraded Betti numbers without computing the minimal free resolution of the ideal. In this section also Mayer-Vietoris ideals are introduced; these are those ideals for which the bounds on the Betti numbers are sharp, and are divided in three types. Each Mayer-Vietoris tree of a monomial ideal gives us a free resolution of it, and therefore, an expression of the multigraded Hilbert function of the ideal. We also give explicit ways to compute resolutions and Koszul homology generators of monomial ideals from their Mayer-Vietoris trees. Mayer-Vietoris trees and some of their main properties were presented in \cite{S06a}.

The second section analyzes several special types of Mayer-Vietoris trees that will be applied in chapter four for the analysis of several types of ideals. This section uses Mayer-Vietoris trees as descriptions of the ideal and shows how the analysis of these trees can be used to obtain the main homological invariants of the ideal allowing a more detailed analysis of it. The use of these features of Mayer-Vietoris trees to analyze the homological structure of monomial ideals is presented in \cite{SC07d,SC07c}.

Finally, the third section treats Mayer-Vietoris tree from an algorithmic point of view. The basic algorithm for the construction of Mayer-Vietoris trees is exposed here together with some issues about its implementation and different versions. Some results of an implementation using the \cpp library $\cocoalib$ are exposed, together with some reasons for the use of this library and some technical questions. In the last part of the section we show some experiments made on some types of ideals, and some comparisons with other algorithms that compute homological invariants of ideals, which are implemented in the available computer algebra systems. We compare the basic Mayer-Vietoris algorithm with the algorithms used to compute minimal free resolutions and multigraded Hilbert series in $\cocoa$, $\singular$  and $\macaulay$. 

\section{Mayer-Vietoris trees of monomial ideal}
\spanishsection{\'Arboles de Mayer-Vietoris de ideales monomiales}

Using recursively the exact sequences we saw in section \ref{m-v} for every $\ab\in\NN^n$ we could compute the Koszul homology of $I=\langle m_1,\dots,m_r\rangle$. When making use of them to compute $H_*(\KK_{\ab}(I))$, we need $H_*(\KK_{\ab}(\tilde I_r))$ and $H_*(\KK_{\ab}(I_{r-1}))$; and for each of these two computations, one needs the corresponding smaller ideals. Note that the size of the involved ideals and their number of minimal generators decreases until they are generated by only one monomial, in which case the Koszul homology is trivial.

The involved ideals can be displayed as a tree, the root of which is $I$ and every node $J$ has  as \emph{children} $\tilde J$ on the left and $J'$ on the right (if $J$ is generated by $r$ monomials, $\tilde J$ denotes $\tilde J_r$ and $J'$ denotes $J_{r-1}$).
This is what we call a \textbf{Mayer-Vietoris tree} of the monomial ideal $I$, and we will denote it $MVT(I)$.

\begin{Remark}\label{strategies}
Note that the construction of $MVT(I)$ depends on the selection of the {\it distinguished} monomial used to split the ideal (see remark \ref{pivot}) which we call the \emph{pivot} monomial. Thus, for a given ideal $I$ we have several (and eventually different) trees depending on the way we choose this generator.

Basically, the first strategy one thinks of is to select the pivot monomial according to some term order $\tau$. In this case, we speak of the $\tau$-Mayer-Vietoris tree, in which we choose as pivot monomial the biggest one according to some term order. These trees are unique given $\tau$ and $I$.

Another basic strategy will try to keep or reach some property of the ideals of the nodes in the tree, such like genericity, being segments with respect to some term order, etc. Sometimes it is useful to eliminate first some monomials that prevent the ideal from having some property that helps in the computation. For example, consider the ideal $I=\langle x^6y, x^2y^3,xy^4, xz\rangle$, it looks reasonable to eliminate first $xz$, since then $I'$ is an ideal in two variables, hence it is Mayer-Vietoris of type $A$ i.e. we are in the optimal situation (see section \ref{MV-ideals} for the definition of Mayer-Vietoris ideals), and moreover it does not contain any appearance of the variable $z$, so there will be no interaction between the two main branches of the tree, and if there would be some repeated multidegree in a relevant node, it would just be \emph{inside} the subtree hanging from $\tilde I$ (observe that in the example, this is not the case, since $\tilde I$ is again a Mayer-Vietoris ideal of type $A$).

Finally, another relevant strategy selects as pivot monomial $m_s$ a generator having the biggest exponent in some of the variables, then the generators in $\tilde J$ will all have the same exponent in this variable. Keeping the nodes sorted in such a way gives us small trees, as the left branch will have a length of at most the number of variables plus one. This strategy has shown good behaviour in random examples.

Note that in each node we can change the strategy for constructing the tree. Although many times the trees will be built following the same strategy in every node, this flexibility is of use in some applications and examples.
\end{Remark}

\begin{Example}\label{mvt-examples-1}
The monomial ideal $I=\langle xy,xz,yz\rangle$, has the following Mayer-Vietoris tree, if we change the order of the generating monomials, we obtain isomorphic trees (i.e. isomorphic as graphs and with isomorphic ideals in the corresponding nodes)

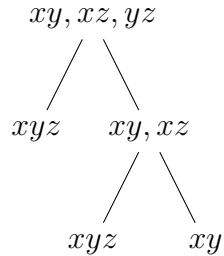
\begin{center}
\begin{tikzpicture}
\node{$xy,xz,yz$}
child{node{$xyz$}}
child{node{$xy,xz$}
child{node{$xyz$}}
child{node{$xy$}}};
\end{tikzpicture}
\end{center}

On the other side, for $I=\langle x^2,y^2,xy\rangle$ we obtain two different trees according to the ordering of the monomials:

\begin{center}

$\begin{array}{cc}

\begin{tikzpicture}[scale=1]
\tikzstyle{level 1}=[sibling distance=3cm]
\tikzstyle{level 2}=[sibling distance=2cm]
\node{$x^2,y^2,xy$}
child{ node {$x^2y,xy^2$}
	child{node{$x^2y^2$}}
	child{node{$x^2y$}}}
child{ node {$x^2,y^2$}
	child{node{$x^2y^2$}}
	child{node{$x^2$}}};
\end{tikzpicture} 

&

\begin{tikzpicture}[scale=1]
\tikzstyle{level 1}=[sibling distance=3cm]
\tikzstyle{level 2}=[sibling distance=2cm]
\node{$xy,x^2,y^2$}
child{ node {$xy^2$}}
child{ node {$xy,x^2$}
	child{node{$x^2y$}}
	child{node{$xy$}}};
\end{tikzpicture} 
\end{array}$
\end{center}

\end{Example}
\begin{Remark}
Since by definition, in $MVT(I)$ every \emph{father} has exactly two \emph{children}, we can assign position indices to every node, in the following way: $I$ has position $1$ and if $J$ has position $p$ then $\tilde J$ has position $2p$ and $J'$ has position $2p+1$. We will denote this $MVT_1(I)=I,\,MVT_{p}(I)=J,\,MVT_{2p}(I)=\tilde J,\, MVT_{2p+1}(I)=J'$. These indices will be very useful for efficiently reading the information hidden in $MVT(I)$.
\end{Remark}

\subsection{$MVT(I)$ and Koszul homology computations}
\spanishsubsection{$MVT(I)$ y c\'alculos de homolog\'ia de Koszul}

Let $I$ be a monomial ideal, we show that all the multidegrees $\ab\in\NN^n$ such that $\beta_{i,\ab}(I)\neq 0$ for some $i$, are present as exponents of generators in some node $J$ of any Mayer-Vietoris tree of $I$. In the following lemmas $J$ is a node of $MVT(I)$ minimally generated by $\{m_1,\dots,m_s\}$, and $m_s$ is the pivot monomial.

\begin{Lemma}
$H_0(\KK(J))=H_0(\KK(J'))\oplus\langle m_s\rangle$
\end{Lemma}

\noindent {\bf Proof: } In homological degree $0$ we have the following exact sequences in homology for each $\ab\in \NN^n$:
$$\cdots\longrightarrow H_{0}(\KK_{\ab}(\tilde J))\longrightarrow H_0(\KK_{\ab}(J')\oplus \KK_{\ab}(\langle  m_s\rangle  ))\longrightarrow H_0(\KK_{\ab}(J))\longrightarrow 0$$

but if $x^\ab$ is such that $\beta_{0,\ab}(J'\oplus\langle m_s\rangle)\neq 0$ then $x^\ab=m_j$ for some $j$. Now, if $\beta_{0,\ab}(\tilde J)\neq 0$ then $x^\ab=m_{i,s}$ for some $i$. Then $m_j=m_{i,s}$, and we have that $m_j|m_i$, $m_j|m_s$ which is a contradiction, since $\{m_1,\dots,m_s\}$ is a minimal generating set of $J$. Then the exact sequence in homology is of the form
$$0\longrightarrow H_0(\KK_{\ab}(J')\oplus \KK_{\ab}(\langle  m_s\rangle  ))\longrightarrow H_0(\KK_{\ab}(J))\longrightarrow 0$$
and hence the result. $\square$
\begin{Lemma}\label{nonzero}Let $i>0$, and $\ab\in\NN^n$ such that $\beta_{i,\ab}(J)\neq0$; then $\beta_{i-1,\ab}(\tilde J)\neq0$ or $\beta_{i,\ab}(J')\neq0$
\end{Lemma}

\noindent {\bf Proof: } Let $\ab\in \NN_0^n$ such that $\beta_{i-1,\ab}(\tilde J)=0$ and $\beta_{i,\ab}(J')=0$ then the exact sequence in homology

$$\cdots\longrightarrow H_{i,\ab}(\KK(J')){\longrightarrow}H_{i,\ab}(\KK(J))\stackrel{\Delta}\longrightarrow H_{i-1,\ab}(\KK(\tilde J))\longrightarrow \cdots$$

is of the form

$$0\longrightarrow H_{i,\ab}(\KK(J))\longrightarrow 0$$

and then $H_{i,\ab}(\KK(J))=0$. $\square$

Bringing together these two lemmas, we have the following 
\begin{Proposition}\label{necesary-nodes}
If $H_{i,\ab}(\KK(I))\neq 0$ for some $i$, then $x^{\ab}$ is a generator of some node $J$ in any Mayer-Vietoris tree $MVT(I)$.
\end{Proposition}

From this proposition we have that all the multidegrees of Koszul generators (equivalently Betti numbers) of $I$ appear in $MVT(I)$. For a sufficient condition, we need the following notation: among the nodes in $MVT(I)$ we call {\it relevant nodes} those in even position or in position 1. This is because for all other nodes in odd position, the corresponding generators have already appeared in a {\it relevant node} thus, no new multidegree appears in the non-relevant nodes. Therefore, we can restrict to relevant nodes $J$ in proposition \ref{necesary-nodes}. 

\begin{Proposition}
If $x^{\ab}$ appears only once as a generator of a relevant node $J$ in $MVT(I)$ then there exists exactly one generator in $H_*(\KK(I))$ which has multidegree $\ab$.
\end{Proposition}

\noindent{\bf Proof: }
Let $x^\ab$ be a generator of $I$ (i.e. of the relevant node in position 1). Then $x^\ab$ appears only in this relevant node in $MVT(I)$ and it is clear that there is a generator of $H_0(\KK(I))$ that has multidegree $\ab$, namely $x^\ab$ itself.

If $x^{\ab}$ appears in the tree only as a generator of the node $J$ in even position $p$, then there exists $L\in MVT(I)$ such that $J=\tilde L$. Then we have the following exact sequence in homology 

$$\cdots\longrightarrow H_{1,\ab}(\KK(L'))\longrightarrow H_{1,\ab}(\KK(L))\longrightarrow H_{0,\ab}(\KK(J))\longrightarrow  H_{0,\ab}(\KK(L'))\longrightarrow\cdots$$
but $ H_{0,\ab}(\KK(L'))=0$, since the set of monomials generating $L'$ is included in the set of generators of $L$, and we know that $\ab$ is not the exponent of any of them. On the other hand, $ H_{1,\ab}(\KK(L'))=0$ because if it was different from $0$, then either $(L')'$ or $\widetilde{L'}$ would have $x^\ab$ among their generators (see lemma \ref{nonzero}), but it is impossible for $(L')'$ since its set of generators is included in the set of generators of $L$; and it would be a contradiction for $\widetilde{L'}$, since it is a relevant node. Therefore, we have
$$0\longrightarrow H_{1,\ab}(\KK(L))\longrightarrow H_{0,\ab}(\KK(J))\longrightarrow  0$$

so we have that the connecting morphism is an isomorphism and we have exactly one generator in $H_1(\KK(L))$ with multidegree $\ab$. When iterating the process to compute $H(\KK(I))$ we will always be in one of the following situations:
\begin{itemize}
\item If there exists $M\in MVT(I)$ such that $M'=L$ then we have the following exact sequence
$$\cdots \longrightarrow H_{1,\ab}(\KK(\tilde M))\longrightarrow H_{1,\ab}(\KK(L)\oplus \KK(\langle  m_s\rangle  ))\longrightarrow H_{1,\ab}(\KK(M))\longrightarrow H_{0,\ab}(\KK(\tilde M))\longrightarrow 0$$ but in this case $\tilde M$ is a relevant node, and then the rightmost element of the sequence is just a $0$; moreover, if $H_{1,\ab}(\KK(\tilde M))\neq 0$ then there is some relevant node in the tree hanging from it that has $x^\ab$ as a generator. Since the relevant nodes of this subtree are relevant nodes in the original one, then the above sequence is just $0 \longrightarrow H_{1,\ab}(\KK(L)\oplus \KK(\langle  m_s\rangle  ))\longrightarrow H_{1,\ab}(\KK(M))\longrightarrow  0$

\item If there exists $M\in MVT(I)$ such that $\tilde M=L$ then we have the following exact sequence
$$\cdots \longrightarrow  H_{2,\ab}(\KK(M')\oplus \KK(\langle  m_s\rangle  ))\longrightarrow H_{2,\ab}(\KK(M))\longrightarrow H_{1,\ab}(\KK(L))\longrightarrow H_{1,\ab}(\KK(M')\oplus \KK(\langle  m_s\rangle  ))$$ but in this case the sequence is just 
$$0 \longrightarrow H_{2,\ab}(\KK(M))\longrightarrow H_{1,\ab}(\KK(L))\longrightarrow 0$$
because if $H_{1,\ab}(\KK(M')\oplus \KK(\langle  m_s\rangle  ))$ is different from zero, then there is some relevant node in the subtree hanging from $M'$ such that $x^\ab$ appears in it as a generator, and the relevant nodes of this subtree are also relevant for the bigger tree.
\end{itemize}
 
Therefore, from this iterative process we have that there is one generator in some $H_i(\KK(I))$ in multidegree $\ab$. $\square$

We know that the only multidegrees relevant to the homology computation of $I$ are those in position $1$ of $MVT(I)$, from which we obtain $H_0(\KK(I))$ and those in even position in $MVT(I)$. The degree of the homology to which they contribute, can also be read from their position in the tree. Assign a {\it dimension} to every node in $MVT(I)$ in the following way: $dim(MVT_1(I))=0$; and if $dim(MVT_p(I))=d$ then $dim(MVT_{2p}(I))=d+1$, $dim(MVT_{2p+1}(I))=d$. Note that the dimension $d$ of a node in dimension $p$ is just the number of zeros in the binary expansion of $p$. The generators of each relevant node contribute to the homology modules in the homological degree given by the dimension of the node. To verify this just consider the Mayer-Vietoris sequence at each multidegree in the tree, and the fact that at the bottoms of the tree, one has always a node with two generators and two children with one generator each, the right child contributes with its homology just in degree zero (it is a generator of its father) and the left child contributes with one generator in degree one, the recursive construction of the tree and the sequences yield the correspondence between degrees.

\begin{Example}\label{mvt-examples-2}
Let us consider the ideal $I=\langle xy^2,xyz^3,y^5,z^6\rangle\subseteq\kb[x,y,z]$. A Mayer-Vietoris tree of this ideal is shown in the picture:
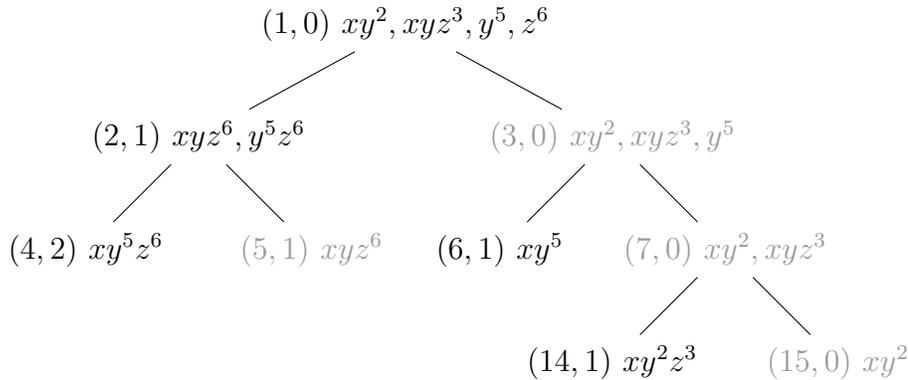
\begin{center}
\begin{tikzpicture}[scale=1]
\tikzstyle{level 1}=[sibling distance=5.5cm]
\tikzstyle{level 2}=[sibling distance=3cm]
\node{$(1,0)$ $xy^2,xyz^3,y^5,z^6$}
child{ node{$(2,1)$ $xyz^6,y^5z^6$}
	child{node{$(4,2)$ $xy^5z^6$}}
	child{node[color=black!40!white]{$(5,1)$ $xyz^6$}}}
child{ node[color=black!40!white]{$(3,0)$ $xy^2,xyz^3,y^5$}
	child{node{$(6,1)$ $xy^5$}}
	child{node[color=black!40!white]{$(7,0)$ $xy^2,xyz^3$}
		child{node{$(14,1)$ $xy^2z^3$}}
		child{node[color=black!40!white]{$(15,0)$ $xy^2$}}}};
\end{tikzpicture}
\end{center} 
Every node is given by a triple {\tt(position,\, dimension)\, ideal} and the relevant nodes are the ones in strong black color. Observe that this tree has no repeated multidegree in the relevant nodes, therefore the multigraded Betti numbers of $I$ are just read from the tree. In this case we have $\beta_0(I)=4$, $\beta_1(I)=4$ and $\beta_2(I)=1$.
\end{Example}

\subsection{Mayer-Vietoris trees and homological computations}\label{computations}
\spanishsubsection{\'Arboles de Mayer-Vietoris y c\'alculos homol\'ogicos}

Now that we have described the Mayer-Vietoris trees, we would like to use them as a tool to efficiently read homological information about monomial ideals. Depending on the goals, the work with the tree will be different. We briefly describe here how to obtain minimal free resolutions from the Mayer-Vietoris trees, and how to use them to obtain the Koszul homology of the ideal.

\subsubsection{Minimal free resolutions}
\spanishsection{Resoluciones libres m\'inimas}

The main object expressing the homological structure of a monomial ideal is its minimal free resolution. As we have seen in sections \ref{resolutions} and \ref{examples-resolutions}, the explicit computation of minimal free resolutions is not a trivial task, and many different approaches have been given in the literature.
This section explores the way between the Mayer-Vietoris tree and the computation of free resolutions of monomial ideal. The question is whether we can explicitly compute a resolution of $I$ based on the data given by $MVT(I)$. In section \ref{m-v} we have seen that using Myer-Vietoris sequences and iterated mapping cones we can compute a (nonminimal) resolution of any monomial ideal. Recalling the construction of that section, we note that the mapping cone of the resolutions of the two children of a given node in $MVT(I)$ is equivalent to a resolution of their father. In particular, the multidegrees of the multigraded pieces of  the resolution of the father are given by the multidegrees of the multigraded pieces of its children's resolutions:

\begin{Lemma}
Let $J$ be a node of $MVT(I)$ and $m_s$ its pivot monomial, $(\Pc_{J'},d_{J'})$ and $(\Pc_{\tilde{J}},d_{\tilde{J}})$ resolutions of its children, then there is a resolution $(\Pc_{J'},d_J')$ of $J$ such that
$$(\Pc_J)_0=(\Pc_{J'})_0\oplus R(-md(m_s))\quad(\Pc_J)_i=(\Pc_{J'})_i\oplus (\Pc_{\tilde{J}})_{i-1}$$
\end{Lemma}

This lemma is an easy consequence of the mapping cone construction given in theorem \ref{ses3} applied to Mayer-Vietoris sequences. As an immediate outcome of this lemma we have the following result:

\begin{Theorem}\label{MVRes}
For every monomial ideal $I$ and every Mayer-Vietoris tree of it $MVT(I)$ there is a resolution of $I$ supported on the relevant nodes of $MVT(I)$. We call this resolution a \emph{Mayer-Vietoris resolution of $I$}. 
\end{Theorem}

In general, Mayer-Vietoris resolutions are nonminimal. Of course, if there are no repeated multidegrees in the generators of the relevant nodes of the tree, the resulting resolution is minimal, but it is not the only case, as will be seen in the next section. Moreover, coming from iterated mapping cones, Mayer-Vietoris resolutions are not difficult to minimize. Recall  from remark \ref{remark-mapping-cones} that if we keep minimality at each step, we know that the only possible part of the matrix which can be reduced is that corresponding to the inclusion, and the minimization process is improved. Moreover, as we keep track of the multidegrees involved, only when the same multidegree appears in the resolutions of both $\tilde J$ and $J'\oplus\langle m_s\rangle$ at the same homological degree we can have some non-minimality on the resolution of $J$.

We illustrate the construction of a Mayer-Vietoris resolution with a detailed example, to show that all the maps in the above constructions have explicit expressions:

\begin{Example}
Let us consider the ideal $I=\langle x^2,y^2,xy\rangle\subseteq R=\kb[x,y]$. We will use the following Mayer-Vietoris tree to build a resolution:
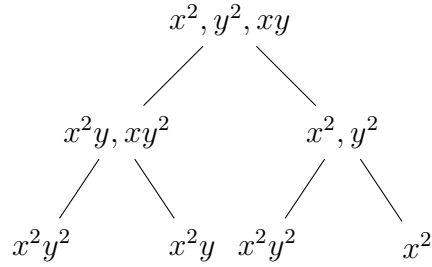
\begin{center}
\begin{tikzpicture}[scale=1]
\tikzstyle{level 1}=[sibling distance=3cm]
\tikzstyle{level 2}=[sibling distance=2cm]
\node{$x^2,y^2,xy$}
child{ node {$x^2y,xy^2$}
	child{node{$x^2y^2$}}
	child{node{$x^2y$}}}
child{ node {$x^2,y^2$}
	child{node{$x^2y^2$}}
	child{node{$x^2$}}};
\end{tikzpicture}
\end{center}
Our starting effective short exact sequence of modules is given by
\begin{center}
\begin{tikzpicture}
\draw[->] (0,0)node[left]{$0$}--(1,0)node[right]{$\langle x^2y^2\rangle$};
\draw (4.5,0) node{$\langle x^2\rangle \oplus \langle y^2\rangle$};
\draw[->] (8.5,0)node[left]{$\langle x^2,y^2\rangle$}--(9.5,0) node[right]{$0$};
\draw[->] (2.3,0.1)--(3.5,0.1)node[above,pos=0.5]{$\mi$};
\draw[->] (3.5,-0.1)--(2.3,-0.1)node[below,pos=0.5]{$\rho$};
\draw[->] (5.5,0.1)--(7,0.1)node[above,pos=0.5]{$\mj$};
\draw[->] (7,-0.1)--(5.5,-0.1)node[below,pos=0.5]{$\sigma$};
\end{tikzpicture}
\end{center}
which gives rise to an effective short exact sequence of chain complexes (only the part we are interested in is depicted):

\begin{center}
\begin{tikzpicture}[scale=1.1]
\draw[->] (-0.5,0)node[left]{$0$}--(0.5,0)node[right]{$\langle x^2y^2\rangle$};
\draw (6.7,0) node{$\langle x^2\rangle\quad\, \oplus\quad\, \langle y^2\rangle$};
\draw[->] (11.5,0)node[left]{$\langle x^2,y^2\rangle$}--(12.5,0) node[right]{$0$};
\draw[->] (2.5,0.1)--(4.8,0.1)node[above,pos=0.5]{$\mi$};
\draw[->] (4.8,-0.1)--(2.5,-0.1)node[below,pos=0.5]{$\rho$};
\draw[->] (8.5,0.1)--(10,0.1)node[above,pos=0.5]{$\mj$};
\draw[->] (10,-0.1)--(8.5,-0.1)node[below,pos=0.5]{$\sigma$};

\draw(1.2,-2)node{$R(-22)$};
\draw (6.7,-2) node{$R(-20) \,\,\oplus\,\, R(-02)$};
\draw[->] (2.3,-1.9)--(4.7,-1.9)node[above,pos=0.5]{$\mi_1$};
\draw[->] (4.7,-2.1)--(2.3,-2.1)node[below,pos=0.5]{$\rho_1$};

\draw[->] (1.3,-0.4)--(1.3,-1.5)node[left,pos=0.5]{$d_1$};
\draw[->] (1.5,-1.5)--(1.5,-0.4)node[right,pos=0.5]{$h_1$};

\draw[->] (5.7,-0.4)--(5.7,-1.5)node[left,pos=0.5]{$d'_1$};
\draw[->] (5.9,-1.5)--(5.9,-0.4)node[right,pos=0.5]{$h'_1$};
\draw[->] (7.7,-0.4)--(7.7,-1.5)node[left,pos=0.5]{$d''_1$};
\draw[->] (7.9,-1.5)--(7.9,-0.4)node[right,pos=0.5]{$h''_1$};

\draw(1.4,-3.9)node{$0$};
\draw (5.7,-3.9) node{$0$};
\draw (7.7,-3.9) node{$0$};

\draw[->] (1.4,-3.5)--(1.4,-2.4);
\draw[->] (5.7,-3.5)--(5.7,-2.4);
\draw[->] (7.7,-3.5)--(7.7,-2.4);

\end{tikzpicture}
\end{center}

The relevant maps $\mi$ and $\rho$ are given as in proposition \ref{mv}. $d,d'$ and $d''$ are the evident maps in the minimal free resolutions of the corresponding ideals, and $h,h'$ and $h''$ are their contracting homotopies. Since Taylor resolution is minimal for ideals generated by one monomial, we can use Fr\"obergs \cite{F78} contracting homotopy. Then, the following maps are induced (assume the module $R(-22)$ in the resolution of $\langle x^2y^2\rangle$ is generated by $\eb_1$ and $R(-20)$, $R(-02)$ are generated respectively by $\eb'_1$ and $\eb'_2$)
$$\mi_1(\eb_1)=h\mi d(\eb_1)=h\mi(x^2y^2)=h(x^2y^2,-x^2y^2)=(y^2\eb'_1,-x^2\eb''_1)$$
$$\rho_1=h\rho(d'\oplus d'')$$
Observe that we have explicit expressions or formulas for every arrow in the diagram, therefore the \emph{effective} cone can be explicitly computed. The cone of $\mi$ is then given by
\begin{center}
\begin{tikzpicture}[scale=1.1]
\draw[->] (0,0)node[left]{$0$}--(1,0)node[right]{$R(-22)$};
\draw (6.5,0) node{$R(-20) \oplus\RR(-02)\oplus R$};
\draw[->] (11,0)node[left]{$R\oplus R$}--(12,0) node[right]{$0$};
\draw[->] (2.6,0.1)--(4.5,0.1)node[above,pos=0.5]{{\scriptsize $\left ( \begin{array}{c} y^2\\-x^2\\-x^2y^2 \end{array}\right )$}};
\draw[->] (4.5,-0.1)--(2.6,-0.1)node[below,pos=0.5]{$\rho_1$};
\draw[->] (8.5,0.1)--(9.7,0.1)node[above,pos=0.5]{{\scriptsize $\left ( \begin{array}{ccc} x^2&0&1\\0&y^2&-1
 \end{array}\right )$}};
\draw[->] (9.7,-0.1)--(8.5,-0.1)node[below,pos=0.5]{$\rho$};
\end{tikzpicture}
\end{center}
It is clear that the last generator of the central module can be always eliminated, form the way we have constructed the cone, therefore, the last column of the rightmost matrix and the last row of the leftmost matrix disappear. Moreover, the equivalence between the resulting complex and the resolution of $\langle x^2,y^2\rangle$ is made explicit in a ``contraction'' of the last matrix into one row; note that this can also always be performed for it comes from the construction of our sequence, and will only happen at the end of the resolution. Therefore, after this first step we have the minimal free resolution of the node $MVT(I)_3$ from the resolutions of its children:

\begin{center}
\begin{tikzpicture}
\draw[->] (0,0)node[left]{$0$}--(1,0)node[right]{$R(-22)$};
\draw (6.5,0) node{$R(-20) \oplus R(-02)$};
\draw[->] (11,0)node[left]{$R$}--(12,0) node[right]{$0$};
\draw[->] (2.6,0.1)--(4.5,0.1)node[above,pos=0.5]{{\scriptsize $\left ( \begin{array}{c} y^2\\-x^2\end{array}\right )$}};
\draw[->] (4.5,-0.1)--(2.6,-0.1)node[below,pos=0.5]{$\rho_1$};
\draw[->] (8.5,0.1)--(9.7,0.1)node[above,pos=0.5]{{\scriptsize $\left ( \begin{array}{cc} x^2&y^2
 \end{array}\right )$}};
\draw[->] (9.7,-0.1)--(8.5,-0.1)node[below,pos=0.5]{$\bar\rho$};
\end{tikzpicture}
\end{center}

Exactly in the same way we produce the minimal free resolution of $MVT(I)_2$ from its children, obtaining
\begin{center}
\begin{tikzpicture}
\draw[->] (0,0)node[left]{$0$}--(1,0)node[right]{$R(-22)$};
\draw (6.5,0) node{$R(-21) \oplus R(-12)$};
\draw[->] (11,0)node[left]{$R$}--(12,0) node[right]{$0$};
\draw[->] (2.6,0.1)--(4.5,0.1)node[above,pos=0.5]{{\scriptsize $\left ( \begin{array}{c} y\\-x\end{array}\right )$}};
\draw[->] (4.5,-0.1)--(2.6,-0.1)node[below,pos=0.5]{$\rho_1$};
\draw[->] (8.5,0.1)--(9.7,0.1)node[above,pos=0.5]{{\scriptsize $\left ( \begin{array}{cc} x^2y&xy^2
 \end{array}\right )$}};
\draw[->] (9.7,-0.1)--(8.5,-0.1)node[below,pos=0.5]{$\bar\rho$};
\end{tikzpicture}
\end{center}

To obtain the minimal free resolution of $I=MVT(I)$ we just plug the obtained resolutions into the same schema and proceed as before:
\begin{center}
\begin{tikzpicture}[scale=1.1]
\draw[->] (-0.5,0)node[left]{$0$}--(0.5,0)node[right]{$\langle x^2y,xy^2\rangle$};
\draw (6.7,0) node{$\langle x^2,y^2\rangle\quad\, \oplus\quad\, \langle xy\rangle$};
\draw[->] (11.5,0)node[left]{$\langle x^2,y^2,xy\rangle$}--(12.5,0) node[right]{$0$};
\draw[->] (2.5,0.1)--(4.8,0.1)node[above,pos=0.5]{$\mi$};
\draw[->] (4.8,-0.1)--(2.5,-0.1)node[below,pos=0.5]{$\rho$};
\draw[->] (8.5,0.1)--(9.5,0.1)node[above,pos=0.5]{$\mj$};
\draw[->] (9.5,-0.1)--(8.5,-0.1)node[below,pos=0.5]{$\sigma$};

\draw(1.5,-2)node{$R(-21)\oplus R(-12)$};
\draw (6.5,-2) node{$R(-20) \oplus R(-02)\oplus R(-11)$};
\draw[->] (3.3,-1.9)--(4,-1.9)node[above,pos=0.5]{$\mi_1$};
\draw[->] (4,-2.1)--(3.3,-2.1)node[below,pos=0.5]{$\rho_1$};

\draw[->] (1.6,-0.4)--(1.6,-1.5)node[left,pos=0.5]{$d_1$};
\draw[->] (1.8,-1.5)--(1.8,-0.4)node[right,pos=0.5]{$h_1$};

\draw[->] (5.5,-0.4)--(5.5,-1.5)node[left,pos=0.5]{$d'_1$};
\draw[->] (5.7,-1.5)--(5.7,-0.4)node[right,pos=0.5]{$h'_1$};

\draw[->] (8,-0.4)--(8,-1.5)node[left,pos=0.5]{$d''_1$};
\draw[->] (8.2,-1.5)--(8.2,-0.4)node[right,pos=0.5]{$h''_1$};

\draw(1.5,-3.9)node{$R(-22)$};
\draw (5.5,-3.9) node{$R(-22)$};\draw (8.2,-3.9) node{$0$};
\draw[->] (2.7,-3.8)--(4.5,-3.8)node[above,pos=0.5]{$\mi_2$};
\draw[->] (4.5,-4)--(2.7,-4)node[below,pos=0.5]{$\rho_2$};

\draw[->] (1.6,-3.5)--(1.6,-2.4)node[left,pos=0.5]{$d_2$};
\draw[->] (1.8,-2.4)--(1.8,-3.5)node[right,pos=0.5]{$h_2$};
\draw[->] (5.5,-3.5)--(5.5,-2.4)node[left,pos=0.5]{$d'_2$};
\draw[->] (5.7,-2.4)--(5.7,-3.5)node[right,pos=0.5]{$h'_2$};;
\draw[->] (8.2,-3.5)--(8.2,-2.4);

\draw(1.5,-5.8)node{$0$};
\draw (5.5,-5.8) node{$0$};
\draw[->] (1.5,-5.4)--(1.5,-4.3);
\draw[->] (5.5,-5.4)--(5.5,-4.3);

\end{tikzpicture}
\end{center}
Observe that again all arrows are explicit, and therefore we can completely compute the mapping cone of $\mi$:

\begin{center}
\begin{tikzpicture}[xscale=0.75,yscale=0.8, transform shape]
\draw[->] (0,0)node[left]{$0$}--(1,0)node[right]{$R(-22)$};
\draw (6.5,0) node{$R(-22) \oplus R(-21)\oplus R(-12)$};
\draw(17,0)node[left]{$R(-20)\oplus R(-02)\oplus R(-11)\oplus R$};
\draw[->] (20,0)node[left]{$R\oplus R$}--(21,0) node[right]{$0$};

\draw[->] (2.6,0.1)--(3.8,0.1)node[above,pos=0.5]{{\scriptsize $\left ( \begin{array}{c} 1\\y\\-x\end{array}\right )$}};
\draw[->] (3.8,-0.1)--(2.6,-0.1)node[below,pos=0.5]{$\rho_2$};
\draw[->] (9.1,0.1)--(10.8,0.1)node[above,pos=0.5]{{\scriptsize $\left ( \begin{array}{ccc} y^2&y&0\\-x^2&0&x\\0&-x&-y\\0&x^2y&xy^2
 \end{array}\right )$}};
\draw[->] (10.8,-0.1)--(9.1,-0.1)node[below,pos=0.5]{$\rho_1$};
\draw[->] (17,0.1)--(18.7,0.1)node[above,pos=0.5]{{\scriptsize $\left ( \begin{array}{cccc} x^2&y^2&0&1\\
0&0&xy&-1
 \end{array}\right )$}};
\draw[->] (18.7,-0.1)--(17,-0.1)node[below,pos=0.5]{$\rho$};
\end{tikzpicture}
\end{center}
After performing the same operations as in the first iterations, we arrive into the Mayer-Vietoris resolution of $I$ corresponding to the used tree:
\begin{center}
\begin{tikzpicture}[xscale=0.75,yscale=0.8, transform shape]
\draw[->] (0,0)node[left]{$0$}--(1,0)node[right]{$R(-22)$};
\draw (6.5,0) node{$R(-22) \oplus R(-21)\oplus R(-12)$};
\draw(16,0)node[left]{$R(-20)\oplus R(-02)\oplus R(-11)$};
\draw[->] (18.5,0)node[left]{$R$}--(19.5,0) node[right]{$0$};

\draw[->] (2.6,0.1)--(3.8,0.1)node[above,pos=0.5]{{\scriptsize $\left ( \begin{array}{c} 1\\y\\-x\end{array}\right )$}};
\draw[->] (3.8,-0.1)--(2.6,-0.1)node[below,pos=0.5]{$\rho_2$};
\draw[->] (9.2,0.1)--(10.7,0.1)node[above,pos=0.5]{{\scriptsize $\left ( \begin{array}{ccc} y^2&y&0\\-x^2&0&x\\0&-x&-y
 \end{array}\right )$}};
\draw[->] (10.7,-0.1)--(9.2,-0.1)node[below,pos=0.5]{$\rho_1$};
\draw[->] (16,0.1)--(17.7,0.1)node[above,pos=0.5]{{\scriptsize $\left ( \begin{array}{cccc} x^2&y^2&xy
 \end{array}\right )$}};
\draw[->] (17.7,-0.1)--(16,-0.1)node[below,pos=0.5]{$\rho$};
\end{tikzpicture}
\end{center}
Observe that the scalar $1$ at the top of the leftmost matrix indicates that this is not a minimal resolution. The existence of such element was pointed by the map $\mi_2$ between modules generated at the same multidegree, namely the multidegree that is repeated in the relevant nodes of our tree, see remark \ref{remark-mapping-cones}.
\end{Example}

\subsubsection{Koszul homology}
\spanishsubsubsection{Homolog\'ia de Koszul}

Another important issue is the computation of the Koszul homology of the ideal $I$. This means not only obtaining the ranks of the homology modules, i.e. the (multigraded) Betti numbers which has been already seen in the preceding paragraphs, but giving explicit sets of generators for them. For those multidegrees that are not repeated in the relevant nodes of the tree, the generating sets can be obtained automatically from the Mayer-Vietoris tree. For the others, simplicial techniques are an alternative.

\subsubsection*{Non-repeated relevant multidegrees}
\spanishsubsubsection{Multigrados relevantes no repetidos}

First of all, recall the trivial case is that of the node in position one, i.e.generators of $I$, here
$$H_i(\KK(\langle m_1\rangle))\simeq
\begin{cases}
\kb&i=0,\\
0&otherwise.\\
\end{cases}$$
and the generator of $H_0(\KK(\langle m_1\rangle))$ can be identified with $m_1$ itself.

Now we can begin with a recursive process starting with the node in which the multidegree we are interested in appears. The ingredient we need for this recursive process is a explicit formula to obtain preimages of the connecting morphism $\Delta$. This procedure has been described in section \ref{mvs-koszul}. In the case of non repeated multidegrees in the relevant nodes we know that the connecting morphism $\Delta$ is an isomorphism. Therefore, the preimages obtained of generators of the relevant node are generators of the corresponding homology groups of the father. Therefore, we start from the leaves of the tree and apply the formulas in section \ref{mvs-koszul} recursively to obtain an explicit generator in some particular multidegree of the corresponding homology module of $I$, the root node.

We illustrate this procedure with the following example:

\begin{Example}
Let us take the ideal and tree from example \ref{mvt-examples-2}, for which we know that the tree had no repeated multidegrees in the relevant nodes.
Let us start with $MVT(I)_4=\langle xy^5z^6\rangle$. The exact sequence we are using is
$$ 0\rightarrow \langle xy^5z^6\rangle\rightarrow \langle xyz^6\rangle\oplus\langle y^5z^6\rangle\rightarrow\langle xyz^6,y^5z^6\rangle\rightarrow 0$$
We have that $\mi(xy^5z^6)=(xy^5z^6,-xy^5z^6)$, now, applying the Spencer and Koszul differential to it, we obtain $\partial\delta(xy^5z^6,-xy^5z^6)=(5xy^5z^6,-xy^5z^6)$ therefore, we consider the element $(xy^4z^6\otimes y,-y^5z^6\otimes x)$ applying $\mj$ we obtain that the class of $\Delta^{-1}(xy^5z^6)=xy^4z^6\otimes y-y^5z^6\otimes x$ is a generator of $H_1(\KK(\langle xyz^6,y^5z^6\rangle))$. We move now to $MVT(I)_1=I$ and consider the connecting morphism between $H_1(\KK(MVT(I)_2))$ and $H_2(\KK(I))$, and apply the same procedure to $xy^4z^6\otimes y-y^5z^6\otimes x$, obtaining the cycle $\frac{-6}{7}(y^4z^6\otimes xy+xy^4z^5\otimes yz-y^5z^5\otimes xz)$, the class of which is a generator of $H_2(\KK(I))$.

Proceeding in the same way with the other generators of the relevant nodes we arrive to the following set of generators of $H_*(\KK(I))$:

$$\begin{array}{lll}
g^0_1=xy^2&g^1_1=y^5z^5\otimes z-y^4z^6\otimes y&g^2_1=y^4z^6\otimes xy+xy^4z^5\otimes yz-y^5z^5\otimes xz\\
g^0_2=xyz^3&g^1_2=xy^4\otimes y-y^5\otimes x&\\
g^0_3=y^5&g^1_3=xy^2z^2\otimes z-xyz^3\otimes y&\\
g^0_4=z^6&g^1_4=2xyz^5\otimes z-yz^6\otimes x-xz^6\otimes y&\\
\end{array}$$
\end{Example}

\subsubsection*{Repeated multidegrees}
\spanishsubsubsection{Multigrados repetidos}

In the case of repeated multidegrees one issue is to know whether they give rise to a generator of the Koszul homology or not. In these cases one can proceed by using the Koszul simplicial complexes of section \ref{simpl-koszul-complex}. Recall the main result concerning these complexes is theorem \ref{simp-kos-hom}, which states an isomorphism between the reduced homology of the (upper) Koszul simplicial complex of $I$ at a multidegree $\ab$ and the Koszul homology of $I$ at that multidegree. This isomorphism is based on the equality between the subjacent chain complexes. This isomorphism can be made explicit: the face $\tau=\{\tau_1,\dots,\tau_j)$ of $\Delta^I_\ab$ corresponds to $\frac{\ab}{\xb^\tau}\otimes x^{\tau_1}\wedge \cdots \wedge x^{\tau_j}$ all differentials are immediately translated from one complex to the other.
\begin{Example}
Consider for example the ideal of example \ref{mvt-examples-2}, and the multidegree $(1,1,6)$, the (upper) Koszul simplicial complex at $(1,1,6)$ is given by $\{\{x\},\{y\},\{z\},\{x,y\}\}$; it has homology at degree $0$ and a generator is given by $[2z-x-y]$, translating this via the isomorphism between the reduced homology of $\Delta^I_{(1,1,6)}$ and $H_{1,(1,1,6)}(\KK(I))$ it corresponds to $[2xyz^5\otimes z-yz^6\otimes x-xz^6\otimes y]$, which is indeed a generator of $H_{1,(1,1,6)}(\KK(I))$.
\end{Example}
 
Computing simplicial homology is in general very expensive, so one should avoid unnecessary computations, i.e. at those multidegrees that are not Betti multidegrees. For this, we have two alternatives: First, compute the minimal resolution using the procedure just seen, and compute reduced homology of the Koszul complex at the repeated multidegrees, the other to be computed using the formulas in the preceding paragraph. The second alternative is valid in cases in which we are only interested in some particular multidegree or if we are detecting whether it has zero homology or not is particularly easy, and consists in directly computing the reduced simplicial homology.
However, for these homology computations we have some advantages that make them easier. First of all, we have coefficients on the field $\kb$ and thus, no torsion will appear, thus simplicial homology computations can be significantly improved. The second advantage is that, once we have the data coming either from the tree or from the minimal resolution, we are not interested in the full reduced homology modules of $\Delta^I_{\ab}$ but only in some given dimension $j$. Then we can use a subcomplex of $\Delta^I_{\ab}$, namely the one that has as facets the $j+1$ faces of it. Moreover, we will only use part of the chain complex of this subcomplex, the one at homological dimension $j+1$ and $j$. This has been discussed in section \ref{simplicial_Koszul_comp}.

These and other techniques to detect homology at a given multidegree, like Koszul ideals (section \ref{stanley-alexander}), will be used in this chapter to perform computations with Mayer-Vietoris trees (see later).

\subsection{Mayer-Vietoris ideals}\label{MV-ideals}
\spanishsubsection{Ideales de Mayer-Vietoris}

Collecting the multidegrees that are not repeated in the relevant nodes of $MVT(I)$, we have a subset of the set of multidegrees $\ab\in\NN^n$ such that $\beta_{i,\ab}\neq 0$ for some $i$. On the other hand, collecting all the multidegrees (repeated and non-repeated) in the relevant nodes, we have a superset. Of course if we only have non-repeated generators in the relevant nodes of $MVT(I)$ then we obtain the exact set of multidegrees with non zero Betti numbers. More explicitly:

Let $I$ be a monomial ideal and $MVT(I)$ a Mayer-Vietoris tree of $I$. Let $\ab\in\NN^n$; let $\overline {\beta_{i,\ab}}(I)=1$ if $\ab$ is the multidegree of a generator of a relevant node of dimension $i$ in $MVT(I)$ which does not appear as a generator of any other relevant node, and $\overline {\beta_{i,\ab}}(I)=0$ in other case. Let $\widehat {\beta_{i,\ab}}(I)$ be the number of times $\ab$ appears as the multidegree of some generator of dimension $i$ in some relevant node in $MVT(I)$. Then for all $\ab\in\NN^n$ we have 
$$\overline{\beta_{i,\ab}}(I)\leq\beta_{i,\ab}(I)\leq\widehat {\beta_{i,\ab}}(I)$$
It is clear that if all the generators of the relevant nodes of $MVT(I)$ are different, we have equalities. These considerations lead us to the following definitions:

\begin{Definition}\label{mv-ideals}
Let $I$ be a monomial ideal.
\begin{itemize}
\item If there exists a Mayer-Vietoris tree of $I$ such that there is no repeated generator in the ideals of the relevant nodes, then we say that $I$ is a {\it Mayer-Vietoris ideal of type A}. In this case, $\overline{\beta_{i,\ab}}(I)=\beta_{i,\ab}(I)=\widehat {\beta_{i,\ab}}(I)$ $\forall i\in\NN,\,\ab\in\NN^n$.
\item If $\overline{\beta_{i,\ab}}(I)=\beta_{i,\ab}(I)$ for all $i\in\NN$, $\ab\in\NN^n$ then we say that $I$ is a {\it Mayer-Vietoris ideal of type B1}.
\item If $\widehat{\beta_{i,\ab}}(I)=\beta_{i,\ab}(I)$ for all $i\in\NN$, $\ab\in\NN^n$ then we say that $I$ is a {\it Mayer-Vietoris ideal of type B2}.
\end{itemize}
\end{Definition}
We will denote the set of Mayer-Vietoris ideals of types $A$, $B1$ and $B2$ as $MV_A$, $MV_{B1}$ and $MV_{B2}$ respectively. Observe that a Mayer-Vietoris ideal of type $A$ is of course Mayer-Vietoris of types $B1$ and $B2$. We give here some examples of each of these three families of Mayer-Vietoris ideals. A more detailed and complete treatment of the following and other examples will be given in section \ref{special-MVT} and chapter \ref{Applications}.
\subsubsection*{Mayer-Vietoris ideals of type A}
\spanishsubsubsection{Ideales de Mayer-Vietoris de tipo A}

Mayer-Vietoris ideals of type $A$ include \emph{Ferrers} ideals (see section \ref{families}) and \emph{consecutive k-out-of-n} ideals (see section \ref{reliability}), which have importance in the reliability theory of coherent systems. These two examples and other Mayer-Vietoris ideals of type $A$ are treated in detail in chapter  \ref{Applications}. This type of ideals have a minimal Mayer-Vietoris resolution, and for some of them, the inspection of the minimal Mayer-Vietoris tree provides actual formulas for their (multigraded) Betti numbers.

A first observation is that every ideal in two variables is Mayer-Vietoris of type $A$. To see this just take an ideal $I\subseteq \kb[x,y]$ and build a Mayer-Vietoris tree using the following strategy: first observe that no two minimal generators of $I$ have the same exponent in any of the variables. Take one of the variables, say $x$, and choose as pivot monomial the generator having biggest exponent in $x$, say $j$. The ideal in the left-hand child is then generated by one generator which has $j$ exponent in $x$ and in $y$ the minimum of the $y$-exponent among the other generators. Go on with this strategy and it is easy to see that there are no repeated generators, since all the generators in dimension $1$ of the tree have different exponents in $x$, and none of them is equal to any generator in dimension $0$, which are just the generators of the ideal.

\subsubsection*{Mayer-Vietoris ideals of type B1}
\spanishsubsubsection{Ideales de Mayer-Vietoris de tipo B1}

This type of ideals includes all that are minimally resolved by their \emph{Scarf complex}, in particular, generic monomial ideals (see section \ref{families}). The Mayer-Vietoris resolution of these ideals is in general non-minimal, but the multigraded Betti multidegrees of them can be immediately read off their Mayer-Vietoris trees.

\subsubsection*{Mayer-Vietoris ideals of type B2}
\spanishsubsubsection{Ideales de Mayer-Vietoris de tipo B2}

The first class of Mayer-Vietoris ideals of type $B2$ are those minimally resolved by the Taylor resolution, as it is easy to see:
\begin {Proposition}
If $I$ is a monomial ideal minimally resolved by its Taylor resolution, then it is a Mayer-Vietoris ideal of type $B2$.
\end{Proposition}
\noindent{\bf Proof:}
If the Taylor resolution of $I$ is minimal, we know that we have a Koszul generator corresponding to every subset of minimal generators of $I$. These are exactly the generators appearing in the relevant nodes of any $MVT(I)$ in this case, since no divisibility will apply in any of the nodes. $\square$

Being minimally resolved by the Taylor resolution is a very restrictive property, however, \cite{F78,B02,HHMT06} gives some criteria for a monomial ideal to have such minimal resolution. There is a much bigger number of Mayer-Vietoris ideal of type $B2$. Among then, we meet another family of ideals corresponding to coherent systems in reliability theory. In this case, we deal with {\it k-out-of-n} systems. These are represented by monomial ideals generated by products of any $k$ variables in the polynomial ring of $n$ variables (see section \ref{reliability}).

Mayer-Vietoris ideals of type $B2$ are also minimally resolved by their Mayer-Vietoris resolution. Some other examples of this type of ideals appear in chapter \ref{Applications}.

\section{Some special Mayer-Vietoris trees}\label{special-MVT}
\spanishsection{Algunos \'arboles de Mayer-Vietoris especiales}
This section is concerned with some particular Mayer-Vietoris trees and ideals, the properties of which will be used in the next sections. We include some characterizations of certain Mayer-Vietoris ideals and introduce some techniques to work with Mayer-Vietoris trees. In section \ref{families}  these techniques will be used to compute Betti numbers of certain ideals, and in section \ref{reliability} we will use these techniques in actual applications.
\subsection{Pure Mayer-Vietoris trees}
\spanishsubsection{\'Arboles de Mayer-Vietoris puros}

\begin{Definition}
Let $I\subseteq R=\kb[x_1\dots,x_n]$ be a graded ideal, and let 
$$ 0\rightarrow \bigoplus^{\beta_l}_{i=1}R(-d_{l_i})\rightarrow\cdots\rightarrow \bigoplus^{\beta_1}_{i=1}R(-d_{1_i})\rightarrow R\rightarrow R/I\rightarrow 0$$
be the minimal graded free resolution of $R/I$. The ring $R/I$ (equivalently the ideal $I$) has a \emph{pure resolution} if there exist integers $d_1,\dots,d_l$ such that
$$d_{1_i}=d_1\,\forall i,\quad \dots,d_{l_i}=d_l\,\forall i$$
If in addition $d_i=d_1+i-1$ for $2\leq i \leq l$ the resolution is said to be $d_1$-linear.

For a pure resolution, the sequence $(d_1,d_2,\dots,d_l)$ is called the \emph{shift type}. The \emph{degree type} of a pure resolution is given by the sequence of differences $(d_l-d_{l-1},\dots,d_2-d_1,d_1-0)$.
\end{Definition}

\begin{Example}
Take $I=\langle xy,yz,zt,xt \rangle$ in $R=\kb[x,y,z,t]$ which is the edge ideal of the square $C_4$, see the figure below. The minimal resolution of $I$ has the form

$$0 \rightarrow R(-4) \rightarrow R^4(-3) \rightarrow R^4(-2)\rightarrow R \rightarrow R/I\rightarrow 0$$

which is $2$-linear.

Let us consider the ideal $I=\langle xy,yz,zt,tu,xu \rangle$ in $R=\kb[x,y,z,t,u]$ which is the edge ideal of the pentagon $C_5$. The minimal resolution of $I$ has the form
$$ 0 \rightarrow R(-5) \rightarrow R^5(-3) \rightarrow R^5(-2) \rightarrow R \rightarrow R/I\rightarrow 0$$
and is therefore not linear, but pure with shift type $(2,3,5)$ and degree type $(2,1,2)$.

If instead we take $I=\langle xy,yz,zt,tu,uv,xv \rangle$ in $R=\kb[x,y,z,t,u,v]$ which is the edge ideal of the hexagon $C_6$. The resolution has the form
$$0 \rightarrow R^2(-6) \rightarrow R^6(-5) \rightarrow R^6(-3)\oplus R^3(-4) \rightarrow R^6(-2) \rightarrow R\rightarrow R/I\rightarrow 0$$
which is not pure.

\begin{center}
$\begin{array}{c@{\hspace{0.5in}}c@{\hspace{0.5in}}c}

 \begin{tikzpicture}[x=1.5cm,y=1.5cm]
\draw(0,0)--(1.5,0)--(1.5,1.5)--(0,1.5)--cycle;

	\foreach \position in {(0,0),(1.5,0),(1.5,1.5),(0,1.5)}
		{
		\fill \position circle(2pt);
		}
	\draw (0,0) node[left] {x};
	\draw (1.5,0) node[right] {y};
	\draw (1.5,1.5) node[right] {z};
	\draw (0,1.5) node[left] {t};
;\end{tikzpicture}

&

 \begin{tikzpicture}[x=1.5cm,y=1.5cm]
\draw (xyz polar cs: angle=72,radius=1) -- (xyz polar cs: angle=144,radius=1)--(xyz polar cs: angle=216,radius=1)
--(xyz polar cs: angle=288,radius=1)--(xyz polar cs: angle=0,radius=1)--(xyz polar cs: angle=72,radius=1);
\foreach \an in {72,144,216,288,0}
		{
		\fill (xyz polar cs: angle=\an,radius=1) circle(2pt);
		}
\draw (xyz polar cs: angle=72,radius=1) node[right] {x};
\draw (xyz polar cs: angle=144,radius=1) node[left] {y};
\draw (xyz polar cs: angle=216,radius=1) node[left] {z};
\draw (xyz polar cs: angle=288,radius=1) node[below] {t};
\draw (xyz polar cs: angle=0,radius=1) node[right] {u};
\end{tikzpicture}

&
 \begin{tikzpicture}[x=1.5cm,y=1.5cm]
\draw (xyz polar cs: angle=60,radius=1) -- (xyz polar cs: angle=120,radius=1)--(xyz polar cs: angle=180,radius=1)
--(xyz polar cs: angle=240,radius=1)--(xyz polar cs: angle=300,radius=1)--(xyz polar cs: angle=0,radius=1)--(xyz polar cs: angle=60,radius=1);
\foreach \an in {0,60,120,180,240,300}
		{
		\fill (xyz polar cs: angle=\an,radius=1) circle(2pt);
		}
\draw (xyz polar cs: angle=60,radius=1) node[above] {x};
\draw (xyz polar cs: angle=120,radius=1) node[above] {y};
\draw (xyz polar cs: angle=180,radius=1) node[left] {z};
\draw (xyz polar cs: angle=240,radius=1) node[below] {t};
\draw (xyz polar cs: angle=300,radius=1) node[below] {u};
\draw (xyz polar cs: angle=0,radius=1) node[right] {v};
\end{tikzpicture}

 \\ [0.4cm]
C_4 & C_5 & C_6
\end{array}$
\end{center}

\end{Example}

\begin{Remark}
Observe that the Betti diagram of an ideal with a linear resolution consists only on one row, namely $d_1$, which gives the regularity of the ideal. In the case of pure resolutions, the Betti diagram might have more rows, but each column of the diagram consists only of one number.
\end{Remark}

In a similar way, we introduce \emph{pure} Mayer-Vietoris trees, and show that any ideal having a pure Mayer-Vietoris tree is a Mayer-Vietoris ideal of type $B2$ and hence its corresponding Mayer-Vietoris resolution is minimal. In chapter \ref{Applications} we will give some families of ideals which have unmixed and pure Mayer-Vietoris trees:

\begin{Definition}
Let $I$ be a monomial ideal and $MVT(I)$ a Mayer-Vietoris ideal of it. Let $x^\mu$ and $x^\nu$ two generators of some relevant nodes $MVT_p(I)$ and $MVT_q(I)$ respectively; let $d_p$ and $d_q$ be the dimensions of these nodes in $MVT(I)$. If $deg(x^\mu)\neq deg(x^\nu)$ whenever $d_p\neq d_q$ then we say that $MVT(I)$ is an \emph{unmixed} Mayer-Vietoris tree.
If for each dimension $d$ there is an integer $g_d$ such that $deg(x^\mu)=g_d$ for every generator of a node of dimension $d$, then we say that the Mayer-Vietoris tree is \emph{pure}. If in addition $g_d=d+g_0$ the we say that the Mayer-Vietoris tree is \emph{$g_0$-linear}.
For a pure tree, the sequence $(g_0,g_1,\dots)$ is called the \emph{shift type} of the tree. The \emph{degree type} of a pure tree is given by the sequence of differences $(\dots,g_1-g_0,g_0-0)$.
\end{Definition}

\begin{Proposition}\label{repeated-relevant-multidegrees}
Let $\mu\in \NN^n$ a multidegree such that all $k$ relevant nodes in $MVT(I)$ of which $\xb^\mu$ is a minimal generator have the same dimension $i$ then $dim(H_{i,\mu}(\KK(I))=k$.
\end{Proposition}

\noindent{\bf Proof: }
Let $x^\mu$ be a generator of two relevant nodes of the tree $MVT(I)$ with the conditions of the statement. Let $p$ be the position of a common ancestor of these nodes and $d_p$ its dimension. Let call $J_p$ the ideal in $MVT_p(I)$. Then we have the following part of the exact sequence in Koszul homology of $J_p$:

\begin{equation*}
\begin{split}
\cdots\longrightarrow H_{i-d_p+1,\mu}(K(J_p))\stackrel{\Delta}{\longrightarrow} H_{i-d_p,\mu}(K(\tilde J_p))\longrightarrow H_{i-d_p,\mu}(K(J_p'))\\
\longrightarrow H_{i-d_p,\mu}(K(J_p))\stackrel{\Delta}{\longrightarrow}H_{i-d_p-1,\mu}(K(\tilde J_p))\longrightarrow\cdots
\end{split}
\end{equation*}

Since $x^\mu$ is a generator of relevant nodes only in dimension $i$, the left most and rightmost modules are $0$, and then we have that $$dim(H_{i-d_p,\mu}(K(J_p)))=dim(H_{i-d_p,\mu}(K(J'_p)))+dim(H_{i-d_p,\mu}(K(\tilde J_p)))$$

this is true for any common ancestor of the nodes in which $x^\mu$ is a generator, in particular for $I$, the only relevant node of dimension $0$. $\square$

\begin{Corollary}\label{unmixed}

If $I$ has an unmixed Mayer-Vietoris tree, then $I$ is a Mayer-Vietoris ideal of type $B2$.
\end{Corollary}

\begin{Corollary}
Every ideal having a pure or linear Mayer-Vietoris tree is Mayer-Vietoris of type $B2$.
\end{Corollary}

\begin{Example}
An unmixed Mayer-Vietoris tree of the edge ideal of the square $C_4$ is given by
 
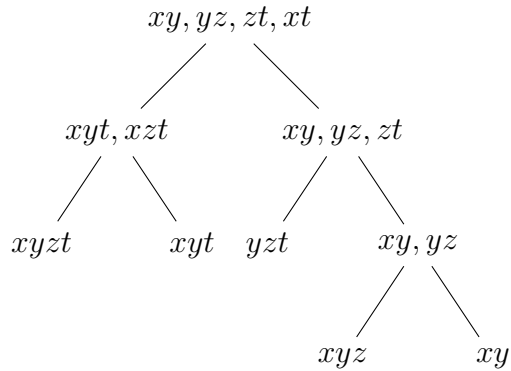
\begin{center}
\begin{tikzpicture}[scale=1]
\tikzstyle{level 1}=[sibling distance=3cm]
\tikzstyle{level 2}=[sibling distance=2cm]
\node{$xy,yz,zt,xt$}
child{ node {$xyt,xzt$}
	child{node{$xyzt$}}
	child{node{$xyt$}}}
child{ node {$xy,yz,zt$}
	child{node{$yzt$}}
	child{node{$xy,yz$}
		child{node{$xyz$}}
		child{node{$xy$}}}};
\end{tikzpicture} 
\end{center}
Although the edge ideal of the pentagon $C_5$ has a pure resolution, it has no unmixed Mayer-Vietoris tree; this ideal is not even Mayer-Vietoris of any type.
\end{Example}

\subsection{Separable Mayer-Vietoris trees}\label{separable}
\spanishsubsection{\'Arboles de Mayer-Vietoris separables}

\begin{Definition}
We say that a monomial ideal $I\subseteq\kb[x_1,\dots,x_n]$ is \emph{separable by addition} if there exist two subsets $A=\{x_{A_1},\dots,x_{A_k}\}$ and $B=\{x_{B_1},\dots,x_{B_l}\}$, $A\sqcup B=\{x_1,\dots,x_n\}$; and two ideals $J\subseteq\kb[A]$, $K\subseteq\kb[B]$ such that $I=J+K$. If $I=J\cap K$ we say that $I$ is \emph{separable by intersection}.
\end{Definition}

In the next paragraphs we will use Mayer-Vietoris trees to compute the (multigraded) Betti numbers of a separable ideal $I$ from the Betti numbers of its pieces $J$ and $K$, provided at least one of this pieces is a Mayer-Vietoris ideal of type $A$ or $B2$. For this we will use a certain strategy for constructing the Mayer-Vietoris tree of $I$. We start with ideals that are separable by addition.

Let $I=J+K$ a separable ideal such that $K$ is Mayer-Vietoris of type $A$ or $B2$. In order to obtain the Betti numbers of $I$ using the Betti numbers of $J\subseteq\kb[A]$ and $K\subseteq\kb[B]$ we will construct a particular Mayer-Vietoris tree of $I$ which we will call a \emph{separable} Mayer-Vietoris tree. The construction of this tree is as follows:
Since $K$ is of type $A$ or $B2$, there is a strategy to construct a Mayer-Vietoris tree of it in which the multidegrees of all generators of the relevant nodes have non-zero homology. To construct the separable tree of $I$ select as pivot monomial one in the variables of $K$, following the same strategy as in the construction of the Mayer-Vietoris tree of $K$. Once this is done everywhere, we have a set of nodes such that the ideals in them are isomorphic to $J$. This means that the ideals in these nodes are generated by the minimal generators of $J$ each of them multiplied by the same monomial in $\kb[B]$. From here on, we can follow any convenient strategy.
Note that for each separable ideal we have a family of separable Mayer-Vietoris trees, depending on what strategy we follow in the nodes isomorphic to $J$.

With this construction, we have a procedure to compute the list of (multigraded) Betti numbers of $I$ in terms of those of $J$ and $K$. We can do it in several steps:
\begin{itemize}
\item First of all, all the Betti numbers of $K$ are in the list of Betti numbers of $I$.
\item Second, for each relevant node of $MVT(K)$ of dimension $i$ we add $r_J$ new elements to the list, where $r_J$ is the number of minimal generators of $J$. We denote by $R_i(MVT(K))$ the number of relevant nodes of $MVT(K)$ of dimension $i$. The multidegrees of these elements are the product of the pivot monomial of $K$ that was used to reach this relevant node multiplied by each of the generators of $J$.
\item Finally, for each of the leaves of $MVT(K)$ we have two nodes in $MVT(I)$. If the corresponding leaf of $MVT(K)$ has dimension $j$, each of these two nodes has an ideal isomorphic to $J$, one in an even position and dimension $j$, one in odd position and dimension $j-1$. Then, for each of these nodes in even position, we have to add to the list of $\beta_i(I)$ all the $\beta_{i-j}(J)$ if $i\geq j$ and for each of the nodes in odd position we add $\beta_{i-j}(J)$ if $i > j$. We denote by $L_j(MVT(K))$ the number of leaves of dimension $j$ in $MVT(K)$.
\end{itemize}

Sumarizing, we have that
\begin{equation}\label{separable_addition}
\beta_i(I)=\beta_i(K)+r_J\cdot R_i(MVT(K))+\sum_{i> j}L_j(MVT(K))\cdot(\beta_{i-j}(J)+\beta_{i-j-1}(J))
\end{equation}

\begin{Example}\label{example-separable}
Let $I=\langle x_1x_3,x_1x_4,x_2x_3,x_2x_4,x_5x_6,x_5x_7,x_6x_7\rangle$. $I$ is separable with $J=\langle x_1x_3,x_1x_4,x_2x_3,x_2x_4\rangle$ and $K=\langle x_5x_6,x_5x_7,x_6x_7\rangle$. Then $A=\{x_1,x_2,x_3,x_4\}$ and $B=\{x_5,x_6,x_7\}$. We have that $K$ is Mayer-Vietoris of type $B2$. The Betti numbers of $J$ and $K$ are
$$\beta_0(J)=3,\,\beta_1(J)=3,\,\beta_2(J)=1\quad \mbox{ and }\quad \beta_0(K)=3,\,\beta_1(K)=2$$
The minimal Mayer-Vietoris tree of $K$ is
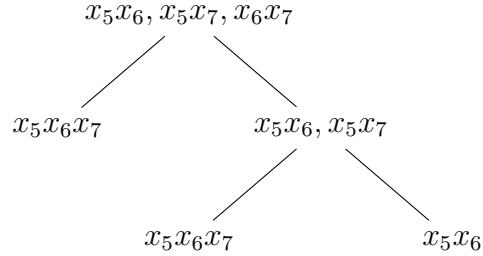
\begin{center}
\begin{tikzpicture}
\tikzstyle{level 1}=[sibling distance=3.5cm]
\node{$x_5x_6,x_5x_7,x_6x_7$}
child{node{$x_5x_6x_7$}}
child{node{$x_5x_6,x_5x_7$}
child{node{$x_5x_6x_7$}}
child{node{$x_5x_6$}}};
\end{tikzpicture}
\end{center}
It has one relevant node in dimension $0$, two relevant nodes in dimension $1$, and three leaves, two in dimension $1$ and one in dimension $0$, i.e. $R_0(MVT(K))=1,\,R_1(MVT(K))=2$ and $L_0(MVT(K))=1,\,L_1(MVT(K))=2$. In the second step of our procedure, we add for each of the relevant nodes $3$ elements to the list of the $\beta_1(I)$. And in the third step, from the two copies of $J$ hanging from the final leave in dimension $0$ we add $\beta_0(J)$ and $\beta_1(J)$ to $\beta_1(I)$; $\beta_1(J)$ and $\beta_2(J)$ to $\beta_2(I)$ and $\beta_2(J)$ to $\beta_3(I)$. Finally, from each of the two final leaves in dimension $1$ we have two copies of $J$ hanging, from which we add $\beta_0(J)$ and $\beta_1(J)$ to $\beta_2(I)$; $\beta_1(J)$ and $\beta_2(J)$ to $\beta_3(I)$ and finally $\beta_2(J)$ to $\beta_4(I)$. And we have
$$\beta_0(I)=\beta_0(K)+R_0(MVT(K))\beta_0(K)=3+3=6$$
$$\beta_1(I)=\beta_1(K)+R_1(MVT(K))\cdot r_J+L_0(MVT(K))(\beta_0(J)+\beta_1(J))=2+6+3+3=14$$
$$\beta_2(I)=L_0(MVT(K))(\beta_2(J)+\beta_1(J))+L_1(MVT(K))(\beta_1(J)+\beta_0(J))=(3+1)+2(3+3)=16$$
$$\beta_3(I)=L_0(MVT(K))\beta_2(J)+L_1(MVT(K))(\beta_2(J)+\beta_1(J))=1+2(3+1)=9$$
$$\beta_4(I)=L_1(MVT(K))\beta_2(J)=2\cdot 1=2$$
Here we show the first part of separable Mayer-Vietoris tree of $I$ where the three steps of the procedure can be clearly seen:
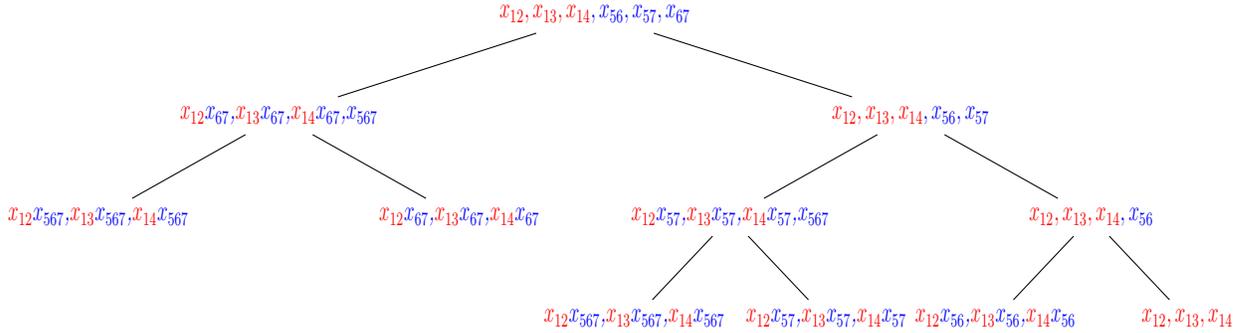
\begin{center}
\begin{tikzpicture}[xscale=0.6,yscale=0.9, transform shape]
\tikzstyle{everynode}=[font=\tiny]
\tikzstyle{level 1}=[sibling distance=14cm]
\tikzstyle{level 2}=[sibling distance=8cm]
\tikzstyle{level 3}=[sibling distance=4.25cm]

\node{\color{red}{$x_{12},x_{13},x_{14}$}\color{blue}{$,x_{56},x_{57},x_{67}$}}
child{node{\color{red}{$x_{12}$}\color{blue}{$x_{67}$},\color{red}{$x_{13}$}\color{blue}{$x_{67}$},\color{red}{$x_{14}$}\color{blue}{$x_{67}$},\color{blue}{$x_{567}$}}
child{node{\color{red}{$x_{12}$}\color{blue}{$x_{567}$},\color{red}{$x_{13}$}\color{blue}{$x_{567}$},\color{red}{$x_{14}$}\color{blue}{$x_{567}$}}}
child{node{\color{red}{$x_{12}$}\color{blue}{$x_{67}$},\color{red}{$x_{13}$}\color{blue}{$x_{67}$},\color{red}{$x_{14}$}\color{blue}{$x_{67}$}}}
}
child{node{\color{red}{$x_{12},x_{13},x_{14}$}\color{blue}{$,x_{56},x_{57}$}}
child{node{\color{red}{$x_{12}$}\color{blue}{$x_{57}$},\color{red}{$x_{13}$}\color{blue}{$x_{57}$},\color{red}{$x_{14}$}\color{blue}{$x_{57}$},\color{blue}{$x_{567}$}}
child{node{\color{red}{$x_{12}$}\color{blue}{$x_{567}$},\color{red}{$x_{13}$}\color{blue}{$x_{567}$},\color{red}{$x_{14}$}\color{blue}{$x_{567}$}}}
child{node{\color{red}{$x_{12}$}\color{blue}{$x_{57}$},\color{red}{$x_{13}$}\color{blue}{$x_{57}$},\color{red}{$x_{14}$}\color{blue}{$x_{57}$}}}
}
child{node{\color{red}{$x_{12},x_{13},x_{14}$}\color{blue}{$,x_{56}$}}
child{node{\color{red}{$x_{12}$}\color{blue}{$x_{56}$},\color{red}{$x_{13}$}\color{blue}{$x_{56}$},\color{red}{$x_{14}$}\color{blue}{$x_{56}$}}}
child{node{\color{red}{$x_{12},x_{13},x_{14}$}}}
}};
\end{tikzpicture}
\end{center}
In this tree the monomial $x_{i_1}\cdots x_{i_k}$ is represented by $x_{i_1\dots i_k}$. The nodes that have some monomial only in the variables of $B$ correspond to the nodes of the Mayer-Vietoris tree of $K$. In the part of the tree depicted, all nodes except the root have a copy of $J$, given by some monomials with mixed variables of $A$ and $B$. At each of the leaves of this part of the tree, we have a copy of $J$ multiplied by some monomial in $K$. To finish the tree, we just build six copies of the tree of $J$ and add the corresponding Betti numbers to those of $I$.
\end{Example}

We move now to ideals that are \emph{separable by intersection}. Note that all separable ideals by addition are splittable, in the sense of Eliahou and Kervaire \cite{EK90} (see also section \ref{reliability}). In particular, $I=J+K$ as above splits into $J$ and $K$. Eliahou and Kervaire provide the following formula for the Betti numbers of splittable ideals:
$$\beta_i(J+K)=\beta_i(J)+\beta_i(K)+\beta_{i-1}(J\cap K)$$
This formula can be used together with the considerations above. Therefore, using equation (\ref{separable_addition}) we have that if $I=J\cap K$ is separable by intersection
\begin{eqnarray*}\label{separable_intersection}
\beta_i(I)&=&\beta_{i+1}(J+K)-\beta_{i+1}(J)-\beta_{i+1}(K)\\
&=&r_J\cdot R_{i+1}(MVT(K))+\sum_{i+1> j}L_j(MVT(K))\cdot(\beta_{i+1-j}(J)+\beta_{i-j}(J))\,-\beta_{i+1}(J)
\end{eqnarray*}

Therefore, when dealing with separable ideals, Eliahou-Kervaire's formula can be improved, since the Betti numbers of either $J+K$ or $J\cap K$ can be computed using only the Betti numbers of $J$ and $K$. In example \ref{example-separable}, the use of Eliahou-Kervaire's formula to compute $\beta_i(J+K)$ needs the computation of the Betti numbers of the ideal $J\cap K$, which has $9$ minimal generators while $I=J+K$, which has just $6$ generators. Now, using the separable Mayer-Vietoris tree, as seen in the example, we have the Betti numbers of $J+K$ and together with the Eliahou-Kervaire formula, we arrive to the above equation to compute the Betti numbers of the ideal $J\cap K$, obtaining
$$\beta_0(J\cap K)=9,\quad \beta_1(J\cap K)=15,\quad \beta_2(J\cap K)=9,\quad \beta_3(J\cap K)=2$$ 

Finally, all these considerations can be summarized in the following
\begin{Proposition}\label{separable-formulas}
Let $I,I'\subseteq\kb[x_1,\dots,x_n]$ be monomial ideals such that $I=J+K$ is separable by addition and $I'=J'\cap K'$ is separable by intersection, with $K$ and $K'$ Mayer-Vietoris of type $A$ or $B2$. Then their Betti numbers are given by
\begin{equation*}
\beta_i(I)=\beta_i(K)+r_J\cdot R_i(MVT(K))+\sum_{i> j}L_j(MVT(K))\cdot(\beta_{i-j}(J)+\beta_{i-j-1}(J))
\end{equation*}
\begin{eqnarray*}
\beta_i(I')=r_J\cdot R_{i+1}(MVT(K'))+\sum_{i+1> j}L_j(MVT(K'))\cdot(\beta_{i+1-j}(J')+\beta_{i-j}(J'))\,-\beta_{i+1}(J')
\end{eqnarray*}
\end{Proposition}

In the most favorable situations, we can even give closed form expressions for the numbers $R_i(MVT(K))$ and $L_j(MVT(K))$ and $R_i(MVT(K'))$ and $L_j(MVT(K'))$. One such situation is given by the ideals coming from series-parallel systems. In section \ref{reliability} this type of systems are studied and actual formulas for the computation of their reliability are based on our discussion above.

\subsection{Mayer-Vietoris trees of powers of prime monomial ideals}\label{powers-primes}
\spanishsubsection{\'Arboles de Mayer-Vietoris de potencias de ideales monomiales primos}

We finish this section with some considerations on the Mayer-Vietoris trees of powers of prime monomial ideals. These ideals are Mayer-Vietoris of type $B2$ and appear also in applications since they are useful to study other ideals (see section \ref{families}) and allow us to generalize certain systems in reliability theory (section \ref{reliability}).
\begin{Lemma}\label{MVT-prime}
let $I=\langle x_i|x\in S\subseteq\{x_1,\dots,x_n\}\rangle$ an ideal generated by a subset of the variables, i.e. a prime monomial ideal, then $I$ is a Mayer-Vietoris ideal of type $A$.
\end{Lemma}
\noindent {\bf Proof: }
Assume for simplicity that $I=\langle x_1,\dots,x_k\rangle$ for some $k\leq n$. If $k=1$ the result is trivial. If $k\geq1$, consider the lexicographic Mayer-Vietoris tree of $I$. The nodes at positions $2$ and $3$ are respectively generated by $\langle x_ix_k\vert i<k\rangle$, and $\langle x_1,\dots,x_{k-1}\rangle$ the nodes of the tree hanging from one of the nodes have no common generator with the nodes of the other tree, since in one case all generators are multiples of $x_k$ and in the other case, no $x_k$ appears. Moreover, both trees are isomorphic, being the right tree $MVT(\langle x_1,\dots,x_{k-1}\rangle)$. An inductive argument yields the result.$\square$

Observe that these trees are $1$-linear, hence the regularity of these ideals is trivially $1$. Also, the Mayer-Vietoris resolution of these ideals is minimal. Moreover, easy formulas for their (multigraded) Betti numbers are obtained:
\begin{Lemma}\label{Betti-prime}
If $I$ is a prime monomial ideal generated by a set of $n$ variables then 
$$\beta_i(I)={{n}\choose{i+1}}\quad \forall i=\{0,\dots,n-1\}$$
and the multidegrees in which $\beta_{i,\alpha}\neq0$ are those $\alpha$ that are the product of exactly $i$ different variables.
\end{Lemma}

If we consider the powers of prime ideals, we move to the class of Mayer-Vietoris ideals of type $B2$, in which the Mayer-Vietoris resolution is still minimal:

\begin{Proposition}\label{J^k_is_MVB2}
Let $J$ be an ideal generated by a subset of the variables, then
$$J^k\in MV_{B2}\quad \forall k\geq 1$$
\end{Proposition}

\noindent{\bf Proof: } The proof will consist on finding a pure Mayer-Vietoris tree for $J^k$. For this, we first sort the generators of $J^k$ in the following way: For each integer $i=k,k-1,k-2,\dots$ in descending order, and for each variable $x_j$, $j=1,\dots,n$ taken in ascending order (this will be called `distinguished variable'), take successively the monomials $x^\mu$ such that $\mu_j=i$, and to the rest of the variables, put the exponents of the generators of $J_j^{k-i}$, where $J_j=\langle x_1,\dots,\hat{x_j},\dots,x_n\rangle$, with the following restrictions: $\mu_l\leq\mu_j\quad \forall l>j$ and $\mu_l<\mu_j\quad \forall l<j$, see the example below.

Consider now that we are in a node of dimension $0$. And consider as pivot monomial the first one according to the order we have just described. The right-hand subtree has as generator just a subideal of $J^k$ generated by a tail of its generators with respect to the just described ordering, therefore the same arguments apply to it. The left-hand subtree has as root the ideal generated by the lcms of the pivot monomial $x^\mu$ and the generators after it. Let $j$ be the distinguished variable of $x^\mu$ and $\mu_j=i$; from the way we have sorted the generators, we see that the monomials of the form $\mu+1_l$ with $l\in\{1,\dots,\widehat{j},\dots,n\}$  such that $\mu_l+1_l\leq i$ and $\mu_l+1_l< i$ if $l<j$, are present in the set of $lcm$'s we are using. The rest of the $lcm$'s are just multiples of these. Therefore, the left-hand tree has as root the ideal generated by a subset of the variables, multiplied by a certain monomial $x^\mu$ and is then isomorphic to the Mayer-Vietoris tree of a prime ideal. Therefore, all the generators of the ideal of a node of dimension $d$ have the same degree, namely $k+d$. Therefore, our tree is not only pure, but $k$-linear. Thus, $J^k$ is  a Mayer-Vietoris ideal of type $B2$.$\square$

\begin{Corollary}\label{tail-of-power-of-prime}
A subideal of $J^k$ generated by a final segment, or `tail'  with respect to the above ordering, of the generators of $J^k$ is a Mayer-Vietoris ideal of type $B2$ 
\end{Corollary}

\begin{Remark}
Some special `tails' are those formed by all generators listed from a given $i$ downwards. We can denote the ideals generated by this list as $J_{[n,i]}^k$ which is just the ideal in $n$ variables generated by monomials of degree $k$ in which the variables involved have exponents at most $i$. For example, $J_{[n,k]}^k$ is just $\mm^k$ and the ideals $J^k_{[n,1]}$ play an important role in reliability theory, as the ideals of $k$-out-of-$n$ systems (see section \ref{reliability}), hence ideals of the form $J^k_{[n,i]}$ allow us to study generalized or multistate $k$-out-of-$n$ systems. Of course, as we have seen, all ideals $J_{[n,i]}^k$ are Mayer-Vietoris of type $B2$ and therefore their Mayer-Vietoris resolution is minimal.

It is easy to see that powers of prime monomial ideals are stable in the sense of \cite{EK90}, and we will see later in section \ref{families} that stable ideals are Mayer-Vietoris of type $B2$. However, the above proof was given to introduce the sorting of the generators that allows us to see that the ideals generated by final segments with respect to this sorting are also Mayer-Vietoris of type $B2$. These ideals need not to be stable, in particular ideals of the form $J^k_{[n,i]}$ are not stable in general.
\end{Remark}

\begin{Example}
Let $J=\langle x,y,z\rangle$ and $k=4$. The table shows the way in which we sort the generators:
\begin{center}
\begin{tabular}{|c|c|c|c|c|c|c|c|c|c|}
\hline
\multicolumn{3}{|c|}{$i=4$}&\multicolumn{3}{c|}{$i=3$}&\multicolumn{3}{c|}{$i=2$}\\
\hline
\scriptsize{$j=1$}&\scriptsize{$j=2$}&\scriptsize{$j=3$}&\scriptsize{$j=1$}&\scriptsize{$j=2$}&\scriptsize{$j=3$}&\scriptsize{$j=1$}&\scriptsize{$j=2$}&\scriptsize{$j=3$}\\
\hline
$x^4$&$y^4$&$z^4$&$x^3y,x^3z$&$xy^3,y^3z$&$xz^3,yz^3$&$x^2y^2,x^2z^2,x^2yz$&$y^2z^2,xy^2z$&$xyz^2$\\
\hline
\end{tabular}
\end{center}
We now take each of the generators in order, from left to right, and consider the ideal generated by the $lcm$'s of it and the generators located at the right of it in the table. These ideals are of the form $x^\mu\cdot \hat J$, where $\hat J$ is an ideal generated by some smaller subset of the variables. For instance, if $x^\mu$ is $x^3y$, we have that $\hat J$ is $\langle y,z\rangle$; if $x^\mu$ is $xy^2$, $\hat J$ is just $\langle z\rangle$. From this point, we just use the Mayer-Vietoris tree of the corresponding $\hat J$s, multiplied by $x^\mu$.
\begin{center}
\begin{tabular}{|c|c|c|c|}
\hline
$x^\mu$&$\hat J$&$x^\mu\hat J$&$\widetilde{x^\mu\hat J}$\\
\hline
$x^4$&$\langle y,z\rangle$&$\langle x^4y,x^4z\rangle$&$\langle x^4yz\rangle$\\
$y^4$&$\langle x,z\rangle$&$\langle xy^4,y^4z\rangle$&$\langle xy^4z\rangle$\\
$z^4$&$\langle x,y\rangle$&$\langle xz^4,yz^4\rangle$&$\langle xyz^4\rangle$\\
$x^3y$&$\langle y,z\rangle$&$\langle x^3y^2,x^3yz\rangle$&$\langle x^3y^2z\rangle$\\
$x^3z$&$\langle y,z\rangle$&$\langle x^3z^2,x^3yz\rangle$&$\langle x^3yz^2\rangle$\\
$xy^3$&$\langle x,z\rangle$&$\langle x^2y^3,xy^3z\rangle$&$\langle x^2y^3z\rangle$\\
$y^3z$&$\langle x,z\rangle$&$\langle xy^3z,y^3z^2\rangle$&$\langle xy^3z^2\rangle$\\
$xz^3$&$\langle x,y\rangle$&$\langle x^2z^3,xyz^3\rangle$&$\langle x^2yz^3\rangle$\\
$yz^3$&$\langle x,y\rangle$&$\langle xyz^3,y^2z^3\rangle$&$\langle xy^2z^3\rangle$\\
$x^2y^2$&$\langle z\rangle$&$\langle x^2y^2z\rangle$&$-$\\
$x^2z^2$&$\langle y\rangle$&$\langle x^2yz^2\rangle$&$-$\\
$x^2yz$&$\langle y,z\rangle$&$\langle x^2y^2z,x^2yz^2\rangle$&$\langle x^2y^2z^2\rangle$\\
$y^2z^2$&$\langle x\rangle$&$\langle xy^2z^2\rangle$&$-$\\
$xy^2z$&$\langle z\rangle$&$\langle xy^2z^2\rangle$&$-$\\
$xyz^2$&$-$&$-$&$-$\\
\hline
\end{tabular}
\end{center}
The first column of this table gives the multidegrees of the $0$-th Betti numbers, and the third and fourth columns give the multidegrees of the first and second Betti numbers respectively. In the example, $\beta_0(I)=15,\,\beta_1(I)=24$ and $\beta_2(I)=10$.

From this example, we see that 

$$\begin{array}{l}
J_{[3,4]}^3=\langle x^3y,x^3z,xy^3,y^3z,xz^3,yz^3,x^2y^2,x^2z^2,x^2yz,y^2z^2,xy^2z,xyz^2\rangle\\
\beta_0(J^3_{[3,4]})=12\quad \beta_1(J^3_{[3,4]})=18\quad\beta_2(J^3_{[3,4]})=7\\
\\
J_{[3,4]}^2=\langle x^2y^2,x^2z^2,x^2yz,y^2z^2,xy^2z,xyz^2 \rangle\\
\beta_0(J^2_{[3,4]})=6\quad \beta_1(J^2_{[3,4]})=6\quad\beta_2(J^3_{[3,4]})=1
\end{array}
$$

\end{Example}

\begin{Remark}
The number of generators of the ideal $J^k$ whit $J=\langle x_1,\dots,x_n\rangle$ is the number of combinations with repetition of $n$ elements choose $k$ which equals $CR(n,k)={{n+k-1}\choose{k}}$. To compute the number of generators of $J_{[n,k]}^i$ we have to substract $n$ times all monomials in which some variable has a value bigger than $l$, and less or equal $k$, this is $\sum_{j=l+1}^k CR(n-1,k-j)$, therefore, the total number of generators of $J_{[n,k]}^i$ is $CR(n,k)-n\cdot\sum_{j=l+1}^k CR(n-1,k-j)={{n+k-1}\choose{k}}-n\cdot\sum_{j=l+1}^k{{n+k-j-2}\choose{k-j}}$.

One can also read the number of generators that have the same $i$ and $j$ in the ordering of the generators of $J^k$ described above. One way to do it is the following: Denote by $B_{n,k,i,l}$ the number we are looking for. It is equal to the number of monomials in $n$ variables, of total degree $k$ such that the $i-1$ first variables have degree less than $l$ the $i$-th variable has degree $l$ and the rest of the variables have degrees less than or equal $l$. Using generating functions, we find that $B_{n,k,i,l}$ is the coefficient of $x^{k-l}$ in the expansion of $$(1+x+\cdots+x^{l-1})^{i-1}(1+x+\cdots+x^l)^{n-i}=(\frac{1-x^l}{1-x})^{i-1}(\frac{1-x^{l+1}}{1-x})^{n-i}$$
$$=(1-x^l)^{i-1}(1-x^{l+1})^{n-i}\sum_{j\geq 0} CR(n-1,j)x^j$$.
Counting these generators and using lemma \ref{Betti-prime} we can obtain formulas for the Betti numbers of these ideals.
\end{Remark}

\section{Algorithm}
\spanishsection{Algoritmo}

In this section we present some details about algorithms that perform computations on monomial ideals using Mayer-Vietoris trees. First of all we discuss the different options that are available depending on what computations we want to perform and on the knw properties of our input ideals. Then, we present the basic algorithm to construct Mayer-Vietoris trees and then enter in some implementation details. Finally, we present the results of some experiments and comparisons with other algorithms that perform similar computations.

\subsection{The Basic Mayer-Vietoris tree algorithm}
\spanishsubsection{El algoritmo b\'asico de \'arboles de Mayer-Vietoris}

The algorithm given in table \ref{algorithm} performs the computation a Mayer-Vietoris tree of a monomial ideal $I$ in a very simple way. It is basically a recursive computation. The pseudo-code in table \ref{algorithm} presents this computation in the form of a loop. The main procedures involved are the computations of the children of a node, which are called \emph{tilde(ideal)} and \emph{ideal'}.

Every node in the tree is given by its position, dimension and generators. Note that the complexity of this algorithm depends strongly on the number of generators of $I$ and has a somehow weaker dependence on the number of variables (the necessary divisibility tests depend on the number of variables); it is (almost) independent on the degrees of the generators involved, provided the exponents fit in the limits of our system. The basic operations we use more often are divisibility tests among monomials, used to obtain the minimal generating sets of the ideals given by the $tilde(ideal)$ procedure, and also taking $lcm$s of monomials. Therefore, the importance of optimizing these two operations is crucial for the performance of the algorithm, as we will see later.

\begin{table}[!htb]
\centering
\begin{tabular}{|p{13cm}|}
\hline
$$\begin{array}{ll} Algorithm:\, \mbox{Mayer-Vietoris Tree of a Monomial ideal }$I$\\
\hline \hline\\
\ Input: \, \mbox{Minimal generating set of } I=\langle m_1,\dots m_r \rangle\\
\ Output: \, MVT(I) \mbox{ as a list of pairs (position, ideal)} \\
\\
\hline\\
\ 1 \quad {\bf if }\, r=1 \,{\bf then\,return }\, \{(\{1,0,\{m_1\}\})\}\\
\ 2 \quad \quad {\bf else } \\
\ 3 \quad \quad tree:=\{(1,0,\{m_1,\dots ,m_r\})\}\\
\ 4 \quad \quad undone:=tree\\
\ 5 \quad \quad {\bf while }\, undone \neq \emptyset\, {\bf do}\\
\ 6 \quad \quad \quad node:= first(undone)\\
\ 7 \quad \quad \quad undone:= tail(undone)\\
\ 8 \quad \quad \quad ideal:= ideal(node)\\
\ 9 \quad \quad \quad pos:= position(node),\, dim:=dimension(node)\\
\ 10 \quad \quad \quad append(tree,(2*pos,dim+1,tilde(ideal)))\\
\ 11 \quad \quad \quad append(tree,(1+2*pos,dim,ideal'))\\
\ 12 \quad \quad \quad {\bf if } \mbox{ number of generators}(tilde(ideal))>1\\
\ \qquad\qquad\qquad {\bf then }\, append(undone,(2*p,dim+1,tilde(ideal)))\\
\ 13 \quad \quad \quad {\bf if } \mbox{ number of generators}(ideal')>1\\
\ \qquad\qquad\qquad {\bf then }\, append(undone,(1+2*p,dim,ideal'))\\
\ 14 \quad \quad {\bf end\,while }\\
\ 15 \quad \quad {\bf return }\, tree\\
\ 16 \quad {\bf end\,if}\\
\end{array}
$$
\\
\hline
\end{tabular}\caption{Basic Mayer-Vietoris tree algorithm}\label{algorithm}
\end{table}

On steps $10$ and $11$ the children of a given node are computed. The procedures $tilde(ideal)$ and $ideal'$ compute the corresponding new ideals. These procedures have implemented in them the selection strategy used to choose the pivot monomial in each node. Many different strategies could be applied to construct the tree and these strategies could vary in each node. The strategy we use in our implementation consists in choosing as pivot monomial one generator whose exponent in some variable is maximal among the exponents in this variable of all generators . This strategy has advantages in relation to the size of the resulting tree and also in relation to the efficiency of the implementation. This will be briefly discussed below. Also, strategies based on term orderings have been implemented.

The output of this basic version of the algorithm is a list of the nodes of the tree. This will not in general be the version used in applications, because sometimes we need fewer information (just the relevant nodes, for instance) and sometimes we will need more information, like knowing which monomial are repeated or not as generators of relevant nodes. For most computations we will need to process the information given by the output of this basic algorithm, therefore it encodes just the first step of the actual Mayer-Vietoris based algorithms. However, an implementation of this basic algorithm is useful to test the times and performance of the actual computation of Mayer-Vietoris trees and compare it with other available algorithms used to make similar computations.

\subsection{Different versions of the algorithm}
\spanishsubsection{Diferentes versiones del algoritmo}

Mayer-Vietoris trees can be used for different computations. The main use of them is to compute (bounds of) multigraded Betti numbers of monomial ideals, but as we have already seen, other computations like Hilbert series, free resolutions, etc. can be achieved using Mayer-Vietoris trees. Also, the information coming from Mayer-Vietoris trees could be used to compute combinatorial and irreducible decompositions among others. Then, different versions should be available depending on what computations we want to perform and also on what we know in advance about our given ideal. Basically, the different options are the following:
\begin{enumerate}
\item Based on the type of ideal we deal with:
	\begin{enumerate}
	\item If we know the ideal is Mayer-Vietoris of type $A$ or $B2$.

		In this case, to obtain the multigraded Betti numbers we just need to store the generators of the relevant nodes. We can construct a different list for each dimension, and the Betti multidegrees will be given in an adequate form. If we are interested in the computation of the Mayer-Vietoris resolutions, since we know that it is minimal in this case, we can take advantage of the fact that we know that no further minimization will be necessary later.
	\item If we know the ideal is Mayer-Vietoris of type $B1$.

		The multigraded Betti numbers of Mayer-Vietoris ideals of type $B1$ are just the list of non-repeated monomials in the relevant nodes. Therefore, to obtain them, we need to store only the non-repeated generators of the ideals in the nodes. The checking and eventual deletion of repeated generators can be done either on the fly, checking the new monomial against the so far computed list of multigraded Betti numbers, or after the computation of the complete tree is already done. The Mayer-Vietoris resolution is not minimal in general in this case, so we have to take into account that if we want the minimal one from the resolution obtained from the tree, we will probably need to perform minimizations. This influences the way we store the necessary data. 
	\item The general case, in which we don't know whether the ideal is Mayer-Vietoris of any type.

		In this case we want to have all the information that the tree provides us about the generators of the relevant nodes. Knowing which of them are repeated and which are not, we can build the upper and lower bounds for the multigraded Betti numbers. Moreover, to decide which of the repeated multidegrees will actually be part of the final list of Betti numbers, we will use further criteria that will need the positions and dimensions of them, to establish relations among them and use the available criteria or procedures to make the decisions. Since the Mayer-Vietoris resolution will be non-minimal in general, we need to store the necessary data in a similar way as it was done in the case of Mayer-Vietoris ideals of type $B1$.
	\end{enumerate}
\item Based on the required output
	\begin{enumerate}
	\item Bounds for the Betti multidegrees.

		These bounds can be directly read from the relevant nodes of the Mayer-Vietoris tree. Then, the construction of the tree is almost all we need. We only need to know which of the generators of the relevant nodes are repeated, and which are not. The output consists then of two lists, one with the non-repeated generators of the relevant nodes, and another one with the repeated multidegrees.
	\item Mayer-Vietoris resolution.

		The construction of the Mayer-Vietoris resolution can be done in a recursive way, at the same time we compute the tree. We need to construct the matrices of the resolution while we compute the tree, so procedures for combining the matrices corresponding to the resolutions of the children into the matrices of the resolution of the father are needed. These procedures are based on the effective mapping cone construction given in sections \ref{resolutions} and \ref{computations}.
	\item Generating sets of the Koszul homology modules.

		To obtain generators of the Koszul homology modules we need procedures that implement the contracting homotopies for the Koszul differential and construct the new generators in a recursive way, form generators of the elements below in the tree.

	\item If we want \emph{minimal} results, i.e. the actual Betti numbers or resolutions.

	When we want to compute the actual sets of multigraded Betti numbers and /or the minimal free resolution, the construction of the Mayer-Vietoris tree is just a first step. Some minimization process must be performed afterwards. In the case of the Betti numbers, we need to select which of the repeated monomials in the relevant nodes do have homology and which do not. For this we need a good way to access them, see the relations betwen them and apply criteria and procedures to make the decision (see later). In the case of minimal free resolutions, we need to minimize the Mayer-Vietoris resolution. Since this has a very special form, the minimization procedure can be optimized, and for this we need a convenient storage of the Mayer-Vietoris resolution.

	\item Further computations, such as combinatorial or irreducible decompositions.

	Since Mayer-Vietoris trees are used to compute multigraded Betti numbers, they can be used to make the computations that were described in chapter \ref{structure}. Therefore, we can use the output of Mayer-Vietoris trees to compute irreducible decompositions, combinatorial decompositions, etc. An adequate storage of the sets of generators of the relevant nodes in the tree will optimize these algorithms. As was seen in sections \ref{irreducible_decomposition} and \ref{primary-decomposition}, the $(n-1)$-st Koszul homology is of particular interest for these decomposition. In the cases in which we are interested in obtaining only the multidegrees of the $(n-1)$-st homology generators, we can take the output of the basic Mayer-Vietoris algorithm and take the nonrepeated nodes of dimension $n-1$, then, for the repeated multidegrees $\mu$ in that dimension we keep those for which $\underline{\mu}$ is not in $I$ and $\underline{\mu}\cdot x_i$ is in $I$ for all $i$.
	\end{enumerate}
\end{enumerate}

\begin{Remark}
Some of the differences among the algorithms required for the different options are based on the knowledge of the ideal being Mayer-Vietoris. To be precise, we should say that we know that the ideal is Mayer-Vietoris, and we know how to construct a Mayer-Vietoris tree that gives us the desired result. Not all trees are equivalent even for Mayer-Vietoris ideals. The basic example is $I=\langle x^2,y^2,xy\rangle$. It is Mayer-Vietoris of type $A$, since it is an ideal in two variables, but not all Mayer-Vietoris trees of it will give us immediately the correct multigraded Betti numbers, see example \ref{mvt-examples-1}.
\end{Remark}

\subsection{Some implementation issues}\label{implementation}
\spanishsubsection{Algunas cuestiones sobre la implementaci\'on}

A full implementation and description of all the versions and features of these algorithms is beyond the scope of this thesis. We will briefly discuss in this section some of the most remarkable aspects to give an idea of the different characteristics that have been taken into account.
\subsubsection*{Data structures.}
\spanishsubsubsection{Estructuras de datos}

There are two levels in which the data structures used are most relevant. First of all, the structures used to represent monomials and monomial ideals. Second, the structures used to represent Mayer-Vietoris trees and their output.

With respect to monomials and monomial ideals, there are several data structures that have been used in the literature. A very effective one is the so called \emph{monomial trees} proposed by R.A. Milowski \cite{M04}, which was also used in \cite{R07}. It is a tree data structure which allows efficient algorithms for minimization, ideal intersection, etc. A different data structure to deal with monomial ideals is given in \cite{J04}. Since we want to use our algorithms in combination with other algebraic tools, it looks reasonable to make our algorithms as compatible as possible with other computer algebra systems, but without loosing efficiency. A good balance of both advantages can be found in the \cpp library $\cocoalib$ \cite{cocoalib}. We will use the structures implemented in $\cocoalib$ for power products and ideals, a short description of them is given later.

For Mayer-Vietoris trees the STL \cpp \emph{vector class} is a good alternative since it provides easy access to the components of the tree. Depending on the version of the algorithm we want to use, we will need different structures for the output of the Mayer-Vietoris tree algorithm. 
\begin{itemize}
\item[-] In algorithms where the output is known to be minimal in advance, like those for Mayer-Vietoris ideals of type $A$ or $B2$ where we just want the Betti numbers, the output will be a simple list of lists, giving the needed multidegrees by homological degree. Since no searches will be needed, the structure contains all the output in an optimal way.
\item[-] In the versions of the algorithm where the output is nonminimal, or at least it is not known to be minimal, we need to know where and which are the repeated multidegrees of the relevant nodes, we need a structure that makes these searches for repeated elements easier. Since the elements of our nodes can be sorted with respect to the monomial ordering we have, and since the operations we need are searches and insertions, binary trees represent a good option.
\item[-] In those versions in which the computation of Mayer-Vietoris trees is a preprocessing of the ideal for further computations, like free resolutions, irreducible decompositions, ... We also need rich structures, on one side because we need to store complicated data (matrices, for instance) and because we will need to perform further operations on them, like searches, comparisons and sortings. 
\end{itemize}

\subsubsection*{The selection strategy.}
\spanishsubsubsection{La estrategia de selecci\'on}

Mayer-Vietoris trees are highly customizable, in the sense that at each step of their construction, one can choose the pivot monomial with complete freedom. When we follow a common pattern at each step, we speak of a \emph{selection strategy} in the construction of the tree. There are three natural selection strategies (see remark \ref{strategies}), two of which can be easily implemented in the basic version of the algorithm, in particular when we have term orderings easily available.

The first strategy used in our implementation of the basic Mayer-Vietoris tree algorithm consists in choosing as pivot monomial one having biggest exponent in some of the variables, we will call this the \emph{maximum} strategy. This has some advantages: Let $j_i$ the exponent of the variable $x_j$ in the generator $m_i$ of the monomial ideal $I$. Assume the chosen pivot monomial $m_s$ is such that for some variable $x_j$ we have that $j_i\leq j_s$ for all $i\neq s$. Then, all generators of $\tilde I$ have the same exponent in $x_j$, namely $j_s$. Thus, keeping this strategy, after at most $n-1$ steps, all generators in the corresponding node have the same exponent in all variables, i.e. there's only one generator, and the corresponding branch of the tree stops here. Thus, for this strategy, every `left' branch of the tree is bounded above by the number of variables minus one. Another good reason to follow this strategy is the following: The basic operation on which the whole algorithm is based is divisibility among monomials. The most frequent operation in the algorithm are divisibility tests in which the exponents of each variable in the monomials are compared. If we know that some monomials have the same exponent in some of the variables, we do not need to compare these variables to check divisibility. Thus, following the strategy just mentioned allows us to improve the implementation of the algorithm avoiding many comparisons between exponents of the monomials involved. In few words, once a variable is used, the left branch of the subtree can be considered as an ideal in one variable less, and so on.

Also strategies based on term orderings have been implemented. Using the term orderings available in $\cocoalib$, we are able to compute \lex, \deglex and \degrevlex Mayer-Vietoris trees. Here, two options are available for each term order, namely, select as pivot monomial respectively the first or the last monomial according to the selected term ordering. Some experiments showing the performance of these strategies can be seen in section \ref{experiments}.
\subsubsection*{Implementation using $\cocoalib$}
\spanishsubsubsection{Implementaci\'on usando $\cocoalib$}

Our algorithm is implemented using the \cpp library $\cocoalib$ \cite{cocoalib}. There are several reasons to do it so:

First of all, it provides a good balance between integration with other computer algebra software, in particular with \cocoa, and the programming capabilities of \cpp. Since we look for competitive speed and at the same time we want to transfer the results to systems capable of computing with algebraic and symbolic objects, $\cocoalib$ looks like a good choice to reach both purposes.

On the other hand, $\cocoalib$ provides us with the natural algebraic structures we deal with, in a very natural way. It has rings, ideals, power products, etc. already implemented in an efficient way. Also term orderings are available, and we can make direct use of them in the implementation of the selection strategies for our algorithm.

Finally, there are some issues in $\cocoalib$ that make it appropriate for the Mayer-Vietoris algorithms. In particular there are good implementations of procedures that are very often used in our algorithms, like obtaining the minimal generating set of a monomial ideal or fast divisibility tests between monomials.

The main features of the $\cocoalib$ classes we use are the following:

The class {\tt PPMonoid} is used to represent the monoid of the power products in our ideals. Basically all operations between monomials are performed considering them as members of this class. When creating a {\tt PPMonoid}, the grading and ordering must be specified, and also the names of the indeterminates. Available orderings are \lex, \deglex and \degrevlex. Also available is a function to create an ordering given a matrix. Although these matrix-based orderings have not been used in the basic implementation of the Mayer-Vietoris algorithms, they could be of use when implementing non-standard strategies to construct the trees. There is a complete catalog of functions available for members of the {\tt PPMonoid} class. Among them we use comparison operators with respect to the specified ordering, product and quotients, divisibility tests and least common multiples. In particular, a fast divisibility test can be used for power products, using masks (see the $\cocoalib$ documentation file {\tt DivMask.txt} and see also the discussion about performing divisibility tests of power products in \cite{B97}).

{\tt SparsePolyRing} is the class used for representing rings of polynomials. A polynomial ring with a specified term ordering coming from our {\tt PPMonoid} can be constructed with the functions {\tt NewPolyRing(CoeffRing, IndetNames, ord)} or {\tt NewPolyRing(CoeffRing, PPM)}; the  power product monoid specifies how many indeterminates are we using, their names, and the term ordering. Considering our monomials as elements in these polynomial rings, we can build monomial ideals and perform computations with them, since so far there is not a specific class for monomial ideals. Some of the available operations for ideals that we use are the addition, or intersection of ideals, membership tests, etc. Of particular importance is the computation of a minimal generating set of a monomial ideal. The $\cocoalib$ function that performs this is {\tt TidyGens} which is very fast, and therefore very important for the actual performance of our implementation.

\subsubsection*{Repeated multidegrees.}
\spanishsubsubsection{Multigrados repetidos}

When the ideal we are dealing with is not of any Mayer-Vietoris type, our algorithm provides bounds for the multigraded Betti numbers. In order to obtain the correct Betti numbers, we need to know which of the multidegrees that are repeated in the relevant nodes of the tree contribute to the homology of the ideal, and which do not. There are several criteria and procedures that can be used to solve this question. We can summarize them as follows:
\begin{itemize}
\item Algebraic criteria:

	The algebraic properties of the ideal can help determining whether a given multidegree contributes or not to the homology. For instance, the number of variables being $n$ obviously tells us that every monomial appearing in a node of dimension bigger than $n$ will have no homology at that degree. Knowing in advance some other properties of the ideal like the dimension, regularity, etc. help us to perform further eliminations.
	
	On the other hand, since the Mayer-Vietoris resolutions corresponding to the given tree is supported on the generators of the relevant nodes, we can perform the usual minimization process to it and hence the surviving multidegrees will have nonzero Koszul homology. This can always be performed.

\item Simplicial procedures:

	One way to decide whether there is Koszul homology for a given ideal at a given multidegree at some homological degree is by means of the Koszul simplicial complexes (see section \ref{simpl-koszul-complex}). This gives a definitive answer to this question. Being this procedure too expensive in many cases, since it amounts to compute simplicial homology, it is sometimes very useful, either if the complexes are easy to determine and compute, or if we are interested in some particular multidegree, for example when computing regularity.
	The Koszul simplicial complex can be used in combination with Stanley-Reisner theory (see section \ref{stanley-alexander}) and again Mayer-Vietoris trees in the following way: If we want to know if $H_{i,\ab}(\KK(I))$ is zero or not, we can construct the Koszul simplicial complex $\Delta_I^\ab$ (or its dual). Then, instead of directly computing its reduced homology, we can consider its Stanley-Reisner ideal $I_{\Delta_I^\ab}$ and compute the Mayer-Vietoris tree of it. These ideals are the upper and lower Koszul ideals of $I$ at $\ab$ and are described in section \ref{stanley-alexander}. Looking at the Mayer-Vietoris tree of these ideals in the dimension we are looking for can help obtaining the answer. This procedure has shown to be useful in examples, see section \ref{experiments}.

\item Homological criteria:
	The machinery under Mayer-Vietoris trees consists basically on having an exact sequence in homology for each multidegree. Most of the times these sequences consists on very few terms. In the case of non-repeated monomials in relevant nodes these sequences consist on just two terms and many other cases considerations on these sequences give us the answer about repeated monomials contributing to homology (see for instance section \ref{special-MVT}, or proposition \ref{repeated-relevant-multidegrees}). In any case, the rank nullity theorem provides us relations, in form of equations, among the appearances of a multidegree in the relevant nodes of a tree. these relations together with the above criteria and procedures help us in detecting homology at given multidegrees.
	
\end{itemize}

Note that whenever we find that one multidegree $\mu$ in the list of repeated multidegrees of the relevant nodes does not contribute to the $i$-th Koszul homology of $I$ we can automatically delete one of the further appearances of $\mu$, since by the chain complex reduction algorithm we have seen that the corresponding element in the resolution forms a reduction pair with another element. Proceeding in descending order with respect to homological degree will allow us to use this feature.

Some of these criteria are difficult to implement and some of them have a more theoretical use as can be seen in chapter \ref{Applications}. Of course, also the computation of different trees using different strategies for a given ideal will provide further knowledge about the ideal. In any case, either by computing the homology of the simplicial complexes or by minimizing the Mayer-Vietoris resolution, we have always one way to obtain the homological description of the monomial ideal using Mayer-Vietoris trees.

We have implemented an easy procedure to eliminate from $MVT(I)$ most of the multidegrees that do not contribute to the Koszul homology of $I$. A table with the performance of this procedure can be seen in section \ref{experiments}. It consists on three steps or tests:
\begin{itemize}
\item [-] The first test deals with multidegrees that are particularly easy to eliminate due to some special criteria. First of all, we eliminate all multidegrees of dimension bigger or equal $n$, since they cannot contribute to the homology. This is done by simply non storing them. Second, we delete all elements which are inside the ideal, i.e. such that $\pi(\mu)\in I$ (see definition \ref{boundary} and remark \ref{boundary-roune}), since only the elements of the boundary can contribute to homology. And finally, for all repeated multidegrees of dimension $n-1$ we check wether they are maximal corners or not, which is easy to examine and is equivalent to have $(n-1)$-st homology (see proposition \ref{closed-maximal}).
\item[-] The second test uses Koszul ideals. For those multidegrees in the boundary, we build the corresponding Koszul ideal(s) and use them to determine whether they have homology or not at our given multidegree (proposition \ref{koszul-stanley-betti}). For this, we use again Mayer-Vietoris trees of the obtained Koszul ideals.
\item[-] Some of the multidegrees that are undecided by the two precedent tests will be the $\mu\in\NN^n$ such that there exists some $i$ with $\beta_{i,\mu}(I)>1$. Therefore, we can still run a third test based on proposition \ref{repeated-relevant-multidegrees}, namely we look for repeated multidegrees that appear only at a given dimension in the tree, and we store them with the corresponding multiplicity. 
\end{itemize}
\begin{Remark}
Some remarks are necessary about this procedure:
\begin{itemize}
\item The output of this procedure for each multidegree can be that the given multidegree has homology, that it has not, or that it is undecided. In this third case, further theoretical criteria can be applied, or finally one can use the Koszul simplicial complex or the minimization of the corresponding Mayer-Vietoris resolution. Observe that although this procedure has a good performance in detecting multidegrees which do not contribute to the Koszul homology (see table \ref{experiment2}), it will leave undecided at least those repeated multidegrees $\mu$ which have homology in different dimensions, i.e. such that there exist $i\neq j$ with $\beta_{i,\mu}\neq 0$ and $\beta_{j,\mu}\neq 0$
\item When using the above described test we will use the Mayer-Vietoris trees of the corresponding Koszul ideals to determine whether they have homology at the given multidegree in the appropriate dimension. Note that although the upper and lower Koszul ideals are Alexander duals, and therefore have equivalent multigraded homology, the Mayer-Vietoris trees of them can be quite different and therefore multidegrees that are undecided by one tree may be decided by the other. Consider for example the following ideal: $I=\langle xyz,yzt,ztu,uvx,vxy\rangle$ which is called the \emph{cyclic $3$-out-of-$6$} ideal. When deciding whether the multidegree $xyztv$ has homology or not, we may consider the lower Koszul ideal of $I$, $KI^{xyztv}=\langle yzt,xyv,xyz\rangle$, depending on the strategy followed, a Mayer-Vietoris tree of $KI^{xyztv}$ leaves $xyztv$ undecided. On the other hand, the upper Koszul ideal $KI_{xyztv}$ is $\langle y,xz,xt,zv\rangle$ which, with the same strategy shows that $H_{i,xyztv}(\KK(I))$ is zero for all $i$.
\item Observe that these tests produce tighter bounds for the Betti numbers of the ideal that come from the corresponding Mayer-Vietoris tree. In particular, when eliminating multidegrees we lower the upper bound, which is done by tests 1 and 2. The third test states that certain repeated multidegrees do have homology, therefore it raises the lower bound.
\end{itemize}
\end{Remark}
Figure \ref{algorithms-versions} shows the relations between the different versions of the algorithm, the components they need and the different outputs we obtain.

\begin{center}
\begin{figure}
\begin{tikzpicture}[auto,>=latex',scale=0.6,transform shape]

\tikzstyle{every node}=[text centered, font=\Large, font=\bf]
\tikzstyle{algorithm} = [draw, thin, draw=blue, fill=blue!20, text width=3.5cm]
\tikzstyle{result} = [ellipse, draw, thin, draw=red,fill=red!40, minimum height=2.5cm,text width=2.5cm]
\tikzstyle{comment} = [draw,double, thin, draw=orange,fill=yellow!40, minimum height=2.5em]


    \draw(0,45) node[result] (c-m) {C-M regularity};
    \draw (5,45) node[result] (betti) {$\beta_{i,\mu}$};
    \draw (15,45) node[result] (min-res) {Min. Res.};
    \draw (25,45) node[result] (koszul) {Koszul basis};

    \draw (5,42) node[algorithm] (simplicial) {Simplicial Computation};
    \draw (5,38) node[algorithm] (tests) {Tests};

    \draw (10,40) node[result] (hilbert) {Mult. Hilbert}; 
    \draw (15,40) node[algorithm] (minimize) {Minimize Resolution};
    \draw (25,40) node[algorithm] (minimize-koszul) {Minimize Set};

    \draw (0,35) node[result] (c-m-bounds) {C-M Reg. bounds};
    \draw (15,35) node[result] (mvr) {MVR};
    \draw (25,35) node[result] (koszul-gens) {Koszul generators};

    \draw (5,30) node[result] (betti-lower) {$\overline{\beta_{i,\mu}}$};
    \draw (10,30) node[result] (betti-upper) {$\widehat{\beta_{i,\mu}}$};

    \draw (22,30) node[algorithm] (koszul-computation) {Computation of Koszul cycles and gens.};
    \draw (15,25) node[algorithm] (tree-res) {Resolution from the tree};
    \draw (5,25) node[algorithm] (repeated) {Repeated and non-repeated};

    \draw (10,15) node[algorithm] (basic) {MVT Basic};

	\path[->] (betti) edge (c-m);
	\path[->] (min-res) edge (betti);
	\path[->] (min-res) edge[bend right](c-m);
	\path[dashed,->] (min-res) edge [bend left] node {[Herzog,Sergeraert]} (koszul);
	\path[dashed,->] (koszul) edge [bend left] node {[Herzog,Sergeraert]} (min-res);
	\path[->] (koszul) edge[out=90,in=90] (hilbert);

	\path[->] (betti) edge (hilbert);
	\path[->] (simplicial) edge (hilbert);
	\path[->] (simplicial) edge (c-m);
	\path[->] (tests) edge (c-m);
	\path[->] (tests) edge (hilbert);
	\path[->] (min-res) edge (hilbert);

	\path[dashed,->] (minimize-koszul) edge (koszul);
	\path[dashed] (koszul-gens) edge (minimize-koszul);
	\path[->] (koszul-gens) edge[out=270,in=270] (hilbert);

	\path (mvr) edge (minimize);
	\path[->] (minimize) edge (min-res);
	\path[dashed,->] (mvr) edge (koszul-gens);

	\path[->] (c-m-bounds) edge (tests);

	\path[->] (betti-lower) edge (c-m-bounds);
	\path[->] (betti-upper) edge (c-m-bounds);
	\path[->] (betti-upper) edge (hilbert);
	\path[->] (betti-upper) edge (tests);
	\path[->] (repeated) edge (c-m-bounds);
	\path[->] (repeated) edge (betti-lower);
	\path[->] (repeated) edge (betti-lower);
	\path[->] (repeated) edge (koszul-computation);
	\path[->] (repeated) edge (betti-upper);
	\path[->] (koszul-computation) edge (koszul-gens);

	\path[->] (tree-res) edge (mvr);

	\path[->] (basic) edge (repeated);
	\path[->] (basic) edge (betti-upper);
	\path[->] (basic) edge (tree-res);
	\path[->] (basic) edge (koszul-computation);

\end{tikzpicture}
\caption{Versions, components and outputs of Mayer-Vietoris algorithms}\label{algorithms-versions}
\end{figure}
\end{center}

\subsection{Experiments}\label{experiments}
\spanishsubsection{Experimentos}

We finish this section giving the results of some experiments using the described algorithms.
First of all we show the importance of choosing a good strategy when building Mayer-Vietoris trees. This experiment constructs seven different trees for several random ideals. The trees are constructed according to the following strategies: $\lex_l$, $\deglex_l$ and $\degrevlex_l$ use as pivot the \emph{last} monomial in the generators of each node, with respect to the term orders \lex, $\deglex$ and $\degrevlex$ respectively; $\lex_f$, $\deglex_f$ and $\degrevlex_f$ take as pivot monomial the \emph{first} monomial with respect to the corresponding term order. Finally $max$ is the \emph{maximum} strategy explained in remark \ref{strategies}. Table \ref{experiment1} shows the number $n$ of variables of the ideal, the number $g$ of minimal generators, the size $S$ of the minimal resolution (sum of all the Betti numbers) and the size of the Mayer-Vietoris resolution from the given tree. The examples are random in the sense that every variable in each generator has a random exponent between $0$ and $50$. The table shows that the choice of a good strategy is crucial for the performance of the algorithm. Also, the best strategies are those which use the last monomial with respect to some ordering as pivot monomial. In particular, $\lex$ has shown a good behaviour for random examples. The maximum strategy shows good average behaviour too. Observe that in the last examples, as the number of variable grows with respect to the number of generators ,the divisibility relations among least common multiples of the generators become more infrequent and therefore Taylor resolution is closer to minimal. In the last example it is indeed the case and therefore every Mayer-Vietoris tree gives the minimal resolution, which has size $2^g-1$.

\begin {table}[h]
\begin {center}

    \begin{tabular}{|c|c|c|c|c|c|c|c|c|c|}
    	\hline
	$n$&$g$&$S$&$\lex_f$&$\lex_l$&$\deglex_f$&$\deglex_l$&$\degrevlex_f$&$\degrevlex_l$&max\\
	\hline
    	$10$&10&227&399&227&685&233&685&233&245\\
	$10$&20&4361&36425&4435&152245&5249&146957&5225&5969\\
	$10$&30&5875&20961&5983&?&6673&?&6701&9743\\
	$15$&10&503&766&503&1023&503&1023&503&593\\
	$15$&20&11993&66763&12127&?&13055&?&13071&27535\\
	$15$&30& & &611637&?&1025935&?&1027233&1312199\\
	$20$&15&9859&18187&9971&31301&10419&31183&10419&17375\\
	$40$&15&24719&28991&25911&32755&26503&32759&26503&28861\\
	$60$&15&32767&32767&32767&32767&32767&32767&32767&32767\\
    	\hline
    \end{tabular}

\caption{Strategies for building Mayer-Vietoris trees of random examples}
\label{experiment1}
\end{center}
\end{table}
Our second table (table \ref{experiment2}) shows the number of repeated multidegrees that can be treated using the procedure shown in section \ref{implementation}. Note that in the cases in which the left branches of the tree are bounded by the number of variables (minus one) like in the \emph{maximum} strategy, the first part of the procedure will not be used. In this case, the column $\vert MVR(I)\vert$ indicates the size of the corresponding Mayer-Vietoris resolution, the column \emph{Repeated} shows the number of multidegrees that are repeated as generators of relevant ideals and the column \emph{Undecided} gives the number of those multidegrees that remain undecided after applying the procedure described in section \ref{implementation}. Observe that for random examples these procedures decide all multidegrees. This is not to be expected in the general case, although it is not easy to find examples in which the procedure leaves undecided multidegrees.

\begin {table}[h]
\begin {center}

    \begin{tabular}{|c|c|c|c|c|c|}
    	\hline
	$n$&$g$&$S$&$\vert MVR(I)\vert$&Repeated&Undecided\\
	\hline
	$5$&$30$&285&345&32&0\\
	$10$&$10$&227&233&3&0\\
	$10$&$20$&4361&5171&400&0\\
	$10$&$30$&5875&6675&393&0 \\
    	\hline
    \end{tabular}
\caption{Number of undecided multidegrees in random examples.}
\label{experiment2}
\end{center}
\end{table}

In table \ref{table3} we show the times and results of computations\footnote{All computations were made on a Pentium IV processor (2.5GHz) running CoCoA 4.7, Singular 3.0.1 and Macaulay2 0.9.8 under Linux (Mandriva 2007).} in some examples. In the table, the example $r_{n,g}$ corresponds to a random ideal in $n$ variables with $g$ generators, generators are given by a product of the $n$ variables each of which has a random exponent between $0$ and $50$. The ideals $v(n,a,b)$ correspond to the ideals described by Valla in \cite{V04}, where $n$ is the number of variables and $a>b$ two parameters that define the ideal (see section \ref{valla}). The ideals $I_{k,n}$ correspond to \emph{$k$-out-of-$n$} ideals, $\bar{I}_{k,n}$ to \emph{consecutive $k$-out-of-$n$} and $\tilde{I}_{k,n}$ to \emph{cyclic $k$-out-of-$n$}; these are described in section \ref{reliability} and refer to some important systems in reliability theory. Finally the ideals $I_{\{4,4,3,2\}}$ and $I_{\{5,6,3,3\}}$ describe two particular series-parallel networks, and are also described in section \ref{reliability}; the number of generators of $I_{\{n_1,\dots,n_k\}}$ is $\prod_i n_i$, and the number of variables is $\sum_i n_i$. Note that for some of these ideals, actual formulas for their Betti numbers are available, therefore the purpose of including some Mayer-Vietoris tree computations with them is only to illustrate the performance of the algorithms with different types of ideals.

The column $gens$ gives the number of minimal generators of the example. $S$ is the size of the minimal resolution; the column CoCoA shows the time (in seconds) used by the CoCoA \cite{cocoa} command {\sf BettiDiagram} in computing the Betti Diagram of $I$ (note that this command computes the minimal free resolution of the ideal), question mark means that computations stopped with no result. The columns \emph{mres},  \emph{hres}, \emph{sres} and \emph{lres} show the time in which the respective Singular \cite{Singular} commands computed the resolutions of the given ideals using different methods. Commands {\sf mres} and {\sf hres} (the second one using Hilbert driven algorithm) compute the minimal free resolution, the other ones compute nonminimal free resolutions using respectively Schreyer and La Scala algorithms. The column M2 contains the times in which the Macaulay2 \cite{M2} command {\sf res} computes the (minimal) resolution of the ideal. Finally, the column tree shows running times of our implementation of the basic version of the Mayer-Vietoris tree algorithm. The algorithm was implemented using \cocoalib. The strategies used are those in table \ref{experiment1}, in particular $\lex_l$ is used in most of the examples. This basic version computes the upper and lower bounds for the actual multigraded Betti numbers. Thus, it is a partial computation with respect to the other columns. To complete this computations one should read the computed tree to obtain the resolutions and homology generators of section \ref{computations}. However, in a variety of examples and applications, the data obtained by the basic version suffice, in particular if we know in advance that the ideal is Mayer-Vietoris of some type, or if the output of the algorithm tells us. These examples are marked with an asterisk in the table.

\begin {table}[h]
\begin {center}
    \begin{tabular}{|c|c|c|c|c|c|c|c|c|c|}
    	\hline
	Example&gens&$S$&CoCoA&mres&hres&sres&lres&M2&Tree\\
	\hline
    	$r_{12,19}$&19&9311&980'71&89'53&202'30&2'52&0'75&2'68&0'15\\
	$r_{22,17}$&17&45431&?&2250'98&?&27'12&15'11&16'72&0'94\\
	$r_{55,17}$&17&106761&?&3839'89&?&119'49&170'42&?&2'60\\
	$r_{65,18}$&18&212514&?&?&?&?&660'11&?&4'9\\
	$^*v(10,4,2)$&715&178177&?&5487'85&?&277'19&301'59&?&5'6\\
	$^*v(11,4,2)$&1001&471041&?&?&?&2268'79&2773'21&?&12'83\\
	$^*I_{5,10}$&252&5503&10'55&1'33&1'17&0'11&0'14&0'49&0'49\\
	$\tilde{I}_{5,20}$&17&439&0'26&0'1&0'1&0'1&0'1&0'02&0'01\\
	$^*\bar{I}_{5,40}$&36&247551&?&?&?&1097'65&1271'99&?&7'38\\
	$^*I_{\{4,4,3,2\}}$&96&4725&4'13&0'91&0'44&0'1&0'13&0'48&0'14\\
	$^*I_{\{5,6,3,3\}}$&270&95697&?&1907'57&3141'12&96'52&95'3&15'02&2'31\\
    	\hline
    \end{tabular}
\caption{Times of algorithms computing free resolutions of different ideals}
\label{table3}
\end{center}
\end{table}

Observe in the CoCoA and $hres$ column that the system is not able to perform computations as the size of the resolutions grow, this is probably due to the fact that the algorithms implemented for these commands are not specially targeted to monomial ideals but to general polynomial ideals; moreover it looks that the strategy implemented in these commands is particularly inneficient when applied to monomial ideals although could be particularly efficient when applied to other types of polynomial ideals. The Macaulay2 column shows that even if the computations with relatively small resolutions are quite efficient, as soon as the size of them grow, some memory problems appear and the system cannot perform the computations.  The fastest commands in Singular, compute nonminimal free resolutions, the size of which is compared to the size of the Mayer-Vietoris resolution in table \ref{table2} for some of the examples. The size of a resolution is considered as the sum of the ranks of all modules in it, i.e. the size of the minimal resolution is the sum of all Betti numbers. Column {\it sres} shows the size of the resolution computed using Schreyer algorithm in Singular, and column {\it lres} shows the size using La Scala algorithm, observe that in these examples La Scala computes the minimal resolution, although in general it is not the case (and the $\singular$ output of the algorithm warns this). We can see that Mayer-Vietoris trees constitute an alternative also to the computation of minimal free resolutions by usual algorithms. Observe that it has good behavior with respect to the growing of the number of variables when compared to other algorithms.

\begin {table}[h]
\begin {center}
    \begin{tabular}{|c|c|c|c|c|c|c|c|c|c|}
    	\hline
	Example&gens&minimal&sres&lres&Mayer-Vietoris\\
	\hline
    	$random(12)$&19&9311&25915&11667&10129\\
	$random(22)$&17&45431&63418&45431&46223\\
	$random(55)$&17&106761&106761&106761&106761\\
	$random(65)$&18&212514&?&212514&221639\\
	$valla(10,4,2)$&715&178177&178177&178177&178177\\
	$valla(11,4,2)$&1001&471041&471041&471041&471041\\
    	\hline
    \end{tabular}
\caption{Size of resolutions computed by Singular and MVT}
\label{table2}
\end{center}
\end{table}

Every Mayer-Vietoris tree of a monomial ideal $I$ gives a resolution of $I$ supported on the generators of the relevant nodes of the tree. Therefore, the alternating sum of the relevant generators of every Mayer-Vietoris tree gives an expression of the multigraded Hilbert series of $I$. Then, the basic Mayer-Vietoris tree algorithm can be used to compute these multigraded series. The computer algebra system $\cocoa$ has a procedure to compute multigraded Hilbert series; tables \ref{table4} and \ref{table5} show comparisons of the performance of both algorithms in some examples. Table \ref{table4} uses some Valla ideals as examples. The size of the minimal free resolutions of these ideals is big with respect to the number of variables. Table \ref{table5} uses \emph{consecutive $k$-out-of-$n$} ideals (see section \ref{reliability}) because the size of the resolution in terms of the number of variables and generators has a reasonable growth. The first number in each cell of this table \ref{table5} shows the time taken by $\cocoa$ and the second one is the time taken by the Mayer-Vietoris tree algorithm.

\begin {table}[h]
\begin {center}
    \begin{tabular}{|c|c|c|}
    	\hline
	Example&$\cocoa$&MVT\\
	\hline
    	valla(6,4,2)&0'20&0'13\\
	valla(8,4,2)&0'41&1'05\\
	valla(10,4,2)&1'95&5'6\\
	valla(12,4,2)&7'64&28'89\\
    	\hline
    \end{tabular}
\caption{{\sf HilbertSeriesMultiDeg} and MVT times for some Valla ideals}
\label{table4}
\end{center}
\end{table}

\begin {table}[h]
\begin {center}
    \begin{tabular}{|c|c|c|c|}
    	\hline
	$n$ &$I_{5,n}$&$I_{10,n}$&$I_{15,n}$\\
	\hline
    	20&0'05::0'02&0'02::0&0'02::0\\
	25&0'1::0'07&0'04::0&0'02::0\\
	30&0'57::0'32&0'09::0'02&0'03::0\\
	35&4'09::1'73&0'12::0'03&0'04::0'01\\
	40&62'39::7'38&0'18::0'09&0'07::0'01\\
	45&?::46'29&0'5::0'3&0'14::0'03\\
    	\hline
    \end{tabular}
\caption{{\sf HilbertSeriesMultiDeg} and MVT times for consecutive $k$-out-of-$n$ ideals}
\label{table5}
\end{center}
\end{table}

We can see in the tables that Mayer-Vietoris trees are an efficient alternative when the number of variables grow, while the $\cocoa$ implementation of multigraded Hilbert series is very efficient when the number of variables is small, even if the size of the resolution is big. In other types of examples, like random or network ideals, the behaviour is the same.             
\fancyhead[ER]{\itshape \leftmark} 
\chapter{Applications}\label{Applications}
\spanishchapter{Aplicaciones}
One of the favorable characteristics of monomial ideals is that they are suitable for applications inside and outside mathematics. The transformation of some problems relative to polynomial algebra in terms of monomial ideals is present in the Gr\"obner basis theory. Some problems in other areas can also be treated in a monomial framework, and monomial ideals and modules are very useful to model problems in actual applications outside mathematics. 
In this chapter we will give an overview of several of these applications, and use the techniques described in the precedent chapters. We only describe some of the applications in which those techniques have a relevant significance, but the field is wide open for further developments.

\section{Families of monomial ideals}\label{families}
\spanishsection{Familias de ideales monomiales}

Combinatorial commutative algebra has been an active area of research in the recent years, and monomial ideals are a main object in it. One of the directions in which progress have been made concerning monomial ideals is to consider some special types of them. Some of the more interesting ones are stable, generic and cogeneric, minimal, Borel-fixed... \cite{EK90,BPS98,MS04,GPW99,MSY00,P05}. In the following sections we will use Mayer-Vietoris trees to study different classes of monomial ideals. Some of these classes, such as Borel-fixed or stable monomial ideals are interesting because of the role they play in the study of Hilbert functions and Gr\"obner basis theory. Scarf monomial ideals, and in particular generic monomial ideals have a big importance in combinatorial commutative algebra and satisfy many interesting properties \cite{MSY00,Y99}. Some of the ideals studied arise from applications in different branches of mathematics, like geometry (Valla ideals) or graph theory (Ferrers ideals). Some of the results in this section were presented in \cite{SC07d}.

\subsection{Borel-fixed, stable and segment monomial ideals}\label{borel-fixed}
\spanishsubsection{Ideales monomials fijos de Borel, estables y segmentos}

Borel-fixed monomial ideals constitute a very important family, not only because of their connection to algebraic geometry, but also because of the relation to an important object: the generic initial ideal (\gin). The basic concepts involved are the following:
\begin{Definition}
The \emph{Borel group} is the subgroup $B_n(\kb)\subset GL_n(\kb)$ consisting of the upper triangular matrices.

For any invertible matrix $g=(g_{i,j})\in GL_n(\kb)$ we define the following action on a polynomial $f\in\kb[x_1\dots,x_n]$: 
$$g\cdot f=f(gx_1,\dots,gx_n),\qquad \mbox{where}\qquad gx_i=\sum_{j=1}^n g_{ij}x_j$$
given an ideal $I$, we get a new ideal by applying the action of $g$:
$$g\cdot I=\langle g\cdot f\vert f\in I\rangle$$
The ideals fixed by the action of the Borel group on them are called \emph{Borel-fixed} ideals.
\end{Definition}
\begin{Proposition}[\cite{MS04}, Proposition 2.3]
The following are equivalent for a monomial ideal $I$ when $char \kb=0$:
\begin{enumerate}
\item $I$ is Borel-fixed
\item If $m\in I$ is any monomial divisible by $x_j$, then $m\frac{x_j}{x_i}\in I$ for all $i<j$
\end{enumerate}
\end{Proposition}
\begin{Remark}
Ideals satisfying condition $2$ on the preceding proposition are called \emph{strongly stable}. The notion of stability for monomial ideals was introduced by \cite{EK90} as those ideals such that for every monomial $m\in I$, if we let $max(m)$ be the largest index of a variable dividing $m$ then $x_im/x_{max(m)}$ is in $I$ for all $i<max(m)$. When checking whether a monomial ideal is (strongly) stable, it is enough to verify the corresponding condition for the minimal generators of the ideal. In characteristic $0$ (strongly) stable and Borel-fixed are equivalent, in positive characteristic one only has that strongly stable implies Borel-fixed (for example $I=\langle x_1^p,x_2^p\rangle$ is Borel-fixed in characteristic $p$, but it is not stable). Minimal free resolutions of stable monomial ideals were described by Eliahou and Kervaire \cite{EK90}. 
\end{Remark}

Strongly related to Borel-fixed and stable ideals are segment ideals:
\begin{Definition}
Let $\tau$ be a term order. A monomial ideal is s $\tau$-segment if when $m_1,m_2$ are monomials of the same degree and $m_1>_{\tau}m_2$ and $m_2\in I$ then $m_1\in I$.
\end{Definition}

To check this condition is not enough to consider only the minimal generators of the ideal, although it is enough if $\tau$ is the $\lex$ term order.

Segment ideals are strongly stable, and in $\kb[x_1,x_2]$ all strongly stable ideals are segments, but if the number of variables is greater than $2$, strongly stable ideals are, in general, not segments. In relation with Mayer-Vietoris ideals we have that
\begin{Proposition}\label{stable-B2}
Stable ideals are Mayer-Vietoris of type $B2$
\end{Proposition}
\noindent {\bf Proof: } Let $I$ be a stable monomial ideal. We will construct a Mayer-Vietoris tree of $I$ with the following pivot selection strategy: Select as pivot monomial in the root node $I$ a generator $m_s$ having the biggest $max(m)$ among all minimal generators of $I$. This strategy will be followed in any node of dimension $0$. In the rest of the nodes the strategy is irrelevant as we will soon see.
Assume we are in a node $J$ of $MVT(I)$ such that $dim(J)=0$ and let $x^\mu$ be the pivot monomial of $J$ and $j=max(x^\mu)$, the biggest index of a variable in the support of $x^\mu$. Now, $\tilde J$ is the ideal $\langle x_i\cdot x^\mu \vert i< max(x^\mu)\rangle$. To see this, note that if $i<max(x^\mu)=j$ then, $x^\mu\frac{x_i}{x_j}$ is in $I$, since the ideal is stable. Then, there exists some generator $m^{\nu}$ of $J$ such that $m^\nu$ divides $x^\mu\frac{x_i}{x_j}$, i.e. $\nu_i\leq\mu_i+1$, $\nu_k\leq\mu_k\,\forall k\neq i,j$, and $\nu_j\leq\mu_j-1$. But we must then have $\nu_i=\mu_i+1$, otherwise $x^\nu$ would divide $x^\mu$ which is a contradiction. Then, $lcm(x^\nu,x^\mu)=x_ix^\mu$. The same situation holds for every $i<j$. The $lcm$ of $m^\mu$ with any other generator of $J$ is clearly divisible by some of the preceding ones, so $\tilde J=\langle x_i\cdot x^\mu \vert i< max(x^\mu)\rangle$. Now, the part of the tree hanging from $\tilde J$ is Mayer-Vietoris of type $B2$, whichever strategy we follow to select the pivot monomials, see lemma \ref{MVT-prime} which applies here since $\tilde J$ is isomorphic to a prime monomial ideal.

Now, assume there is a multidegree $x^\mu$ that appears in different nodes of $MVT(I)$. Then all relevant nodes in which $x^\mu$ appears, have the same dimension. To see this, assume on the contrary that there is some relevant node $MVT_p(I)$ of dimension $k$ and some other relevant node $MVT_q(I)$ of dimension $l>k$ such that $x^\mu$ appears in both of them. Then, by the construction of the tree, there exist two minimal generators of $x^\nu$ and $x^{\nu'}$ of $I$ such that $x^\mu=x_1\cdots x_{i_k}x^\nu$, $x^\mu=x_1\cdots x_{i_l}x^{\nu'}$ therefore $x^{\nu'}$ divides $x^\nu$ which is a contradiction. Then , since all of the appearances of $x^\nu$ appear at a given dimension, $I$ is Mayer-Vietoris of type $B2$ (see the proof of proposition \ref{repeated-relevant-multidegrees}).$\square$

\begin{Corollary}
Let $I=\langle m_1,\dots,m_r\rangle$ a stable ideal such that $max(m_i)=j_i$ then
$$\beta_0(I)=r,\qquad \beta_i(I)=\sum_{k=1}^r {{j_k -1}\choose{i}}\quad \forall i>0$$
Moreover, the multidegrees of the nonzero Betti numbers in degree $i$ are
$$\{x_{k_1}\cdots x_{k_i}\cdot m_s\vert x_{k_1},\dots,x_{k_i}\in \{1,\dots,j_s\}\}$$
\end{Corollary}
\noindent{\bf Proof: } The result comes from the observation of the Mayer-Vietoris trees of subsets of the variables, see lemma \ref{Betti-prime}. This result is already present in \cite{EK90}. $\square$

\begin{Remark}
A monomial ideal $I$ generated by squarefree monomials is called a \emph{squarefree stable ideal} \cite{AHH98} if for all $x^\mu\in I$, we have that $(x_jx^\mu)/x_{max(x^\mu)}$ is in $I$ for each $j<max(x^\mu)$ such that $x_j$ does not divide $x^\mu$. Squarefree stable monomial ideals are Mayer-Vietoris of type $B2$. The proof of this fact is just a direct translation of proposition \ref{stable-B2}. The corresponding formula for the Betti numbers of a squarefree monomial ideal $I=\langle m_1,\dots,m_r\rangle$ such that $max(m_i)=j_i$ is
$$\beta_0(I)=r,\qquad \beta_i(I)=\sum_{k=1}^r {{j_k -\vert supp(m_j)\vert}\choose{i}}\quad \forall i>0$$

\end{Remark}

\subsection{Generic monomial ideals}\label{generic}
\spanishsubsection{Ideales monomiales gen\'ericos}

Generic Monomial Ideals have been studied by E.Miller, B.Sturmfels and others, see \cite{MS04,MSY00} for example. They have nice combinatorial properties which make them interesting from a computational point of view and useful for applications \cite{GW04}. The definition is the following:
\begin{Definition}
A monomial ideal $I=\langle m_1,\dots, m_r\rangle$ is called  \emph{generic} if whenever two distinct minimal generators $m_i$ and $m_j$ have the same positive degree in some variable, there is a third generator $m_k$ which strictly divides $lcm(m_i,m_j)$. $I$ is called {\it strongly generic} if no two generators have the same nonzero degree in any variable.
\end{Definition}

For every monomial ideal $I=\langle m_1,\dots,m_r\rangle$ its {\it Scarf Complex}, $\Delta_I$ is the collection of all subsets of $\{m_1,\dots,m_r\}$ whose least common multiple is unique:
$$\Delta_I=\{\sigma\subseteq \{1,\dots, r\}\vert m_\sigma=m_\tau\Rightarrow\sigma=\tau\}$$ 
Let us call {\it Scarf monomial ideals} to those ideals $I$ such that the chain complex supported on the Scarf complex $\Delta_I$ minimally resolves $R/I$.
\begin{Remark}
The Scarf complex $\Delta_I$ is a simplicial complex of dimension at most $n-1$. It may be disconnected. Every free resolution of $R/I$ contains the chain complex supported on $\Delta_I$.
\end{Remark}
The important feature of the Scarf complex in relation with generic monomial ideals is that it provides their minimal free resolution:
\begin{Proposition}
If $I$ is a generic monomial ideal, then the chain complex supported on the Scarf Complex $\Delta_I$ minimally resolves $R/I$.
\end{Proposition}
With respect to Mayer-Vietoris trees, we have the following
\begin {Proposition}
If $I$ is a Scarf monomial ideal, then $I$ is a Mayer-Vietoris ideal of type $B1$.
\end{Proposition}
\noindent{\bf Proof:} 
Let $c$ be a generator of the Koszul homology of $I$, let $\ab$ be its multidegree. Then, $\ab$ is the multidegree of a generator of a relevant node in $MVT(I)$ and thus it is the least common multiple of a set of minimal generators of $I$, lets call this set $S$; we have then that $m_S=x^\ab$. Assume $\ab$ appears in some other relevant node, then there exists another set $T$ of minimal generators of $I$ such that $m_T=x^\ab$. If $I$ is Scarf, we know that the Scarf complex minimally resolves $I$ and then the multidegree of each generator of the Koszul homology of $I$ corresponds to the least common multiple of exactly one set of minimal generators of $I$. Hence, we have a contradiction, and $x^\ab$ appears only once in a relevant node of $MVT(I)$ $\square$

\begin{Corollary}
All generic monomial ideals are Mayer-Vietoris ideals of type $B1$ 
\end{Corollary}

Note that even if we can immediately read the multidegrees in which the Betti numbers of Scarf ideals are nonzero, the Mayer-Vietoris resolution of these ideals is in general non-minimal. On the other hand, if we are interested in the (multigraded) Betti numbers of generic ideals, Mayer-Vietoris trees are an alternative to the construction of their Scarf complex.

%
%



\subsection{Valla ideals}\label{valla}
\spanishsubsection{Ideales de Valla}

In his paper \cite{V04}, G. Valla studies the Betti numbers of some monomial ideals that define algebras that are the artinian reduction of the homogeneous coordinate ring of the scheme of two general fat points in $\PP^n$. These ideals are not stable and in \cite{V04}, formulas for their Betti numbers are given. If we want to study the multigraded Betti numbers of these ideals or their Koszul homology, we can use their Mayer-Vietoris trees. These ideals constitute an example of how to make assertions on a particular type of ideals from the examination of their Mayer-Vietoris trees. The ideals Valla studies can be described as follows:

$$I_{a,b}=\langle x_1^{a+b-2j}J^j,x_1^{a-t}J^t\rangle,\quad j=0,\dots,b-1,\quad t=b,\dots a$$
where $a$ and $b$ are positive integers such that $a\geq b$, and $J=\langle x_2,\dots, x_n\rangle$.

Valla ideals are zero dimensional, therefore their height is $n$ and their projective dimension is $n-1$. The size of these ideals grow very rapidly when $n$, $a$ and $b$ grow. In table \ref{growth-valla}, $n$ is the number of variables, $a$ and $b$ define the corresponding Valla ideal, $g$ is the number if generators and $Min$ is the size of the minimal resolution.

\begin {table}[h]
\begin {center}
    \begin{tabular}{|c|c|c|c|c|}
    	\hline
	$n$&$a$&$b$&$g$&Min\\
	\hline
    	3&3&2&10&31\\
	3&6&4&28&97\\
    	4&3&2&20&111\\
	4&6&4&84&545\\
	6&3&2&56&1023\\
	6&6&4&462&10625\\
	8&3&2&120&7423\\
	10&4&2&715&178177\\
	11&4&2&1001&471041\\
    	\hline
    \end{tabular}\caption{Size of Valla ideals and their minimal free resolutions.}\label{growth-valla}
\end{center}
\end{table}

The size of these ideals when the number of variables grows, makes the computation of their minimal resolutions and multigraded Hilbert series very hard to perform. Their (graded) Betti numbers are given by the formulas in \cite{V04}, one can also compute them using Mayer-Vietoris trees. In fact, these ideals have some Mayer-Vietoris trees that gives us the complete collection of multigraded Betti numbers. Moreover, the corresponding Mayer-Vietoris resolution is minimal, and the inspection of the corresponding Mayer-Vietoris trees can give us some insight on the structure of these ideals:

\begin{Proposition}\label{Valla_is_MVB2}
Valla ideals are Mayer-Vietoris ideals of type $B2$.
\end{Proposition}

\noindent{\bf Proof: } Valla ideals are of the form
$$I_{a,b}=\langle x_1^{a+b-2j}J^j,x_1^{a-t}J^t\rangle,\quad j=0,\dots,b-1,\quad t=b,\dots, a$$
for $a\geq b$.
We will sort the monomials generating $I_{a,b}$ in a way such that the construction of their Mayer-Vietoris trees using as pivot monomial always the first monomial, will prove that $I_{a,b}$ is a monomial ideal of type $B2$. First, we see that the listing of the minimal generators of $I_{a,b}$ consists in two blocks:
\begin{itemize}
\item In the first block, for each $j=0,\dots,b$, we have the generators $x_1^{a+b-2j}J^j$. We take these first, in the same order, using the ordering we saw in proposition \ref{J^k_is_MVB2} to sort the elements that have the same exponent in $x_1$. Note that the degree of each of these generators is $a+b-j$.
\item The second block has all the elements of the form $x_1^{a-t}J^t$ for $t=b,\dots, a$. We sort them in the same way, after the generators of the first block. These generators have all degree $a$.
\end{itemize}
Consider now a generator of  the first block. It has the form $x_1^{a+b-2j}x^\mu$ for some $x^\mu\in \kb[x_2,\dots,x_n]$ and $\vert\mu\vert=j$. The left-hand tree hanging from this node has as generators all monomials of the form  $x_1^{a-t}x^\mu x_k$ for $k=2,\dots,n$. This is because the exponent of $x_1$ is maximal among the remaining generators, and thus it is unchanged. And from the sorting of the generators we see that for the generators of the first block, there is always an element below them with exponent $\mu_j+1$ in the variable $x_j$, $j=2,\dots,n$. All these subtrees are unmixed, and the generators in them do not appear in any other subtree with a different exponent in $x_1$, moreover, their appearances in different subtrees with the same exponent in $x_1$ occur all in the same dimension, for they have the same degree. This applies to the generators in the second block, for the values of $t$ up to $a-1$. These generators contribute with ${{n-1}\choose{d}} \cdot {{n+j-2}\choose{j}}$ to the $d$th Betti number in the case of the generators of the first block, and with ${{n-1}\choose{d}} \cdot  {{n+t-2}\choose{t}}$ in the case of the considered generators of the second block.

For the generators of the form $J^a$ with $J=\langle x_2,\dots,x_n\rangle$, we follow the same argument we saw in proposition \ref{J^k_is_MVB2}. The generators in the relevant nodes of these subtrees do not interfere with the previous ones, again because of the zero exponent in the first variable. These generators contribute with $\beta_d(J^a)$ to the $d$-th Betti number of $I_{a,b}$.

All these consideration prove proposition \ref{Valla_is_MVB2}. Therefore, when studying these ideals, Mayer-Vietoris trees provide an efficient way to compute their multigraded Betti numbers without computing their minimal free resolutions. Moreover, the Mayer-Vietoris resolution of a Valla ideal is minimal. $\square$
\begin{Corollary}
The (graded) Betti numbers of the Valla ideal $I_{a,b}\subseteq\kb[x_1,\dots,x_n]$ are given by:
$$\beta_i(I_{a,b})=\sum_{j=0}^{a-1}{{n-1}\choose{i}}{{n+j-2}\choose{j}}+\beta_i(J^a)$$
$$\beta_{i,j}(I_{a,b})={{n-1}\choose{i}}{{n+i+a+b-j-2}\choose{i+a+b-j}}+\beta_{i,j}(J^a)\qquad a+i\leq j\leq a+b+i$$
\end{Corollary}

\noindent{\bf Proof: }The inspection of the Mayer-Vietoris tree constructed in proposition \ref{Valla_is_MVB2} easily gives the result. Observe that $\beta_{i,j}(J^a)$ is zero if $j\neq a+i$ and is ${{n+a-1}\choose{a+i}}{{a+i-1}\choose{i}}$ if $j=a+i$, therefore our formulas are equivalent to the ones in \cite{V04}. $\square$

\begin{Corollary}
An irredundant irreducible decomposition of the Valla ideal $I_{a,b}\subseteq\kb[x_1,\dots,x_n]$ is given by
$$I_{a,b}=\bigcap_{\substack{j=0,\dots,b-1\\ \mu_2+\cdots+\mu_n=j+2}}\langle x_1^{a+b-2j},x_2^{\mu_2},\cdots ,x_n^{\mu_n}\rangle$$
\end{Corollary}

\noindent{\bf Proof: } Valla ideals are zero dimensional, therefore, their irredundant irreducible decomposition is given by $\mm^\mu$ with $\mu$ the multidegree of a maximal corner (see section \ref{irreducible_decomposition}). The part of the Mayer-Vietoris tree corresponding to the second block of the ideal has multidegrees in the variables $x_2,\dots,x_n$ and then no multidegree of this part corresponds to a maximal corner. Thus, the only multidegrees we are interested in are those corresponding to $(n-1)$-st homology. The inspection of the tree constructed in proposition \ref{Valla_is_MVB2} gives the result. $\square$

Valla ideals constitute an example of how Mayer-Vietoris trees may help to perform a close examination of a monomial ideal in order to read its Betti numbers and other properties. In this case we have shown a simple way to obtain the results given in \cite{V04} and even more. We finish with an example

\begin{Example}\label{example-valla}
Consider the ideal $I_{3,2}\subseteq\kb[x,y,z]$, its minimal generating set is given by
$$I_{3,2}=\langle x^5,x^3y,x^3z,xy^2,xyz,xz^2,y^3,y^2z,yz^2,z^3\rangle$$
The list of relevant nodes of $MVT(I_{3,2})$ is the following (the node at position $0$ is given by just the minimal generating set of the ideal):
\begin{itemize}
\item[-] Nodes of dimension $1$:
	\begin{itemize}
	\item[]First block:

	$\begin{array}{ll}
	MVT_{2}(I_{3,2})=\langle x^5y,x^5z \rangle&MVT_{6}(I_{3,2})=\langle x^3y^2,x^3yz \rangle\\
	MVT_{14}(I_{3,2})=\langle x^3z^2,x^3yz \rangle&MVT_{30}(I_{3,2})=\langle xy^3,xy^2z \rangle\\
	MVT_{62}(I_{3,2})=\langle xz^3,xyz \rangle&MVT_{126}(I_{3,2})=\langle xy^2z,xyz^2 \rangle\\
	\end{array}$

	\item[]Second block:

	$\begin{array}{ll}
	MVT_{254}(I_{3,2})=\langle y^3z \rangle&MVT_{510}(I_{3,2})=\langle yz^3 \rangle\\
	MVT_{1022}(I_{3,2})=\langle y^2z^3 \rangle&\\
	\end{array}$
	\end{itemize}
\item[-] Nodes of dimension $2$:
	\begin{itemize}
	\item[]First block:

	$\begin{array}{ll}
	MVT_{4}(I_{3,2})=\langle x^5yz \rangle&MVT_{12}(I_{3,2})=\langle x^3y^2z \rangle\\
	MVT_{28}(I_{3,2})=\langle x^3yz^2 \rangle&MVT_{60}(I_{3,2})=\langle xy^3z \rangle\\
	MVT_{124}(I_{3,2})=\langle xyz^3 \rangle&MVT_{252}(I_{3,2})=\langle xy^2z^2 \rangle\\
	\end{array}$
	\end{itemize}
\end{itemize}
Then from the tree, we have that
$$\beta_0(I_{3,2})=10,\qquad \beta_1(I_{3,2})=15,\qquad \beta_2(I_{3,2})=6$$
and the corresponding graded version:
$$\beta_{0,5}(I_{3,2})=1,\qquad \beta_{0,4}(I_{3,2})=2,\qquad \beta_{0,3}(I_{3,2})=7$$
$$\beta_{1,6}(I_{3,2})=2,\qquad \beta_{1,5}(I_{3,2})=4,\qquad \beta_{1,4}(I_{3,2})=9$$
$$\beta_{2,7}(I_{3,2})=1,\qquad \beta_{2,6}(I_{3,2})=2,\qquad \beta_{2,5}(I_{3,2})=3$$
Finally, the irredundant irreducible decomposition of $I_{3,2}$ is given by
$$I_{3,2}=\langle x^5,y,z\rangle\cap\langle x^3,y^2,z\rangle\cap\langle x^3,y,z^2\rangle\cap\langle x,y^3,z\rangle\cap\langle x,y,z^3\rangle\cap\langle x,y^2,z^2\rangle$$
\end{Example}

\subsection{Ferrers ideals}\label{ferrers}
\spanishsubsection{Ideales de Ferrers}

A \emph{Ferrers graph} $G$ is a bipartite graph on two distinct vertex sets $X=\{x_1,\dots,x_n\}$ and $Y=\{y_1,\dots,y_m\}$ such that if $(x_i,y_j)$ is an edge of $G$, then so is $(x_p,y_q)$ for $1\leq p\leq i$, $1\leq q\leq j$. In addition, $(x_1,y_m)$ and $(x_n,y_1)$ are required to be edges of $G$. A \emph{Ferrers ideal} is the edge ideal associated with a Ferrers graph. These ideals are studied in \cite{CN06, CN07}.
\begin{Remark}
For any Ferrers graph $G$ there is an associated sequence of non-negative integers $\lambda=(\lambda_1,\dots,\lambda_n)$, where $\lambda_i$ is the degree of the vertex $x_i$; since $\lambda_1=m\geq\lambda_2\geq\cdots\geq\lambda_n\geq1$, then $\lambda$ is a partition.
We can also associate to a Ferrers graph, a \emph{Ferrers tableau}, i.e. an array of $n$ rows of cells with $\lambda_i$ adjacent cells, left justified, in the $i$-th row.
Ferrers graphs and tableaux have been studied in relation to chromatic polynomials, Schubert varieties, permutation statistics, quantum mechanical operators, rook theory, etc \cite{B04,D97,EW04}.
\end{Remark}
In \cite{CN06,CN07}, resolutions of Ferrers ideals are studied, and more general ideals are studied in terms of Ferrers ideals. The authors describe in \cite{CN06} a cellular minimal free resolution for Ferrers ideals. In \cite{CN07} monomial ideals generated in degree two are studied using the results about Ferrers ideals. Here, we show that Mayer-Vietoris ideals give an alternative way for the computation of the multigraded Betti numbers of  such ideals and their minimal free resolutions. In particular, we proof the folowing

\begin{Proposition}\label{Ferrers-MV}
Ferrers ideals are Mayer-Vietoris ideals of type $A$
\end{Proposition}

\noindent {\bf Proof of proposition \ref{Ferrers-MV}:} 
Let $I_\lambda$ the Ferrers ideal associated to a partition $\lambda=(\lambda_1,\dots,\lambda_n)$, i.e. $I=\langle x_1y_1,\dots, x_1y_{\lambda_1},\dots,x_ny_1,\dots,x_ny_{\lambda_n}\rangle$. Consider the lexicographic Mayer-Vietoris tree $MVT_{\lex}(I_\lambda)$ of it i.e. we take as pivot monomial in each node the last generator according to \lex order.
We have that $MVT_2(I_\lambda)=I_{\lambda'}$ where $\lambda'=(\lambda_1,\dots,\lambda_{n-1},\lambda_n-1)$ if $\lambda_n>1$ or $(\lambda_1,\dots,\lambda_{n-1})$ if $\lambda_n=1$. On the other hand, it is easy to see that $MVT_3(I_\lambda)$ is the ideal $\langle x_1x_ny_{\lambda_n},x_2x_ny_{\lambda_n}\dots,x_{n-1}y_{\lambda_n},x_ny_1y_{\lambda_n},\dots,x_ny_{\lambda_{n-1}}y_{\lambda_n}\rangle$. This ideal is just the ideal generated by the subset of the variables $x_1,\dots,x_{n-1},y_1,\dots,y_{\lambda_{n-1}}$, with each generator multiplied by $x_ny_{\lambda_n}$; therefore, this ideal is Mayer-Vietoris of type $A$, i.e. this subtree has no repeated generators in the relevant nodes; moreover, it is an unmixed subtree.

We assume without loss of generality that $\lambda_n>1$. Every generator in the left-hand tree, is a multiple of $x_ny_{\lambda_n}$ as we have just seen. The right hand tree is also unmixed. All the nodes of dimension one in the right hand tree are of the form $\langle x_1x_my_{p},x_2x_my_{p}\dots,x_{m-1}x_my_{p},x_my_1y_{p},\dots,x_my_{p-1}y_{p}\rangle$ with $(m,p)<(n,\lambda_n)$ in the usual sense. Therefore, none of them is a multiple of $x_ny_{\lambda_n}$ and thus no generator of the relevant nodes in the lefthand tree is present as a generator of a relevant node in the righthand tree. Proceeding recursively yields the result.$\square$
\begin{Example}
Let $\lambda=(3,2,2)$, we have the following tree
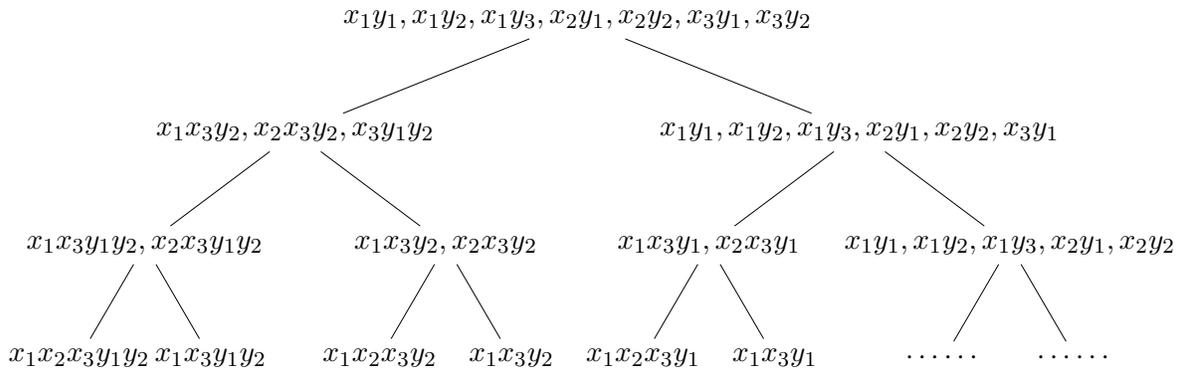
\begin{small}
\begin{center}
\begin{tikzpicture}[scale=1,very thin]
\tikzstyle{level 1}=[sibling distance=7.5cm]
\tikzstyle{level 2}=[sibling distance=4cm]
\tikzstyle{level 3}=[sibling distance=1.75cm]
\node{$x_1y_1,x_1y_2,x_1y_3,x_2y_1,x_2y_2,x_3y_1,x_3y_2$}
	child{node{$x_1x_3y_2,x_2x_3y_2,x_3y_1y_2$}
		child{node{$x_1x_3y_1y_2,x_2x_3y_1y_2$}
			child{node{$x_1x_2x_3y_1y_2$}}
			child{node{$x_1x_3y_1y_2$}}}
	child{node{$x_1x_3y_2,x_2x_3y_2$}
		child{node{$x_1x_2x_3y_2$}}
		child{node{$x_1x_3y_2$}}}}
child{node{$x_1y_1,x_1y_2,x_1y_3,x_2y_1,x_2y_2,x_3y_1$}
	child{node{$x_1x_3y_1,x_2x_3y_1$}
		child{node{$x_1x_2x_3y_1$}}
		child{node{$x_1x_3y_1$}}}
child{node{$x_1y_1,x_1y_2,x_1y_3,x_2y_1,x_2y_2$}
	child{node{$\cdots\cdots$}}
child{node{$\cdots\cdots$}}}};
\end{tikzpicture}
\end{center}
\end{small}
\end{Example}

As we have just seen, \lex-Mayer Vietoris trees of Ferrers ideals are not only unmixed but $2$-linear. Therefore, Ferrers ideals have a $2$-linear minimal free resolution. These trees provide us with a nice formula for the Betti numbers of Ferrers ideals in the following way: 
Let $I_\lambda$ be a Ferrers ideal. In the construction of the Mayer-Vietoris tree, we have that the nodes in dimension $0$ are $I_{(1)}=I_{\lambda^1}<\dots<I_{\lambda^r}=I_\lambda$ , where each $\lambda^s$ is obtained from $\lambda^{s+1}$ by subtracting $1$ from the rightmost integer. Each of them has a lefthand tree (therefore of dimension $1$) generated by $(j-1)+(i-1)$ generators, where $j$ ranges from $1$ to the number of nonzero integers $\lambda^s_i$, and $i$ ranges from $1$ to the rightmost $\lambda^s_i$. Since the ideals in these subtrees are all isomorphic to ideals generated by $(j-1)+(i-1)$ different variables, the $k$-th Betti number is ${j-1+i-1}\choose {k-1}$, and because of the shift in dimension, each of these ideals contribute with ${j-1+i-1}\choose {k}$ to the $k$-th Betti number of $I_\lambda$. Therefore:

$$\beta_k(I_\lambda)=\sum_{j=1}^n\sum_{i=1}^{\lambda_j}{{j+i-2}\choose{k}}\quad \mbox{for }0\leq k\leq max\{\lambda_i+i-2\}$$

Observe that all the $k$-th Betti multidegrees have degree $k+2$, and therefore $\beta_{i,j}(I_\lambda)=\beta_i(I_\lambda)$ if $j=i+2$ and zero otherwise. This formula is equivalent to the one given in \cite{CN06} and similarly simple, although the one we propose does not contain any minus sign. From this formula, we can derive the relevant invariants that can be read off the Betti numbers:
\begin{Proposition}[see \cite{CN06}, Corollary 2.2]
The projective dimension of a Ferrers ideal $I_\lambda$ is $max_j\{\lambda_j+j-2\}$ and the Castelnuovo-Mumford regularity $reg(I_\lambda)$ is equal to $2$.
\end{Proposition}

Now we turn to finding explicit formulas for irreducible and primary decompositions of Ferrers ideals. For a Ferrers ideal $I_\lambda$ let us consider the following set of the variables: $\{i \in\{1,\dots,n\} \vert \lambda_i\neq\lambda_{i-1}\}$. Then we have the following decomposition

\begin{equation}\label{ferrers-decom}
I_\lambda=\langle x_1,\dots,x_n\rangle \bigcap_{\{i \in\{1,\dots,n\} \vert \lambda_i\neq\lambda_{i-1}\}} \langle x_1,\dots,x_{i-1},y_1,\dots,y_{\lambda_i}\rangle
\end{equation}

\begin{Proposition}[see \cite{CN06}, Corollary 2.2 and Corollary 2.5]\label{ferrers-primary-irreducible}
The decomposition (\ref{ferrers-decom}) of the Ferrers ideal $I_\lambda$ is both a primary irredundant decomposition and an irreducible irredundant decomposition of $I_{\lambda}$.
The height of the edge ideal $I_{\lambda}$ of a Ferrers graph $G$ is $min\{min_j\{\lambda_j+j-1\},n\}$.
\end{Proposition}

\noindent{\bf Proof:}
An easy proof is based on the free directions at the monomial $x_1\cdots x_ny_1\cdots y_{\lambda_n}$: Since our ideals are squarefree, every multidegree $\mu$ of a Koszul generator (i.e. every Betti multidegree) is also squarefree, therefore, $\overline{\mu}$ is $1$ for every such multidegree. Considering the maximal cones of free directions at $1$ (equivalently the facets of the lower Koszul simplicial complex $\Delta^I_{x_1\cdots x_ny_1\cdots y_{\lambda_n}}$) we have the result. In fact, such maximal corners include $[x_1,\dots,x_n]$ and also $[y_1,\dots,y_n]$. For the other cones, observe that if $\lambda_i=j$ then $[x_i,\dots,x_n,y_{j+1},\dots,y_{\lambda_1}]$ is a cone of true free directions. Among these, the maximal ones correspond to those $i$ such that $\lambda_i\neq\lambda_{i-1}$, then the $\underline{\mu}$ corresponding to these maximal cones are exactly the factors in decomposition (\ref{ferrers-decom}). $\square$
\begin{Remark}
The observation of the Mayer-Vietoris tree of the artinian closure of $I_\lambda$ gives us the form of the maximal corners of it. Deleting all indices which have an exponent two in these corners gives us an alternative proof of proposition \ref{ferrers-primary-irreducible}.
\end{Remark}

\begin{Example}
Let $\lambda=(5,3,3,3,2,1)$. The set $\{i \in\{1,\dots,n\} \vert \lambda_i\neq\lambda_{i-1}\}$ is $\{1,2,5,6\}$ and therefore we have that the following decomposition:

\begin{eqnarray*}
I_\lambda=&\langle x_1,x_2,x_3,x_4,x_5,x_6\rangle\cap\langle y_1,y_2,y_3,y_4,y_5\rangle \cap\langle x_1,y_1,y_2,y_3\rangle\\
&\cap\langle x_1,x_2,x_3,x_4,y_1,y_2\rangle\cap\langle x_1,x_2,x_3,x_4,x_5,y_1\rangle
\end{eqnarray*}

is an irredundant irreducible and primary decomposition.
\end{Example}

\subsection{Quasi-stable ideals}
\spanishsubsection{Ideales quasi-estables}

Quasi-stable ideals appear in the context of involutive bases, in particular in relation with Pommaret bases. A complete treatment of this class of ideals is given in \cite{S02a,S02b}. In this section we briefly introduce Pommaret bases and quasi-stable ideals, together with a rough exposition of their context and main properties. Finally we show how the knowledge of (part of) the Koszul homology of a quasi-stable ideal $I$ allows us to complete the minimal generating set of $I$ to a Pommaret basis.

\subsubsection{Pommaret bases}
\spanishsubsubsection{Bases de Pommaret}

Pommaret bases are a type of involutive bases, which are some special non-reduced Gr\"obner bases which posses additional combinatorial properties. Their origin lies in the Janet-Riquier theory of differential systems and they play a very interesting role in a certain homological approach to the formal theory of differential systems (see \cite{S07} and section \ref{formal_theory} below).

Involutive division is introduced in the abelian monoid $(\NN^n,+)$ of multi indices. For a multi index $\nu\in\NN^n$, we call $\Cc(\nu)$ to the cone of $\nu$, $\Cc(\nu)=\nu+\NN^n$, and given a finite subset $\Nc\subset\NN^n$, we define its span $\langle \Nc\rangle=\cup_{\nu\in\Nc}\Cc(\nu)$. We also write for a subset $N$ of $\{1,\dots,n\}$,  $\NN^n_N=\{\nu\in\NN^n\vert \forall j\notin N, \nu_j=0\}$.

\begin{Definition}
An \emph{involutive division} $L$ is defined on $(\NN^n,+)$ if for any finite set $\Nc\subset\NN^n$, a subset $N_{L,\Nc}(\nu)\subseteq\{1,\dots,n\}$ of \emph{multiplicative indices} is associated to every multi index $\nu\in\Nc$ such that the following two conditions on the \emph{involutive cones} $\Cc_{L,\Nc}(\nu)=\nu+\NN^n_{N_{L,\Nc}(\nu)}$ are satisfied:
\begin{enumerate}
\item If there exist two elements $\mu,\nu\in\Nc$ with $\Cc_{L,\Nc}(\mu)\cap\Cc_{L,\Nc}(\nu)\neq\emptyset$, either $\Cc_{L,\Nc}(\mu)\subseteq\Cc_{L,\Nc}(\nu)$ or $\Cc_{L,\Nc}(\nu)\subseteq\Cc_{L,\Nc}(\mu)$ holds.
\item If $\Nc'\subset\Nc$ then $N_{L,\Nc}(\nu)\subseteq N_{L,\Nc'}(\nu)$ for all $\nu\in\Nc'$.

An arbitrary multi index $\mu\in\NN^n$ is \emph{involutively divisible} by $\nu\in\Nc$, written $\nu\vert_{L,\Nc}\mu$ if $\mu\in\Cc_{L,\Nc}(\nu)$.
\end{enumerate}
\end{Definition}

Involutive divisions introduce restrictions of the cone of a multi index in order to obtain the span of a subset of $\NN^n$ as a \emph{disjoint} union of such restricted cones. Two particular involutive divisions play an important role in the literature, namely \emph{Janet} division \cite{J29} and \emph{Pommaret} division \cite{J20}. The first one has been studied in \cite{GB98,GBY01} etc, and the second one in \cite{A95,M98,S02a,S02b,S07b} among others.

Pommaret division $P$ assigns the multiplicative indices in the following simple way: If $1\leq k\leq n$ is the smallest index such that $\nu_k>0$ for some $\nu\in\NN^n$, then we call $k$ the \emph{class} of $\nu$, $cls(\nu)$, and set $N_{P,\Nc}(\nu)=\{1,\dots,k\}$. We also define $N_{P,\Nc}([0,\dots,0])=\{1,\dots,n\}$. Pommaret division is \emph{globally defined} in the sense that the assignment of the multiplicative indices is independent of the set $\Nc$.

\begin{Definition}
The involutive span with respect to an involutive division $L$ of a finite set $\Nc\subset \NN^n$ is
\begin{equation}\label{involutive-span}
\langle \Nc\rangle_L=\bigcup_{\nu\in\Nc}\Cc_{N,L}(\nu)
\end{equation}
The set $\Nc$ is a \emph{weak involutive basis} of the monoid ideal $\langle \Nc\rangle$ if
$\langle \Nc\rangle_L=\langle \Nc\rangle$ (in general only $\langle \Nc\rangle_L\subseteq\langle \Nc\rangle$ holds). A weak involutive basis is an involutive basis if the union of the right hand side of (\ref{involutive-span}) is disjoint. We call any finite set $\Nc\subseteq\overline{\Nc}\subset\NN^n$ such that $\langle \overline{\Nc}\rangle_L=\langle \Nc\rangle$ an \emph{involutive completion} of $\Nc$.
\end{Definition}

\begin{Definition}
Let $I\subseteq\kb[x_1,\dots,x_n]$ a polynomial ideal. A finite subset $\Hc\subset I$ is a \emph{weak involutive basis} of $I$ for an involutive division $L$ on $\NN^n$ if its leading exponents $le_<\Hc$ form a weak involutive basis of the monoid ideal $le_<I$. The subset $\Hc$ is an involutive basis of $I$ if $le_<\Hc$ is an involutive basis of $le_<I$ and no two distinct elements of $\Hc$ have the same leading exponents.
\end{Definition}

Observe that this definition implies that any (weak) involutive basis is a Gr\"obner basis.
\begin{Example}
The following two diagrams represent the difference between the usual and involutive cones (with respect to Pommaret division) of the same ideal. In the second diagram the arrows represent the Pommaret multiplicative directions:

\begin{figure}[h]
\begin{center}
$\begin{array}{l@{\hspace{1in}}r}
\begin{tikzpicture}[scale=1.5, fill opacity=0.6, text opacity=1]
\fill[blue] (0,1.5)rectangle(2.4,2.4);
\fill[blue] (1,1)rectangle(2.4,1.5);
\fill[blue] (1.5,0.5)rectangle(2.4,2.4);
\draw [step=.5cm, very thin](0,0) grid (2.4,2.4);
\fill (0,1.5) circle(1pt) node[left]{$y^3$};
\fill (1.5,0.5) circle(1pt) node[below]{$x^3y$};
\fill (1,1) circle(1pt) node[above left]{$x^2y^2$};
\end{tikzpicture}

&

\begin{tikzpicture}[scale=1.5, fill opacity=0.6, text opacity=1]
\fill[blue] (0,1.5)rectangle(2.4,2.4);
\fill[blue] (1,1)rectangle(2.4,1.1);
\fill[blue] (1.5,0.5)rectangle(2.4,0.6);
\draw [step=.5cm, very thin](0,0) grid (2.4,2.4);
\fill (0,1.5) circle(1pt) node[left]{$y^3$};
\fill (1.5,0.5) circle(1pt) node[below]{$x^3y$};
\fill (1,1) circle(1pt) node[above left]{$x^2y^2$};

\draw[->, line width = 1.5pt](0,1.5)--(0.25,1.5);\draw[->,line width = 1.5pt](0,1.5)--(0,1.75);
\draw[->,line width = 1.5pt](1.5,0.5)--(1.75,0.5);
\draw[->,line width = 1.5pt](1,1)--(1.25,1);
\end{tikzpicture}

\end{array}$
\end{center}
\caption{Usual and involutive cones for $I=\langle x^3y,x^2y^2,y^3)$}
\label{example-pommaret}
\end{figure}
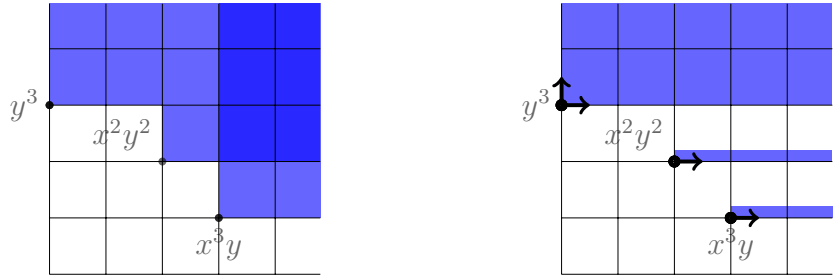

In figure \ref{example-pommaret-2} we show the difference between a quasi-stable ideal and a non quasi-stable ideal. Again, the arrows represent the Pommaret multiplicative directions. The white dots represent elements of the Pommaret basis that are not minimal generators of the ideal. Observe that in the second case, every element of the form $xy^k$ for $k\geq 3$ is an element of the Pommaret basis.

\begin{figure}[h]
\begin{center}
$\begin{array}{l@{\hspace{1in}}r}
\begin{tikzpicture}[scale=1.5, fill opacity=0.6, text opacity=1]
\fill[blue] (0,1.5)rectangle(2.4,2.9);
\fill[blue] (1.5,1)rectangle(2.4,1.1);
\fill[blue] (1.5,0.5)rectangle(2.4,0.6);
\draw [step=.5cm, very thin](0,0) grid (2.4,2.9);
\fill (0,1.5) circle(1pt) node[left]{$y^3$};
\fill (1.5,0.5) circle(1pt) node[below]{$x^3y$};
\filldraw[draw=black, fill=white] (1.5,1) circle(1.3pt) node[above left, text=black!50]{$x^3y^2$};

\draw[->, line width = 1.5pt](0,1.5)--(0.25,1.5);\draw[->,line width = 1.5pt](0,1.5)--(0,1.75);
\draw[->,line width = 1.5pt](1.5,0.5)--(1.75,0.5);
\draw[->,line width = 1.5pt](1.5,1)--(1.75,1);
\end{tikzpicture}
&

\begin{tikzpicture}[scale=1.5, fill opacity=0.6, text opacity=1]
\fill[blue] (0.5,1.5)rectangle(2.4,1.6);
\fill[blue] (0.5,2)rectangle(2.4,2.1);
\fill[blue] (0.5,2.5)rectangle(2.4,2.6);
\fill[blue] (1.5,1)rectangle(2.4,1.1);
\fill[blue] (1.5,0.5)rectangle(2.4,0.6);
\draw [step=.5cm, very thin](0,0) grid (2.4,2.9);
\fill (0.5,1.5) circle(1pt) node[left]{$xy^3$};
\fill (1.5,0.5) circle(1pt) node[below]{$x^3y$};
\filldraw[draw=black, fill=white] (1.5,1) circle(1.3pt) node[above left, text=black!50]{$x^3y^2$};
\filldraw[draw=black, fill=white] (0.5,2) circle(1.3pt) node[above left, text=black!50]{$xy^4$};
\filldraw[draw=black, fill=white] (0.5,2.5) circle(1.3pt) node[above left, text=black!50]{$xy^5$};

\draw[->, line width = 1.5pt](0.5,2)--(0.75,2);
\draw[->, line width = 1.5pt](0.5,1.5)--(0.75,1.5);
\draw[->, line width = 1.5pt](0.5,2.5)--(0.75,2.5);
\draw[->,line width = 1.5pt](1.5,0.5)--(1.75,0.5);
\draw[->,line width = 1.5pt](1.5,1)--(1.75,1);
\end{tikzpicture}

\end{array}$
\end{center}
\caption{Ideals with and without finite Pommaret bases.}
\label{example-pommaret-2}
\end{figure}
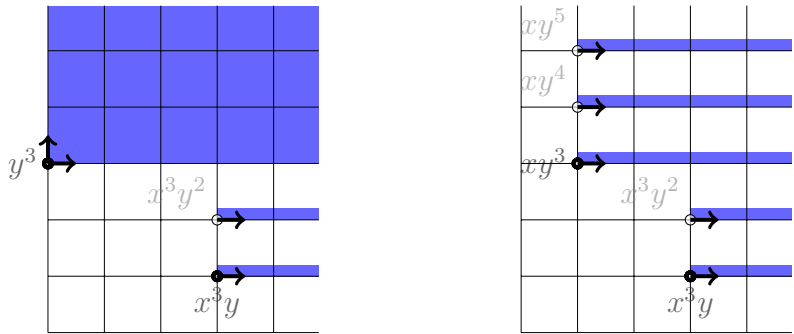

\end{Example}

\subsubsection{Quasi-stable ideals}\label{quasi-stable}
\spanishsubsubsection{Ideales quasi-estables}

Not every monoid ideal in $\NN^n$ possesses a finite Pommaret basis, this implies that there are polynomial ideals without a finite Pommaret basis for a given term order. The coordinates in which a Pommaret basis exists are called \emph{$\delta$-regular}. Since generic coordinates are $\delta$-regular, the problem of having a finite Pommaret basis at the level of polynomial ideals can be solved using generic changes of coordinates. In the monomial case those ideals having a finite Pommaret basis deserve an special name:

\begin{Definition}
A monomial ideal $I\subseteq \kb[x_1,\dots,x_n]$ is called \emph{quasi-stable} if it possesses a finite Pommaret basis.
\end{Definition}

\begin{Remark}
The reason for this terminology is the particular relation of stable monomial ideals (or monomial submodules in general) with Pommaret bases: A monomial ideal is stable if and only if its minimal generating set is simultaneously a Pommaret basis. Quasi-stable ideals are those such that their minimal generating set is not yet a Pommaret basis, but can be completed to a Pommaret basis of the ideal.
\end{Remark}

There are several algebraic characterizations of quasi-stable ideals, which are independent of the theory of Pommaret bases:

\begin{Proposition}[\cite{S02b}, Prop. 4.3]
Let $I$ be a monomial ideal with $dim(R/I)=D$. Then the following six statements are equivalent:
\begin{enumerate}
\item[i)] $I$ is quasi-stable.
\item[ii)] The variable $x_1$ is not a zero divisor for $R/I^{sat}$ and for all $1\leq j<D$ the variable $x_{j-1}$ is not a zero divisor for $R/\langle I,x_j,\dots,x_n\rangle^{sat}$ ($I^{sat}$ is the saturation of $I$).
\item[iii)] We have $I:x_n^\infty\subseteq I:x_{n-1}^\infty\subseteq\cdots\subseteq I:x_D^\infty$ and for all $n<j\leq D$ an exponent $k_j\geq 1$ exists such that $x_j^{k_j}\in I$.
\item[iv)] For all $1\leq j\leq n$ the equality $I:x_j^\infty=I:\langle x_1,\dots,x_j\rangle^\infty$ holds.
\item[v)] For every associated prime ideal $\mp\in Ass(R/I)$ an integer $1\leq j \leq n$ exists such that $\mp=\langle x_1,\dots,x_j\rangle$.
\item[vi)] If $x^\mu\in I$ and $\mu_i>0$ for some $1\leq i < n$, then for each $0<r\leq \mu_i$ and $i<j\leq n$ an integer $s\geq 0$ exists such that $x^{\mu-r_i+s_j}\in I$.
\end{enumerate}

\end{Proposition}
Ideals satisfying conditions $ii$ and $iv$ are called \emph{monomial ideals of nested type} in \cite{BG06} where the equivalence between the two statements is proved, Bermejo and Gimenez proved in fact the equivalence of $ii$, $iii$, $iv$ and $v$. Ideals satisfying conditions $iv$ and $vi$ are called \emph{ideals of Borel type} in \cite{HPV03} where the corresponding equivalence is proved.

\begin{Remark}
Given a quasi-stable monomial ideal and a Pommaret basis of it, a free resolution of $I$ is constructed in \cite{S02b}. This resolution is minimal if and only if the ideal is stable, and it is identical to the one given in \cite{EK90}.
\end{Remark}

\subsubsection{Monomial completion via Koszul homology}
\spanishsubsubsection{Compleci\'on monomial via homolog\'ia de Koszul}

Given a monomial ideal $I$ the question of deciding whether it is quasi-stable and in the affirmative case completing its minimal generating set to a Pommaret basis can be answered using Koszul homology. The opposite problem, namely obtaining the Koszul homology of a quasi-stable ideal $I$ given a Pommaret basis of it has been treated in \cite{Sa07}.

Given a monomial ideal and the set of multidegrees $\mu\in\NN^n$ such that $H_{i,\mu}(\KK(I))\neq 0$ for some $i$, we want to construct a Pommaret basis of $I$. We follow the following steps:
\begin{itemize}
\item Let $\hat{I}$ be the artinian closure of $I$ and take the maximal and closed corners of $\hat{I}$. Recall that a multidegree $\mu$ is a \emph{closed corner} if $\vert supp(\mu)\vert=k$ and $H_{k-1,\mu}(\KK(\hat{I}))\neq 0$, in case $k=n$ it is a \emph{maximal corner}.
\item From each such multidegree, take one step back in each nonmultiplicative directions, i.e. take $\nu$ such that $\nu\cdot \vec{k}=\mu$, where $\vec{k}$ is the set of nonmultiplicative directions of $\nu$. Now, we put into the Pommaret basis of $\hat{I}$ all nonmultiplicative divisors of these $\nu$. These form a Pommaret basis of $\hat{I}$ (every artinian ideal is quasi-stable). To obtain a Pommaret basis of $I$ we just delete those elements of the Pommaret basis of $\hat{I}$ having an exponent in $x_i$ bigger than $\lambda_i$ for some $i$.
\item For those such multidegrees $\nu$ obtained in the preceding step which are in the boundary of $I$ and have free directions in $I$, we look whether there is some nonmultiplicative free directions for them, then we are in $\delta$-singular coordinates, i.e. $I$ is not quasi-stable. Otherwise, the ideal is quasi-stable.
\end{itemize}

\begin{Proposition}
Let $I$ be a monomial ideal. The above procedure decides whether $I$ is quasi-stable, and in the affirmative case it gives a Pommaret basis of $I$. 
\end{Proposition}

\noindent{\bf Proof:} All elements of the Pommaret basis of $I$ lie on the boundary of the ideal. Then starting form an element $h$ of a Pommaret basis of $I$ and moving along the boundary we always end up at a maximal or closed corner. If $cls(h)=1$ then we arrive at a maximal corner, because we have $n-1$ non-multiplicative directions, and if we cannot go further in the boundary we are at an $lcm$ with some generator of lower $x_1$ degree, hence we are at a maximal corner. If $cls(h)>1$ then we are in a situation with less variables and arrive to a closed corner following the argument of the case $cls(h)=1$.

For the decision of being quasi-stable, assume we are in a multidegree $\mu$ such that $x^\mu$ is an element of the Pommaret basis of $\hat{I}$ given by the precedent procedure such that the exponent of $x_i$ is less or equal $\lambda_i$ for all $i$, when multiplying $x^\mu$ by some nonmultiplicative variables, we must reach one of the  $x^\nu$ obtained in the procedure, but if any of these has nonmultiplicative free directions, we can always multiply by some monomial in these free directions without getting into an involutive cone, and then the ideal is not quasi-stable. If no $x^\nu$ has nonmultiplicative free directions, then we have obtained a Pommaret basis for $I$ and therefore it is quasi-stable.$\square$

\begin{Example}
Let us consider the following ideal: $I=\langle xy^2z,y^5,z^3\rangle$. The artinian closure of $I$ is $\hat{I}=\langle x^2,xy^2z,y^5,z^3\rangle$. The set of closed and maximal corners of $\hat{I}$ is $\{x^2y^5, y^5z^3, x^2z^3, xy^5z^3,x^2y^5z, x^2y^2z^3\rangle$. The elements $x^\nu$ obtained in the above procedure are $\{x^2y^4,y^5z^2,x^2z^2,xy^4z^2,x^2yz^2\}$.  A Pommaret basis of $\hat{I}$ is given by
$$
\begin{array}{rl}
\Hc_{\hat{I}}=&\{x^2,x^2z,x^2z^2,x^2y,x^2yz,x^2yz^2,x^2y^2,x^2y^3,x^2y^4,xy^2z,\\
&y^5,z^3,xy^2z^2,xy^3z,xy^3z^2,xy^4z,xy^4z^2,y^5z,y^5z^2\}
\end{array}
$$
None of the elements $\nu$ in the boundary of $I$ have nonmultiplicative directions, and then a Pommaret basis of $I$ is 
$$\Hc_I=\{xy^2z,y^5,z^3,xy^2z^2,xy^3z,xy^3z^2,xy^4z,xy^4z^2,y^5z,y^5z^2\}$$

Consider now $J=\langle x^2,xy^2z,z^3\rangle$. Then the artinian closure of $J$ is the same as the artinian closure of $I$. Now, the element $xy^4z^2$ is one of the $x^\nu$ that is in the boundary of $J$, and it has $y$ as a free directions, which is nonmultiplicative, therefore $J$ is not quasi-stable.
\end{Example}
 
\section{Formal theory of differential systems}\label{formal_theory}
\spanishsection{Teor\'ia formal de sistemas diferenciales}

\subsection{Geometric approach to PDEs and algebraic analysis of them}
\spanishsubsection{Acercamiento geom\'etrico a las EDP's y su an\'alisis algebraico}

In this section we follow very closely \cite{S02d,S07,S07b}, in which the author demonstrates how differential equations can be analyzed using algebraic and homological techniques. For this analysis it is necessary to provide  first a differential geometric framework, which comes from the \emph{formal geometry of differential equations} \cite{KLV86,S02d,KL07,P78}.
\subsubsection{Geometric framework.}
\spanishsubsubsection{Marco geom\'etrico}

The first key concept is that of \emph{Jet bundles}. Let $\pi:\Ec\rightarrow \Xc$ be a fibred manifold. For our purposes we consider $\Xc=\RR^n$, $\Ec=\RR^{n+m}$ and $\pi$ the projection on the first $n$ components. We provide this manifold with local coordinates $\xb=(x^1,\dots,x^n)$ on $\Xc$ and fibre coordinates $\ub=(u^1,\dots,u^m)$ on $\Ec$. A \emph{section} is a map $\sigma:\Xc\rightarrow\Ec$ such that $\pi\circ\sigma=id_\Xc$; in local coordinates $\sigma$ corresponds to a smooth function $\ub=s(\xb)$.
\begin{Definition}
Consider the following equivalence relation between sections: $\sigma_1\sim\sigma_2$ if their graphs have at the point $\sigma_i(\xb_0)$ a contact of order $q$, i.e. if their Taylor expansion at $\xb_0$ coincide up to order $q$.

A \emph{$q$-jet} is an equivalence class $[\sigma]^{(q)}_{\xb_0}$ under this relation. The \emph{$q$-th order jet bundle} $J_q\pi$ is the set of all $q$-jets. $J_q\pi$ is a $(n+m{{n+q}\choose{q}})$-dimensional manifold. Projection on the point $\xb_0\in\Xc$ defines a fibration $\pi^q:J_q\pi\rightarrow \Xc$ with fibre coordinates for $[\sigma]^{(q)}_{\xb_0}$ given by $\ub^{(q)}=(u^\alpha_\mu)$ with $1\leq\alpha\leq n$, $0\leq\vert\mu\vert\leq q$ where $u^\alpha_\mu$ is the value of $\frac{\partial^{\vert\mu\vert}s^\alpha}{\partial x^\mu}$ at $\xb_0$.
\end{Definition}

\begin{Definition}
A \emph{differential equation} is a fibred submanifold $\Rc_q\subseteq J_q\pi$ which will be locally described as the zero set of a function $\Phi:J_q\pi\rightarrow \RR^t$, i.e. $\Phi(\xb,\ub^{(q)})=0$.
\end{Definition}

\begin{Definition}
If $q>r$ let us consider the canonical fibrations $\pi^q_r:J_q\pi\rightarrow J_r\pi$. Then we have the following operations on a differential equation $\Rc_q$:

\emph{Projection}: The $r$-fold projection of $\Rc_q\subseteq J_q\pi$ is $\Rc^{(r)}_{q-r}=\pi^q_{q-r}(\Rc_q)\subseteq J_{q-r}\pi$.

\emph{Prolongation}: Given an equation $\Rc_q\subseteq J_q\pi$, we consider the restriction $\hat{\pi}^q:\Rc_q\rightarrow \Xc$ of $\pi^q$. This provides $\Rc_q$ with a fibred manifold structure, and we can construct jet bundles over it. Considering $J_r\hat{\pi}^q$ and $J_{q+r}\pi$ as submanifolds of $J_r\pi^q$ then the $r$-fold prolongation of $\Rc_q$ is $\Rc_{q+r}=J_r\hat{\pi}^q\cap J_{q+r}\pi\subseteq J_{q+r}\pi$.
\end{Definition}

In the local coordinate picture, prolongation requires only \emph{formal derivative}. If we have the function $\Phi:J_q\pi\rightarrow \RR$, then its formal derivative with respect to the variable $x^i$ is a smooth function $D_i\Phi:J_{q+1}\rightarrow \RR$ given by
$$D_i\Phi=\frac{\partial\Phi}{\partial x^i}+\sum_{\alpha=1}^{m}\sum_{0\leq\vert\mu\vert\leq q}\frac{\partial \Phi}{\partial u^\alpha_\mu}u^\alpha_{\mu+1_i}.$$
If $\Rc_q$ is locally described as the zero set of the functions $\Phi^\tau:J_q\pi\rightarrow \RR$, then its first prolongation $\Rc_{q+1}$ is the zero set of the functions $\phi^\tau$ and their formal derivatives $D_i\Phi^\tau$, for $1\leq i\leq n$, $1\leq\tau \leq t$. Projection is in this sense more difficult to describe, since it requires elimination of variables.
One would expect that prolongation and projection are inverse operations in the sense that if one first prolongs $\Rc_q$ to $\Rc_{q+r}$ and then projects back, the obtained equation $\Rc^{(r)}_q$ coincides with $\Rc_q$. However, this is not true in general, as integrability conditions may arise. It is true in general that $\Rc^{(r)}_q\subseteq\Rc_q$. Integrability conditions are not additional restrictions on the solution space of $\Rc_q$, but conditions implicitly contained in $\Rc_q$ that are made visible by a sequence of prolongations and projections. They can be considered as obstructions for the order by order construction of formal power series solutions. Therefore, the question of the existence of such formal solutions leads to the definition of  \emph{formally integrable systems}:

\begin{Definition}
The differential equation $\Rc_q\subseteq J_q\pi$ is called \emph{formally integrable} if for all $r\geq 0$, the equality $\Rc^{(1)}_{q+r}=\Rc_{q+r}$ holds. 
\end{Definition}
 
\begin{Example}\label{example-differential-1}
Let us consider the following simple example from \cite{S02d} of a differential equation $\Rc_1$ in one dependent variable $u$ and three independent variables $x,y,z$ given by
\begin{equation}\label{example-integrability-conditions}
 \left\{
\begin{array}{rl}
 u_{3}+yu_{1}=&0\\
u_{2}=&0
\end{array}
\right.
\end{equation}
We will illustrate the operations of prolongation and projection, as well as the concept of integrability condition. To obtain the prolongation $\Rc_2$ of our equation, we have to formally differentiate the two equations of the system \ref{example-integrability-conditions} with respect to all the variables and add the result to $\Rc_1$. We obtain
\begin{equation}\label{example-prolonged}
 \left\{
\begin{array}{rl}
u_{3}+yu_{1}=&0\\
u_{2}=&0\\
u_{13}+yu_{11}=&0\\
u_{23}+yu_{12}+u_{1}=&0\\
u_{33}+yu_{13}=&0\\
u_{12}=&0\\
u_{22}=&0\\
u_{23}=&0
\end{array}
\right.
\end{equation}
The projection operation is not so easy to perform effectively as it requires to find some way to perform elimination of the variables. However, since we deal with a linear system, the situation is more favorable. The first two equations of $\Rc_2$ are of first order, since they come from the original system. In the fourth equation we can eliminate the two second order derivatives using the sixth and eighth equations, yielding the integrability condition $u_{1}=0$. There are no further possible eliminations and we obtain the system
\begin{equation}\label{example-prolonged-restricted}
 \Rc^{(1)}_1=
\left\{
\begin{array}{rl}
 u_{1}=&0\\
u_{2}=&0\\
u_{3}=&0\\
\end{array}
\right.
\end{equation}
To see how integrability conditions really represent conditions implicitly contained in $\Rc_1$ and not additional restrictions on its solution space, consider a formal power series solution of $\Rc_1$ of the form
\begin{equation}\label{example-power_series}
u(\xb)=\sum^\infty_{\vert\mu\vert=0}\frac{a_\mu}{\mu!}(\xb-\xb_0)^\mu
\end{equation}
From $\Rc_1$ we obtain that the coefficients $a_{1}$ and $a_{3}$ must satisfy $a_{3}+y_0 a_{1}=0$, which can be done in many ways. But if we look at $\Rc^{(1)}_1$ we see that $a_{1}$ must be zero, and therefore $a_1=a_3=0$ for any formal power series solution of $\Rc_1$.
\end{Example}

The verification of the formal integrability of a differential system using the definition requires the verification of an infinite number of conditions and thus it cannot be done effectively. A finite criterion for formal integrability can be formulated via algebraic and homological methods \cite{G65,G69,S02d,S07b}. The key for the application of these methods lies in a natural polynomial structure hidden in the jet bundle hierarchy, called \emph{fundamental identification}. A main result towards this identification is the following:
\begin{Proposition}
The jet bundle $J_q\pi$ of order $q$ is affine over the jet bundle $J_{q-1}\pi$ of order $q-1$.
\end{Proposition}
Let $[\sigma]^{(q)}_\xb$ and $[\sigma']^{(q)}_\xb$ be two points in $J_q\pi$ such that $[\sigma]^{(q-1)}_\xb$=$[\sigma']^{(q-1)}_\xb$, i.e. the two points belong to the same fibre with respect to the fibration $\pi^q_{q-1}$. Thus, $[\sigma]^{(q)}_\xb$ and $[\sigma']^{(q)}_\xb$ correspond to two Taylor series truncated at degree $q$ which coincide up to degree $q-1$. Then, their difference consists of one homogeneous polynomial of degree $q$ for each variable $u^\alpha$. This implies that the difference between two points in $J_q\pi$ can be interpreted as a vector. To identify the vector space, consider $\rho=[\sigma]^{(q)}_\xb\in J_q\pi$ and $\bar{\rho}=[\sigma]^{(q-1)}_\xb=\pi^q_{q-1}(\rho)$ its projection to $J_{q-1}\pi$. Also, set $\xi=\sigma(\xb)=\pi_0^q(p)\in\Ec$. Then, the fibre $(\pi^q_{q-1})^{-1}(\bar{\rho})$ is an affine space modelled on the vector space $S_q(T^*_\xb\Xc)\otimes V_\xi\pi$, where $S_q$ denotes the $q$-fold symmetric product, and $V_\xi\pi\subset T_\xi\pi$ is the vertical bundle defined as the kernel of the tangent map $T_\xi\pi$.

The fact that the jet bundle $J_q\pi$ is an affine bundle over $J_{q-1}\pi$ implies that the tangent space to the affine space $(\pi^q_{q-1})^{-1}(\bar{\rho})$ at the point $\rho\in J_q\pi$ is canonically isomorphic to the vector space $S_q(T^*_\xb\Xc)\otimes V_\xi\pi$. This isomorphism is called the \emph{fundamental identification}, and allows the algebraic analysis of differential equations below. An expression for the fundamental identification is the map $\epsilon_q:V_\rho\pi^q_{q-1}\rightarrow S_q(T^*_\xb\Xc)\otimes V_\xi\pi$ given by
$$\epsilon_q(\partial_{u^\alpha_\mu})=\frac{1}{\mu!}dx^\mu\otimes\partial_{u^\alpha}$$
here, $V_\rho\pi^q_{q-1}$ is the vertical space defined as the kernel of the tangent map $T_\rho\pi^q_{q-1}$. It is spanned by all the vectors $\partial_{u^\alpha_{\mu}}$ with $\vert\mu\vert=q$.

\subsubsection{Algebraic Analysis}
\spanishsubsubsection{An\'alisis algebraico}

To each differential equation we can associate a family of vector spaces called the \emph{geometric symbol} of the equation. This classical construction can be transformed to obtain the so called principal symbol, which is taken as the starting point for the algebraic analysis of differential systems using the geometric framework described above.

\begin{Definition}
Let $\Rc_q\subseteq J_q\pi$ be a differential equation. The (geometric) \emph{symbol} $\Nc_q$ of $\Rc_q$ is a family of vector spaces over $\Rc_q$ where the value at $\rho\in\Rc_q$ is given by
$$(\Nc_q)_\rho=T_\rho\Rc_q\cap V_\rho\pi^q_{q-1}=V_p(\pi^q_{q-1}\vert_{\Rc_q})$$
\end{Definition}

Thus, the symbol is the vertical part of the tangent space of the submanifold $\Rc_q$ with respect to the filtration $\pi^q_{q-1}$. Locally, $\Nc_q$ can be described as the solution space of the system of linear equations:
$$(\Nc_q)_\rho=\left\{ \sum_{\substack{1\leq \alpha\leq m\\ \vert\mu\vert=q}} \frac{\partial\Phi^\tau}{\partial u^\alpha_\mu}\dot{u}^\alpha_\mu=0\qquad \tau=1\dots t\right .$$ 
Since the derivatives $\frac{\partial\Phi^\tau}{\partial u^\alpha_\mu}$ are evaluated at $\rho$, it is a system with real coefficients, the matrix of which is called the \emph{symbol matrix}, $\Mc_q(\rho)$.  In the non linear case, we perform a linearisation at the point $\rho$ to obtain the symbol at this point.

Every prolongation $\Rc_{q+r}\subseteq J_{q+r}\pi$ of $\Rc_q$ possesses a symbol $\Nc_{q+r}\subseteq T(J_{q+r}\pi)$. We can derive a local representation of $\Nc_{q+r}$ from a local representation of $\Nc_q$. A local representation of $\Nc_{q+1}$ is
$$(\Nc_{q+1})_\rho=\left\{ \sum_{\substack{1\leq \alpha\leq m\\ \vert\mu\vert=q}} \frac{\partial\Phi^\tau}{\partial u^\alpha_\mu}\dot{u}^\alpha_{\mu+1_i}=0\qquad
\begin{array}{c}
\tau=1\dots t\\
i=1\dots n
\end{array}
\right.$$
\begin{Proposition}\label{symbol-dimension}
If $\Nc_{q+1}$ is a vector bundle, then $dim\Rc^{(1)}_q=dim\Rc_{q+1}-dim \Nc_{q+1}$
\end{Proposition}

\begin{Remark}
The condition of $\Nc_{q+1}$ being a vector bundle and not just a family of vector spaces is not satisfied in general even by symbols of linear differential equations. Of course, $dim(\Nc_q)_{\rho}$ might vary with $\rho$. Only if the dimension remains constant over $\Rc_q$ the symbol is a vector bundle. We will always assume that the considered symbols are vector bundles.
\end{Remark}
The existence of integrability conditions is signaled by a dimension inequality, namely $dim\Rc^{(1)}_q<dim\Rc_q$. Proposition \ref{symbol-dimension} relates the dimension of $\Rc^{(1)}_q$ to the dimension of $\Nc_{q+1}$; hence, analysing the prolonged symbol matrix gives information about potential integrability conditions.

Another notion of symbol can be defined as follows: Given a one-form $\chi\in T^*\Xc$, it induces for every $q>0$ a map $\iota_{\chi,q}:V\pi\rightarrow V\pi^q_{q-1}$ defined by $\iota_{\chi,q}(v)=\epsilon^{-1}_q(\chi^q\otimes v)$. In local coordinates, $\chi=\chi_idx^i$ and we have $\iota_{\chi,q}:v^\alpha\partial_{u^\alpha}\mapsto\chi_{\mu}v^\alpha\partial_{u^\alpha_\mu}$ where $\mu$ runs over all multi indices of length $q$, and $\chi_\mu=\chi^{\mu_1}_1\cdots\chi^{\mu_n}_n$. Let $\sigma$ be the symbol map of the differential equation $\Rc_q$ globally described by the map $\Phi:J_q\pi\rightarrow \Ec'$. Then, the \emph{principal symbol of $\Rc_q$} is the linear map $\tau_\chi:V\pi\rightarrow T\Ec'$ given by $\tau_\chi=\sigma\circ \iota_{\chi,q}$. Locally, we can associate a matrix $T[\chi]$ with $\tau_\chi$:
$$T^\tau_\alpha[\chi]=\sum_{\vert\mu\vert=q}\frac{\partial \Phi^\tau}{\partial u^\alpha_\mu}\chi^\mu$$
We may think of $T[\chi]$ as a contraction of the (geometric) symbol matrix $M_q$. Both matrices have the same number of rows ($p=dim\Ec'$). $T[\chi]$ has $m=dim\Ec$ columns, and the column with index $\alpha$ is a linear combination of all columns in $M_q$ corresponding to a variable $\dot{u}^\alpha_\mu$, with the coefficients given by $\chi_\mu$.

\begin{Remark}\label{principal-symbol-syzygies}
Using $T[\chi]$ we can relate the construction of integrability conditions with syzygy computations. If the functions $\Phi^\tau$ locally representing $\Rc_q$ lie in a differential field $\FF$ then the entries of $T[\chi]$ are polynomials in $\Pc=\FF[\chi_1,\dots,\chi_n]$. The rows of $T[\chi]$ may be considered as elements of $\Pc^m$ and generate a module $\Mc\subseteq\Pc^m$. Lets $\Sb\in Syz(\Mc)$ be a syzygy of the rows of $T[\chi]$. The substitution $\chi_i\rightarrow D_i$ transforms each component $S_\tau$ of $\Sb$ into a differential operator $\hat{S}_\tau$. By construction, $\Psi=\sum_{r-1}^t \hat{S}_\tau \Phi^\tau$ is a linear combination of differential consequences of $\Rc_q$ in which the highest order terms cancel. This represents the formulation of taking cross derivatives.
\end{Remark}

\begin{Example}\label{example-differential-2}
Let us consider an example due to Janet, the partial differential equation $\Rc_2$ in one dependent variable $u$ and three independent variables $x,y,z$ given by
\begin{equation}\label{Janet}
 \left\{
\begin{array}{rl}
 u_{33}+yu_{11}=&0\\
u_{22}=&0
\end{array}
\right.
\end{equation}
Then, we have $n=3,\, m=1,\,t=2,\,q=2$. At a given point $\rho\in\Rc_2$ the (geometric) symbol matrix has $t=2$ rows and $m{{n+q-1}\choose{n-1}}=1{{4}\choose{2}}=6$ columns. Each of the rows corresponds to one of the equations, and each column is labelled by one of the second order derivatives of the (in this case unique) dependent variables:
$$
M_2=
\left(
\begin{array}{cccccc}
1&0&0&0&0&y\\
0&0&1&0&0&0\\
\end{array}
\right)
$$
The columns of this matrices are labelled respectively by $u_{33},u_{23},u_{22},u_{13},u_{12}$ and $u_{11}$.
The principal symbol matrix has then $2$ rows and just one column, and is given by
$$T[\chi]=
\left(
\begin{array}{c}
\chi^2_3+y\chi^2_1\\
\chi^2_2\\
\end{array}
\right)
$$
Considering the entries of $T[\chi]$ as polynomials in $\Pc=\FF[\chi_1,\chi_2,\chi_3]$ as in remark \ref{principal-symbol-syzygies} above, the rows of $T[\chi]$ can be considered as elements of $\Pc$ and generate a module (in this case an ideal, because we have just one dependent variable) $\Mc=\langle \chi^2_3+y\chi^2_1,\chi^2_2\rangle$. The syzygy module of this ideal is generated by
$$S_1=\chi_2^2\eb_1-(\chi^2_3+y\chi^2_1)\eb_2$$
and using the correspondence $\chi_i\mapsto D_i$ we obtain the differential operator
$$\hat{S}_1=D_2^2\eb_1-(D^2_3+yD^2_1)\eb_2$$
We apply this operator to equation \emph{(\ref{Janet})}, and obtain the integrability condition $u_{112}=0$.
\end{Example}

\subsection{The role of Koszul homology}
\spanishsubsection{El papel de la homolog\'ia de Koszul}

\subsubsection{Involution and formal integrability}
\spanishsubsubsection{Involuci\'on e integrabilidad formal}

The concept of formal integrability is sometimes not sufficient for the study of differential systems. The notion of involution achieves a better understanding of these systems: ``whatever one intends to do with a differential equation, one should first render it involutive: \emph{involution is the central principle in the theory of under- and overdetermined systems}'' \cite[p.VIII]{S02d}. The theory of involution will make any subsequent analysis of under- and overdetermined systems significantly easier. Involution is often confused with formal integrability, but it is a stronger notion. Involution can be used to give rigorous definitions of under- and overdetermined systems or even prove results on the existence and uniqueness of solutions. Involutive differential equations have many applications, in particular in mathematical physics and control theory, numerical
analysis, etc; see \cite{S02d,P88,P01a,P01b}.

The homological treatment of involution uses as a main tool Spencer homology (see section \ref{spencer}). In this context, the finiteness of Spencer cohomology has important consequences when studying the degree of involution of symbolic systems, see for example \cite{LS02,S07b}. The algebraic definitions of symbolic systems and involution have been given in section \ref{spencer}, now we transfer them to the formal theory of differential equations, showing first that every differential equation defines a symbolic system, which will allow us to define the degree of involution of a differential equation.

\begin{Proposition}
Let $\Rc_q\subseteq J_q\pi$ be a differential equation, and $(\rho_r\in\Rc_r)_{r\geq q}$ be a sequence of points such that $\pi^r_q(\rho_r)=\rho_q$, and set $\xi=\pi_0^q(\rho_q)$ and $\xb=\pi^q(\rho_q)$. If we set $\Nc_r=\mathfrak{S}_r(T^*_\xb \chi)\otimes V_\xi\pi$ for $0\leq r<q$, then the sequence $((\Nc_r)_{\rho_r})_{r\in \NN}$ defines a symbolic system in $\mathfrak{S}(T^*_\xb\chi)\otimes V_\xi\pi$ which satisfies $\Nc_{r+1}=\Nc_{r,1}$ for all $r\geq q$. 
\end{Proposition}

Observe that even if in this proposition a sequence of points $\rho_r\in\Rc_r$ were used to consider the symbols $(\Nc_r)_{\rho}$, the obtained symbolic system is independent of the choice of these points (the sequence of which might even not exist). Hence, at each point $\rho\in\Rc_q$ the symbol $(\Nc_q)_{\rho}$ induces a symbolic system which can be considered as a submodule $\Nc[\rho]\subseteq \mathfrak{S}(T^*_\xb\chi)\otimes V_\xi\pi$ (lemma \ref{symbolic-system-comodule}). We then speak of the \emph{symbol comodule} of $\Rc_q$ at the point $\rho$.

\begin{Definition}
The symbol $\Nc_q$ of the differential equation $\Rc_q\subseteq J_q\pi$ of order $q$ is \emph{involutive} at the point $\rho\in\Rc_q$ if the symbol comodule $\Nc[\rho]$ is involutive at degree $q$.
\end{Definition} 

The annihilator $\Nc^0\subseteq \Sc(T_\xb\chi)\otimes V_\xi\pi$ of the comodule $\Nc$ is an $\Sc(T_\xb\chi)$-submodule which will be called the \emph{symbol module} of $\Rc_q$ at $\rho$. In the chosen coordinates and bases, $\Nc^0$ is generated by
$$\sum_{\substack{1\leq \alpha\leq m\\ \vert\mu\vert=q}} \frac{\partial\Phi^\tau}{\partial u^\alpha_\mu}\partial^\mu_{\xb}\otimes \partial_{\mu^\alpha}$$

Identifying $\Sc(T_\xb\chi)$ with the polynomial ring $\Pc=\RR[\partial_\xb^1,\dots,\partial_\xb^n]$ one sees that $\Nc^0$ is the polynomial module generated by the rows of the matrix $T[\chi]$ of the principal symbol.
\begin{Remark}
There are several approaches to the detection of the degree of involution of a symbolic system, and explicit criteria for a comodule to be involutive have been given. One of them is the Cartan test \cite{M53}, which can only be applied in $\delta$-regular bases (see definition in \cite{S07b}). Generic bases are $\delta$-regular, hence we can avoid the problem of $\delta$-regularity performing a random coordinate transformation to our basis; however this is computationally very expensive in general. A dual Cartan test, of a more homological nature is based on a letter of J-P. Serre appended to \cite{GS64}. To face the problem of $\delta$-regularity, an efficient solution has been given based on Pommaret bases, a special kind of involutive bases thoroughly studied in \cite{S02a,S02b} which have a direct application to the study of differential systems \cite{S02d,S07,S07b}.

From the different ways to compute $Tor(\kb,\Mc)$, we have seen in section \ref{tor} that the dimensions of the Koszul homology modules of an $R$-module $\Mc$ are the Betti numbers of it. Hence, if the Castelnuovo-Mumford regularity of an ideal $I\subseteq R$, $reg(I)$ is $q$, then all the homology modules $H_{r,p}(\Mc)=H_{r,p}(R/I)$ with $r\geq q$ vanish. Hence, the Castelnuovo-Mumford regularity of $I$ is the same as the degree of involution of $\Mc$, see \cite{S02d,M03} which report this fact. 
\end{Remark}

\begin{Definition}
The differential equation $\Rc_q$ is called \emph{involutive} if it is formally integrable and if its symbol is involutive.
\end{Definition}
\begin{Theorem}\label{involutive-equation-symbol}
$\Rc_q$ is an involutive differential equation if and only if its symbol is involutive and $\Rc_q^{(1)}=\Rc_q$.
\end{Theorem}

The fact that every (regular) differential equation can be completed to involution is known as the \emph{Cartan-Kuranishi} theorem, which states the following:
\begin{Theorem}[Cartan-Kuranishi]
Let $\Rc_q\subseteq J_q\pi$ be a regular differential equation. Then two integers $r,s\geq 0$ exist such that $\Rc_{q+r}^{(s)}$ is involutive.
\end{Theorem}
A constructive proof of this theorem consists on an algorithm that performs such completion to involution. The algorithm consists of two nested loops. The inner loop prolongs the equation until an involutive symbol is reached. The outer loop checks whether one further prolongation and subsequent projection yields integrability conditions. It can be shown that this algorithm terminates, producing an involutive equation of the form $\Rc_{q+r}^{(s)}$.

Theorem \ref{involutive-equation-symbol} is a consequence of the following theorem which gives a criterion for formal integrability based on Spencer and Koszul (co)homologies. We reproduce the proof given in \cite{S07b} for it provides us with a good example of the interactions of techniques from homological and commutative algebra in the theory of differential equations. In particular it demonstrates that a generating set of the Koszul homology module $H_1(\Nc^0)$ is relevant for the analysis of the involution of a differential system, since it shows us exactly which cross derivatives may produce integrability conditions:

\begin{Theorem}[\cite{S07b},Theorem 7.15]
The differential equation $\Rc_q$ is formally integrable if and only if an integer $r>0$ exists such that the symbolic system $\Nc$ defined by the symbol $\Nc_q$ and all its prolongations is $2$-acyclic at degree $q+r$ and the equality $\Rc^{(1)}_{q+r'}=\Rc_{q+r'}$ holds for all values $0\leq r'\leq r$.
\end{Theorem}

\noindent{\bf Proof: } One direction is trivial. For a formally integrable equation $\Rc_q$ we even have $\Rc^{(1)}_{q+r'}=\Rc_{q+r'}$ for all $r'\geq 0$, and we know that every symbolic system $\Nc$ must become acyclic at some degree $q+r$.

For the converse, we first note that, because of lemma \ref{homology-involutive-system} $\Nc$ is trivially $1$-acyclic at degree $q$. Our assumption says that in addition the Spencer cohomology modules $H^{q+s,2}(\Nc)$ vanish for all $s\geq r$, this implies dually that the Koszul homology modules $H_{q+s,1}(\Nc^0)$ vanish for all $s>r$.

The Koszul homology corresponds to a minimal free resolution of $\Nc^0$ and hence our assumption tells us that the maximal degree of a minimal generator in the first syzygy module $Syz(\Nc^0)$ is $q+r$. We have seen that the syzygies of $\Nc^0$ are related to those integrability conditions arising from generalised cross-derivatives between the highest-order equations. If we know the equality $\Rc^{(1)}_{q+r'}=\Rc_{q+r'}$ holds for all $0\leq r'\leq r$, then none of these cross-derivatives can produce an integrability condition. Furthermore, no integrability conditions can arise from lower-order equations. Hence $\Rc_q$ is formally integrable.$\square$ 

\begin{Example}\label{example-differential-3}
Let us continue example \ref{example-differential-2}, let $\Rc_2$ be the Janet equation \emph{(\ref{Janet})}. As we have already seen, we have that $\Rc^{(1)}_3$ is given by  the integrability condition $u_{112}=0$ obtained by the analysis of the syzygy module of the ideal corresponding to the principal symbol of $\Rc_2$, the original equations in \ref{Janet}, and their formal derivatives.
The principal symbol of it is given by
$$T[\chi]=
\left(
\begin{array}{c}
\chi^2_3+y\chi^2_1\\
\chi^2_2\\
\chi^2_1\chi_2
\end{array}
\right)
$$
then, the ideal generated by its rows is $I_2=\langle\chi^2_3+y\chi^2_1,\chi^2_2,\chi^2_1\chi_2\rangle$. A generating set of its syzygy module is given by $S_1=\chi^2_2\eb_1-(\chi^2_3+y\chi^2_1)\eb_2$ (which we already had), $S_2=\chi^2_1\eb_2-\chi_2\eb_3$ and $S_3=\chi^2_1\chi_2\eb_1-y\chi^2_2\eb_3$. If we apply the corresponding differential operators to $\Rc^{(1)}_3$ we only obtain one new integrability condition, which comes from $S_3$: $u_{1111}=0$.

Adding this new condition to the system, together with the original equations and their prolongations up to order $4$, we obtain $\Rc^{(2)}_4$. The ideal corresponding to the principal symbol is given by $I_3=\langle\chi^2_3+y\chi^2_1,\chi^2_2,\chi^2_1\chi_2,\chi_1^4\rangle$. If we compute its syzygy module and apply the corresponding differential operators, none of them leads to a new integrability condition, so we have that $\Rc^{(2)}_4$ is formally integrable.

However, $\Rc^{(2)}_4$ is not involutive. We can see it in two ways: One consists in noting that the resolution of $I_4$ is not linear, and since we know that we can always obtain a linear resolution from an involutive symbol \cite{S07b} then $\Rc^{(2)}_4$ is not involutive. The second way consists in choosing a point $\rho \in \Rc^{(2)}_4$, where $y=a$ for some constant. Taking the ideal $I=\langle\chi^2_3+a\chi^2_1,\chi^2_2,\chi^2_1\chi_2,\chi_1^4\rangle$ corresponding to $I_3$ at this point, one observes that the truncated ideal $I_{\geq 4}$ is the annihilator of the symbol $\Nc$ at $\rho$. The regularity of this ideal is $reg(I_{\geq 4})=5$ so the degree of involution of $\Nc$ is $5$, not $4$.
\end{Example}

\subsubsection{Initial value problems}
\spanishsubsubsection{Problemas de valor inicial}

For general differential systems the notion of a well posed initial value problem is not completely clear. We will consider that in an initial value problem one prescribes certain derivatives $u^\alpha_{\mu}(\xb)$ on planes of the form $x^{i_1}=\cdots=x^{i_k}=0$ with possibly different codimensions $k$ for different derivatives. Since we are considering formal power series solutions, the initial data will have the form of formal power series.
\begin{Definition}
An initial value problem for a differential equation $\Rc_q$ is \emph{formally well-posed} if it possesses a unique formal power series solution for arbitrary formal power series as initial data.
\end{Definition}
A formally well-posed initial value problem has then a solution, and it is unique. Thus, in a formally well-posed initial value problem the Taylor coefficients of the initial data are in a one-to-one correspondence with the parametric derivatives of the differential equation.

We want to determine in a systematic way formally well-posed initial value problems for (involutive) differential systems. First we analyze linear systems and restrict to the case of one dependent variable (in the case of more dependent variables one performs the same analysis for each of them separately). These systems correspond to an ideal $I\subseteq \Dc=\FF[\partial_1,\dots,\partial_n]$.

When solving a differential system using power series solutions, one needs to use some of the derivatives as principal, and the others as parametric. In order to construct formally well-posed initial value problems, we need an algebraic approach to the distinction into principal and parametric derivatives:

\begin{Definition}
Let $I\subseteq \Dc$ be a (left) ideal. For a given term order $<$, we call the derivatives $\partial^\mu$ with $\mu\in\Delta_<(I)=le_<I$ \emph{principal}. All remaining derivatives are \emph{parametric} and their exponents form the complementary set $\Pi_<(I)=\NN^n\smallsetminus\Delta_<(I)$.
\end{Definition}
The construction of formally well-posed initial value problems relies on the construction of a disjoint decomposition of the complementary set $\Pi_<(I)$ which is the complement of a monoid ideal, namely $\Delta_<(I)$.

To find such a decomposition, choose a point $\xb_0\in\chi$ and assume that it corresponds to the origin. Consider a subset $N\subseteq\{1,\dots,n\}$ and consider on one hand $\yb=N(\xb)$ where $y^i=x^i$ if $i\in N$ and $y^i=0$ if $i\notin N$, and in the other hand consider $\xb_N=(x^{i_1},\dots,x^{i_k})$ if $N=\{i_1,\dots,i_k\}$. We know that the set $\Pi_<(I)$ has a disjoint decomposition of the form given in proposition \ref{stanley-monoid}. This decomposition is defined by a set $\bar{N}\subseteq \Pi_<(I)$ and for each $\nu\in\bar{N}$, a set $N_\nu\subseteq\{1,\dots,n\}$. Based on these data, we can prescribe the following initial conditions for our system in a neighbourhood of $\xb_0$:
\begin{equation}\label{initial-conditions}
 D_\nu(N_\nu(\xb))=f_\nu(\xb_{N_\nu})\qquad \forall \nu\in\bar(N)
\end{equation}
where the $f_\nu$ are arbitrary formal power series in their arguments and if $N_\nu=\emptyset$, $f_\nu$ is an arbitrary constant. The result we are looking for is then the following:
\begin{Theorem}[\cite{S07}]
The initial value problem for $I$ defined by the initial conditions (\ref{initial-conditions}) above is formally well posed.
\end{Theorem}
Therefore, the construction of a formally well-posed initial value problem for an involutive system $\Fc\subset\Dc$ is equivalent to the determination of a Stanley decomposition of $\Dc/\Fc$, which is equivalent to find a decomposition of the monoid ideal $le_<\Fc$ of the form given by proposition \ref{stanley-monoid}.
\begin{Remark}
If the ideal $I\subseteq\Dc$ is monomial then any complementary decomposition of $\Dc/I$ leads to a closed-form representation of the general solution of the system \cite{S07}: Let $I$ be a monomial ideal and consider a disjoint decomposition of $\Pi(I)$ defined by $\bar N\subseteq \Pi(I)$ and the associated sets $N_\nu\subseteq\{1,\dots\}$. then, every smooth solution $u(\xb)$ of the differential system corresponding to $I$ may uniquely be written in the form
$$u(\xb)=\sum_{\nu\in\bar{N}}x^\nu f_\nu(\xb_{N_\nu})$$
with smooth functions $f_\nu$. Conversely, every expression of this form is a smooth solution of the differential system. 
\end{Remark}

As we have seen in section \ref{stanley}, the knowledge of the Koszul homology of the monomial ideal $I$ corresponding to $le_<\Fc$ gives us a procedure to obtain a Stanley decomposition of the latter. While in the previous paragraph we have seen that the relevant piece of the Koszul homology that is needed for the detection of involution is $H_1(\KK(I))$, it turns out that the relevant part of the Koszul homology of $I$ that provides the necessary information for the Stanley decomposition is $H_{n-1}(\KK(I))$. The procedure for the construction of such a decomposition is given in section \ref{stanley}, and the relevant results are proposition \ref{Stanley-artinian} and \ref{Stanley-nonartinian}. To illustrate this application of the Koszul homology to differential systems we end this section with two examples: First we just continue with example \ref{example-differential-3} to complete the analysis of the Janet equation, and then we treat the $U(1)$ Yang-Mills equation in two dimensions:
\begin{Example}\label{example-differential-4}
We have seen in example \ref{example-differential-3} that the Janet equation $\Rc_2$ defined by $u_{33}+yu_{11}=u_{22}=0$ can be completed to an involutive equation. This needs several prolongations and we obtained that $\Rc^{(2)}_5$ is involutive, and already $\Rc^{(2)}_4$ is formally integrable. The annihilator ideal of the symbol comodule $\Nc$ at a point $\rho$ where $y=a$ is given by
$$I=\langle z^2+ax^2,y^2,x^2y,x^4\rangle$$
The leading ideal $le_<I$ is artinian, and the multidegrees if the generators of $H_2(\KK(le_<I))$ are $(4,1,2)$ and $(2,2,2)$. Therefore, following proposition \ref{Stanley-artinian}, a Stanley decomposition of $\RR[x,y,z]/le_<I$ is given by
\begin{eqnarray*}
\RR[x,y,z]/le_<I&=&1\oplus x\oplus x^2\oplus x^3\oplus y\oplus xy\oplus z\\
& &\oplus xz\oplus y^2z\oplus x^3z\oplus yz\oplus xyz
\end{eqnarray*}
Thus, for a well-posed initial value problem, we need to prescribe
$$\begin{array}{cccc}
u(\xb_0)=f_{[0,0,0]}& u_{1}(\xb_0)=f_{[1,0,0]}& u_{11}(\xb_0)=f_{[2,0,0]}&u_{111}(\xb_0)=f_{[3,0,0]}\\ u_{2}(\xb_0)=f_{[0,1,0]}&u_{12}(\xb_0)=f_{[1,1,0]}&u_{3}(\xb_0)=f_{[0,0,1]}&u_{13}(\xb_0)=f_{[1,0,1]}\\
u_{113}(\xb_0)=f_{[2,0,1]}&u_{1113}(\xb_0)=f_{[3,0,1]}&u_{23}(\xb_0)=f_{[0,1,1]}&u_{123}(\xb_0)=f_{[1,1,1]}
\end{array}
$$
where all the $f_\nu$ are constants.
\end{Example}

\begin{Example}
The $U(1)$ Yang-Mills equation has two independent variables and two dependent variables in two dimensions, and is given by
\begin{equation}\label{Yang-Mills}
 \left\{
\begin{array}{rl}
 u_{22}-v_{12}=&0\\
u_{12}-v_{11}=&0
\end{array}
\right.
\end{equation}
We must analyze each dependent variable separately. For the variable $v$ we have that the ideal we have to consider is just the zero ideal, since none of the equations has a leading term on the variable $v$ and therefore, $v(x,t)$ can be chosen arbitrarily.

For the second dependent variable $u$ we have the ideal $I=\langle xt,t^2\rangle$. Note that this ideal is not artinian. Following proposition \ref{Stanley-nonartinian}, we easily find the Stanley decomposition
$$R/\langle xt,t^2\rangle=t\oplus1\cdot\RR[x]$$
therefore, the conditions for a well-posed initial value problem are
\begin{eqnarray*}
v(x,t)&=&f(x,t)\\
u(x,0)&=&g(x)\\
u_t(0,0)&=&h
\end{eqnarray*}
where $f,g$ are arbitrary functions, and $h$ an arbitrary constant.
Note in particular that this is an example in which the fact of being underdetermined is not captured by the usual ``counting rules'' since it has as many equations as dependent variables and still it is undetermined as we have just seen.
\end{Example}

\section{Reliability of coherent systems}\label{reliability}
\spanishsection{Fiabilidad de sistemas coherentes}

Reliability theory studies systems of components subject to failure. A main goal is to evaluate the performance of the system computing its probability of failure. This probability is computed in terms of the events that cause failure and the probability of failure of each component in the system. Typical examples of such systems are networks, although there are a number of important non-network systems, such as \emph{k-out-of-n} systems. In general, each component may assume a finite number of states that correspond to different efficiency levels. If the components of a system can only admit two states, i.e. failure and non failure, we say the system is \emph{binary}; otherwise we say we have a \emph{multistate} system. Under some assumptions, the continuous case, in which the components may have continuous states, can be mapped into this discrete setting. 

In the case of coherent systems, the evaluation of its reliability benefits from tools coming from commutative algebra. B. Giglio and H.P. Wynn established in \cite{GW03,GW04} the correspondence between coherent systems and monomial ideals. Based on their work, we provide with our techniques new bounds for the reliability of coherent systems, and show in some relevant examples that Koszul homology and Mayer-Vietoris trees constitute excellent tools for this goal.

In this section we begin providing the necessary definitions from reliability theory; then we show connection with the theory of monomial ideals and finally give some relevant examples in which the computation of Koszul homology and the application of Mayer-Vietoris trees provide good tools for reliability evaluation of the different systems. The main results of this part are presented in \cite{SW07}.

\subsection{Reliability theory of coherent systems}
\spanishsubsection{Teor\'ia de fiabilidad de sistemas coherentes}

We consider the concepts we need from reliability theory as defined in \cite{GW03} and \cite{GW04}.
\begin{Definition}
A \emph{system} is a set $\Sc$ of $n$ components, each of which has different increasing efficiency levels coded by the integers $\{0,1,\dots \}$.

An \emph{outcome} is a nonnegative integer vector of length $n$ that describes the state of the set of components. We call $\mD$ the set of all possible outcomes.

A \emph{failure outcome} is an outcome that leads to failure of $\Sc$. The \emph{failure set} $\mF$ is the set of all failure outcomes. The \emph{nonfailure set} $\bar{\mF}$ is the complement of $\mF$ in $\mD$. Another way to describe these sets is through a binary function $\Phi$ defined on $\Sc$ that takes the value $0$ on the outcomes that lead to failure on the system, and $1$ otherwise. Then
$$\mF=\{\alpha \in \Sc\vert \Phi(\alpha)=0\}$$
$$\bar{\mF}=\{\alpha \in \Sc\vert \Phi(\alpha)=1\}$$
The function $\Phi$ is called the \emph{structure function} of the system $\Sc$.
\end{Definition}
A particularly interesting class of systems are the so called \emph{coherent systems}, in which increasing the efficiency level of any component will not lead to a deterioration of the system, or equivalently, deceasing the efficiency level of any component does not improve the performance of the system:
\begin{Definition}
A coherent system is  a multistate system $\Sc$ with a monotone structure function $\Phi$:
$$\alpha\preceq\beta\Rightarrow \Phi(\alpha)\leq\Phi(\beta),\, \forall \alpha,\beta\in \Sc$$
where $\alpha\preceq\beta$ if $\alpha_i\leq\beta_i$ for all $1\leq n$.
\end{Definition}
Note that the monotonicity of $\Phi$ in coherent systems gives the failure and nonfailure sets the following respective properties:
$$\alpha\in\mF\Rightarrow\beta\in\mF\quad \forall \beta\preceq\alpha$$
$$\alpha\in\bar{\mF}\Rightarrow\beta\in\bar{\mF}\quad \forall \alpha\preceq\beta$$
These properties enable us to define minimal failure and nonfailure sets, the knowledge of which is sufficient to describe the operation of the system.
\begin{Definition}
The \emph{minimal failure set} of a coherent system $\Sc$
are the maximal points in $\mF$ with respect to $\preceq$:
$$\mF^{*}=\{\alpha\in \mF\vert \nexists\beta\in\mF\mbox{ s.th. }\alpha\preceq\beta\}$$
similarly, the \emph{minimal nonfailure set} collects the minimal points in $\bar{\mF}$ with respect to $\preceq$:
$$\bar{\mF}^{*}=\{\alpha\in \bar{\mF}\vert \nexists\beta\in\mF\mbox{ s.th. }\beta\preceq\alpha\}$$
\end{Definition}
\begin{Remark}
The minimal failure (resp. nonfailure) set gives a minimal representation of the failure (resp. nonfailure) set of the system. It is minimal in the sense that they are the minimal set of components whose failure causes the failure of the system, and symmetrically with the minimal nonfailure set.
\end{Remark}

The problem reliability theory of coherent systems faces is to calculate the \emph{reliability} $\Rc$ of a system $\Sc$, defined as the probability that the system is operating. Equivalently, the \emph{unreliability} $\Uc$ of $\Sc$ is the probability that $\Sc$ is not operating:
$$\Rc=Prob(\bar{\mF})=Prob(\Phi=1)$$
$$\Uc=Prob(\mF)=Prob(\Phi=0)=1-\Rc$$

These probabilities can be expressed in terms of \emph{orthants}, i.e. sets of the form $Q_\alpha=\{\beta\in\Sc\vert\alpha\preceq\beta\}$. In particular, since $\bar{\mF}=\bigcup_{\alpha\in\bar{\mF}^*}Q_\alpha$, we have that the reliability of $S$, can be expressed as
$$\Rc=Prob(\bar{\mF})=Prob(\bigcup_{\alpha\in\bar{\mF}^*}Q_\alpha)$$
The classical approach to the computation of such possibilities makes use of the inclusion-exclusion identity:
$$Prob(\bar{\mF})=\sum_{\alpha\in\bar{\mF}^*}Prop(Q_\alpha)-\sum_{\alpha,\alpha'\in\bar{\mF}^*}Prop(Q_\alpha\cap\alpha')+\dots+(-1)^{\vert\mF^*\vert+1}Prob(\bigcap_{\alpha\in\bar{\mF^*}}Q_\alpha)$$
Truncations of this identity give upper and lower bounds for the reliability of the systems. These are known as Bonferroni inequalities:
$$Prob(\bar{\mF})\leq\sum_{\substack{I\in\Pc(\bar{\mF}^*)\\ \vert I\vert\leq r}}(-1)^{\vert I\vert+1}Prob(\bigcap_{\alpha\in I}Q_\alpha)\quad \mbox{(r odd)}$$
$$Prob(\bar{\mF})\geq\sum_{\substack{I\in\Pc(\bar{\mF}^*)\\ \vert I\vert\leq r}}(-1)^{\vert I\vert+1}Prob(\bigcap_{\alpha\in I}Q_\alpha)\quad \mbox{(r even)}$$

The Bonferroni identity and inequalities are computationally very inefficient, in particular when dealing with large number of sets. Therefore, there have been many efforts to give improved versions of these formulas. A particular example of these methods is given by the \emph{abstract tube} theory \cite{D03,NW92,NW97,NW01}, which associates certain simplicial complexes to probability arrangements.

\begin{Example}\label{network_example}
In order to clarify the concepts just introduced, let's consider the network of figure \ref{network}, which appears in \cite{D03}. It has $5$ nodes and $8$ edges connecting them. We have two selected nodes, labeled $I$ and $O$ that represent the \emph{input} and \emph{output} nodes. We say that the network works if we can connect the two selected nodes.
\begin{figure}
\begin{center}
\begin{tikzpicture}[thin, scale=0.75]
\draw (0,0) -- (4,2) node[pos=0.5, above]{\tiny $1$};
\draw (0,0) -- (4,0) node[pos=0.5, above]{\tiny $2$};
\draw (0,0) -- (4,-2) node[pos=0.5, below]{\tiny $3$};
\draw (4,2) -- (4,0) node[pos=0.5, left]{\tiny $4$};
\draw (4,0) -- (4,-2) node[pos=0.5, left]{\tiny $5$};
\draw (4,2) -- (8,0) node[pos=0.5, above]{\tiny $6$};
\draw (4,0) -- (8,0) node[pos=0.5, above]{\tiny $7$};
\draw (4,-2) -- (8,0) node[pos=0.5, below]{\tiny $8$};
\fill[blue!50,draw=blue] (0,0) circle (2ex) node[text=black]{I};
\fill[red!20,draw=red] (4,2) circle (2ex);
\fill[red!20,draw=red] (4,0) circle (2ex);
\fill[red!20,draw=red] (4,-2) circle (2ex);
\fill[blue!50,draw=blue] (8,0) circle (2ex) node[text=black]{O};
\end{tikzpicture}\caption{An example of coherent system: a network}\label{network}
\end{center} 
\end{figure}
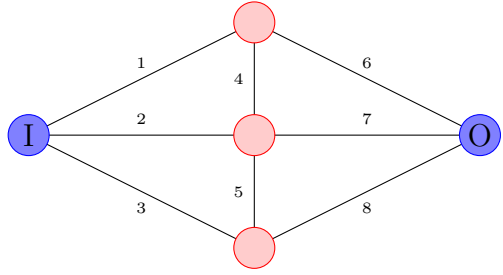
In networks, the minimal nonfailure set is given by the minimal pathes connecting the input and output node, i.e. the leftmost and rightmost node in the figure of the examples. The minimal pathes are $16,147,246,1458,27,3456,258,357$ and $38$. If $0$ encodes the failure state of a path and $1$ its working state. Then we have that for this network the minimal nonfailure set is given by
$$\bar{\mF}^*=\{(1,0,0,0,0,1,0,0),(1,0,0,1,0,0,1,0),(0,1,0,1,0,1,0,0),$$
$$(1,0,0,1,1,0,0,1),(0,1,0,0,0,0,1,0),(0,0,1,1,1,1,0,0),$$
$$(0,1,0,0,1,0,0,1),(0,0,1,0,1,0,1,0),(0,0,1,0,0,0,0,1)\}$$
The intersection of two pathes $\alpha$ and $\alpha'$ is coded by the tuple that has a $1$ in position $i$ if either $\alpha_i$ or $\alpha'_i$ equal $1$. The classical Bonferroni inequality consists in $2^9-1=511$ terms, which is very redundant, since many of the terms are repeated.
If we assume that each path has the same probability $p$ of being in the working state, we have that the reliability of the system is
$$\Rc=3p^2+4p^3-9p^4-10p^5+27p^6-18p^7+4p^8$$
In table \ref{network_bounds} we see the bounds given by the Bonferroni inequalities when truncating the inclusion-exclusion identity of this network for the different values of $r$. The column $s_r$ gives the number of sets used to compute each bound

\begin{table}
\begin{center}
\begin{tabular}{llr}
\hline\\
$r$&bound on $\Rc$&$s_r$\\
\\
\hline\\
$1$&$3p^2+4p^3+2p^4$&$9$\\
$2$&$3p^2+4p^3-9p^4-16p^5-9p^6$&$45$\\
$3$&$3p^2+4p^3-9p^4-8p^5+34p^6+30p^7+3p^8$&$129$\\
$4$&$3p^2+4p^3-9p^4-10p^5+27p^6-50p^7-34p^8$&$255$\\
$5$&$3p^2+4p^3-9p^4-10p^5+27p^6-12p^7+54p^8$&$381$\\
$6$&$3p^2+4p^3-9p^4-10p^5+27p^6-18p^7-24p^8$&$465$\\
$7$&$3p^2+4p^3-9p^4-10p^5+27p^6-18p^7+12p^8$&$501$\\
$8$&$3p^2+4p^3-9p^4-10p^5+27p^6-18p^7+3p^8$&$510$\\
$9$&$3p^2+4p^3-9p^4-10p^5+27p^6-18p^7+4p^8$&$511$\\
\hline\\
\end{tabular}\caption{Bounds for the reliability of the network in figure \ref{network} (table taken from \cite{D03})} \label{network_bounds}
\end{center}
\end{table}
\end{Example}

\subsection{Coherent systems and monomial ideals}
\spanishsubsection{Sistemas coherentes e ideales monomiales}

A main result in \cite{GW03,GW04} establishes the parallelism between monomial ideals and nonfailure sets of coherent systems. 
\begin{Proposition}[\cite{GW04}]
Given a system $\Sc$ of $n$ components with nonfailure set $\bar{\mF}$ and minimal nonfailure set $\bar{\mF}^*$:
\begin{enumerate}
\item The minimal nonfailure points in $\bar{\mF}^*$, seen as the exponent vectors of monomials in $R=\kb[x_1,\dots,x_n]$ are the (exponents of the) minimal generators of a monomial ideal $I_\Sc$
\item The points in $\bar{\mF}$ represent the (exponents of the) monomials belonging to the monomial ideal $I_\Sc$ generated by the minimal failure points.
\item The points in $\mF$ represent the (exponents of the) monomials belonging to the factor ring $R/I$.
\end{enumerate}
\end{Proposition}

From this proposition, we see that the problem of finding the reliability of the system, i.e. $\Rc=Prob(\bar{\mF})=Prob(\bigcup_{\alpha\in\bar{\mF}^*}Q_\alpha)$ can be faced in an algebraic setting. In this algebraic framework, the orthants $Q_\alpha$ are just the principal ideals $\langle x^\alpha\rangle$, and the intersections of orthants correspond to intersection of these ideals. Therefore, the problem of computing the reliability $\Rc$ of $\Sc$ uses the computation of the multigraded Hilbert series of $I_\Sc$, which enumerates the monomials in $I_\Sc$. What we have to do is to count the monomials in $I_\Sc$ which are given by $H_{I_\Sc}(\xb)$, and divide by the monomials in $R$, which is given by $H_R(\xb)$. Since $\frac{H_{I_\Sc}(\xb)}{H_R(\xb)}=\Kc_{I_\Sc}$, we just have to compute the $\Kc$-polynomial of $I_\Sc$.

As we have seen in section \ref{algebra_monomial}, there are several options for this computation. One way is computing the multigraded Betti numbers of $I_\Sc$ and then, the $\Kc$-polynomial of $I_\Sc$ is just the alternating sum of them. These Betti numbers can be obtained via the computation of the Koszul homology or equivalently of the minimal free resolution of $I_\Sc$. This provides the $\Kc$-polynomial with minimal redundancy. Moreover, the alternating sum of the ranks of the modules in any free resolution of $I_\Sc$ provide also $\Kc$-polynomials which allow us to compute the multigraded Hilbert function of the ideal. These computations are normally faster, but they provide in general $\Kc$-polynomials with bigger redundancy. Since Mayer-Vietoris trees support resolutions of their corresponding ideals, we can use $\Kc$-polynomials coming from the results of them.

\begin{Remark}
The Taylor resolution of an ideal $I$ has one generator for every $J\subseteq\{1,\dots,r\}$, where $r$ is the number of generators of $I$. The $\Kc$-polynomial of it corresponds to the classical inclusion-exclusion formula. In \cite{GW04} formulas corresponding to Scarf complexes are used. These provided sharper bounds and better identities, which correspond to the ones obtained using the \emph{abstract tube} theory from \cite{NW92,NW97,NW01}.
\end{Remark}

\begin{Example}
The ideal corresponding to the nonfailure set of the network $\Sc$ in example \ref{network_example} is
$$I_\Sc=\langle x_1x_6,x_1x_4x_7,x_2x_4x_6,x_1x_4x_5x_8,x_2x_7,x_3x_4x_5x_6,x_2x_5x_8,x_3x_5x_7,x_3x_8\rangle$$
In this case, the Mayer-Vietoris tree of $I$ with respect to degree-lexicographic term order, gives us the exact multigraded Betti numbers, so we can avoid the computation of the Koszul homology or minimal free resolution to obtain them. The Betti numbers are $\beta_0=9$, $\beta_1=25$, $\beta_2=31$, $\beta_3=18$ and $\beta_4=4$, so that we have that the $\Kc$-polynomial of $I_\Sc$ has just $87$ summands. This is a big improvement with respect to the formula obtained by the Taylor resolution, corresponding to the classical inclusion-exclusion principle, which has $511$ terms, and improves also the formula in \cite{GW04}, based on the Scarf complex, which has $103$ elements.
\end{Example}

\subsection{Examples}
\spanishsubsection{Ejemplos}

We finish this section with the application of the techniques seen in chapters \ref{structure} and \ref{computation}. We study a very natural type of networks, called series-parallel networks. Then we move to some of the most studied non-network examples, namely $k$-out-of$n$ systems and some of their generalizations and special types; a throughout presentation of these types of systems can be seen in \cite{KZ03}, where some of the systems studied here and many more are treated. 
\subsubsection{Series-parallel networks}
\spanishsubsubsection{Redes series-paralelo}

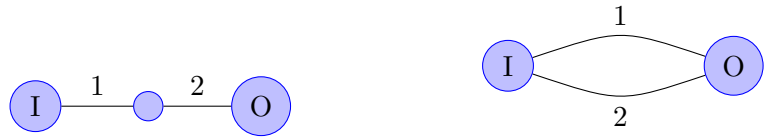
\begin{figure}
\begin{center}
$\begin{array}{l@{\hspace{1in}}r}

\begin{tikzpicture}[scale=.5]
\begin{scope}[shape=circle,minimum size=0.1 cm,fill=blue!25, draw=blue]
\tikzstyle{every node}=[draw,fill]
\node (q_1) at (0,1) {\small I};
\node (q_2) at (3,1) {};
\node (q_3) at (6,1) {\small O};
\end{scope}
\draw (q_1) -- (q_2) node[pos=0.5, above]{\small $1$};
\draw (q_2) -- (q_3) node[pos=0.5, above]{\small $2$};

\end{tikzpicture}

&

\begin{tikzpicture}[scale=.5]
\begin{scope}[shape=circle,minimum size=0.1 cm,fill=blue!25, draw=blue]
\tikzstyle{every node}=[draw,fill]
\node (q_1) at (0,1) {\small I};
\node (q_2) at (6,1) {\small O};
\end{scope}
\draw (q_1) ..controls +(3,1) .. (q_2) node[pos=0.5, above]{\small $1$};
\draw (q_1) ..controls +(3,-1) .. (q_2) node[pos=0.5, below]{\small $2$};

\end{tikzpicture}
\\
\end{array}$
\end{center}\caption{Basic series and parallel networks.}\label{series-parallel}
\end{figure}

Our first example will be a (rather big) class of networks, which we call \emph{series-parallel networks}. Assume we connect two nodes by two edges, we can do it either in parallel: both edges connect both nodes or in series: the first edge connects the first node and the beginning of the second edge, and the second edge ends at the second node (see figure \ref{series-parallel}).

\begin{figure}
\begin{center}
$\begin{array}{c@{\hspace{1in}}c}

\begin{tikzpicture}[scale=.5]
\begin{scope}[shape=circle,minimum size=0.1 cm,fill=blue!25, draw=blue]
\tikzstyle{every node}=[draw,fill]
\node (q_1) at (0,1) {\small I};
\node (q_2) at (3,1) {};
\node (q_3) at (6,1) {\small O};
\end{scope}
\draw (q_1) -- (q_2) node[pos=0.5, above]{\small $1$};
\draw (q_2) -- (q_3) node[pos=0.5, above]{\small $2$};

\end{tikzpicture}

&

\begin{tikzpicture}[scale=.5]
\begin{scope}[shape=circle,minimum size=0.1 cm,fill=blue!25, draw=blue]
\tikzstyle{every node}=[draw,fill]
\node (q_1) at (0,1) {\small I};
\node (q_2) at (6,1) {\small O};
\end{scope}
\draw (q_1) ..controls +(3,1) .. (q_2) node[pos=0.5, above]{\small $1$};
\draw (q_1) ..controls +(3,-1) .. (q_2) node[pos=0.5, below]{\small $2$};

\end{tikzpicture}
 \\ [0.4cm]
N_1 & N_2\\[0.5cm]

\begin{tikzpicture}[scale=.5]
\begin{scope}[shape=circle,minimum size=0.1 cm,fill=blue!25, draw=blue]
\tikzstyle{every node}=[draw,fill]
\node (q_1) at (0,1) {\small I};
\node (q_2) at (3,-4) {};
\node (q_3) at (6,1) {\small O};
\end{scope}
\draw (q_1) -- (q_2) node[pos=0.5, left]{\small $11:=1$};
\draw (q_2) -- (q_3) node[pos=0.5, right]{\small $12:=2$};
\draw (q_1) ..controls +(3,1) .. (q_3) node[pos=0.5, above]{\small $21:=3$};
\draw (q_1) ..controls +(3,-1) .. (q_3) node[pos=0.5, below]{\small $22:=4$};

\end{tikzpicture}

&

\begin{tikzpicture}[scale=.5]
\begin{scope}[shape=circle,minimum size=0.1 cm,fill=blue!25, draw=blue]
\tikzstyle{every node}=[draw,fill]
\node (q_1) at (0,1) {\small I};
\node (q_2) at (4,1) {};
\node (q_3) at (8,1) {\tiny S};
\node (q_4) at (12,1) {\small O};

\end{scope}
\draw (q_1) -- (q_2) node[pos=0.5, above]{\small $11:=1$};
\draw (q_2) -- (q_3) node[pos=0.5, above]{\small $12:=2$};
\draw (q_3) ..controls +(2,1) .. (q_4) node[pos=0.5, above]{\small $21:=3$};
\draw (q_3) ..controls +(2,-1) .. (q_4) node[pos=0.5, below]{\small $22:=4$};

\end{tikzpicture}
 \\ [0.4cm]
N_1\times N_2 & N_1 + N_2\\
\end{array}$
\end{center}\caption{Example of series-parallel network construction.}\label{series-parallel-construction}
\end{figure}
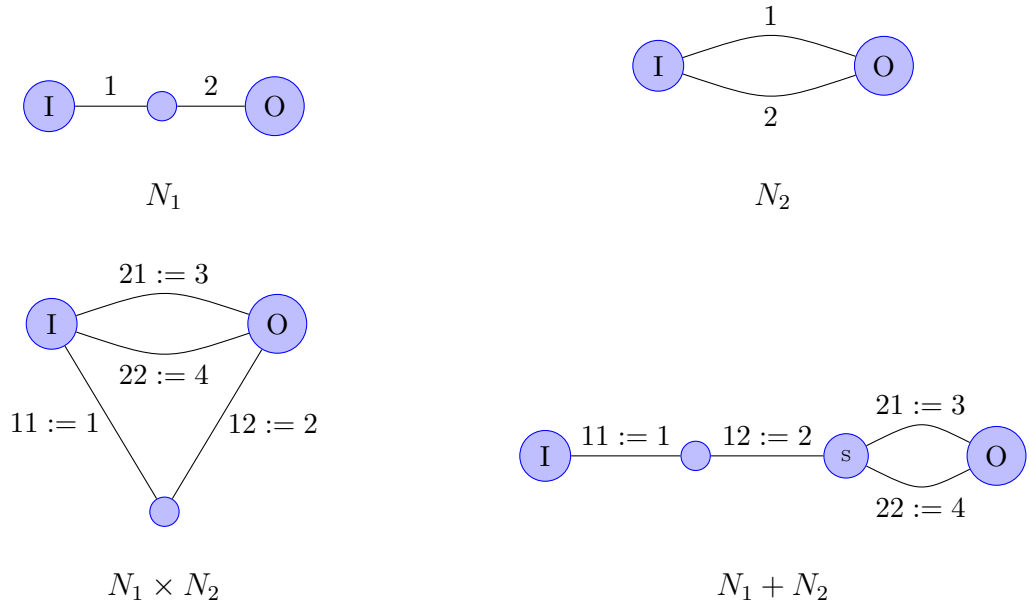

This basic construction can be generalized in the following way: We consider a edge $p$ joining two nodes $I$ and $O$ a series-parallel network. We call such a network a \emph{basic} series-parallel network. Consider now two series-parallel networks $N_1$ and $N_2$. We can connect them in series or in parallel, and the result is a series-parallel network. This is done in the following way: 
\begin{itemize}
\item First, in any case, we rename the edges in each node so that each edge has a different label. If the edge $p_S$ for some (possibly empty) set $S$ of subindices is in network $i$ we can rename it $p_{\{i\}\cup S}$. After this, we can still rename them just by counting them in lexicographic order.
\item If the initial (input) node of $N_i$ is labelled $I_i$ and its final (output) node is labelled $O_i$  for $i=1,2$, then the parallel union of $N_1$ and $N_2$, which we will denote $N=N_1 \times N_2$ identifies $I_1$ and $I_2$ in one node $I$, which will be the initial node of $N$ , and identifies $O_1$ and $O_2$ in one node $O$, which will be its final node.
\item With the same notation as above, the series union of $N_1$ and $N_2$, which we will denote $N=N_1 + N_2$ has as initial node $I_1$, as final node $O_2$, and identifies $O_1$ and $I_2$ in one intermediate node $S$.
\end{itemize}
We just formalize these considerations in the following definition of \emph{series-parallel networks:}
\begin{Definition}
We say that a network $N$ is a \emph{parallel-series network} if either $N$ consists of an input node, an output node and a edge joining them, or if $N=N_1+N_2$ or $N=N_1\times N_2$ with $N_1,N_2$ series-parallel networks.
\end{Definition}
These constructions can be seen in figure \ref{series-parallel-construction}, in which the label of the edge $p_S$ is just $S$.

Given any network $N$ (not necessary a series-parallel one), we associate a monomial ideal to it as follows: Consider one variable $x_S$ for each connection $S$ in $N$. Then, the monomial ideal $I_N$ associated to $N$ is minimally generated by the monomials $x_{S_1}\cdots x_{S_k}$ where $S_1,\dots,S_k$ is a minimal cut in the network $N$ \cite{D03,GW04}. Let us consider now the ideals associated to series-parallel networks. It is clear
that the ideal $I_N$ of a network $N$ with just one edge $p_1$
connecting two nodes $I$ and $O$ is just $I_N=\langle x_1\rangle$.
The construction operations $+$ and $\times$ we have just seen, have
their counterpart in the ideals of the resulting networks:

\begin{Proposition}
Let $N_1$ and $N_2$ be two networks the edges of which are labelled (after renaming as seen above) $p_1,\dots,p_{n_1}$ and $p_{n_1+1},\dots,p_{n_1+n_2}$. Then,
$$I_{N_1+N_2}=I_{N_1}+I_{N_2}\qquad I_{N_1\times N_2}=I_{N_1}\cap I_{N_2}$$
where $I_{N_1+N_2}$ and $I_{N_1\times N_2}$ are ideals in $\kb[x_1,\dots,x_{n_1+n_2}]$
\end{Proposition}

\noindent{\bf Proof: }We have that $$I_N=\langle x_S\vert S=\{s_1,\dots,s_{k_s}\} \mbox{ is a minimal cut in N}\rangle$$

Any minimal cut in $N_1$ or $N_2$ is a minimal cut in $N$, and there is no mixture among them. Then it is clear that the generating set of $I_{N_1+N_2}$ is just the union of the generating sets of $I_{N_1}$ and $I_{N_2}$, each being generated in a different set of variables.

Now, the minimal cuts of $N_1\times N_2$ can be always considered as a combination of one minimal cut in
$N_1$ and one minimal cut in $N_2$ and there are no other minimal pathes. Since there is no
intersection between the set of variables of $I_{N_1}$ and $I_{N_2}$
the combination of minimal cuts simply means product of their variables, and hence the result.$\square$

\begin{Example}
Let's take the networks in figure \ref{series-parallel-construction}. After relabelling, the edges in $N_1$ are $p_1$ and $p_2$, and the edges in $N_2$ are $p_3$ and $p_4$. We have that
$$I_{N_1}=\langle x_1x_2\rangle,\quad I_{N_2}=\langle x_3,x_4\rangle,\quad I_{N_1+N_2}=\langle x_1x_2,x_3,x_4\rangle,\quad I_{N_1\times N_2}=\langle x_1x_2x_3,x_1x_2x_4\rangle$$
\end{Example}

Observe that this result is quite general in the sense that the original networks need not to be series-parallel. For ideals of series-parallel networks, we have that Mayer-Vietoris trees give a good way to compute their multigraded Betti numbers, and hence, the reliability of the corresponding network:

\begin{Proposition}
The ideals associated to series-parallel networks, i. e. \emph{series-parallel ideals}, are Mayer-Vietoris ideals of type A.
\end{Proposition}
\noindent {\bf Proof:}
Recall that an ideal $I$ is Mayer-Vietoris of type A if there is some Mayer-Vietoris tree of $I$ such that there are no repeated generators in the relevant nodes. 
If $N$ is a basic series-parallel network with unique edge $p_1$ then $I_N=\langle x_1\rangle$ which is Mayer-Vietoris of type A.
Now consider two series-parallel networks $N_1$ and $N_2$ whose ideals are Mayer-Vietoris of type A, i.e. there is some strategy for selecting the pivot monomials when constructing a Mayer-Vietoris tree such that it is of type $A$. We have to proof that $I_{N_1}+I_{N_2}$ and $I_{N_1}\cap I_{N_2}$ are Mayer-Vietoris of type A:
\begin{itemize}
\item  The generators of $I_{N_1}+I_{N_2}$ are the union of the generating sets of $I_{N_1}$ and $I_{N_2}$. We sort them so that the generators of $I_{N_2}$ all appear after the generators of $I_{N_1}$. We now proceed taking as pivot monomial always a generator of $I_{N_2}$ following the strategy used to build the minimal Mayer-Vietoris tree of $I_{N_2}$. Doing so, we have that $MVT_p(I_{N_1}+I_{N_2})$ has as generators the generators of $I_{N_1}$ each one multiplied by some product of the variables of $I_{N_2}$ and also the generators of $MVT_p(I_{N_2})$. So far, we have no repeated generators in the relevant nodes: Assume that there is some generator repeated in two relevant nodes at positions $p$ and $q$ then they have the same exponents in the variables of $I_{N_1}$ and the same in the variables of $I_{N_2}$. If the generator has only variables of the second ideal, to be equal would mean that there are equal in $MVT(I_{N_2})$. And since those generators with `mixed variables' all are of the form $m\cdot m'$ with $m$ a minimal generator of $I_{N_1}$, no two of these are repeated.

This procedure takes us to nodes in which no further element only in the variables of $I_{N_2}$ is available. From this moment, on each node we follow the strategy of $MVT(I_{N_1})$. Since these nodes in positions $p$ have as generators all the minimal generators of $I_{N-1}$ times some polynomial $ m'_p$ in the variables of the second ideal. And since the $m'_p$ are different for different $p$, we have that al the trees hanging from these nodes are isomorphic to $MVT(I_{N_1})$, therefore, there's no repeated generator in the relevant nodes in each of them. There is also no repetition among the different `copies' of $MVT(I_{N_1})$ because of each ${m'}_p$ is unique.

\item $I_{N_1}\times I_{N_2}$: Let us denote by $m_1,\dots,m_r$ the generators of $I_{N_1}$, and by $n_1,\dots,n_s$ the generators of $I_{N_2}$. And assume without loss of generality that the minimal Mayer-Vietoris trees of $I_{N_1}$ and $I_{N_2}$ were obtained using always the last generator as the pivot monomial. Then $I_{N_1\times N_2}=I_{N_1}\cap I_{N_2}$ is generated by $\{m_in_j\vert i=1,\dots,r;j=1,\dots,s\}$. To build $MVT(I_{N_1\times N_2})$ consider as pivot monomial $m_rn_s$; then, $MVT_2(I_{N_1\times N_2})$ is generated by

$$m_rn_{1,s},\dots,m_rn_{s-1,s},m_{1,r},\dots,m_{r-1,r}n_s$$

now, select $m_{r-1}n_s$ as pivot monomial in $MVT_3(I_{N_1\times N_2})$ and we obtain that $MVT_6(I_{N_1\times N_2})$ is generated by

$$m_{r-1}n_{1,s},\dots,m_{r-1}n_{s-1,s},m_{1,r-1}n_s,\dots,m_{r-2,r-1}n_s$$

It is clear that since $I_{N_1}$ and $I_{N_2}$ are Mayer-Vietoris of type $A$ and they are generated in disjoint sets of variables, there is no repeated relevant multidegree inside the subtree of $MVT(I_{N_1\times N_2})$ hanging from $MVT_2(I_{N_1\times N_2})$ or inside the subtree hanging from $MVT_6(I_{N_1\times N_2})$.

On the other hand we have that no multidegree appears in a relevant node in both trees. The reason is the following: Every multidegree in $MVT_2(I_{N_1\times N_2})$ is of the form $n_{\sigma}m_\alpha$ with $s\in\sigma$ and $r\in\alpha$; and every element in $MVT_6(I_{N_1\times N_2})$ is of the form $n_{\sigma'}m_{\alpha'}$ with $s\in\sigma$, $\alpha\subseteq\{1,\dots,r-1\}$. If $n_{\sigma}m_\alpha=n_{\sigma'}m_{\alpha'}$ then in particular $m_\alpha=m_{\alpha'}$, but then they would be repeated in $MVT(I_{N_1}$, which is a contradiction.

Following the same argument while taking elements of the form $m_in_s$ as pivot monomials in the nodes of dimension $0$ we obtain the same contradiction, based on the fact that $MVT(I_{N_1})$ has no repeated relevant multidegrees. Then we turn to take pivot monomials of the form $m_in_{s-1}$ ($i\in\{1,\dots,r\}$) and we can follow a symmetric argument, the contradiction coming now from the fact that $MVT(I_{N_2})$ has no repeated relevant multidegrees. A recursive application of these two arguments yields the result.

\end{itemize}

So, both $I_{N_1+N_2}$ and $I_{N_1\times N_2}$ are Mayer-Vietoris of type $A$.$\square$

It is easy to see that the union ideal of two ideals $I_1$ and $I_2$ generated in two disjoint sets of variables is splittable in the sense of \cite{EK90}. In particular this applies to series parallel networks: Let $I(N_1)=\langle m_1,\dots,m_{r_1}\rangle$ and $I(N_2)=\langle m_1',\dots,m_{r_2}'\rangle$ be the ideals of two series-parallel networks with the same conditions as above. Then the ideal $I(N_1)+I(N_2)$ is \emph{splittable}.
The obvious splitting function is
$$G(I(N_1)\cap I(N_2))\longrightarrow G(I(N_1))\times G(I(N_2))$$
$$m_im_{j}'\longmapsto(m_i,m_{j}')$$
where $G(I)$ denotes the set of minimal generators of $I$.
Then we have the following formula:
$$\beta_q(I(N_1)+I(N_2))=\beta_q(I(N_1))+\beta_q(I(N_2))+\beta_{q-1}(I(N_1)\cap I(N_2))$$

Consider now a network $N_1\times N_2$. It has $r_1\cdot r_2$ minimal generators. Then using the above formula, we can compute the Betti numbers of $I(N_1\times N_2)$ from those of $I(N_1)$, $I(N_2)$ and $I(N_1+N_2)$, which have, respectively $r_1$, $r_2$ and $r_1+r_2$ generators, so the complexity of the problem decreases if $r_1\cdot r_2<r_1+r_2$ otherwise, from the intersection ideal we can compute the Betti numbers of the union ideal. When computing the Betti numbers of these kind of ideals, we just compute normally and whenever we find an intersection or a union of networks, we compute its ideal or the corresponding `dual' depending on the number of generators.

In particular, a pure series network $N$, has an ideal of the form $I_N=\langle x_1,\dots,x_r\rangle$ and then it's Betti numbers are given by
$$\beta_i(I_N)={{r}\choose{i+1}}$$
and the multidegrees of the Betti numbers of degree $i$ are all the products of $i+1$ of the variables. The simplest case of application of the above comments is when we have two pure series networks, and we want to compute the parallel connection of them, $N_1\times N_2$. Then, we have that
$$\beta_i(I(N_1)\cap I(N_2))=\beta_{i+1}(I(N_1)+ I(N_2))+\beta_{i+1}(I(N_1))+\beta_{i+1}(I(N_2))$$
Since $I(N_1)+I(N_2)$ is again a pure series network, we have that
$$\beta_{i}(I(N_1)\cap I(N_2))={{r_1+r_2}\choose{i+2}}-{{r_1}\choose{i+2}}-{{r_2}\choose{i+2}}$$
and the multidegrees of the $i$-th Betti numbers are those products of $i+2$ variables such that there are always some variable belonging to the support of each of the original networks.

In a more general setting, the situation is the following: 

Consider the parallel connection of several pure series networks, i.e.
$$N=\{N_{11}+\cdots + N_{1r_1}\}\times\cdots \times\{N_{k1}+\cdots+ N_{kr_k}\}$$
where $N_{ij}$ are basic networks, i.e. the have only one connection. If we associate a variable $x_{ij}$ to each of these basic networks, then the corresponding ideals are of the form $I_N=\langle x_{11},\dots,x_{1r_1}\rangle \cap \cdots \cap \langle x_{k1},\dots,x_{kr_k}\rangle$, which has $\prod_{1}^k r_k$ generators in $\sum_1^k r_k$ variables. We will denote these networks by $N\equiv\{r_1,\dots,r_k\}$ and the correspondent ideal by $I_{\{r_1,\dots,r_k\}}$. These ideals are splittable and separable, and it is easy to see that they have linear resolutions. We can use their separable Mayer-Vietoris trees to obtain a recursive way to compute their Betti numbers. The basic procedure is to use the process given in section \ref{separable} to have the Betti numbers of $I_{\{r_1,\dots,r_k\}}$  form those of $I_{\{r_1,\dots,r_{k-1}\}}$ and $I_{\{r_k\}}$. It is obvious that $I_{\{r_k\}}=\langle x_{k1},\dots,x_{kr_k}\rangle$ is Mayer-Vietoris of type $A$. Moreover, its Betti numbers are given by $\beta_i(I_{\{r_k\}})={{r_k}\choose{i+1}}$. The Mayer-Vietoris tree of this type of ideals is \emph{complete}, in the sense that the children of a node with $j$ generators has $j-1$ generators, hence the tree has $2^{r_k}-1$ nodes. A carefull inspection of the tree reveals that it has $\sum_{j=1}^{r_k-1}{{j-1}\choose{i-1}}$ relevant nodes of dimension $i$. Finally, this tree has $\sum_{j=1}^{min\{r_k,i\}}{{r_k-1}\choose{j-1}}$ final nodes in even position, and $\sum_{j=0}^{min\{r_k,i-1\}}{{r_k-1}\choose{j}}$ final nodes in odd position. Hence, our formula is
$$\beta_i(I_{\{r_1,\dots,r_k\}})={{r_k}\choose{i+1}}+\sum_{j=1}^{r_k-1}{{j-1}\choose{i-1}}\cdot r_{k-1}+\sum_{j=1}^{min\{r_k,i\}}{{r_k-1}\choose{j-1}}\beta_{i-j}(I_{\{r_1,\dots,r_{k-1}\}})$$ $$+\sum_{j=0}^{min\{r_k,i-1\}}{{r_k-1}\choose{j}}\beta_{i-j}(I_{\{r_1,\dots,r_{k-1}\}})$$

\begin{Remark}
The ideals associated to these network grow very rapidly, and they become impracticable for their study using traditional algebraic methods to compute their Hilbert function and/or Betti numbers. For example, the ideal $I_{\{4,2,5,3,4\}}$ is generated by $480$ minimal generators in $18$ variables.
\end{Remark}

The dual situation to the one just presented consists of networks presented as a series combination of parallel blocks:
$$N=\{N_{11}\times\cdots\times N_{1k_1}\}+\cdots +\{N_{r1}\times\cdots\times N_{rk_r}\}$$
The ideal associated to $N$ will be denoted $I_{[r_1,\dots,r_k]}$. These networks are much easier to compute, since the ideal corresponding to any pure-parallel network consists just of one monomial, namely the product of all variables representing connections in the network. These ideals are generated in $\sum_1^k r_k$ variables and have $k$ generators.
The ideal of the series connection of such networks is then isomorphic to a pure-series network and hence its Betti numbers are given by
$$\beta_i(I_{[r_1,\dots,r_k]})={{k}\choose {i+1}}$$
With these considerations we see that we can handle in the same way we have seen above the parallel connection either of pure series networks or of series connections of pure parallel networks. Let us see it with an example that finishes the section:
\begin{Example}
Let us consider the network in figure \ref{series-parallel network}. It is the parallel connection of four series connections of pure-parallel networks. We call each of these four sub-networks $A$, $B$, $C$ and $D$. Network $A$ is generated by $6$ minimal generators in the $12$ variables {\sf a01, a11, a12, a13, a14, a21, a31, a32, a33, a41, a42, a51}. Network $B$ is generated by $6$ minimal generators in the $10$ variables {\sf b01, b11, b12, b13, b21, b31, b41, b42, b43, b51}.Network $C$ is generated by $5$ minimal generators in the $7$ variables {\sf c01, c11, c21, c31, c32, c33, c41}. And finally, network $D$ is generated by $6$ minimal generators in the $9$ variables {\sf d01, d11, d21, d22, d23, d24, d31, d41, d15}. Therefore, our the ideal of our network has $1080$ generators in $36$ variables. Note that each of the inner parallel components, like the one formed by {\sf a11, a12, a13, a14} in network $A$ can be substituted by just a single connection {\sf a1}, and the resulting ideal would be isomorphic to the starting one, via ${\sf a11a12a13a14}\mapsto {\sf a1}$. Therefore, the network in the example is denoted by $N=\{6,6,5,6\}$. It is generated by $1080$ minimal generators in $23$ variables (after relabelling the inner parallel components by single variables). Using the procedure described above, we have that the Betti numbers of the ideal associated to this network are
$$\beta_0=1080,\,\beta_1=10260,\,\beta_2=49410,\,\beta_3=158355,\,\beta_4=375606$$
$$\beta_5=696465,\,\beta_6=1042195,\,\beta_7=1283165,\,\beta_8=1314953,\,\beta_9=1128358$$
$$\beta_{10}=812090,\,\beta_{11}=489090,\,\beta_{12}=244953,\,\beta_{13}=100926,\,\beta_{14}=33648$$
$$\beta_{15}=8855,\,\beta_{16}=1771,\,\beta_{17}=253,\,\beta_{18}=23,\,\beta_{19}=1$$

\begin{figure}
\begin{center}
\begin{tikzpicture}[scale=.5, transform shape]
\begin{scope}[shape=circle,minimum size=0.1 cm,fill=blue!25, draw=blue]
\tikzstyle{every node}=[draw,fill]
\node (input) at (0,0) {\small I};
\node (a_1) at (5,5) {};\node (a_2) at (10,5) {};\node (a_3) at (15,5) {};\node (a_4) at (20,5) {};\node (a_5) at (25,5) {};
\node (b_1) at (5,0) {};\node (b_2) at (10,0) {};\node (b_3) at (15,0) {};\node (b_4) at (20,0) {};\node (b_5) at (25,0) {};
\node (c_1) at (6,-5) {};\node (c_2) at (12,-5) {};\node (c_3) at (18,-5) {};\node (c_4) at (24,-5) {};
\node (d_1) at (5,-10) {};\node (d_2) at (10,-10) {};\node (d_3) at (15,-10) {};\node (d_4) at (20,-10) {};\node (d_5) at (25,-10) {};
\node (output) at (30,0) {\small O};

\end{scope}
\draw (input) -- (a_1) node[pos=0.5, left]{ \textsf{a01}};
\foreach \sl / \lbl in {2/11,1/12,-1/13,-2/14} {\draw (a_1)..controls +(2.5,\sl) .. (a_2) node[pos=0.5, above]{\textsf {a\lbl}};}
\draw (a_2)--(a_3) node[pos=0.5, above]{ \textsf{a21}};
\foreach \sl / \lbl in {1/31,0/32,-1/33} {\draw (a_3)..controls +(2.5,\sl) .. (a_4) node[pos=0.5, above]{\textsf {a\lbl}};}
\foreach \sl / \lbl in {1/41,-1/42} {\draw (a_4)..controls +(2.5,\sl) .. (a_5) node[pos=0.5, above]{\textsf {a\lbl}};}
\draw (a_5)--(output) node[pos=0.5, right]{ \textsf{a51}};

\draw (input) -- (b_1) node[pos=0.5, above]{ \textsf{b01}};
\foreach \sl / \lbl in {1/11,0/12,-1/13} {\draw (b_1)..controls +(2.5,\sl) .. (b_2) node[pos=0.5, above]{\textsf {b\lbl}};}
\draw (b_2)--(b_3) node[pos=0.5, above]{ \textsf{b21}};
\draw (b_3)--(b_4) node[pos=0.5, above]{ \textsf{b31}};
\foreach \sl / \lbl in {1/41,0/42,-1/43} {\draw (b_4)..controls +(2.5,\sl) .. (b_5) node[pos=0.5, above]{\textsf {b\lbl}};}
\draw (b_5)--(output) node[pos=0.5, above]{ \textsf{b51}};

\draw (input) -- (c_1) node[pos=0.5, above]{ \textsf{c01}};
\draw (c_1)--(c_2) node[pos=0.5, above]{ \textsf{c11}};
\draw (c_2)--(c_3) node[pos=0.5, above]{ \textsf{c21}};
\foreach \sl / \lbl in {1/31,0/32,-1/33} {\draw (c_3)..controls +(2.5,\sl) .. (c_4) node[pos=0.5, above]{\textsf {c\lbl}};}
\draw (c_4)--(output) node[pos=0.5, above]{ \textsf{c41}};

\draw (input) -- (d_1) node[pos=0.5, left]{ \textsf{d01}};
\draw (d_1)--(d_2) node[pos=0.5, above]{ \textsf{d11}};
\foreach \sl / \lbl in {2/21,1/22,-1/23,-2/24} {\draw (d_2)..controls +(2.5,\sl) .. (d_3) node[pos=0.5, above]{\textsf {d\lbl}};}
\draw (d_3)--(d_4) node[pos=0.5, above]{ \textsf{d31}};
\draw (d_4)--(d_5) node[pos=0.5, above]{ \textsf{d41}};
\draw (d_5)--(output) node[pos=0.5, right]{ \textsf{d51}};

\end{tikzpicture}
\end{center}\label{series-parallel network}
\end{figure}

\end{Example}

\subsubsection{$k$-out-of-$n$ systems}
\spanishsubsubsection{Sistemas $k$-entre-$n$}

The first non-network examples we will study are $k$-out-of-$n$ systems, introduced by Birnbaum et al. \cite{BES61} (see also \cite{KZ03}). A $k$-out-of-$n$:G (G for `good') has $n$ components, the system operates (is `good') whenever $k$ or more components operate (are `good'). Symmetrically, $k$-out-of-$n$:F (F for `failure') systems are defined. Since G and F $k$-out-of-$n$ systems can be expressed in terms of each other, we will just refer to $k$-out-of-$n$ systems. A $k$-out-of-$n$ system can be modelled by the ideal
$$I_{k,n}=\langle x^\mu \mbox{ where }\mu \mbox{ is a squarefree monomial of degree k in n variables} \rangle$$

for example, $I_{3,5}=\langle xyz,xyu,xyv,xzu,xzv,xuv,yzu,yzv,yuv,zuv\rangle$ is the ideal corresponding to the $3$-out-of-$5$ problem. Observe that $I_{k,n}$ has a minimal generating set which consists of ${n}\choose{k}$ monomials. Using the result pointed in section \ref{simplicial_Koszul_comp}, we know that we have to check the Koszul homology only in the multidegrees that are in the lcm-lattice of $I$, $L_I$. It is easy to see that $L_I$ consists of all squarefree monomials involving a number of variables between $k$ and $n$. The following lemma characterizes the Koszul simplicial complex at each of these multidegrees:

\begin{Lemma}
If $\ab\in L_{I_{k,n}}$ has $k+i$ nonzero indices, $k<k+i\leq n$, the simplicial Koszul complex $\Delta_\ab^I$ consists of all $j$-faces with $0\leq j\leq i-1$.
\end{Lemma}

\noindent {\bf Proof:} Let $\xb^\ab$ be a squarefree monomial consisting of the product of $k+i$ variables, $k<k+i\leq n$. If we divide $\xb^\ab$ by the product of $j$ of these variables then: If $j\leq i$ then the resulting monomial is the product of a set of $k+i-j$ variables, and thus, a $j-1$ face is present in the Koszul simplicial complex. If $j>i$ then the result of the division is the product of $k+i-j$ variables, being $j>i$, $k+i-j<k$ and thus this product is not in $I$, so no $j-1$ face is in the simplicial Koszul complex for $j>i$.$\square$

Thus, the $(i,\ab)$-th Betti number at the multidegree given by any combination of $k+i$ variables is $dim(\tilde{H}_{i-1}(C_{k,i}))$, where $C_{k,i}$ is the subcomplex of the $k+i$ dimensional simplex $\Delta_{k+i}$ having as facets all the $(i-1)$-faces. And then, $\beta_i(I_{k,n})= {{n}\choose{k+i}}\cdot dim(\tilde{H}_{i-1}(C_{k,i}))$, for all $i\in\{0,\dots,n-k\}$.

Our next goal is then to compute the dimension of the reduced homology of the complexes $C_{k,i}$. Since all faces in dimension less or equal $i-1$ are present in the complex, we know that $C_{k,i}$ has zero homology at all dimensions except possibly at dimension $i-1$. The chain complex of $C_{k,i}$ has the following form:
$$0\rightarrow C_{i-1}\stackrel{\delta_{i-1}}{\rightarrow}\cdots\rightarrow C_{1}\stackrel{\delta_{1}}{\rightarrow}C_0\rightarrow 0$$
we have $\tilde{H}_j(C_{k,i})=0\,\forall j<i-1$ thus $ker\,\delta_j/im\,\delta_{j+1}=0$ and $dim(ker\,\delta_j)=dim(im\,\delta_{j+1})$ for all $j<i-1$. On the other hand, we have the usual equality
$$dim(ker\,\delta_j)=dim (C_j)-dim(im\,\delta_j)$$
putting these together we have that
$$dim(\tilde{H}_{i-1})(C_{k,i})=dim(ker\,\delta_{i-1})={{k+i}\choose {i-1}}-{{k+i}\choose {i-2}}+\cdots+(-1)^{i-2}{{k+i}\choose {1}}+(-1)^{i-1}$$
We can use now the following combinatorial identity:
$${{k+i}\choose {i-1}}-{{k+i}\choose {i-2}}+\cdots+(-1)^{i-2}{{k+i}\choose {1}}+(-1)^{i-1}={{i+k-1}\choose{k-1}}$$
and we obtain that for every $\ab\in L_I$ where $\ab$ is the product of $k+i$ variables, we have that
$$\beta_{(i,\ab)}(I_{k,n})={{i+k-1}\choose{k-1}}$$
and since we have ${{n}\choose{k+i}}$ such $\ab$, then
$$\beta_i(I_{k,n})={{n}\choose{k+i}}\cdot{{i+k-1}\choose{k-1}}\quad\forall\quad 0\leq i \leq {n-k}$$.

Finally, these considerations lead us to the following formula for the {\it multigraded Hilbert series} of $I$:
$$\Hc(I_{k,n};\xb)=\frac{\sum_{i}(-1)^i {{i+k-1}\choose{k-1}}\cdot (\sum_{\ab\in[n,k+i]}\cdot\xb^{\ab})}{\prod_i(1-x_i)},$$
where $[n,k+i]$ denotes the set of $(k+i)$-subsets of $\{1,\dots,n\}$.

\begin{Example}
For $I_{3,5}$ we have
$$\Hc(I_{3,5};\xb)=\frac{(xyz+xyu+xyv+xzu+xzv+xuv+yzu+yzv+yuv+zuv)}{(1-x)(1-y)(1-z)(1-u)(1-v)}$$ $$-\frac{3(xyzu+xyzv+xyuv+xzuv+yzuv)}{(1-x)(1-y)(1-z)(1-u)(1-v)}+\frac{6(xyzuv)}{(1-x)(1-y)(1-z)(1-u)(1-v)},$$
the Betti numbers of $I_{3,5}$ are then: $\beta_0=10$, $\beta_1=15$ and $\beta_2=6$.
\end{Example}

\begin{Remark}
It is easy to see that $I_{k,n}$ are Mayer-Vietoris of type $B2$. One way to see it is to observe that one can easily construct a linear Mayer-Vietoris tree. Also, the fact that these ideals are Mayer-Vietoris of type $B2$ is a consequence of corollary \ref{tail-of-power-of-prime}, since $I_{k,n}$ is generated by a tail of $\langle x_1,\dots,x_n\rangle^k$ with respect to the ordering given in proposition \ref{J^k_is_MVB2}.
\end{Remark}
\subsubsection{Multistate $k$-out-of-$n$ systems}
\spanishsubsubsection{Sistemas $k$-entre-$n$ multiestado}

Let us consider now systems of $n$ components in which every component can reach a finite number of states $\{0,1,\dots,i\}$. Assume that such a system fails whenever the sum of the states of its components reaches a level $k$. We call such  a system \emph{$i$-multistate $k$-out-of-$n$} system. The reason for this terminology is that a $1$-multistate $k$-out-of-$n$ system is just the ordinary $k$-out-of-$n$ system studied above (considering the value $1$ indicates failure).

These systems are modelled by ideals of the form $J^k_{[n,i]}$, studied in section \ref{powers-primes}. Recall that the ideal $J^k_{[n,i]}$ is minimally generated by all monomials in $n$ variables of degree $k$ such that each variable has an exponent less than or equal $i$. As we have seen, all ideals $J_{[n,i]}^k$ are Mayer-Vietoris of type $B2$ and therefore their Mayer-Vietoris resolution is minimal. These ideals have even linear resolutions, see the details in section \ref{powers-primes}.

\subsubsection{Consecutive $k$-out-of-$n$ systems}
\spanishsubsubsection{Sistemas $k$-entre-$n$ consecutivos}

Consecutive $k$-out-of-$n$ systems \cite{K80} operate (resp. fail) whenever $k$ or more consecutive components operate (resp. fail). These systems can be modelled by the ideals
$$\bar I_{k,n}=\langle x^\mu \mbox{ s.th }x^\mu \mbox{ is a product of k consecutive variables} \rangle$$

for example, $\bar I_{3,5}=\langle xyz,yzu,zuv\rangle$ is the ideal corresponding to the consecutive $3$-out-of-$5$ system. Observe that the number of generators of such an ideal is $n-k+1$, therefore, the ideal corresponding to the \emph{consecutive} $k$-out-of-$n$ system is much smaller than the ideal of the corresponding ordinary $k$-out-of-$n$ system. But the nice combinatorial properties of $k$-out-of-$n$ systems are not present in their consecutive counterparts, which makes their study more difficult and interesting, in particular when $k$ is small with respect to $n$ ($k<\frac{n}{2}$).

In order to find the multigraded Betti numbers and Hilbert series of $\bar I_{k,n}$ we will use their Mayer-Vietoris trees. The explicit construction of this tree will give us the results we need. For more clearness, we will denote the monomials by their exponents in brackets, e.g. the monomial $x_1x_3x_6$ will be denoted by $[1,3,6]$, since we are dealing with squarefree monomials, this notation suffices.

We sort the generators of $\bar I_{k,n}$ using the lexicographic order. The construction of $MVT(\bar I_{k,n})$ goes as follows:
\begin{itemize}
\item The root node is just $\bar I_{k,n}$, which is minimally generated by $n-k+1$ monomials.
\item The right child of the root, i.e. $MVT(\bar I_{k,n})_3$ is $\bar I_{k,n-1}$, so we hang here the corresponding tree.
\item The left child of the root, $MVT(\bar I_{k,n})_2$, consists of the following $n-2k+1$ monomials:
\subitem $[j,\cdots (j+k-1),(n-k+1),\cdots ,n]$ for $1\leq j\leq n-2k$ which are the least common multiples of each of the first $n-2k$ generators of the root with the last one. These generators have $2k$ variables.
\subitem $[n-k, \cdots ,n]$ which is the lcm of the last two generators of $MVT(\bar I_{k,n})_1$ and divides $[n-k-j, \cdots, n]$ for $1\leq j\leq (k-1)$ and hence these last will not appear as minimal generators of this node. This generator has $k+1$ variables and since we are using lexicographic order, it will appear as the last generator in $MVT(\bar I_{k,n})_2$.
\item The following nodes to consider are $MVT(\bar I_{k,n})_4$ and $MVT(\bar I_{k,n})_5$, but only if $MVT(\bar I_{k,n})_2$ has more than one generator i.e. if $2k< n$, otherwise they are empty. If it is the case, then
\subitem $MVT(\bar I_{k,n})_4$ consists of $n-2k$ generators, namely the lcms of the first $n-2k$ generators of $MVT(\bar I_{k,n})_2$ with the last one. These have the form $[j,\cdots (j+k-1),(n-k),\cdots ,n]$ for $1\leq j\leq n-2k$ and hence, this node is exactly equal to $\bar I_{k,n-k-1}$ with each monomial in it multiplied by $[n-k,\cdots,n]$. Hence, we hang here  a tree isomorphic to $MVT(\bar I_{k,n-k-1})$.
\subitem $MVT(\bar I_{k,n})_5$ is completely analogous to $MVT(\bar I_{k,n})_5$ and hence equal to $\bar I_{k,n-k-1}$ but this time each monomial in it is multiplied by $[n-k+1,\cdots,n]$. Hence, we also hang here  a tree isomorphic to $MVT(\bar I_{k,n-k-1})$.
The trees we have hanged from the corresponding nodes are of the same form, except that they have less variables, i.i. they are of the form $MVT(\bar I_{k,j})$ with $j<n$. Eventually, we will have the situation in which $2k\geq n$ and in this case, the left child of the root has only one generator, namely $[j-k,\dots,j]$, and the right node is the consecutive $k$-out-of-$(j-1)$ tree, so we proceed in this manner until $j=k+1$.
\end{itemize}

\begin{Example}
Here is the tree corresponding to the consecutive 2-out-of-6 system:
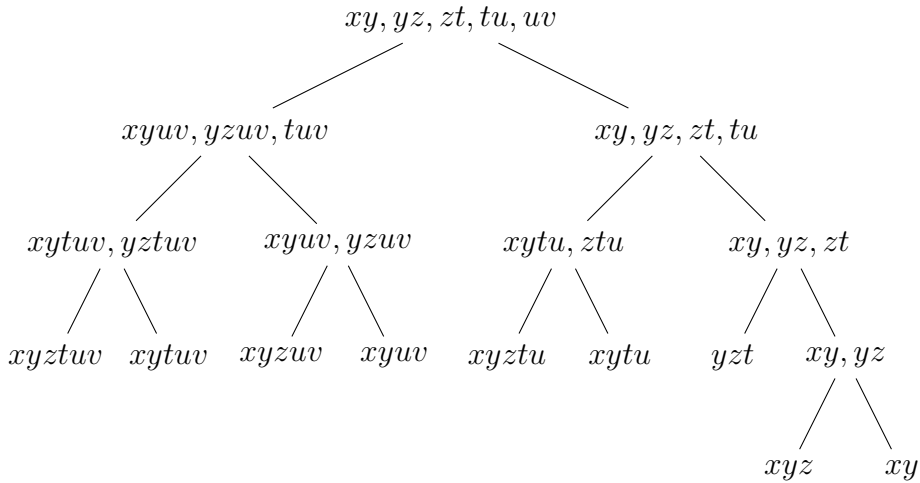
\begin{center}
\begin{tikzpicture}
\tikzstyle{level 1}=[sibling distance=6cm]
\tikzstyle{level 2}=[sibling distance=3cm]
\tikzstyle{level 3}=[sibling distance=1.5cm]
\node{$xy,yz,zt,tu,uv$}
child{node{$xyuv,yzuv,tuv$}
child{node{$xytuv,yztuv$}
child{node{$xyztuv$}}
child{node{$xytuv$}}}
child{node{$xyuv,yzuv$}
child{node{$xyzuv$}}
child{node{$xyuv$}}}}
child{node{$xy,yz,zt,tu$}
child{node{$xytu,ztu$}
child{node{$xyztu$}}
child{node{$xytu$}}}
child{node{$xy,yz,zt$}
child{node{$yzt$}}
child{node{$xy,yz$}
child{node{$xyz$}}
child{node{$xy$}}}}};
\end{tikzpicture}
\end{center}
\end{Example}

Taking into account the properties of the Mayer-Vietoris trees of these ideals, we see that we can read the multigraded Betti numbers directly from the ideal:
\begin{Proposition}\label{consecutive-k-n-tree}
The lexicographic Mayer-Vietoris tree corresponding to the consecutive $k$-out-of-$n$ system has no repeated multidegrees in the relevant nodes.
\end{Proposition}

\noindent {\bf Proof: } Assume we have $\bar I_{k,n}$ as the root of our tree, sorted with respect to lexicographic order, then the variable $n$ appears only in the left child of the root, and it will appear in every multidegree of every node in the tree hanging from this node (see the construction above). Thus, no multidegree of the tree hanging from the left child will appear in the tree hanging from the right child, and vice versa. If $2k\geq n$ then we are done, since the left node has just one generator, and the tree hanging from the right node is the one corresponding to the $k$-out-of-$(n-1)$ system.
If the left child of the root has more than one generator, then we look at its children, $MVT(\bar I_{k,n})_4$ and $MVT(\bar I_{k,n})_5$. The generators of the first one are not present in any node seen so far, and all of them contain the variables $(n-k),\dots,n$; moreover, every generator of the nodes of the tree hanging from it will have these variables. On the other hand, the variable $n-k$ does not appear in the generators of $MVT(\bar I_{k,n})_5$ hence, no multidegree of a generator in the tree hanging from it will appear in the tree hanging from $MVT(\bar I_{k,n})_4$ and vice versa. 
Finally, we have to see that no multidegree appearing in any relevant node of the tree hanging from $MVT(\bar I_{k,n})_5$ is in $MVT(\bar I_{k,n})_2$. We know that $MVT(\bar I_{k,n})_5$ is generated by the generators of $MVT(\bar I_{k,n})_2$ except the last one. Now, every generator of every node in the tree hanging from $MVT(\bar I_{k,n})_5$ will have at least $2k+1$ different variables, $k$ of which will be $(n-k+1),\dots,n$ (see the construction of the tree), and on the other hand, the generators in $MVT(\bar I_{k,n})_2$ have at most $2k$ different variables. $\square$

With this proposition we have that collecting all the generators of the relevant nodes in $MVT(\bar I_{k,n})$ we have the multigraded Betti numbers of $\bar I_{k,n}$ in this case, since no generator in the relevant nodes is repeated, we have that the Betti number at each multidegree is $1$, every multidegree appears only once in the minimal resolution of the ideal.
The description of the tree and its recursive construction give us also means to count how many multidegrees appear in each dimension (i.e. the Betti numbers) and which multidegrees are present. A thorough description of this process would be tedious, but it is not difficult to obtain a complete list of the multidegrees of the Betti numbers, and hence, of the Hilbert series. however, here we only give an idea of the procedure; an algorithm has been implemented by the author to perform this listing.
The main lines of the construction of this list of multidegrees are the following:
\begin{itemize}
\item In dimension $0$ collect all the generators of $\bar I_{k,n}$.
\item In dimension $1$ collect all the multidegrees of the form $[j,\dots,j+k]$ for $1\leq j\leq (n-k)$. \footnote{Note that in the case $2k\geq n$ these are the only ones we have to add, and the corresponding formula is equivalent to the one appearing in \cite{D03}} Moreover, for $k\leq j \leq (n-k)$, add the multidegrees $[1,\cdots,k,(j+1),\cdots,(j+k)],\dots,[(j-k),\dots,(j-1),(j+1),\cdots,(j+k)]$.
\item For every dimension $l$ add the corresponding multidegrees that appear in $\bar I_{k,j-k-1}$ in dimension $(l-2)\geq 0$ multiplied by $[(j-k),\dots,j]$ and the multidegrees that appear in $\bar I_{k,j-k-1}$ in dimension $(l-1)\geq 0$ multiplied by $[(j-k+1),\dots,j]$ for all $(2k+1)\leq j\leq n$
\end{itemize}

From the constructions we have seen we can obtain a recursive relation for the Betti numbers of the \emph{consecutive $k$-out-of-$n$} systems. The construction of the tree in the proof of proposition \ref{consecutive-k-n-tree} gives us the following relations:
$$
\begin{array}{rll}
\beta_0(\bar{I}_{k,n})&=n-k+1&\\
\\
\beta_1(\bar{I}_{k,n})&=n-k& \mbox{if }k\geq n/2\\
\beta_1(\bar{I}_{k,n})&=n-2k+1+\beta_1(\bar{I}_{k,n-1})& \mbox{if }k < n/2\\
\\
\beta_i(\bar{I}_{k,n})&=0&\mbox{if }k\geq n/2,\quad i>1\\
\beta_i(\bar{I}_{k,n})&=\beta_{i-2}(\bar{I}_{k,n-k-1})+\beta_{i-1}(\bar{I}_{k,n-k-1})+\beta_i(\bar{I}_{k,n-1})&\mbox{if }k < n/2,\quad i>1\\
\end{array}
$$     							


\fancyhead[ER]{\itshape Conclusions and further work} 
\fancyhead[OL]{\itshape Conclusions and further work}
		
\chapter{Conclusions}
\spanishchapter{Conclusiones}

With a particular focus on explicit computations and applications of the Koszul homology and Betti numbers of monomial ideals, the main goals of this thesis have been the following:
\begin{itemize}
\item Analyze the Koszul homology of monomial ideals and apply it to describe the structure of monomial ideals.
\item Describe algorithms to perform efficient computations of the homological invariants of monomial ideals, in particular Betti numbers, free resolutions, Koszul homology and Hilbert series.
\item Apply the theory of monomial ideals to problems inside and outside mathematics, in particular making use of the homological invariants of these ideals.
\end{itemize}

The thesis started with a description of the concept of Koszul homology and the properties of monomial ideals, together with some combinatorial techniques that can be applied to the study of these ideals. In the first chapter it was shown that the combinatorial nature of monomial ideals allows us to study their algebraic and homological structure with simple tools. In particular, several tools to study resolutions and Koszul homology were shown, introducing the Mayer-Vietoris sequences associated with a monomial ideal as a new tool.

On the second chapter we have analyzed the algebraic and homological structure of monomial ideals using Koszul homology. We begun with the problem of obtaining a minimal free resolution of a monomial ideal from the knowledge of its Koszul homology. This problem was already solved by Aramova and Herzog \cite{AH95,AH96}, but the explicit constructions were still involved even in some simple cases. In \cite{RS06,RSS06} the constructions of Aramova and Herzog were transformed into algorithms with the help of \emph{effective homology}, therefore, this problem was solved both from a theoretic and from an algorithmic point of view. There is still a lack of simple explicit methods that allow us to read the differentials of a minimal resolution of a monomial ideal from a set of generators of its Koszul homology.

We have also treated in this chapter the problem of finding a combinatorial decomposition of the factor ring $R/I$ where $I\subset R=\kb[x_1,\dots,x_n]$ is a monomial ideal. In particular a procedure has been given to obtain a Stanley decomposition from the knowledge of the Koszul homology of $I$. In this case, it is enough to know the multidegrees of the Koszul generators of $I$. We have seen that the $(n-1)$-st homology generators play a very important role in this problem. This is due to the special relation of these homology generators with the boundary of a monomial ideal, in particular with the so called \emph{maximal corners}, which define the Stanley decomposition in the artinian case. The non-artinian case was studied in terms of the artinian one.

Irredundant irreducible and primary decompositions are also treated in this chapter. Several procedures are given to compute an irredundant irreducible decomposition from the Koszul homology of a monomial ideal or from the Koszul homology of its artinian closure. First of all we treat the artinian case, using $(n-1)$-st Koszul homology and the concept of maximal corners. This approach is different to other ones \cite{R07,MS04} which are also based in concepts equivalent to closed corners. While other approaches are based on the artinian closure for the nonzero dimensional case, we developed another procedure which exclusively uses the Koszul homology of the ideal itself. This is an alternative approach to the existing ones. With respect to irredundant primary decomposition, a simple procedure is introduced to obtain a particular irredundant primary decomposition of a monomial ideal given the $(n-1)$-st Koszul homology of its artinian closure.

The third chapter is focused on computations and introduces one of the central contributions of the thesis, namely Mayer-Vietoris trees. Given a monomial ideal we have introduced the Mayer-Vietoris trees associated to that ideal and the properties of them. Mayer-Vietoris trees are a computational tool based on the Mayer-Vietoris sequences introduced on chapter one. They allow us to perform homological computations on the monomial ideals they are associated to. Every Mayer-Vietoris tree produces a multigraded resolution of the corresponding ideal, and hence an expression of its multigraded Hilbert series. Moreover, they give a subset and a superset of the multidegrees of the Koszul generators of the ideal in each homological degree. We also give conditions under which these subsets and/or supersets provide the actual multigraded Betti numbers of the ideal, and therefore the Mayer-Vietoris trees give its minimal free resolution. The families of ideals for which we can obtain the multigraded Betti numbers directly from the Mayer-Vietoris trees are called \emph{Mayer-Vietoris ideals} and they are divided in several types. Some important families of Mayer-Vietoris ideals are also given.

Mayer-Vietoris trees are not only a computational tool but also a mean to analyze monomial ideals, and their homological structure, including free resolutions, Betti numbers and Hilbert series. In their paper \cite{PS07}, I. Peeva and M. Stillmann propose several conjectures and open problems around syzygies and Hilbert series. In relation with monomial ideals the main goal is stated as 

\noindent{\bf Problem 3.9.1 \cite{PS07}: }Introduce new constructions and ideas on monomial resolutions.

In this context, Mayer-Vietoris trees are simple but computationally efficient tool to analyze monomial resolutions. In particular, they provide a closer look to the divisibility relations inside the $lcm$-lattice of the ideal with a view toward the Koszul homology of it, and therefore toward the free resolutions of the ideal and its Betti numbers. In this sense, several examples of such analysis are provided in this section and in the fourth chapter.

Also in the third chapter some algorithmic issues are presented. A description of a basic algorithm that computes Mayer-Vietoris trees is given together with several implementation details. An implementation of this algorithm has been made using the \cpp library $\cocoalib$. This library, under development by the $\cocoa$ team \cite{cocoa}, brings together the efficiency of \cpp and capabilities to make computations in commutative algebra, and to communicate with other computer algebra systems. Several experiments with this implementation have been made to demonstrate the performance of different strategies used to build Mayer-Vietoris trees, the timings when compared with algorithms to compute resolutions, multigraded Hilbert series and Betti numbers in other systems, and several improvements on the basic algorithm. These experiments show that the computation of the Mayer-Vietoris tree is an efficient approach to compute resolutions and Hilbert series. In particular when the number of variables grow or when the minimal free resolution is big, Mayer-Vietoris trees can give us the information that is needed for most applications, in particular the (multigraded) Betti numbers or at least bounds of them. They are also an efficient alternative to compute multigraded Hilbert series.

The fourth chapter is devoted to applications. We give several applications of the techniques, procedures and algorithms presented in the previous chapters. We apply them to different classes of monomial ideals, to other area of mathematics, namely the formal theory of differential equations, and finally, outside mathematics, to reliability theory.

First of all, several types of ideals are studied. We analyze Borel-fixed, stable and segment ideals with the help of Mayer-Vietoris trees, and demonstrate that they are a good tool to make statements about these types of ideals. We give alternative proofs to several important results about the homological structure of these ideals. The same is applied to Scarf ideals, in particular to generic ideals. All these kind of ideals have some important theoretical features. Second, three types of monomial ideals with applications in other areas of mathematics are studied. Valla ideals \cite{V04} have applications in algebraic geometry, they are proved to be Mayer-Vietoris of type $B2$ and using Mayer-Vietoris trees we give explicit formulas for their (graded) Betti numbers and irredundant irreducible decompositions of them. Ferrers ideals \cite{CN06,CN07} have applications among others in graph theory; they are Mayer-Vietoris of type $A$, and from their Mayer-Vietoris trees we obtain explicit formulas for their Betti numbers and irredundant primary and irreducible decompositions of them, as well as other invariants. Finally, we study quasi-stable ideals, which have strong relations to involutive (in particular Pommaret) bases and also appear in the formal theory of differential systems. Here we use the Koszul homology of them to make a procedure that completes such an ideal to a Pommaret basis and can also decide about the $\delta$-regularity of the given coordinates, giving therefore an alternative approach to those existing in the literature \cite{S02b,HS02}.

The applications to the \emph{formal theory of differential systems} are based on the duality betwen Spencer cohomology and Koszul homology. The framework for this application is a geometric approach to PDEs which makes an algebraic analysis of them possible. The role of Koszul homology in this context is exploited in relation with involution and formal integrability. Knowing the Koszul homology and/or the Betti numbers of the symbol of a differential system gives us a way to detect involution and quasi-regularity of a given differential system. This is based on the correspondence between the degree of involution and the algebraic concept of Castelnuovo-Mumford regularity \cite{M03,S07}. Another application in this area is the formulation of well posed initial value problems. These are strongly related to Stanley decompositions of the corresponding factor ring, and we apply the procedure seen in chapter two to obtain such decompositions from the Koszul homology of the corresponding ideal.

Finally we use monomial ideals in reliability theory. A connection between monomial ideals and coherent systems was established in \cite{GW04}. In particular, certain expressions of the Hilbert series of a monomial ideal provide bounds for the reliability of the corresponding system. We have used several techniques to obtain sharp bounds or even formulas to compute the reliability of several important types of systems. A very natural type of networks, namely \emph{series-parallel} networks is treated using Mayer-Vietoris trees. The trees associated to these networks have been proved to be Mayer-Vietoris of type $A$. Therefore, the Mayer-Vietoris resolution of these ideals is minimal and we can obtain sharp bounds for the relability of such networks. Moreover, we provide explicit procedures to compute the reliability of these networks without computing resolutions, and to compute the Betti numbers of them. We also give some subtypes for which explicit formulas for the Betti numbers of the associated ideals are provided. Some non-network systems are also studied. In particular, several types of $k$-out-of-$n$ systems, namely the usual one, the consecutive and the multistate one. In these cases, either simplicial Koszul homology or Mayer-Vietoris trees are used to produce formulas for the Betti numbers of the associated ideals, and therefore for the reliability of the system. In some cases we give simple proofs of known results and in other cases new results are obtained using our techniques.

The directions explored in this thesis provide us several open ways in which further work has to be done, both from a theoretical and  an applied viewpoint, and also in a computational framework.

The Koszul homology of monomial ideals has a combinatorial nature that has been studied with some simple techniques coming from algebraic topology and basic homological algebra. This direction has to be studied in bigger depth, using more techniques that can give rise to further results. Two obvious directions to follow in this context are \emph{relative homology} and \emph{spectral sequences}, from which new results about combinatorial Koszul homology have to be expected. Also, a study of the topics treated in the thesis from the cohomological point of view using $Ext$ and $Hom$ functors can provide further information.

In the study of monomial ideals, several directions have to be followed. A big area in the theory of monomial ideals is referred to \emph{edge ideals} and \emph{facet ideals}, i.e. ideals associated to graphs and simplicial complexes, see for example \cite{V01}, \cite{FV07,HV07} and \cite{F02}. Some of the results of this thesis can be read in this context, but a rigorous study of the Koszul homology of edge and facet ideals and of the application of Mayer-Vietoris trees to these ideals, has to be made. Also, there are several types of monomial ideals which can be studied with our techniques that have importance from a theoretic point of view, such as \emph{$LPP$-ideals} \cite{FR07}, \emph{pretty clean} monomial ideals \cite{JZ07} etc. Also some concepts like polarization \cite{F05} or splitting \cite{EK90} of monomial ideals might be treated with these techniques, expecting new analysis of them.

With respect to applications, monomial ideals are an ubiquitous concept. The applications presented in the thesis can be completed with a more deep study. The relations between Koszul homology and Pommaret bases, Mayer-Vietoris trees of quasi-stable ideals and different coherent systems in reliability such as multidimensional consecutive $k$-out-of-$n$ systems and others, are the first directions to follow. Also different applications of monomial ideals can be studied using our techniques, in particular those related to the sciences of life, such as the relation of monomial ideals and phylogenetic trees, see for example \cite{GLRRV06}.

Finally, in the algorithmic context, a complete implementation of the procedures to obtain minimal free resolutions of monomial ideals using Mayer-Vietoris trees has to be made. Also, improve and complete the implemented algorithms and add new procedures to detect homology in repeated relevant multidegrees has to be done. This will result in a new alternative algorithm to compute homological invariants of monomial ideals. Together with the effective homology techniques seen in chapter two, this can be applied to more general polynomial ideals.

\appendix
\renewcommand{\chaptermark}[1]{\markboth{Appendix \thechapter \quad #1}{}} 
\fancyhead[ER]{\itshape \leftmark} 
\fancyhead[OL]{\rightmark}

\chapter{Algebra}\label{apal}
\spanishchapter{Algebra}

\section{Homological algebra}
\spanishsection{\'Algebra homol\'ogica}
The main classical reference we have used for this subject is \cite{M95}, also \cite{GM03} is a good reference; books in commutative algebra like \cite{E95} use to include sections on homological methods and some introduction to the main points of homological algebra.
\subsubsection*{Complexes and free resolutions.}
\spanishsubsubsection{Complejos y resoluciones libres}
One of the objectives of homological algebra is the study of complexes and their homology
\begin{Definition}
A {\rm complex} of modules over a ring $R$ is a sequence of modules and morphisms
$$\MM: \dots\rightarrow M_{i+1}\stackrel{\delta_{i+1}}{\rightarrow}M_{i}\stackrel{\delta_{i}}{\rightarrow}M_{i-1}\rightarrow\dots$$
such that  $\delta_i\delta_{i+1}=0$ for each $i\in \ZZ$, i.e. ${\rm im}\delta_{i+1}\subseteq {\rm ker}\delta_{i}$. The {\rm $i$-th homology group} of $\MM$ is defined as
$$H_i(\MM):={\rm ker}\delta_i/{\rm im}\delta_{i+1}$$
If $H_i(\MM)=0$, i.e. ${\rm im}\delta_{i+1}={\rm ker}\delta_{i}$ we say that $\MM$ is {\rm exact} at $i$. If $\MM$ is exact at every $M_i$ we say $\MM$ is {\rm acyclic}.
\end{Definition}

\begin{Definition}
Given an $R$-module $\Mc$, a {\rm resolution} of $\Mc$ is a complex $\MM: \Res M k \delta$ that is exact everywhere except at dimension $0$ in which we have $H_0(\MM)\simeq\Mc$.
\end{Definition}

If the modules in a resolution are flat, projective, free, etc. we say that the resolution is respectively flat, projective, free, etc. Free resolutions have particular interest for us:

\begin{Definition}
Given a finitely generated graded $R$-module $\Mc$, a {\rm free resolution} of $\Mc$ is a long exact sequence $\PP: \Res \Pc k \delta$ where the $\Pc_i$ are finitely generated free modules and the maps $\delta_i$ have degree $0$. We say that a free resolution of $\Mc$ is {\rm minimal} if it is minimal as a complex, i.e. if for each $i$ the image of $\delta_i$ is contained in $(x_1 \cdots x_n)\Pc_{i-1}$, or, more informally, if the differentials $\delta_i$ are represented by matrices with entries in $(x_1\cdots x_n)$ \cite{E95,E04}. This is equivalent to say that for each $i$, the map $\delta_i$ takes a basis of $P_i$ to a minimal set of generators of the image of $\delta_i$. The proof of this equivalence is a consequence of Nakayama's Lemma and can be seen for example in \cite{E04}
\end{Definition}

We observe here that minimal resolutions are unique, up to isomorphism; more exactly (Theorem 20.2 in \cite{E95}): If $\PP$ is a minimal free resolution of $\Mc$, then any free resolution of $\Mc$ is isomorphic to the direct sum of $\PP$ and a trivial complex. 

\begin{Definition} The meaning of "trivial complex" and "isomorphic" is the following:

\begin{itemize}
	\item{A {\rm trivial complex} is a direct sum of complexes of the form $0\rightarrow R\stackrel{1}{\rightarrow}R\rightarrow 0$, which have no homology.}
	\item {Let $\PP :\Res \Pc k \delta$ and $\PP': \Res {\Pc'} k {\delta'}$ be two graded resolutions. An {\rm isomorphism of graded resolutions} is a sequence of graded isomorphisms $\varphi_i:\Pc\rightarrow\Pc'_i$ of degree zero, such that $\delta'_0\cdot\varphi_0=\delta_0$ and $\forall i\geq 1$ the following diagram 	commutes, i.e. $\varphi_{i-1}\cdot\delta_i=\delta'_i\cdot\varphi_i$ :}
\end{itemize}
\end {Definition}

A consequence of the uniqueness of minimal free resolutions is that the number of generators of each degree required for the free modules $\Pc_i$ of the minimal free resolutions depends only on the module $\Mc$. This is a consequence of the definition and basic properties of the ${Tor}$ functor:

\begin{Definition}
If $M$ and $N$ are $R$-modules and $\PP :\Res \Pc k\delta$ is a projective resolution of $N$ as an $R$-module, then $Tor^R_i(M,N)$ is the $i$-th homology of the complex $M\otimes\PP$, i.e. $ker(M\otimes \delta_i)/im(M\otimes \delta_{i+1})$. This homology is (up to isomorphism) independent of the chosen resolution. Moreover, it can also be computed starting with a resolution $\PP'$ of $M$ and computing the homology of $\PP'\otimes N$.
\end{Definition}
In the case of $\PP$ (equiv. $\PP'$) being a minimal free resolution of the $R$-module $\kb$ (respect. $\Mc$) then $dim(Tor^R_i(M,\kb)_j)$ (see definition of $Tor$ below) is exactly the number of degree $j$ generators of $P_i$, see \cite{E04}, Proposition 1.7. These important numbers have a special name:

\begin{Definition}
If $\PP: \Exseq \Pc k \delta$ is a minimal free resolution of $\Mc$, then any minimal set of homogeneous generators of $\Pc_i$ contains exactly ${\rm dim}_\kb Tor_i^R(\kb, \Mc)_j$ generators of degree $j$. This is denoted $\beta_{i,j}(\Mc)$, and we call the  $\beta_{i,j}(\Mc) \forall i,j$, the {\rm graded Betti numbers} of $\Mc$. If we forget about the grading, then the $\beta_i(\Mc)=\sum_j\beta_{i,j}(\Mc)$ are called the {\rm Betti numbers} of $\Mc$.
\end{Definition}

\begin{Definition}
Given a finite free resolution $\PP$, we say that the sum of the ranks of its modules is the {\rm size of $\PP$} and will be denoted $Size(\PP)$.
\end{Definition}
\subsubsection*{Derived functors}
\spanishsubsubsection{Funtores derivados}

The $Tor$ functor is a particular case of a standard tool in homological algebra: Derived Functors, here are the basic definitions ond the two main examples of derived functors: $Tor$ and $Ext$.

\begin{Definition}
Let $\Mc$ be an $R$-module, and let $\PP: \Exseq \Pc k \delta$ be a projective resolution of $\Mc$. If $F$ is a right-exact functor from the category of $R$-modules to the category of abelian groups, then applying $F$ to $\PP$ we obtain a complex $F\PP:0\rightarrow FP_{k}\stackrel{\delta_k}{\rightarrow}FP_{k-1}\stackrel{\delta_{k-1}}{\rightarrow}\cdots FP_1\stackrel{\delta_1}{\rightarrow}FP_0 \rightarrow 0$ which is no longer exact.Then we define the $i$-th {\it left derived functor} to be $L_iF(\Mc)=H_i(F\PP)$.

Similarly, one defines {\it right derived functors} using injective co-resolutions and left-exact functors.
\end{Definition}

Two of the most frequently used derived functors are {\it Tor} and {\it Ext}:

\begin{Definition}
Let $M$ be an $R$-module. The $n$-th left derived functor of the right-exact functor $M\otimes -$ is denoted by $Tor^R_n(M,-)$.

The $n$-th right derived functor of the left-exact functor $Hom(-,M)$ is denoted by $Ext_R^n(-,M)$.
\end{Definition}

\subsubsection*{Spectral sequences}
\spanishsubsubsection{Sucesiones espectrales}

\begin{Definition}
An $R$-\emph{bimodule} (bigraded module) is  a collection $\EE$ of $R$-modules $E_{s,t}$. A differential in a bimodule $d:\EE\rightarrow\EE$ of bidegree $(a,b)$, is a collection of homomorphisms $d_{s,t}:E_{s,t}\rightarrow E_{s+a,t+b}$ such that $d^2=0$. This differential defines the homology bimodule of $\EE$, $H_{s,t}(\EE)=ker(d_{s,t})/Im(d_{s-a,t-b})$.
\end{Definition}
\begin{Definition}
An \emph{espectral sequence} $\EE$ is a sequence $\{\EE^r,d^r\}$ for $r$ greater or equal to a given $k$, such that
\begin{itemize}
\item $\EE^r$ is a bimodule and $d^r$ is a differential of bidegree $(-r,r-1)$ on $E^r$
\item For $r\geq k$ there is an isomorphism $H(\EE^r)\simeq\EE^{r+1}$
\end{itemize}
$\EE^r$ is called the $r$-th \emph{page} of the spectral sequence $\EE$.
\end{Definition}
\begin{Definition}
A spectral sequence $\EE$ is \emph{convergent} if for every $s,t\in \ZZ$, there exists an integer $r=r(s,t)$ such that we have that $d^r_{s,t}=0=d^r_{s+r,t-r+1}$ for all $r\geq r_{s,t}$. If $\EE$ is convergent, we denote $\EE^\infty_{s,t}:=\rm{lim}_{r\rightarrow\infty}\EE^{r}_{s,t}$
\end{Definition}
\begin{Remark}
Very often, spectral sequences lie in some quadrant of the $(s,t)$-plane, in the sense that the differential is null outside that quadrant. We speak then of \emph{first quadrant} spectral sequences, \emph{second quadrant} spectral sequences, and so on. Note that for example a first quadrant spectral sequence is always convergent, which doesn't hold for second quadrant spectral sequences, in general.

Famous examples of spectral sequences are the \emph{Serre} spectral sequence or the \emph{Eilenberg-Moore} spectral sequence.
\end{Remark}
We describe shortly the notion of spectral sequence of a bicomplex, and its application to the computation of $Tor$ modules:
\begin{Definition}
A \emph{bicomplex} $\Cc$ is a bimodule with two differentials, $d'$ and $d''$ of bidegrees $(-1,0)$ and $(0,-1)$ respectively, such that $d'd''+d''d'=0$. First quadrant bicomplexes, second quadrant bicomplexes etc. can be defined in the natural way.

The \emph{totalization} of a bicomplex is a chain complex $(\Tc,d)$ where $T_n=\oplus_{s+t=n}C_{s,t}$ and the differential $d(c)$ with $c\in C_{s,t}\subset T_{s+t}$ is $d(c)=d'(c)\oplus d''(c)$. 
\end{Definition}
Associated with bicomplexes are two natural spectral sequences, coresponding to the vertical and horizontal filtrations of its totalization. Under certain circumstances, these spectral sequences are convergent.
\begin{Theorem}
Let $(C_{s,t},d',d'')$ be a bicomplex. A spectral sequence $_{ver}\EE$ can be defined with $_{ver}\EE^0_{s,t}=C_{s,t}$, $d^0_{s,t}=d''_{s,t}$ and $\EE^1_{s,t}=H''_{s,t}(\Cc)$ is the ``vertical'' homology group of the $s$-column at index $t$. Similarly, a spectral sequence $_{hor}\EE$ can be defined with $_{hor}\EE^0_{s,t}=C_{s,t}$, $d^0_{s,t}=d''_{s,t}$ and $\EE^1_{s,t}=H'_{s,t}(\Cc)$ is the ``horizontal'' homology group of the $t$-row at index $s$. These spectral sequences converge to the homology of the totalization: $E^r_{s,t}\rightarrow H_{s+t}(\Tc)$ under certain circumstances, for example, $_{hor}\EE$ converges if $C_{s,t}=0$ for all $s<0$ or for all $t>0$, and $_{ver}\EE$ converges if $C_{s,t}=0$ for all $s>0$ or for all $t<0$. Both sequences converge if $(C_{s,t},d',d'')$ is a first quadrant or a third quadrant bicomplex.
\end{Theorem}

A nice application of this theorem is the proof that the $Tor$ functor is balanced, i.e. that $Tor^R_\bullet(M,N)$ can be computed using either a resolution of $M$ or of $N$. A sketch of the proof (see details in \cite[Appendix 3]{E95} or in \cite[exercise 1.12]{MS04}) is the following: 
Let 
$$\Pc :\Res P k \phi$$ be a resolution of $M$ and
$$\Qc :\Res Q k \psi$$ be a resolution of $N$.
Using the vertical or horizontal spectral sequences of the bicomplex $(P_s\otimes Q_t,\phi,\psi)$ one proofs that 
$$H(\Pc\otimes N)\simeq H(tot(\Pc\otimes\Qc))\simeq H(M\otimes\Qc)$$
And applying the definition of the $Tor^R_\bullet(M,N)$ modules, we have the result.

\section{Commutative algebra}
\spanishsection{\'Algebra conmutativa}

Reference books for this part are \cite{Sh00} and \cite{E95} which contains almost everything, or \cite{S03} that has a good introduction to the basic concepts with a focus on computation. From a computational point of view \cite{KR00,KR05},  \cite{EGSS02} and \cite{GP02} are good books. Computer algebra systems specialized in commutative algebra are \cocoa \cite{cocoa}, \macaulay \cite{M2} and \singular \cite{Singular}.

The main objects in commutative algebra are modules over polynomial rings. They have been studied from different points of view, being a prominent one its relation with algebraic geometry. Another important point of view in the most recent treatments of commutative algebra is its computational side. In particular, the development of Gr\"obner basis theory has given a big impulse to this area of research and its applications. 

Here we give some definitions of the objects from commutative algebra that are used in the text. Their presentation is a bit random and is just a sequence of the necesary definitions.

\subsubsection*{Grading, multigrading and term orders.}
\spanishsubsubsection{Grados, multigrados y \'ordenes de t\'erminos}

Given a monoid $M$, an $M$-graded ring $R$ is a ring with a direct sum decomposition
$$R = \bigoplus_{i\in M}A_i$$
such that $A_i\cdot A_j \subseteq A_{i \cdot j}$. We say then that $M$ is the index set of the grading. The elements of $A_d$ are known as homogeneous elements of degree $d$. In our considerations, usually the index sets are $\ZZ$ or $\NN$. 

The polynomial ring $R=\kb[x_1,\dots,x_n]$ has a natural $\NN$-grading, considering the degree of the monomials, i.e. $deg(x_i)=1$ for all $i$. Also, $R$ is $\NN^n$-graded or \emph{multigraded} by setting $deg(x_i)=(0,\dots,\stackrel{i}{1},\dots,0)$ i.e. the $i$-th standard vector in $\NN^n$. 

An $R$-module $\Mc$ is $\NN$-graded if it is of the form $\Mc = \bigoplus_{i\in \NN}\Mc_i$ and $R_i\Mc_j\subseteq \Mc_{i+j}$. Similarly multigraded modules are defined. In particular, monomial ideals (i.e. ideals generated by monomials)  are multigraded as ideals in $R$ and as $R$-modules. Monomial ideals are multihomogeneous, therefore they have multigraded resolutions and Betti numbers (see later).

\begin{Definition}
We say that a total ordering $ \preceq$ on $\NN^n$ is a \emph{term order} or \emph{monomial order} if $\mu\preceq \nu$ implies $\mu + \rho \preceq \nu+\rho$ for all $\mu,\nu,\rho\in \NN^n$ and $\preceq$ is a well ordering, i.e. every nonempty subset of $\NN^n$ has a smallest element under $\preceq$. 
\end{Definition}
There are many different term orders, here are the most used:
\begin{itemize}
\item \emph{Lexicographic order} (\lex): Let $\mu=(\mu_1,\dots,\mu_n)$ and $\nu=(\nu_1,\dots,\nu_n)$ two elements in $\NN^n$. Then $\mu>_{\lex}\nu$ if the leftmost nonzero entry of $\mu-\nu\in \NN^n$ is positive.
\item \emph{Degree-lexicographic order} (\deglex): Is the homogeneous version of $\lex$, $\mu>_{\deglex}\nu$ if $deg(\mu)>deg(\nu)$ or if $deg(\mu)=deg(\nu)$ and $\mu>_{\lex}\nu$.
\item \emph{Degree-reverse-lexicographic order} (\degrevlex): In this ordering $\mu>_{\degrevlex}\nu$ if $deg(\mu)>deg(\nu)$ or if $deg(\mu)=deg(\nu)$ and the rightmost nonzero entry of $\mu-\nu\in \NN^n$ is negative.
\end{itemize}
\begin{Remark}
The last two term orders are called homogeneous or degree-compatible. Among all term orders, probably the most important one is $\degrevlex$, note that if we drop the degree-compatible condition in $\degrevlex$, we do not obtain a term order.

A very interesting approach to term orders is that of Robbiano and Kreuzer \cite{KR00}, which use matrices to represent term orders, also called weight orders. Every term order can be represented this way, in particular it is not difficult to find matrices representing $\lex,\deglex$ and $\degrevlex$.
\end{Remark}
\begin{Definition}
If $f$ is a polynomial in $I$, and $>$ a term order, we say that the initial or leading monomial $\init(f)$ is the largest monomial of $f$ with respect to $>$. The leading monomial together with its coefficient is the leading term of the polynomial $\lt(f)$.

Given an ideal $I$ in $R$ and a term order $>$ we define the \emph{initial ideal} of $I$ with respect to $>$ as
$$\init_>(I)=\langle\init(f)\vert f\in I\rangle$$
\end{Definition}

\subsubsection*{Gr\"obner bases}
\spanishsubsubsection{Bases de Gr\"obner}

Many questions in commutative algebra have a constructive flavour. Many problems look for constructive and/or algorithmic solutions and there exist indeed many algorithmical answers to central problems in commutative algebra. At the center of these questions stands the notion of \emph{Gr\"obner basis}, a milestone in the development of computational commutative algebra. Gr\"obner bases are essentialy special sets of generators of ideals in the polynomial ring. They are defined with respect to some term order:
\begin{Definition}
Let $I=\langle g_1,\dots,g_r\rangle$ an ideal in $R$ and $>$ a term order. A set of generators $\{g_1,\dots,g_r\}$ is a \emph{Gr\"obner basis} if $\init(I)=\langle\init(g_1),\dots,\init(g_r)\rangle$
\end{Definition}

Every ideal in $R$ posseses a Gr\"obner basis for every term order, and indeed Gr\"obner bases are bases of the correspondent ideals. In general, an ideal might have different Gr\"obner bases for a fixed term order. However, fixed a term order, any ideal has a unique \emph{reduced} Gr\"obner basis. We say that a Gr\"obner basis $\{g_1,\dots,g_r\}$ is reduced if the coefficient of the initial or leading term of $g_i$ is $1$ for every $i$ and no monomial of the polynomials $g_j,\, j\neq i$ is divisible by $\init(g_i)$.

Gr\"obner bases were introduced by B. Buchberger in his PhD thesis \cite{B06}, and they are very important not only because of their important theoretical properties, but also because of the fact that they can be \emph{algorithmically} computed, i.e., there is a finite algorithm that computes a Gr\"obner basis for a given ideal and a fixed term order, it is known as Buchberger's algorithm \cite{B06}; it is at the core of many algorithms in computational commutative algebra and algebraic geometry. One of the main features of Gr\"obner basis is that they can be used to perform division in the polynomial ring in several variables, and compute ideal membership. This leads to many applications in constructive module thoery. Another set of applications of these bases are in elimination theory, solution of polynomial systems, etc. Good reviews of the applications of Gr\"obner bases are \cite{E95,CLO96,BW98}.

\subsubsection*{Syzygies}
\spanishsubsubsection{Sicigias}

A key role in Buchberger's algorithm for the construction of Gr\"obner bases is played by $S$-polynomials, which are special combinations of polynomials designed to cancell leading terms. Further improvements of this algorithm lead to study cancellation in general. This takes us to the notion of \emph{syzygy}:
\begin{Definition}
Take a set of polynomials $\{f_1,\dots,f_r\}$; a syzygy in the leading terms $\lt(f_1),\dots,\lt(f_r)$ is  a tuple of polynomials $(s_1,\dots,s_r)$ such that
$$\sum_{i=1}^r s_i\cdot \lt(f_i)=0$$
Given a generating set $G_I$ of an ideal $I$, the syzygies of $G_I$ constitute a module, called the \emph{syzygy module} of $I$. This can be generalized from ideals to $R$-modules.
\end{Definition}

Syzygies are a main tool to describe modules of the polynomial ring via generators and relations, and via resolutions.
If we have a set of elements $\{m_i\vert i\in \Ic\}$ that generate $\Mc$ as $R$-module, we define a map from the free module $F_0=R^\Ic$ that maps the $i-th$ generator to $m_i$. Let $\Mc_1$ be the kernel of this map; it is also finitely generated as  $R$-module, and we can define a new map from a free module $F_1$ to $F_0$ with image $\Mc_1$. In this way we can construct a graded free resolution of $\Mc$. The elements of $\Mc_i$ are called the {\it $i$-th syzygies} of $\Mc$, and $\Mc_i$ the {\it $i$-th syzygy module of $\Mc$}. Graded and multigraded versions of this construction are defined in a similar way. Taking homogeneous elements in the generating sets of $\Mc$ and the syzygy modules, we construct graded resolutions and these are in strong relation with Hilbert functions, from which the structure of the module and many important invariants can be read off, see for example \cite{P07}. One main theorem is syzygy theory is Hilbert's syzygy theorem, which states that every finitely generated $R$-module has a finite graded free resolution of length at most the number of variables of the polynomial ring $R$.

Syzygies are not only used in Buchberger's algorithm or to build free resolutions of modules, they have a broad use in geometry \cite{E04}.

An important invariant of homological nature associated to modules over the polynomial ring is the following:

\begin{Definition}
Let $\Mc$ be a finitely generated graded $R$-module and let 
$$\cdots\rightarrow\Pc_i\rightarrow\cdots\rightarrow\Pc_0\rightarrow 0$$
be its minimal graded free resolution. Let $d_i$ be the maximum of the degrees of the minimal generators of $\Pc_i$. We say that $\Mc$ is \emph{$m$-regular} for some integer $m$ if $(d_i-i)\leq m$ for all $i$. We define the \emph{Castelnuovo-Mumford regularity} or just the \emph{regularity} of $\Mc$ to be the smallest integer $m$ for which $\Mc$ is $m$-regular.
\end{Definition}

A compact way to describe free resolutions of modules over $R$ is the so called \emph{Betti diagram}, which consists of a table with columns labeled $0,1,\dots$ corresponding to the modules $\Pc_0,\Pc_1,\dots$ in the resolution. The rows are labelled with integers corresponding to degrees. The entry of the $j-th$ row of the $i-th$ column contains the Betti number $\beta_{i,i+j}$. Betti diagrams appear in many places in this thesis, see for instance example \ref{example-complete}.

\subsubsection*{Regular sequences}
\spanishsubsubsection{Sucesiones regulares}

\begin{Definition}
Let $R$ be a ring and $\Mc$ an $R$-module. A sequence of elements $x_1,\dots,x_n\in R$ is called a {\it regular sequence} on $\Mc$ or {\it $\Mc$-regular sequence} if
\begin{enumerate}
\item $(x_1,\dots,x_n)\Mc\neq\Mc$,
\item For $i=1,\dots,n$, $x_i$ is a nonzerodivisor on $\Mc/(x_1,\dots,x_{i-1})\Mc$
\end{enumerate}
\end{Definition}
Regular Sequences extend the notion of non-zerodivisor and are in strong relation with the Koszul Complex (see \cite{E95}). In general, the fact that a sequence of more than two elements in $R$ form an $R$-regular sequence, depends on their order.
\begin{Definition}
The \emph{depth} of $R$ is defined as the maximum length of a regular $R$-sequence on $R$. More generally, the depth of an $R$-module $M$ is the maximum length of an $R$-regular sequence on $M$.
\end{Definition}

The depth of a module is always at least $0$ and no greater than the dimension of the module. The relation between depth and dimension gives rise the famous \emph{Cohen-Macaulay} condition of a module, one of which caracterizations says that a module $\Mc$ is Cohen-Macaulay if $depth(\Mc)=dim(\Mc)$ (see a list of caracterizations of this condition in \cite[Theorem 13.37]{MS04}).

\subsubsection*{Some basic concepts of Commutative Algebra}
\spanishsubsubsection{Algunos conceptos b\'asicos del \'algebra conmutativa}

We finish this section with the definition of some basic concepts of Commutative Algebra that appear throughout this thesis:

\begin{Definition}
Some relevant types of ideals:
\begin{itemize}
\item A proper ideal $I\subset R$ is \emph{prime} if $f\cdot g\in I$ implies either $f$ or $g$ is in $I$.
\item A proper ideal $I\subset R$ is \emph{primary} if $f\cdot g\in I$ implies either $f$ or $g^m$ is in $I$ for some $m\in \NN$.
\item A proper ideal $I\subset R$ is \emph{irreducible} if there do not exist ideals $J_1, J_2$ such that $I=J_1\cap J_2$, $I\subsetneq J_i$.
\item An ideal $I$ is radical if $f^m\in I$ ($m\in\NN$) implies $f\in\NN$
\end{itemize}
\end{Definition}

\begin{Definition}
Let $I$ be an ideal of $R$. An element $f\in R$ is said to be \emph{integral} over $I$ if it satisfies an equation of the form
$$f^m+a_1f^{m-1}+\cdots+a_m=0,\qquad \mbox{with }a_i\in I^i$$

The \emph{integral closure} $\bar{I}$ of $I$ is the set of the integral elements over $I$.
\end{Definition}

\begin{Definition}
The \emph{Krull} dimension of a ring $R$ is the supremmum of the lengths of chains of distinct prime ideals in $R$.
\end{Definition}
\section{Multilinear algebra}
\spanishsection{\'Algebra multilineal}

In this section we recall the definitions of the basic multilinear structures that are used in the text. The main object is the {\it Tensor algebra} and two of its quotients, namely the {\it Symmetric} and {\it Exterior} algebras.
\subsubsection*{The Tensor algebra}
\spanishsubsubsection{El \'algebra tensorial}

\begin{Definition}
Let $\kb$ a field and $V$ be a finite dimensional $\kb$-vector space, with $n=dim(V)$. The $i$-th tensor power of $V$ is the tensor product of $V$ with itself $i$ times:
$$T^iV:=V\otimes  \stackrel{i}{\cdots} \otimes V$$ 
By convention, $T^0=\kb$.
The {\rm Tensor algebra}, $TV$, of the vector space $V$ is the direct sum of all tensor powers of $V$:
$$T(V):=\bigoplus_{i=0}^{\infty} T^iV$$ 
multiplication is the tensor product, i.e. determined by the canonical isomorphism $T^iV\oplus T^jV\mapsto T^{i+j}V$. This makes $TV$ a graded algebra.
\end{Definition}

\subsubsection*{The Symmetric algebra}
\spanishsubsubsection{El \'algebra sim\'etrica}

\begin{Definition}
Let $I=\langle x\otimes y-y\otimes x\rangle$ be the ideal of $TV$ generated by the differences of all `symmetric products' of elements of $TV$. Then we define the {\it Symmetric algebra}, $SV$, of $V$ as the quotient $TV/I$
\end{Definition}

$SV$ is equivalent to the polynomial ring in $n$ variables $\kb[v_1,\dots,v_n]$ where $\{v_1,\dots,v_n\}$ form a basis of $V$. $SV$ inherits the grading of $TV$.

\subsubsection*{The Exterior algebra}
\spanishsubsubsection{El \'algebra exterior}

\begin{Definition}
Let $J=\langle x\otimes x\rangle$ be the ideal of $TV$ generated by the squares of all of elements of $TV$. Then we define the {\it Exterior algebra}, $\wedge V$ of $V$ as the quotient $TV/J$
\end{Definition}

$SV$ inherits the grading of $TV$. But in this case, we have that $\wedge V^i=0\forall i>n$. Moreover, we have that the dimension of $\wedge^iV$ as a $\kb$-vector space is ${i}\choose{2}$, and thus $dim_\kb(\wedge V)=2^n$.

\chapter{Algebraic Topology}\label{aptop}
\spanishchapter{Topolog\'ia Algebraica}
Algebraic Topology has two main faces: homotopical and homological. We will only be concerned with the second one, in particular with simplicial homology. A basic text for algebraic topology is \cite{M91}, which introduces homology in the cubical framework. For simplicial homology see for example \cite{H02}. A complete introduction that contains almost everything is \cite{S66}. From a computational viewpoint a good introduction to the standard algorithms which includes several chapters on applications is \cite{KMM04}.
\section{Simplicial complexes and homology}
\spanishsection{Complejos simpliciales y homolog\'ia}
Inside this thesis, the combinatorial nature of monomial ideals shows the importance of combinatorial structures to perform homological computations. When speaking about combinatorics and homology, one of the main objects are simplicial complexes. They have been put in relation with monomial ideals many times, see \cite{S96,B96,B02,MS04,P07,R07} among others. Therefore, we need to compute the homology of some simplicial complexes. The main definitions are the following

\begin{Definition}
The standard $n$-simplex is the subset of $\RR^{n+1}$ given by
$$\Delta^n =\{(t_0,\cdots,t_n)\in\RR^{n+1}\vert\sum_{i}{t_i} = 1 \mbox{ and } t_i \geq 0 \mbox{ for all } i\}$$ 
Setting $k$ components to $0$ we obtain the corresponding $n-k$ faces of this simplex, which are just copies of the standard $n-k$-simplex. Thus, the $0$-faces or vertices of $\Delta^n$ are the points $(1,0\dots,0),(0,1,0,\dots,0),\dots,(0,\dots,0,1)$.
Therefore, the standard $n$-simplex is the convex hull of a set of $(n+1)$ affinely independent points in $\RR^{n}$ (the vertices).
\end{Definition}

\begin{Definition}
A simplicial complex $\Delta$ is a set of simplices that satisfies the following conditions:
\begin{itemize}
\item If a simplex $\Pc$ is in $\Delta$, then every face $\Qc$ of $\Pc$ is also in $\Delta$
\item The intersection of any two simplices $\Pc,\Qc$ in $\Delta$ is either empty or is a face of both $\Pc$ and $\Qc$. 
\end{itemize}
\end{Definition}
In our work with simplicial complexes we will only focus on the combinatorial nature of them, without the geometric properties, so we will only use abstract simplicial complexes. Simplicial complexes as defined above are the geometric realizations of them:
\begin{Definition}
An \emph{abstract simplicial} complex is defined as follows:
Let $\{1,\dots,n\}$ be a set, which we will call the {\rm vertex set}; an {\rm abstract simplicial complex  $\Delta$ on the vertex set $\{1,\dots, n\}$} is a collection of subsets closed under taking subsets. The subsets in this collection are called {\rm faces}. The dimension of a face of cardinality $i$ is $i-1$, and the dimension of $\Delta$ is the maximum of the dimensions of its faces. The empty set is the only $(-1)$-dimesional face in any simplicial complex that is not the {\rm void} complex, the dimension of which is $-\infty$.

A maximal face in a (abstract) simplicial complex $\Delta$ is called a \emph{facet}. Note that such a complex is described by the set of its facets.
\end{Definition}
\begin{Example}\label{simplicial_example}
The abstract simplicial complex $\Delta$ on $9$ vertices whose geometric realization is in the picture has facets $\{1,2,3\},\{1,2,4\},\{1,3,4\},\{2,3,4\},\{4,5\},\{4,6\},\{5,6,7\},\{6,9\}$ and $\{7,8\}$, observe that the tethraedron in the vertices $1,2,3$ and $4$ is hollow.
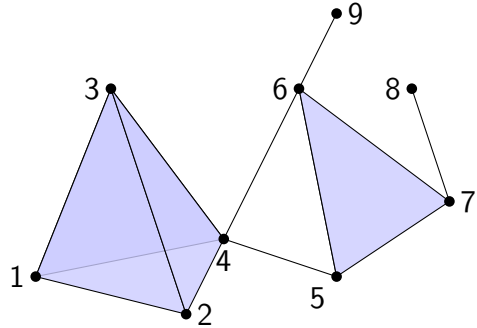
\begin{center}
\begin{tikzpicture}
    \tikzstyle{myfill} = [fill=blue!20,fill opacity=0.8]
	\filldraw[myfill](0,0)--(2.5,.5)--(1,2.5)--cycle;
	\filldraw[myfill](0,0)--(2,-.5)--(1,2.5)--cycle;
	\filldraw[myfill](1,2.5)--(2,-.5)--(2.5,.5)--cycle;
	\draw(2.5,.5)--(4,0)--(3.5,2.5)--cycle;
	\filldraw[myfill](5.5,1)--(4,0)--(3.5,2.5)--cycle;
	\draw(5.5,1)--(5,2.5);
	\draw(3.5,2.5)--(4,3.5);

	\foreach \position in {(0,0),(2,-.5),(1,2.5),(2.5,.5),(4,0),(3.5,2.5),(5.5,1),(5,2.5),(4,3.5)}
		{
		\fill \position circle(2pt);
		}
	\draw (0,0) node[left] {{\sf 1}};
	\draw (2,-.5) node[right] {{\sf 2}};
	\draw (1,2.5) node[left] {{\sf 3}};
	\draw (2.5,.5) node[below] {{\sf 4}};
	\draw (4,0) node[below left] {{\sf 5}};
	\draw (3.5,2.5) node[left] {{\sf 6}};
	\draw (5.5,1) node[right] {{\sf 7}};
	\draw (5,2.5) node[left] {{\sf 8}};
	\draw (4,3.5) node[right] {{\sf 9}};

\end{tikzpicture}
\end{center}
\end{Example}

Simplicial homology can be defined for standard simplicial complexes or for abstract simplicial complexes equivalently. The idea is to associate a chain complex to each simplicial complex, in which the differential combines the faces of each simplex in such a way that the square of the differential vanishes. We use here the formulation in terms of abstract simplicial complexes.
\begin{Definition}
For each $i$, let $F_i(\Delta)$ be the (finite) set of $i$-dimensional faces of $\Delta$, and let $\kb^{F_i(\Delta)}$ a vector space over $\kb$ whose basis elements $e_{\sigma}$ correspond to $i$-faces $\sigma\in F_i(\Delta)$. Then the {\rm reduced chain complex} of $\Delta$ over $\kb$ is the complex $\hat {C_\bullet}(\Delta;\kb): 0\rightarrow\kb^{F_{n-1}(\Delta)}\stackrel{\partial_{n-1}}{\rightarrow}\cdots\kb^{F_i(\Delta)}\stackrel{\partial_{i}}{\rightarrow}\kb^{F_{i-1} (\Delta)}\stackrel{\partial_{i-1}}{\rightarrow}\cdots \stackrel{\partial_{0}}{\rightarrow}\kb^{F_{-1}(\Delta)}\rightarrow 0$ where the differential or {\rm boundary map} is given by
$$\partial_i(e_{\sigma})=\sum_{j\in\sigma}(-1)^{s_j}e_{\sigma-j}$$ being $s_j=r+1$ if $j$ is in position $r$ in $\sigma$.

For each integer $i$, the $\kb$-vector space $\widetilde {H_i}(\Delta;\kb)=ker(\partial_i)/im(\partial_{i+1})$ is the $i$-th {reduced homology} of $\Delta$ over $\kb$.
\end{Definition}

In lower dimensions, simplicial homology has a geometrical meaning: for example, the dimension of the $0$-degree homology is the number of connected components of the simplicial complex.

\begin{Remark}
Computing simplicial homology is in general a difficult task, although recent progress in this area has made simplicial complexes a computational object, in particular new concepts as \emph{persistent homology} and new computer packages have been developed in the recent years. Some contributors to this progress have been Edelsbrunenr et al. \cite{ELZ02}, Dumas et al. \cite{DHSW04}, Carlsson, Perry, de Silva et al. \cite{PS06} or Kacynski et al. \cite{KMM04} among many others. This is an area of very active research which also has many applications inside and outside matheamtics.
\end{Remark}

\section{Homological tools}
\spanishsection{Herramientas homol\'ogicas}

There are a number of mathematical tools which are useful to perform homology computations. Some of them are of a topological or geometrical type and some are of a more algebraic nature. These tools include \emph{relative homology}, subdivisions, the Mayer-Vietoris exact sequence, homology of $CW$-complexes, several types of products, suspensions, dualities, the mapping cone and cylinder, etc. We are interested in the Mayer-Vietoris exact sequence and in mapping cones. These two constructions have a topological origin but are standard algebraic tools frequently used in the literature. When considered under an algebraic viewpoint, we can see these two important tools as instances of the exact homology sequence associated to a short exact sequence of complexes.

\subsection*{The exact homology sequence}
\spanishsubsection{La sucesi\'on exacta en homolog\'ia}

This is one of the most useful basic tools for computing homology. Here we just give the basic result, a proof of it can be seen in almost any book about algebraic topology or homological algebra. We follow the standard source \cite{M95}.

Consider any short exact sequence of chain complexes and chain-complex homomorphisms
\begin{equation}\label{ses_chain_complexes}
0\rightarrow \Kc\stackrel{\alpha}{\rightarrow}\Lc\stackrel{\beta}{\rightarrow}\Mc\rightarrow 0
\end{equation}

\begin{Theorem}[\cite{M95},Theorem 4.1]
For every short exact sequence \ref{ses_chain_complexes} of chain complexes, the corresponding long sequence
\begin{equation}\label{exact_sequence_homology}
\cdots\rightarrow H_{n+1}(\Mc)\stackrel{\Delta}{\rightarrow}H_n(\Kc)\stackrel{\alpha_*}{\rightarrow}H_n(\Lc)\stackrel{\beta_*}{\rightarrow}H_n(\Mc)\stackrel{\Delta}{\rightarrow}H_{n-1}(\Kc)\rightarrow\cdots
\end{equation}

of homology groups, with maps the connecting homomorphism $\Delta$, $\alpha_*=H_n(\alpha)$ and $\beta_*=H_n(\beta)$, is exact. 
\end{Theorem}
Here, $H_n(\alpha)$ denotes the induced morphism in homology. The proof of this theorem is a standard diagram chase exercise.

\subsection*{The Mayer-Vietoris exact sequence}
\spanishsubsection{La sucesi\'on exacta de Mayer-Vietoris}

The Mayer-Vietoris sequence is used to find the homology groups of a space from those of some other spaces, which we already know.
Supose we have two topological spaces $A$ and $B$. We are interested in the relations between the homology modules of $A$,$B$, $X=A^\circ\cup B^\circ$ and $A\cap B$. We can consider the obvious inclusion maps $i:A\cap B\rightarrow A$, $j:A\cap B\rightarrow B$, $k:A\rightarrow X$ and $l:B\rightarrow X$. These inclusions give rise to a short exact sequence of chain complexes
$$0\rightarrow \Cc(A\cap B)\stackrel{i\oplus j}{\rightarrow} \Cc(A)\oplus \Cc(B)\stackrel{k-l}{\rightarrow} \Cc(X)\rightarrow 0.$$
With certain mild assumptions, from this exact sequence, and the morphisms in homology induced by the inclusions we have an exact sequence in homology, called the {\it Mayer-Vietoris sequence}:
$$\cdots\rightarrow H_i(A\cap B)\rightarrow H_i(A)\oplus H_i(B)\rightarrow H_i(X)\stackrel{\Delta}{\rightarrow} H_{i-1}(A\cap B)\rightarrow\cdots$$
Here, $\Delta$ is called the \emph{connecting homomorphism}. There are many versions of this sequence appearing in different areas, and in many cases the proof is a special case of the exact sequence in homology associated to a short exact sequence of chain complexes. This is the approach used in section \ref{m-v} and chapter \ref{computation}.

\subsection*{The mapping cone}
\spanishsubsection{El cono de una aplicaci\'on}

Mapping cones are constructions coming from topology, which have analogues in algebra. Here we define the {\it algebraic Mapping Cone} of chain maps (chain complexes morphisms).
\begin{Definition}
Let $(C,\partial)$ and $(C',\partial')$ chain complexes, and $f:C\rightarrow C'$ a chain map. The algebraic mapping cone of $f$, denoted $(M(f),\delta)$ is a chain complex defined as follows:
$$M(f)_i=C_{i-1}\oplus C_i$$
$$\delta_i(c,c')=(\partial_{i-1}c,\partial'_ic'+(-1)^i f_{i-1}c)$$
for $c\in C_{i-1}$, $c'\in C'_{i}$. The differential $\delta$ verifies $\delta^2=0$.
\end{Definition}
The chain morphisms given by the inclusions $\mi_i:C'\rightarrow M(f)$ and the projections $\mj_i:M(f)\rightarrow C_{i-1}$ defined by $\mi_i(c')=(0,c')$ and $\mj(c,c')=c$, give rise to a short exact sequence of chain complexes
$$0\rightarrow C'\rightarrow M(f)\rightarrow C\rightarrow 0$$
Again, using theorem \ref{exact_sequence_homology} we build the corresponding exact sequence in homology:
$$\cdots\rightarrow H_i(C')\rightarrow H_i(M(f))\rightarrow H_{i-1}(C)\stackrel {\Delta}{\rightarrow}H_{n-1}(C')\rightarrow\cdots$$

where the connecting homomorphism is given by $\Delta=f_*$, the homomorphism in homology induced by $f$. Therefore, we have that every chain complex morphism produces such an exact sequence in homology. Applications of this construction to resolutions in the polynomial ring can be seen in \cite{HT02,S99} and in section \ref{hom_tools}.

\chapter{Effective Homology}\label{apef}
\spanishchapter{Homolog\'ia Efectiva}
\begin{flushright}
\mbox{\parbox{0.5\textwidth}{\footnotesize\it%
Proverb: the difference between effective homology and ordinary homology consists in using the \emph{explicit} homotopy operators. 
\rm\cite[p.62]{RS06}}}
\end{flushright}

Effective homology is a constructive approach to algebraic topology and homological algebra. First introduced by F. Sergeraert \cite{S94} and developed by him and Rubio \cite{RS02} it is founded on one hand on the so-called \emph{Basic Perturbation Lemma} from a theoretical viewpoint, and on the other hand on \emph{Functional Programming} from a computational point of view. Actual computer programs have been developed based on this method, namely EAT \cite{EAT} and the newer Kenzo \cite{Kenzo}. Using effective homology new homotopy and homology groups have been computed; among the first applications was the computation of the homology groups of iterated loop spaces (see for example \cite{RS02}). In recent times, effective homology has been actually applied to $A$-infinity structures \cite{Be06}, spectral sequences \cite{RRS06,Ro07}, or Koszul homology \cite{RS06} producing concrete computer programs.

The constructive nature of effective homology is well expressed in its systematic methods for solving what is called the \emph{homological problem for a chain complex}, which encodes the information one expects to obtain for a chain complex using this constructive approach:.

\begin{Definition}[\cite{RS06}, Definition 41]\label{homological-problem}
Let $R$ be a ground ring and $(\Cc,d)$ a chain complex of $R$-modules. A \emph{solution $S$ of the homological problem for $\Cc$} is a set $S=(\sigma_i)_{1\leq i\leq5}$ of five algorithms:
\begin{enumerate}
\item $\sigma_1:\Cc\rightarrow\{\bot=\mbox{false},\top=\mbox{true}\}$ is a predicate deciding for every $n\in\ZZ$ and every $n$-chain $c\in\Cc_n$ whether $c$ is an $n$-cycle or not, in other words whether $dc=0$ or $dc\neq 0$, whether $c\in Z_n(\Cc)$ or not.
\item  $\sigma_2:\ZZ\rightarrow\{R-modules\}$ associates to every integer $n$ some $R$-module $\sigma_2(n)$ in principle isomporphic to $H_n(\Cc)$. The image $\sigma_2(n)$ will \emph{model} the isomorphism class of $H_n(\Cc)$ in an effective way.
\item The algorithm $\sigma_3$ is indexed by $n\in\ZZ$; for every $n\in \ZZ$ the algorithm $\sigma_{3,n}:\sigma_2(n)\rightarrow Z_n(\Cc)$ associates to every $n$-homology class $\mathfrak h$ coded as an element $\mathfrak h\in\sigma_2(n)$ a cycle $\sigma_{3,n}(\mathfrak h)\in Z_n(\Cc)$ representing this homology class.
\item The algorithm $\sigma_4$ is indexed by $n\in\ZZ$; for every $n\in\ZZ$, the algorithm $\sigma_{4,n}:C_n\supset Z_n(\Cc)\rightarrow \sigma_2(n)$ associates to every $n$-cycle $z\in Z_n(\Cc)$ the homology class of $z$ coded as an element of $\sigma_2(n)$
\item The algorithm $\sigma_5$ is indexed by $n\in\ZZ$; for every $n\in\ZZ$, the algorithm $\sigma_{5,n}:Z_n(\Cc)\rightarrow C_{n+1}$ associates to every $n$-cycle $z\in Z_n(\Cc)$ known as a boundary by the previous algorithm, a boundary preimage $c\in C_{n+1}$ s.th. $dc=z$. In particular $Z_n(\Cc):=ker\sigma_{4,n}$
\end{enumerate}
\end{Definition}

\section{Reductions}
\spanishsection{Reducciones}
A main concept in effective homology is that of a \emph{reduction}, it can be seen as an explicit homology equivalence between two chain complexes, a \emph{big} one and a \emph{small} one. This notion makes more precise the statements in the definition of the \emph{homological problem}.
\begin{Definition}
A reduction $\rho: \widehat{\Cc}\rrdc\Cc$ is a diagram

\begin{center}
\begin{tikzpicture}

\draw[->](0,0) node(x){$\widehat{\Cc}$}	
	(x) ..controls +(1,0.5)and +(1,-0.5) ..(x) node[pos=0.5,right]{$h$}; 
\draw[->](-0.1,-0.3)--(-0.1,-1.5) node[pos=0.5,left]{$f$} node[below]{$\Cc$};
\draw[->](0.1,-1.5)--(0.1,-0.3) node[pos=0.5,right]{$g$};
\end{tikzpicture}
\end{center}

where:
\begin{enumerate}
\item $\widehat{\Cc}$ and $\Cc$ are chain-complexes.
\item $f$ and $g$ are chain-complex morphisms
\item $h$ is a homotopy operator of degree $+1$
\item {The following relations are satisfied
	\begin{enumerate}
	\item $fg=id_\Cc$
	\item $gf+dh+hd=id_{\widehat{\Cc}}$
	\item $fh=hg=hh=0$
	\end{enumerate}
	}
\end{enumerate}
\end{Definition}
The meaning of the homological equivalence  between the two chain complexes that is expressed by such a reduction is explained in the following results from \cite{RS06}:
\begin{Proposition}
Let $\rho: \widehat{\Cc}\rrdc\Cc$ be a reduction. This reduction is equivalent to a decomposition $\widehat{\Cc}=\Ac\oplus\Bc\oplus\Cc'$ such that:
\begin{enumerate}
\item $\widehat{\Cc}\supset \Cc'=im(g)$ is a subcomplex of $\widehat{\Cc}$
\item $\Ac\oplus\Bc=ker(f)$ s a subcomplex of $\widehat{\Cc}$
\item $\widehat{\Cc}\supset\Ac=ker(f)\cap ker(h)$ is not in general a subcomplex of $\widehat{\Cc}$
\item $\widehat{\Cc}\supset\Bc=ker(f)\cap ker(d)$ is a subcomplex of $\widehat{\Cc}$ with null differentials
\item The chain-complex morphisms $f$ and $g$ are inverse isomorphisms between $\Cc'$ and $\Cc$
\item The arrows $d$ and $h$ are module isomorphisms of respective degrees $-1$ and $+1$ between $\Ac$ and $\Bc$
\end{enumerate}
\end{Proposition}

In the terminology of effective homology, the notion of \emph{locally effective} object is applied to objects for which we do not have \emph{global} information, but we can perform particular i.e. \emph{local} computations. In contrast, \emph{effective} objects are those about which we know ``everything'', even \emph{global} properties; see more details in \cite{RS06}. The following theorem applies to \emph{locally effective} objects.

\begin{Theorem}
Let $\rho: \widehat{\Cc}\rrdc\Cc$ be a reduction where the chain-complexes $\widehat{\Cc}$ and $\Cc$ are locally effective. If the homological problem is solved in the small chain-complex $\Cc$ then the reduction $\rho$ induces a solution of the homological problem for the big chain-complex $\widehat{\Cc}$.
\end{Theorem}

Of course, the result is trivial the other way round:

\begin{Proposition}
Let $\rho: \widehat{\Cc}\rrdc\Cc$ be a reduction where the homological problem is solved for $\widehat{\Cc}$. Then the  homological problem is solved for the small chain-complex $\Cc$.
\end{Proposition}

\section{The Basic Perturbation Lemma}
\spanishsection{El lema b\'asico de perturbaci\'on}
The Basic Perturbation Lemma (\cite{S62,B67,RS02}), is the main tool of {\it effective homology}, and can be seen as an algorithm for computations in many different situations where homology appears.
The basics are as follows: Let $\rho=(\widehat C,C,f,g,h)$ be a reduction, which as we have seen is a description of the homology of the big complex $\widehat C$ through the small one $C$. If we modify the differential of $\widehat C$ under certain conditions, then coherent modifications can be applied to the other components of the reduction, so that a new reduction is obtained, and thus, the homology of the {\it new} big chain complex (note that what is new is the differential, not the underlying graded module) can be described again through the homology of the {\it new} small complex. More explicitly:
\begin{Definition}
A {\rm perturbation} of the differential $d$ of a chain complex $C$ is an operator $\delta:C\rightarrow C$ of degree $-1$ such that the sum $d+\delta$ is again a differential i.e. $(d+\delta)\circ(d+\delta)=0$.
\end{Definition}

\begin{Definition}
A perturbation $\widehat\delta$ of the differential $\widehat d$ of the top chain complex $\widehat C$ of a reduction $\rho=(\widehat C,C,f,g,h)$ satisfies the condition of {\it local nilpotency} if for every $x\in\widehat C$, there exists an integer $n$ satisfying the relation $(h\circ\widehat \delta)^n(x)=0$; equivalently, for every $x\in \widehat C$, there exists an integer $n$ such that $(\widehat\delta\circ h)^n(x)=0$
\end{Definition}

Now we can state the Basic Perturbation Lemma, as written in \cite{RS02} following \cite{B67} and \cite{S62}:
\begin{Theorem}
{\bf (Basic Perturbation Lemma)} Let $\mathfrak{R}$ be the type of reductions. An algorithm can be constructed
$${\mathbf {bpl:}}[\mathfrak R\times \mathfrak P]_\chi\longrightarrow \mathfrak {R}$$
where $[\mathfrak R\times \mathfrak P]_\chi$ is the set of coherent pairs $(\rho=(\widehat C,C,f,g,h),\widehat{\delta})$, that is, $\widehat\delta$ is a perturbation of the differential $\widehat d$ of $\widehat C$ satisfying the condition of local nilpotency; the output is a reduction $\rho'=(\widehat C',C',f',g',h')$ where the new top chain complex $\widehat{C'}$ is the old one provided with the new differential $\widehat{d'}=\widehat d+\widehat\delta$; in particular the new bottom chain complex $C'$ is the old one with a new differential $d+\delta$, where $d$ is the old differential of $C$ and $\delta$ is a perturbation determined by the algorithm {\bf bpl}; the same for the new maps $f', g'$ and $h'$.
\end{Theorem}

\begin{Remark}
The basic perturbation lemma provides actual formulas for the morphisms in the new reduction, which can be expressed in the following way:
\begin{description}
\item [$\widehat{\delta'}=$] $\widehat d+\widehat\delta$
\item [$d'=$]$d+\delta$ where $\delta=f\widehat\delta\phi g=f\psi\widehat\delta g$
\item [$f'=$]$f\psi$
\item [$g'=$]$\phi g$
\item [$h'=$]$\phi h=h\psi$

\end{description}
where $\phi$ and $\psi$ are given by
$$\phi=\sum_{i=0}^\infty (-1)^i(h\widehat\delta)^i\qquad \psi=\sum_{i=0}^\infty (-1)^i(\widehat\delta h)^i$$
Note that because of the local nilpotency condition, these series have only a finite number of summands for each $x$.
\end{Remark}
An easy version of the Perturbation Lemma shows up when the perturbation occurs in the {\it small} complex:
\begin{Theorem}
{\bf Easy Basic Perturbation Lemma.} If we have a reduction $\rho=(\widehat C,C,f,g,h)$ and a perturbation $\delta$ of the diferential $d$ of $C$, then the differential $\widehat d$ of $\widehat{C}$ can be perturbed so that we have a new reduction between $\widehat C$ and $C$ with the new differentials.
\end{Theorem}

Many constructions in algebraic topology and homological algebra are suitable to be combined with the basic perturbation lemma in order to perform effective computations. As an example, we give here the combination of the basic perturbation lemma and the cone construction we have seen in Appendix \ref{aptop}, which is widely used in section \ref{hom_tools} chapter \ref{computation}:

\begin{Theorem}[\cite{RS06}, Theorem 62] 
Let $\rho=(f,g,h):\Cc\rrdc\Dc$ and $\rho'=(f',g',h'): \Cc'\rrdc\Dc'$ be two reductions and $\phi:\Cc\rightarrow\Cc'$ a chain-complex morphism. Then these data define a canonical reduction
$$\rho''=(f'',g'',h''): Cone(\phi)\rrdc Cone(f\phi g')$$
\end{Theorem}

\bibliographystyle{alpha}

\addcontentsline{spt}{chapter}{Bibliograf\'ia}


\fancyhead[EL]{\it \thepage} 
\fancyhead[ER]{\itshape Bibliography} 
\fancyhead[OL]{\itshape Bibliography}
\fancyhead[OR]{\it \thepage}  
\fancyfoot{}        

\bibliography{tesis}

\newcommand{\etalchar}[1]{$^{#1}$}
\begin{thebibliography}{HHMT06}

\bibitem[AB58]{AB58}
M.~Auslander and D.~Buchsbaum.
\newblock Codimension and multiplicity.
\newblock {\em The Annals of Mathematics}, 68,3:625--657, 1958.

\bibitem[Agn97]{A97}
G.~Agnarsson.
\newblock On the number of outside corners of monomial ideals.
\newblock {\em Journal of Pure and Applied Algebra}, 117:3--22, 1997.

\bibitem[AH95]{AH95}
A.~Aramova and J.~Herzog.
\newblock Free resolutions and {K}oszul homology.
\newblock {\em Journal of Pure and Applied Algebra}, 105:1--16, 1995.

\bibitem[AH96]{AH96}
A.~Aramova and J.~Herzog.
\newblock Koszul cycles and {E}liahou-{K}ervaire type resolutions.
\newblock {\em Journal of Algebra}, 181:347--370, 1996.

\bibitem[AHH98]{AHH98}
A.~Aramova, J.~Herzog, and T.~Hibi.
\newblock Squarefree lexsegment ideals.
\newblock {\em Matematische Zeitschrift}, 228:353--378, 1998.

\bibitem[Ape95]{A95}
J.~Apel.
\newblock A {G}r\"obner approach to involutive bases.
\newblock {\em Journal of Symbolic Computation}, 19:441--458, 1995.

\bibitem[Bat02]{B02}
E.~Batzies.
\newblock {\em Discrete {M}orse theory for cellular resolutions}.
\newblock PhD thesis, Universit\"at Marburg, 2002.

\bibitem[Bay96]{B96}
D.~Bayer.
\newblock Monomial ideals and duality.
\newblock Lecture notes, Berkeley, 1996.

\bibitem[Ber06]{Be06}
A.~Berciano.
\newblock {\em C\'alculo simb\'olico y t\'ecnicas de control de {A}-infinito
  estructuras}.
\newblock PhD thesis, Universidad de {S}evilla, {S}pain, 2006.

\bibitem[BES61]{BES61}
Z.W. Birnbaum, J.D. Esary, and S.~Saunders.
\newblock Multi-component systems and structures and their reliability.
\newblock {\em Technometrics}, 3:55--77, 1961.

\bibitem[BG06]{BG06}
I.~Bermejo and P.~Gimenez.
\newblock Saturation and {C}astelnuovo-{M}umford regularity.
\newblock {\em Journal of Algebra}, 303:592--617, 2006.

\bibitem[BH97]{BH97}
W.~Bruns and J.~Herzog.
\newblock Semigroup rings and simplicial complexes.
\newblock {\em Journal of Pure and Applied Algebra}, 122:185--208, 1997.

\bibitem[BH98]{BH98}
W.~Bruns and J.~Herzog.
\newblock {\em Cohen-{M}acaulay rings}, volume~39 of {\em Cambridge Studies in
  Advanced Mathematics}.
\newblock Cambridge University Press, Cambridge, revised edition edition, 1998.

\bibitem[Big97]{B97}
A.~M. Bigatti.
\newblock Computation of {H}ilbert-{P}oincar\'e series.
\newblock {\em Journal of Pure and Applied Algebra}, 119(3):237--253, 1997.

\bibitem[BPS98]{BPS98}
D.~Bayer, I.~Peeva, and B.~Sturmfels.
\newblock Monomial resolutions.
\newblock {\em Mathematical Research Letters}, 5:353--378, 1998.

\bibitem[BR64]{BR64}
D.~A. Buchsbaum and D.~S. Rim.
\newblock A generalized {K}oszul complex. {II}. {D}epth and multiplicity.
\newblock {\em Transactions of the American Mathematical Society},
  111(2):197--224, 1964.

\bibitem[BR65]{BR65}
D.~A. Buchsbaum and D.~S. Rim.
\newblock A generalized {K}oszul complex. {III}. {A} remark on generic
  acyclicity.
\newblock {\em Proceedings of the American Mathematical Society},
  16(3):555--558, 1965.

\bibitem[Bro67]{B67}
R.~Brown.
\newblock The twisted {E}ilenberg-{Z}ilberg theorem.
\newblock {\em Celebrazioni Arch. Secolo XX, Simp. Top.}, pages 34--37, 1967.

\bibitem[BS98]{BS98}
D.~Bayer and B.~Sturmfels.
\newblock Cellular resolutions.
\newblock {\em J. Reine Agew. Math.}, 502:123--140, 1998.

\bibitem[Buc64]{B64}
D.~A. Buchsbaum.
\newblock A generalized {K}oszul complex. {I}.
\newblock {\em Transactions of the American Mathematical Society},
  111(2):183--196, 1964.

\bibitem[Buc06]{B06}
B.~Buchberger.
\newblock Bruno {B}uchberger's {P}h{D} thesis 1965: {A}n algorithm for finding
  the basis elements of the residue class ring of a zero dimensional polynomial
  ideal.
\newblock {\em J. Symb. Comput.}, 41(3-4):475--511, 2006.

\bibitem[But04]{B04}
F.~Butler.
\newblock Rook theory and cycle-counting permutation statistics.
\newblock {\em Advances in Applied Mathematics}, 33:655--675, 2004.

\bibitem[BW98]{BW98}
B.~Buchberger and F.~Winkler, editors.
\newblock {\em Gr\"obner bases and applications}, volume 251 of {\em London
  Mathematical Society lecture notes series}.
\newblock Cambridge University press, 1998.

\bibitem[CLO96]{CLO96}
D.~Cox, J.~Little, and D.~O'Shea.
\newblock {\em Ideals, varieties and algorithms}.
\newblock Undergraduate Texts in Mathematics. Springer Verlag, New York, second
  edition, 1996.

\bibitem[CLO98]{CLO98}
D.~Cox, J.~Little, and D.~O'Shea.
\newblock {\em Using algebraic geometry}, volume 185 of {\em Graduate Texts in
  Mathematics}.
\newblock Springer, New York, 1998.

\bibitem[CN06]{CN06}
Alberto Corso and Uwe Nagel.
\newblock Monomial and toric ideals associated to {F}errers graphs.
\newblock Available at
  \url{http://www.citebase.org/abstract?id=oai:arXiv.org:math/0609371}, 2006.

\bibitem[CN07]{CN07}
Alberto Corso and Uwe Nagel.
\newblock Specializations of {F}errers ideals.
\newblock Available at {\tt
  http://www.citebase.org/abstract?id=oai:arXiv.org:math/0703695}, 2007.

\bibitem[{CoC}a]{cocoa}
{CoCoA}Team.
\newblock {{\hbox{\rm C\kern-.13em o\kern-.07em C\kern-.13em o\kern-.15em A}}}:
  a system for doing {C}omputations in {C}ommutative {A}lgebra.
\newblock Available at \url{http://cocoa.dima.unige.it}.

\bibitem[{CoC}b]{cocoalib}
{CoCoA}Team.
\newblock {{\hbox{\rm C\kern-.13em o\kern-.07em C\kern-.13em o\kern-.15em
  A\kern-.13em L\kern-.13em i\kern-.13em b}}}:a {GPL} {C++} library for doing
  {C}omputations in {C}ommutative {A}lgebra.
\newblock Available at \/ {\tt http://cocoa.dima.unige.it}.

\bibitem[Cri02]{C02}
V.~Crispin.
\newblock Integrally closed monomial ideals and powers of ideals.
\newblock Research Reports in Mathematics~7, Stockholm University, 2002.

\bibitem[DER97]{DER97}
I.S. Duff, A.M. Erisman, and J.K. Reid.
\newblock {\em Direct methods for sparse matrices}.
\newblock Monographs on Numerical Analysis. Oxford University press, 1997.

\bibitem[DHSW04]{DHSW04}
J-G. Dumas, F.~Heckenbach, B.~D. Saunders, and V.~Welker.
\newblock Simplicial homology, a (proposed) {GAP} share package.
\newblock Available at \url{http://www.cis.udel.edu/~dumas/Homology/}, 2004.

\bibitem[Din97]{D97}
K.~Ding.
\newblock Rook placements and cellular decomposition of partition varieties.
\newblock {\em Discrete Mathematics}, 170:107--151, 1997.

\bibitem[Doh03]{D03}
K.~Dohmen.
\newblock {\em Improved {B}onferroni inequalities via abstract tubes}, volume
  1826 of {\em Lecture Notes in Mathematics}.
\newblock Springer, Berlin, 2003.

\bibitem[DRSS]{Kenzo}
X.~Dousson, J.~Rubio, F.~Sergeraert, and Y.~Siret.
\newblock The {K}enzo program.
\newblock Available at \url{\tt
  http://www-fourier.ujf-grenoble.fr/~sergerar/Kenzo/}.

\bibitem[DSV01]{DSV01}
J-G. Dumas, B.D. Saunders, and G.~Villard.
\newblock On efficient sparse integer matrix {S}mith normal form computations.
\newblock {\em Journal of Symbolic Computation}, 32:71--99, 2001.

\bibitem[EGSS02]{EGSS02}
D.~Eisenbud, D.~R. Grayson, M.~Stillman, and B.~Sturmfels, editors.
\newblock {\em Computing in algebraic geometry with {M}acaulay2}, volume~8 of
  {\em Algorithms and computation in mathematics}.
\newblock Springer Verlag, 2002.

\bibitem[Eis95]{E95}
D.~Eisenbud.
\newblock {\em Commutative algebra with a view towards algebraic geometry},
  volume 150 of {\em Graduate Texts in Mathematics}.
\newblock Springer, New York, 1995.

\bibitem[Eis04]{E04}
D.~Eisenbud.
\newblock {\em The geometry of syzygies}, volume 229 of {\em Graduate Texts in
  Mathematics}.
\newblock Springer, New York, 2004.

\bibitem[EK90]{EK90}
S.~Eliahou and M.~Kervaire.
\newblock Minimal resolutions of some monomial ideals.
\newblock {\em Journal of Algebra}, 129:1--25, 1990.

\bibitem[ELZ02]{ELZ02}
H.~Edelsbrunner, D.~Letscher, and A.~Zomorodian.
\newblock Topological persistence and simplification.
\newblock {\em Discrete Comput. Geom.}, 28:511--533, 2002.

\bibitem[EM66]{EM66}
S.~Eilenberg and J.C. Moore.
\newblock Homology and fibrations {I}: {C}oalgebras, cotensor product ans its
  derived functors.
\newblock {\em Comment. Math. Helv.}, 40:199--236, 1966.

\bibitem[ER98]{ER98}
J.~A. Eagon and V.~Reiner.
\newblock Resolutions of {S}tanley-{R}eisner rings and {A}lexander duality.
\newblock {\em Journal of Pure and Applied Algebra}, 130:265--275, 1998.

\bibitem[EvW04]{EW04}
R.~Ehrenborg and S.~van Willigenburg.
\newblock Enumerative properties of {F}errers graphs.
\newblock {\em Discrete Computational Geometry}, 32:481--492, 2004.

\bibitem[Far02]{F02}
S.~Faridi.
\newblock The facet ideal of a simplicial complex.
\newblock {\em Manuscripta Mathematica}, 109:159--174, 2002.

\bibitem[Far05]{F05}
S.~Faridi.
\newblock {\em Commutative Algebra: Geometric, Homological, Combinatorial and
  Computational Aspects}, volume 244 of {\em Lecture notes in pure and applied
  mathematics}, chapter Monomial ideals via square-free monomial ideals, pages
  85--114.
\newblock Taylor \& Francis, CRC press, 2005.

\bibitem[FR07]{FR07}
C.A. Francisco and B.P. Richert.
\newblock {\em Syzygies and Hilbert Functions}, volume 254 of {\em Lecture
  Notes in Pure and Applied Mathematics}, chapter Lex-{P}lus-{P}owers ideals,
  pages 113--144.
\newblock Chapman and Hall/CRC, Boca Raton, 2007.

\bibitem[Fro78]{F78}
R.~Froeberg.
\newblock {\em S\'eminaire d'{A}lg\`ebre paul {D}ubreil 31\`eme ann\'ee
  ({P}aris 1977-1978)}, chapter Some complex constructions with appliactions to
  {P}oincar\'e series, pages 272--284.
\newblock Springer Verlag, 1978.

\bibitem[Fro99]{F99}
R.~Froeberg.
\newblock {\em Advances in Commutative Ring Theory. Proceedings of the 3rd
  International Conference, Fez}, volume 205 of {\em Lecture Notes in Pure and
  Applied Mathematics}, chapter Koszul Algebras, pages 337--350.
\newblock Marcel Dekker, 1999.

\bibitem[FV07]{FV07}
C.~A. Francisco and A.~VanTuyl.
\newblock Sequentally {C}ohen-{M}acaulay edge ideals.
\newblock {\em Proceedings of the American Mathematical Society},
  135:2327--2337, 2007.

\bibitem[GB98]{GB98}
V.~P. Gerdt and Y.~A. Blinkov.
\newblock Minimal involutive bases.
\newblock {\em Math. Comp. Simul.}, 45:543--560, 1998.

\bibitem[GBY01]{GBY01}
V.~P. Gerdt, Y.~A. Blinkov, and D.~A. Yanovich.
\newblock {\em Computer Algebra in Scientific Computing-CASC 2001}, chapter
  Construction of {J}anet bases {I}: {M}onomial bases, pages 233--247.
\newblock Springer, Berlin, 2001.

\bibitem[GLR{\etalchar{+}}06]{GLRRV06}
T.~Gonzalez, M.~Llabres, J.~Rocha, F.~Rosello, and G.~Valiente.
\newblock Algebraic distances for phylogenetic trees.
\newblock In {\em Proceedings of the Xth EACA}, pages 107--110, Sevilla
  (Spain), 2006.

\bibitem[GM03]{GM03}
S.~I. Gelfand and Y.~I. Manin.
\newblock {\em Methods of homological algebra}.
\newblock Springer monographs in Mathematics. Springer, 2nd edition, 2003.

\bibitem[Gol65]{G65}
H.L. Goldschmidt.
\newblock Existence theorems for analytic linear partial differential
  equations.
\newblock {\em Ann. Math.}, 82:246--270, 1965.

\bibitem[Gol69]{G69}
H.L. Goldschmidt.
\newblock Integrability criteria for systems of non-linear partial differential
  equations.
\newblock {\em J. Diff. Geom.}, 1:269--307, 1969.

\bibitem[GP02]{GP02}
G.M. Greuel and G.~Pfister.
\newblock {\em A Singular introduction to commutative algebra}.
\newblock Springer, 2002.

\bibitem[GPS05]{Singular}
G.-M. Greuel, G.~Pfister, and H.~Sch\"onemann.
\newblock {\sc Singular} 3.0.
\newblock {A Computer Algebra System for Polynomial Computations}, Centre for
  Computer Algebra, University of Kaiserslautern, 2005.
\newblock {\tt http://www.singular.uni-kl.de}.

\bibitem[GPW99]{GPW99}
V.~Gasharov, I.~Peeva, and V.~Welker.
\newblock The lcm-lattice in monomial resolutions.
\newblock {\em Mathematical Research Letters}, 6:521--532, 1999.

\bibitem[Gre84]{G84}
M.~L. Green.
\newblock Koszul cohomology and the geometry of projective varieties.
\newblock {\em Journal of Differential Geometry}, 19:125--171, 1984.

\bibitem[Gre87]{G87}
M.~Green.
\newblock {\em Lectures on {R}iemann surfaces}, chapter Koszul cohomology and
  geometry, pages 177--200.
\newblock World Scientific Press, 1987.

\bibitem[GS]{M2}
Daniel~R. Grayson and Michael~E. Stillman.
\newblock Macaulay 2, a software system for research in algebraic geometry.
\newblock Available at http://www.math.uiuc.edu/Macaulay2/.

\bibitem[GS64]{GS64}
V.W. Guillemin and S.~Stenberg.
\newblock An algebraic model of transitive differential geometry.
\newblock {\em Bull. Amer. Math. Soc.}, 70:16--47, 1964.

\bibitem[GW03]{GW03}
B.~Giglio and H.P. Wynn.
\newblock Alexander duality and moments in reliability modelling.
\newblock {\em Applicable Algebra in Engineering, Communication and Computing},
  14:153--174, 2003.

\bibitem[GW04]{GW04}
B.~Giglio and H.P. Wynn.
\newblock Monomial idelas and the {S}carf complex for coherent systems in
  reliability theory.
\newblock {\em Annals of Statistics}, 32:1289--1311, 2004.

\bibitem[Hal87]{H87}
S.~Halperin.
\newblock Le complexe de {K}oszul en algebre et topologie.
\newblock {\em Annales de L'Institut Fourier}, 37(4):77--97, 1987.

\bibitem[Hat02]{H02}
A.~Hatcher.
\newblock {\em Algebraic {T}opology}.
\newblock Cambridge University Press, 2002.

\bibitem[Her92]{H92}
J.~Herzog.
\newblock Canonical {K}oszul cycles.
\newblock {\em Aportaciones matem\'aticas, {N}otas de investigaci\'on},
  6:33--41, 1992.

\bibitem[Her04]{H04}
J.~Herzog.
\newblock Finite free resolutions.
\newblock School on Commutative Algebra and Interactions with Algebraic
  Geometry and Combinatorics, Trieste, 2004. The Abdus Salam International
  Centre for Theoretical Physics.
\newblock Lecture notes.

\bibitem[HHMT06]{HHMT06}
J.~Herzog, T.~Hibi, S.~Murai, and Y.~Takayama.
\newblock Componentwise linear ideals with minimal or maximal {B}etti numbers.
\newblock Available at
  \url{http://www.citebase.org/abstract?id=oai:arXiv.org:math/0610470}, 2006.

\bibitem[Hoc77]{H77}
M.~Hochster.
\newblock Cohen-{M}acaulay rings, combinatorics and simplicial complexes.
\newblock In B.R.~Mc Donald and R.~Morris, editors, {\em Ring Theory II, Proc.
  Second Conf, Univ. Oklahoma, Norman, Okla. 1975}, volume~26 of {\em Lecture
  Notes in Pure and Applied Mathematics}. Marcel Dekker, New York, 1977.

\bibitem[HPV03]{HPV03}
J.~Herzog, D.~Popescu, and M.~Vladoiu.
\newblock {\em Commutative Algebra}, volume 331 of {\em Contemp. Math.},
  chapter On the {E}xt-modules of ideals of {B}orel type, pages 171--186.
\newblock American Mathematical Society, Providence, 2003.

\bibitem[HS01]{HS01}
S.~Hosten and G.G. Smith.
\newblock {\em Computations in algebraic geometry with {M}acaulay2}, volume~8
  of {\em Algorithms and Computations in Mathematics}, chapter Monomial ideals,
  pages 73--100.
\newblock Springer, New York, 2001.

\bibitem[HS02]{HS02}
M.~Hausdorf and W.~M. Seiler.
\newblock An efficient algebraic algorithm for the geometric completion to
  involution.
\newblock {\em Applicable Algebra in Engineering Communication and Computing},
  13:163--207, 2002.

\bibitem[HT99]{HT99}
S.~Hosten and R.R. Thomas.
\newblock Standard pairs and group relaxations in integer programming.
\newblock {\em Journal of Pure and Applied Algebra}, 139:133--157, 1999.

\bibitem[HT02]{HT02}
J.~Herzog and Y.~Takayama.
\newblock Resolutions by mapping cones.
\newblock {\em Homology, Homotopy and Applications}, 4(2):277--294, 2002.

\bibitem[HV07]{HV07}
H.~T. Ha and A.~VanTuyl.
\newblock Resolutions of square-free monomial ideals via facet ideals: a
  survey.
\newblock {\em Contemporary Math}, 2007.
\newblock To appear.

\bibitem[Jan20]{J20}
M.~Janet.
\newblock Sur les syst\`emes d'\'equations aux d\'eriv\'ees partielles.
\newblock {\em J. Math. Pure Appl.}, 3:65--151, 1920.

\bibitem[Jan29]{J29}
M.~Janet.
\newblock {\em Le\c{c}ons sur les Syst\`emes d'\'Equations aux D\'eriv\'ees
  Partielles}, volume Fascicule IV of {\em Cahiers Scientifiques}.
\newblock Gauthier-Villars, Paris, 1929.

\bibitem[JLS02]{JLS02}
L.~Johansson, L.Lambe, and E.~Sk\"oldberg.
\newblock On constructing resolutions over the polynomial algebra.
\newblock {\em Homology, Homotopy and Applications}, 4:315--336, 2002.

\bibitem[Joh04]{J04}
M.~Johansson.
\newblock Computation of denominator polynomials for {P}oincare series in
  monomial rings.
\newblock Master's thesis, Stokholm University, 2004.

\bibitem[JZ07]{JZ07}
A.~S. Jahan and X.~Zheng.
\newblock Pretty clean monomial ideals and linear quotients.
\newblock Available at
  \url{http://www.citebase.org/abstract?id=oai:arXiv.org:math/0707.2914}, 2007.

\bibitem[KL06]{KL06}
B.~Kugrligov and V.~Lychagin.
\newblock Dimension of the solutions space of {PDE}s.
\newblock In J.~Calmet, W.~M. Seiler, and R.~W. Tucker, editors, {\em Global
  integrability of field theories}, pages 5--26, 2006.

\bibitem[KL07]{KL07}
B.~S. Kruglikov and V.~V. Lychagin.
\newblock Geometry of differential equations.
\newblock Preprint IHES/M/07/04, Institut des Hautes \'Etudes Scientifiques,
  Bures-sur-Yvette, 2007.

\bibitem[KLV86]{KLV86}
I.~S. Krasilshchik, V.~V. Lychagin, and A.~M. Vinogradov.
\newblock {\em Geometry of Jet Spaces and Nonlinear Partial Differential
  Equations}.
\newblock Gordon \& Breach, New York, 1986.

\bibitem[KMM04]{KMM04}
T.~Kaczynski, M.~Mischaikow, and M.~Mrozek.
\newblock {\em Computational Homology}, volume 157 of {\em Applied mathematical
  Sciences}.
\newblock Springer Verlag, 2004.

\bibitem[Kon80]{K80}
J.M. Kontoleon.
\newblock Reliability determination of a r-sucesive-out-of-n:f system.
\newblock {\em IEEE Transactions in Reliability}, 29:315--336, 1980.

\bibitem[Kos50a]{K50a}
J-L. Koszul.
\newblock Homologie et cohomologie des algebres de {L}ie.
\newblock {\em Bulletin de la S. M. F.}, 78:65--127, 1950.

\bibitem[Kos50b]{K50b}
J-L. Koszul.
\newblock Sur un type d'alg\`ebres differ\'entielles en rapport avec la
  transgression.
\newblock {\em Colloque de Topologie, Bruxelles}, pages 73--81, 1950.

\bibitem[KR00]{KR00}
M.~Kreuzer and L.~Robbiano.
\newblock {\em Computational commutative algebra 1}.
\newblock Springer, Berlin, 2000.

\bibitem[KR05]{KR05}
M.~Kreuzer and L.~Robbiano.
\newblock {\em Computational commutative algebra 2}.
\newblock Springer, Berlin, 2005.

\bibitem[KZ03]{KZ03}
W.~Kou and M.J. Zuo.
\newblock {\em Optimal reliability modelling}.
\newblock Wiley and Sons, New Jersey, 2003.

\bibitem[LS02]{LS02}
L.~Lambe and W.M. Seiler.
\newblock Differential equations, {S}pencer cohomology and computing
  resolutions.
\newblock {\em Georgian Mathematical Journal}, 9:723--772, 2002.

\bibitem[Lyu98]{L98}
G.~Lyubeznik.
\newblock A new explicit finite free resolution of ideals generated by
  monomials in an {R}-sequence.
\newblock {\em Journal of Pure and Applied Algebra}, 51(1-2):193--195, 1998.

\bibitem[Mac95]{M95}
S.~MacLane.
\newblock {\em Homology}.
\newblock Classics in Mathematics. Springer, New York, 1995.
\newblock Reprint of the 1974 edition.

\bibitem[Mal98]{M98}
D.~Mall.
\newblock On the relation between {G}r\"obner and pommaret bases.
\newblock {\em Applicable Algebra in Engineering, Communication and Computing},
  9:117--123, 1998.

\bibitem[Mal03]{M03}
B.~Malgrange.
\newblock Cartan involutiveness = {M}umford regularity.
\newblock {\em Contemp. Math.}, 331:193--205, 2003.

\bibitem[Mas91]{M91}
W.S. Massey.
\newblock {\em A basic course in algebraic topology}, volume 127 of {\em
  Graduate Texts in Mathematics}.
\newblock Springer, New York, 1991.

\bibitem[Mat53]{M53}
Y.~Matsushima.
\newblock On a theorem concerning the prolongation of a differential system.
\newblock {\em Nagoya Math. J.}, 6:1--16, 1953.

\bibitem[Mil00]{M00}
E.~Miller.
\newblock {\em Resolutions and duality for monomial ideals}.
\newblock PhD thesis, University of California, Berkeley, 2000.

\bibitem[Mil04]{M04}
R.~A. Milowski.
\newblock Computing irredundant irreducible decompositions and the {S}carf
  complex of large scale monomial ideals.
\newblock Master's thesis, San Francisco State university, 2004.

\bibitem[MP01]{MP01}
E.~Miller and D.~Perkinson.
\newblock Eight lectures on monomial ideals.
\newblock In A.~V. Geramita, editor, {\em COCOA VI: proceedings of the
  international school, Villa Gualino, May-June 1999}, number 120 in Queens
  Papers in Pure and Applied Mathematics, pages 3--105. Queens Univerity,
  Kingston, Ontario, 2001.

\bibitem[MS04]{MS04}
E.~Miller and B.~Sturmfels.
\newblock {\em Combinatorial commutative algebra}, volume 227 of {\em Graduate
  texts in Mathematics}.
\newblock Springer, New York, 2004.

\bibitem[MSY00]{MSY00}
E.~Miller, B.~Sturmfels, and K.~Yanagawa.
\newblock Generic and cogeneric monomial ideals.
\newblock {\em Journal of Symbolic Computation}, 29:691--708, 2000.

\bibitem[NW92]{NW92}
D.Q. Naimann and H.P. Wynn.
\newblock Inclusion-exclusion-{B}onferroni identities and inequalities for
  discrete tube-like problems via {E}uler characteristics.
\newblock {\em Annals of Statistics}, 20:43--76, 1992.

\bibitem[NW97]{NW97}
D.Q. Naiman and H.P. Wynn.
\newblock Abstract {T}ubes, improved inclusion-exclusion identities and
  inequalities and importance sampling.
\newblock {\em Annals of Statistics}, 25:1954--1983, 1997.

\bibitem[NW01]{NW01}
D.Q. Naiman and H.P. Wynn.
\newblock Improved inclusion-exclusion-{B}onferroni inequalities for simplex
  and orthants arrangements.
\newblock {\em Journal of Inequalities in Pure and Applied Mathematics},
  2:article 18, 2001.

\bibitem[OW07]{OW07}
P.~Orlik and V.~Welker.
\newblock {\em Algebraic Combinatorics}.
\newblock Universitext. Springer, Berlin Heidelberg, 2007.
\newblock Lectures on Arrangements and Cellular Resolutions at a Summer School
  in {N}ordfjordeid, {N}orway, {J}une 2003.

\bibitem[PdS06]{PS06}
P.~Perry and V.~de~Silva.
\newblock Plex: Simplicial complexes in {MATLAB}.
\newblock Available at \url{http://comptop.stanford.edu/programs/}, 2006.

\bibitem[Pee02]{P02}
I.~Peeva.
\newblock Resolutions and lattices.
\newblock {\em Homology, Homotopy and Applications}, 4:427--437, 2002.

\bibitem[Pee07]{P07}
I.~Peeva, editor.
\newblock {\em Syzygies and {H}ilbert functions}, volume 254 of {\em Lecture
  Notes in Pure and Applied Mathematics}.
\newblock Chapman and Hall/CRC, Boca Raton, 2007.

\bibitem[Pha05]{P05}
Jeffry Phan.
\newblock Minimal monomial ideals and linear resolutions.
\newblock Available at
  \url{http://www.citebase.org/abstract?id=oai:arXiv.org:math/0511032}, 2005.

\bibitem[Pom78]{P78}
J.F. Pommaret.
\newblock {\em Systems of Partial Differential Equations and Lie Pseudogroups}.
\newblock Gordon \& Breach, London, 1978.

\bibitem[Pom88]{P88}
J.F. Pommaret.
\newblock {\em Lie Pseudogroups and Mechanics}.
\newblock Gordon \& Breach, London, 1988.

\bibitem[Pom01a]{P01a}
J.~F. Pommaret.
\newblock {\em Partial Differential Control Theory I: Mathematical Tools},
  volume 530 of {\em Mathematics and its Applications}.
\newblock Kluwer, Dordrecht, 2001.

\bibitem[Pom01b]{P01b}
J.~F. Pommaret.
\newblock {\em Partial Differential Control Theory II: Control Systems}, volume
  530 of {\em Mathematics and its Applications}.
\newblock Kluwer, Dordrecht, 2001.

\bibitem[PS07]{PS07}
I.~Peeva and M.~Stillmann.
\newblock Syzigies and {H}ilbert functions.
\newblock available at
  \url{http://www.math.cornell.edu/~irena/papers/overview.pdf}, 2007.

\bibitem[Rom07]{Ro07}
A.~Romero.
\newblock {\em Effective homology and spectral sequences}.
\newblock PhD thesis, Universidad de {L}a {R}ioja, {S}pain, 2007.

\bibitem[Rou07]{R07}
B.~H. Roune.
\newblock The label algorithm for irreducible decomposition of monomial ideals.
\newblock Available at \url{http://arxiv.org/abs/0705.44833}, 2007.

\bibitem[RRS06]{RRS06}
A.~Romero, J.~Rubio, and F.~Sergeraert.
\newblock Computing spectral sequences.
\newblock {\em Journal of symbolic computation}, 41(10):1059--1079, 2006.

\bibitem[RS02]{RS02}
J.~Rubio and F.~Sergeraert.
\newblock Constructive algebraic topology.
\newblock {\em Bulletin des Sciences Math\'emathiques}, 126:389--412, 2002.

\bibitem[RS06]{RS06}
J.~Rubio and F.~Sergeraert.
\newblock Constructive homological algebra and applications.
\newblock Available at
  \url{http://www-fourier.ujf-grenoble.fr/~sergerar/Papers/Genova-Lecture-Note%
s.pdf}, 2006.
\newblock lecture notes MAP summer School.

\bibitem[RSS]{EAT}
J.~Rubio, F.~Sergeraert, and Y.~Siret.
\newblock {EAT}: {S}ymbolic software for effective homology computation.
\newblock Available at \/ {\tt ftp://ftp-fourier.ujf-grenoble.fr/pub/EAT/}.

\bibitem[RSS06]{RSS06}
J.~Rubio, E.~{S\'aenz de Cabez\'on}, and F.~Sergeraert.
\newblock Minimal resolutions.
\newblock SIGSAM Bulletin, submitted, 2006.

\bibitem[RV01]{RV01}
V~Reiner and W.~Velker.
\newblock Linear syzygies of {S}tanley-{R}eisner ideals.
\newblock {\em Mathematica Scandinavica}, 89(1):117--132, 2001.

\bibitem[S06a]{S06b}
E.~{S\'aenz-de-Cabez\'on}.
\newblock Computing {K}oszul homology for monomial ideals.
\newblock In {\em Proceedings of the {R}hine {W}orkshop on {C}omputer
  {A}lgebra, Basel}, 2006.
\newblock Available at http://arxiv.org/abs/math.AC/0605325.

\bibitem[S06b]{S06a}
E.~{S\'aenz-de-Cabez\'on}.
\newblock {M}ayer- {V}ietoris trees of monomial ideals.
\newblock In J.~Calmet, W.~M. Seiler, and R.~W. Tucker, editors, {\em Global
  integrability of field theories}, pages 311--334, 2006.

\bibitem[S07a]{SC07d}
E.~{S\'aenz-de-Cabez\'on}.
\newblock Homological description of monomial ideals using {M}ayer-{V}ietoris
  trees.
\newblock In {\em Proceedings of {MEGA} 2007, Strobl am Wolfgangsee (Austria)},
  2007.

\bibitem[S07b]{SC07c}
E.~{S\'aenz-de-Cabez\'on}.
\newblock Multigraded {B}etti numbers without computing minimal free
  resolutions.
\newblock Submitted, 2007.

\bibitem[Sah07]{Sa07}
M.~Sahbi.
\newblock Pommaret bases and the computation of the koszul homology in the
  monomial case.
\newblock Master's thesis, Fakult\"at f\"ur Informatik, Universit\"at
  Karlsruhe, 2007.

\bibitem[Sch80]{S80}
F.~O. Schreyer.
\newblock Die berechnung von {S}yzygien mit dem verallgemeinerten
  {W}eierstrasschen {D}ivisionssatz.
\newblock Master's thesis, Fakult\"at f\"ur Mathematik, Universit\"at Hamburg,
  1980.

\bibitem[Sch03]{S03}
H.~Schenk.
\newblock {\em Computational algebraic geometry}, volume~58 of {\em London
  Mathematical Society student texts}.
\newblock Cambridge University press, 2003.

\bibitem[Sei02a]{S02d}
W.~M. Seiler.
\newblock Involution - the formal theory of differential equations and its
  applications in computer algebra and numerical analysis.
\newblock Habilitation thesis, Dept. of Mathematics and Computer Science,
  Universit\"at Mannheim, 2002.

\bibitem[Sei02b]{S02c}
W.~M. Seiler.
\newblock Taylor and {L}yubeznik resolutions via {G}r\"obner bases.
\newblock {\em Journal of Symbolic Computation}, 34:597--608, 2002.

\bibitem[Sei07a]{S02a}
W.~M. Seiler.
\newblock A combinatorial approach to involution and $\delta$-regularity {I}:
  Involutive bases in polynomials algebras of solvable type.
\newblock Preprint Universit\"at Kassel, 2007.

\bibitem[Sei07b]{S02b}
W.~M. Seiler.
\newblock A combinatorial approach to involution and $\delta$-regularity {II}:
  Structural analysis of polynomial modules with pommaret bases.
\newblock Preprint Universit\"at Kassel, 2007.

\bibitem[Sei07c]{S07b}
W.~M. Seiler.
\newblock Spencer cohomology, differential equations and {P}ommaret bases.
\newblock In M.~Rosenkranz and D.~M. Wang, editors, {\em Gr\"obner basis in
  symbolic analysis}, volume~2 of {\em Radon series on Computational and
  Applied Mathematics}, pages 171--219. de Gruyter, Berlin, 2007.

\bibitem[Sei07d]{S07}
W.M. Seiler.
\newblock {\em Involution - The formal theory of differential equations and its
  applications in computer algebra}.
\newblock Algorithms and {C}omputation in {M}athematics. Springer, Berlin,
  2007.
\newblock To Appear.

\bibitem[Ser94]{S94}
F.~Sergeraert.
\newblock The computability problem in algebraic topology.
\newblock {\em Advances in Mathematics}, 104:1--29, 1994.

\bibitem[Ser06]{S06}
F.~Sergeraert.
\newblock Effective homology of {K}oszul complexes.
\newblock Preprint available at
  \url{http://www-fourier.ujf-grenoble.fr/~sergerar/Papers/Koszul.pdf}, 2006.

\bibitem[Sha00]{Sh00}
R.~Y. Sharp.
\newblock {\em Steps in commutative algebra}, volume~51 of {\em London
  Mathematical Society Student Texts}.
\newblock Cambridge University Press, second edition, 2000.

\bibitem[Shi62]{S62}
W.~Shih.
\newblock Homologie des espaces fibr\'es.
\newblock {\em Publications Math\'emathiques de l'I.H.E.S.}, 13, 1962.

\bibitem[Sid07]{Si07}
J.~Sidman.
\newblock {\em Syzygies and Hilbert Functions}, volume 254 of {\em Lecture
  Notes in Pure and Applied Mathematics}, chapter Resolutions and Subspace
  arrangements, pages 249--265.
\newblock Chapman and Hall/CRC, Boca Raton, 2007.

\bibitem[Sie99]{S99}
T.~Siebert.
\newblock Algorithms for the computation of free resolutions.
\newblock In B.H. Matzat, G.~M. Greuel, and G.~Hiss, editors, {\em Algorithmic
  algebra and number theory. Selected papers from a conference; Heidelberg,
  Germany, October 1997}, pages 295--310. Springer Verlag, 1999.

\bibitem[Spa66]{S66}
E.~H. Spanier.
\newblock {\em Algebraic topology}.
\newblock Springer, New York, 1966.

\bibitem[Spe69]{S69}
D.C. Spencer.
\newblock Overdetermined systems of linear partial differential equations.
\newblock {\em Bull. Amer. Math. Soc.}, 75:179--239, 1969.

\bibitem[SS98]{LS98}
R.~La Scala and M.~Stillman.
\newblock Strategies for computing minimal free resolutions.
\newblock {\em Journal of Symbolic Computation}, 26,4:409--431, 1998.

\bibitem[Sta78]{S78}
R.~P. Stanley.
\newblock Hilbert functions of graded algebras.
\newblock {\em Advances in Mathematics}, 28:175--193, 1978.

\bibitem[Sta96]{S96}
R.~P. Stanley.
\newblock {\em Combinatorics and commutative algebra}, volume~41 of {\em
  Progress in Mathematics}.
\newblock Birkh\"auser, Boston, second edition, 1996.

\bibitem[SW91]{SW91}
B.~Sturmfels and N.~White.
\newblock Computing combinatorial decompositions of rings.
\newblock {\em Combinatorica}, 11:275--293, 1991.

\bibitem[SW07]{SW07}
E.~{S\'aenz-de-Cabez\'on} and H.~P. Wynn.
\newblock Betti numbers and minimal free resolutions for multi-state system
  reliability bounds.
\newblock Submitted. Appears in Proceedings of {MEGA} 2007, Strobl am
  Wolfgangsee (Austria), 2007.

\bibitem[Tay60]{T60}
D.~Taylor.
\newblock {\em Ideals generated by monomials in an $R$-sequence}.
\newblock PhD thesis, University of Chicago, 1960.

\bibitem[Val04]{V04}
G.~Valla.
\newblock Betti numbers of some monomial ideals.
\newblock {\em Proceedings of the AMS}, 133-1:57--63, 2004.

\bibitem[Vas05]{V05}
W.~Vasconcelos.
\newblock {\em Integral closure. Rees algebras, multiplicities, algorithms}.
\newblock Springer Monographs in Mathematics. Springer, Berlin, 2005.

\bibitem[Vil01]{V01}
R.~Villarreal.
\newblock {\em Monomial algebras}, volume 238 of {\em Pure and Applied
  Mathematics}.
\newblock Marcel Dekker inc., New York, 2001.

\bibitem[Yan99]{Y99}
K.~Yanagawa.
\newblock ${F}_{\Delta}$-type free resolutions of monomial ideals.
\newblock {\em Proceedings of the {AMS}}, 127(2):377--383, 1999.

\end{thebibliography}



\newpage \thispagestyle{empty}
\backmatter


\begin{otherlanguage}{spanish}

\renewcommand{\headrulewidth}{0.0pt} 
\fancyhf{} 

\begin{titlepage}
  \vspace*{1cm}
  \begin{center}%
    {\Huge {\bf Homolog\'ia de Koszul Combinatoria, \\ C\'alculo y Aplicaciones \par}}%
    \vspace*{1.5cm}%
    {\Large \bf Eduardo S\'aenz-de-Cabez\'on Irigaray \par}
    \vspace*{1.5cm}
    {\Large Memoria presentada para la \\ obtenci\'on del grado de Doctor\par}
    \vspace*{0.5cm}%
        
      \par
      \vskip 5em%
\end{center}
\begin{flushright}
   \begin{large}
      \begin{tabular}{l l}
           Directores: & Prof. Dr. Luis Javier Hern\'andez Paricio \\
                       & Prof. Dr. Werner M. Seiler
      \end{tabular}
   \end{large}
\end{flushright}
\vspace*{3.5cm}
\begin{center}

{\Large Universidad de La Rioja \\}
{\large Departamento de Matem\'aticas y Computaci\'on \par}

\vspace*{1cm}

{\large Logro\~no, Diciembre de 2007}
\end{center}\end{titlepage}
\newpage
\thispagestyle{empty} \vspace*{19cm}

\noindent Este trabajo ha sido parcialmente subvencionado por los proyectos CALCULEMUS y GIFT (NEST Contract No 5006) de la Uni\'on Europea, \mbox{project ANGI-2005/10} de la Comunidad Aut\'onoma de La Rioja, y \mbox{becas ATUR-04/51}, \mbox{ATUR-05/46}, \mbox{ATUR-06/31} de la Universidad de La Rioja.

\pagenumbering{roman}

\fancyhf{}
\fancyfoot[C]{\thepage} 

\renewcommand{\headrulewidth}{0.1pt} 

\fancyhf{}
\fancyhead[EL]{\it \thepage} 
\fancyhead[ER]{\itshape \'Indice general} 
\fancyhead[OL]{\itshape \'Indice general}
\fancyhead[OR]{\it \thepage}  
\fancyfoot{}        

\spanishTOC 
\newpage


\mainmatter

\fancyhf{}
\fancyhead[EL]{\it \thepage} 
\fancyhead[ER]{\itshape Introducci\'on} 
\fancyhead[OL]{\itshape Introducci\'on}
\fancyhead[OR]{\it \thepage}  
\fancyfoot{}        

\chapter*{Introducci\'on}

Los ideales monomiales son un tipo especial de ideales del anillo de polinomios que tienen una naturaleza combinatoria. Juegan un papel muy importante en el \'algebra conmutativa, ya que algunos problemas sobre ideales o m\'odulos en el ideal de polinomios se pueden reducir a problemas sobre ideales monomiales, en particular en el contexto de las t\'ecnicas basadas en bases de Gr\"obner. Por otro lado, algunas propiedades de ciertas clases de ideales monomiales son importantes en la teor\'ia de sizigias y funciones de Hilbert. Adem\'as, la naturaleza combinatoria de los ideales de monomios los hace adecuados para aplicaciones en otras \'areas de las matem\'aticas y tambi\'en fuera de las matem\'aticas, desde la teor\'ia de grafos o los sistemas diferenciales hasta la teor\'ia de fiabilidad. En los \'ultimos a\~nos ha habido un gran inter\'es en este tipo de ideales, y se han convertido en \'area de investigaci\'on muy activa.
%

En esta tesis nos ocupamos de las propiedades homol\'ogicas de los ideales de monomios. Normalmente la descripci\'on homol\'ogica de un ideal monomial viene dada por su resoluci\'on libre minimal, a partir de la que se calculan los principales invariantes homol\'ogicos del ideal. Sin embargo, la obtenci\'on de una descripci\'on expl\'icita de forma compacta de la resoluci\'on minimal de un ideal de monomios es un problema abierto, aunque ha habido trabajos muy interesantes sobre este problema en a\~nos pasados. Aqu\'i usaremos la homolog\'ia de Koszul para dar la descripci\'on homol\'ogica de un ideal monomial. Veremos que esta homolog\'ia nos proporciona un buen modo de describir las propiedades homol\'ogicas del ideal as\'i como otras de sus propiedades estructurales. Ambos acercamientos son esencialmente equivalentes, ya que representan dos modos diferentes de calcular ciertos m\'odulos $Tor$ del ideal.


La naturaleza combinatoria de los ideales monomiales introduce un modo combinatorio de trabajar la homolog\'ia de Koszul, por eso en este contexto hablamos de {\bf homolog\'ia de Koszul combinatoria}, que da t\'itulo a esta tesis.


En el primer cap\'itulo se presentan los protagonistas de la historia: la homolog\'ia de Koszul y los ideales monomiales. La homolog\'ia de Koszul es la homolog\'ia de un complejo introducido originalmente por J-L. Koszul en un contexto geom\'etrico \cite{K50a,K50b}. Ha sido objeto de inter\'es para el \'algebra conmutativa desde hace a\~nos, y est\'a tambi\'en en la confluencia de problemas importantes en la teor\'ia formal de sistemas diferenciales con el \'algebra conmutativa, debido a su relaci\'on con la cohomolog\'ia de Spencer \cite{S69}. Por otro lado, ha habido grandes trabajos por aprte de los algebristas sobre ideales monomiales al ser objetos importantes, en particular en la teor\'ia de bases de Gr\"obner, lo que nos permite reducir muchos problemas relacionados con idelaes polinomiales a problemas sobre ideales monomiales, m\'as sencillos de manejar debido a su anturaleza combinatoria. Este cap\'itulo estar\'a dedicado a presentar las notaciones y nociones b\'asicas y las propiedades principales de estos dos conceptos.


El segundo cap\'itulo est\'a dedicado a la descripci\'on de las propiedades homol\'ogicas y estructurales de los ideales de monomios que pueden ser le\'idas a partir de la homolog\'ia de Koszul. Primero nos centramos en invariantes y propiedades homol\'ogicos, que son el principal objetivo del cap\'itulo. Despu\'es tratamos algunas propiedades algebraicas de estos ideales. \'Estas incluyen descomposiciones de Stanley y descomposiciones irreducibles y primarias irredundantes. Tambi\'en transferimos los resultados sobre homolog\'ia de ideales monomiales a ideales polinomiales. Para ello, necesitamos la teor\'ia de perturbaci\'on homol\'ogica junto con la de bases de Gr\"obner. Esto hace a los m\'etodos descritos en el cap\'itulo segunso aplicables an contextos m\'as amplios y nos permite seguir un programa similar al usado en la teor\'ia de bases de Gr\"obner.


El tercer cap\'itulo est\'a dedicado a los c\'alculos. En \'el damos un algoritmo para calcular la homolog\'ia de Koszul de ideales monomiales basado en diferentes t\'ecnicas. Usamos t\'ecnicas homol\'ogicas y combinatorias e introducimos los \'arboles de Mayer-Vietoris, que nos permiten no s\'olo realizar c\'alculos homol\'ogicos sobre ideales monomiales sino que constituyen una nueva t\'ecnica para analizar la estructura de estos ideales. En este contexto se analizan distintos tipos de \'arboles de Mayer-Vietoris. Otra herramienta usada en este cap\'itulo es la homolog\'ia simplicial; para esto, se obtienen mejoras provinientes de dualidades y teor\'ia de Morse y en particular de la aplicaci\'on de la teor\'ia de Stanley-Reisner a complejos simpliciales de Koszul. Se proporciona un estudio de nuestro algoritmo junto a algunos puntos sobre su implementaci\'on y experimentos y comparaciones con otros algoritmos de prop\'osito similar. ëstos muestran que los \'arboles de Mayer-Vietoris son una alternativa eficiente para la realizaci\'on de c\'alculos homol\'ogicos en ideales monomiales.


El cap\'itulo cuarto est\'a dedicado a las aplicaciones. Se muestran diversas aplicaciones: Aplicaciones de los \'arboles de Mayer-Vietoris a diferentes tipos de ideales monomiales que son aplicados ya sea dentro del \'algebra conmutativa (tipo Borel-fijo, ideales estables, segmento, gen\'ericos) o en otras \'areas (Valla, Ferrers, cuasi-estables). Tambi\'en se desarrollan aplicaciones en otros campos, como la teor\'ia formal de sistemas diferenciales y la teor\'ia de fiabilidad. Estas aplicaciones usan las propiedades de la homolog\'ia de Koszul descitas en el segundo cap\'itulo y las heramientas computacionales mostradas en el tercero.


Hay conceptos de distintas \'areas de las matem\'aticas que aparecen en diversos lugares de la tesis. Algunos lectores tendr\'an probablemente familiariedad a algunos de ellos y no con otros, o viceversa. Por esta raz\'on y por ayudar a la facilidad de lectura, inclu\'imos varios ap\'endices en los que se dan las definiciones m\'as relevantes. La intenci\'on de estos ap\'endices es servir como referencia para los conceptos principales, no como introducciones o explicaciones de las diferentes teor\'ias implicadas.


\fancyhead[ER]{\itshape Resumen de los cap\'itulos} 
\fancyhead[OL]{\itshape Resumen de los cap\'itulos}

\chapter*{Resumen de los cap\'itulos}

Presentamos a continuaci\'on un breve resumen en castellano de los cap\'itulos de esta memoria
\section*{Cap\'itulo 1: Homolog\'ia de Koszul e ideales monomiales}

En este cap\'itulo se introducen los principales conceptos sobre cuya interacci\'on versa la presente tesis: la homolog\'ia de Koszul y los ideales monomiales. El cap\'itulo est\'a dividido en tres secciones.

En la primera secci\'on se da la definici\'on, origen y principales propiedades de la homolog\'ia de Koszul. Se hace especial hincapi\'e en su versi\'on graduada y multigraduada. Se explica tambi\'en la dualidad entre la homolog\'ia de Koszul y la cohomolog\'ia de spencer, origen de la aplicaci\'on de la primera en el estudio de sistemas diferenciales. Las referencias principales para esta secci\'on son \cite{K50a,K50b,S69,S07}.

En la segunda secci\'on se da la definici\'on y principales propiedades de los ideales monomiales. Bas\'andonos en las propiedades combinatorias de estos objetos, se pueden obtener caracterizaciones de sus principales invariantes algebraicos y homol\'ogicos, as\'i como algoritmos para su c\'alculo. Se presta en esta secci\'on especial atenci\'on a las resoluciones de estos ideales. Las referencias principales son \cite{V01,MS04}.

En la tercera secci\'on se presenta un cat\'alogo de t\'ecnicas de origen topol\'ogico y homol\'ogico aplicables al estudio de la homolog\'ia de Koszul de los ideales monomiales. Las principales herramientas introducidas son los complejos simpliciales de Koszul, los ideales de Stanley-Reisner y la dualidad de Alexander, as\'i como algoritmos basados en ellos. Adem\'as se aplican las sucesiones de Mayer-Vietoris y el cono de una aplicaci\'on, para el c\'alculo de resoluciones y de la homolog\'ia de Koszul de ideales de monomios. Algunas herramientas y algoritmos son presentados aqu\'i por vez primera. Las referencias principales son \cite{MS04,B96,S96}.

\section*{Cap\'itulo 2: Homolog\'ia de Koszul y estructura de ideales monomiales}

En este cap\'itulo la homolog\'ia de Koszul se utiliza para describir la estructura de los ideales de monomios. El objetivo principal es obtener descripciones de un ideal de monomios partiendo del conocimiento de su homolog\'is de Koszul. Estas descripciones son de distinto tipo: resoluci\'on m\'inima, descomposiciones combinatorias y descomposiciones en ideales irreducibles o primarios. Finalmente una \'ultima secci\'on trata de la obtenci\'on de resoluciones de ideales polinomiales a partir de resoluciones de sus ideales iniciales, que son de tipo monomial.

En la primera secci\'on se describe la obtenci\'on de los n\'umeros de Betti y de resoluciones m\'inimas de ideales de monomios a partir de la homolog\'ia de Koszul. El trabajo se basa en dos bicomplejos, el bicomplejo asociado al funtor $Tor$ y el bicomplejo de Aramova y Herzog \cite{AH95}. Usando t\'ecnicas de perturbaci\'on homol\'ogica (ap\'endice \ref{apef}) se obtienen m\'etodos expl\'icitos para el c\'alculo de resoluciones m\'inimas a partir de generadores de la homolog\'ia de Koszul. las referencias principales de esta secci\'on son \cite{AH95, AH96} y \cite{RS06,RSS06}.

En la segunda secci\'on, dado un ideal monomial $I\subseteq R=\kb[x_1\dots,x_n]$ y su homolog\'ia de Koszul, se obtiene una descomposici\'on de Stanley de $R/I$. Se trata por un lado el caso de ideales cero-dimensionales y basado en \'este, el caso general. En este tipo de descomposiciones juega un papel destacado la homolog\'ia de grado $n-1$ y su relaci\'on con el borde del ideal.

En la tercera secci\'on se trata el problema de obtener una descomposici\'on m\'inima en ideales irreducibles de un ideal de monomios a partir de su homolog\'ia de Koszul. Al igual que en la secci\'on anterior, se distingue el caso de ideales artinianos, y en el caso de los no artinianos se obtiene una descomposici\'on a partir del caso artiniano y se da otro proceso para obtenerla de modo independiente. Estos m\'etodos suponen una alternativa a otras formas de obtener descomposiciones de ideales de monomios en ideales irreducibles, un problema que ya ah sido tratado en la literatura \cite{MS04,M04,R07}.

En la cuarta secci\'on se da un sencillo procedimiento para obtener una descomposici\'on m\'inima en ideales priamrios de un ideal de monomios a partir de su homolog\'is de Koszul. Obteniendo esta descomposici\'on se obtine la altura del ideal as\'i como sus primos asociados. De nuevo, tiene una importancia capital la homolog\'ia de Koszul de grado $n-1$.

Finalmente, la quinta secci\'on se ocupa de ideales polinomiales. El ideal inicial de un ideal $I$ generado por polinomios es un ideal monomial al que podemos aplicar las t\'ecnicas que conocemos para calcular su homolog\'ia de koszul y/o su resoluci\'on minimal. A partir de \'esto, usando t\'ecnicas de perturbaci\'on homol\'ogica y bases de Gr\"obner, obtenemos la homolog\'ia de Koszul y/o resoluciones del ideal $I$. Esta t\'ecnica ha sido utilizada en la literatura de dos formas diferentes, ambas son expuestas aqu\'i. Por un lado Lambe y otros \cite{JLS02,LS02} utilizan la perturbaci\'on de las resoluciones e Taylor y Lyubeznik. Por otro lado, Sergeraert y otros \cite{RS06,S06,RSS06} utilizan el cono efectivo de una aplicaci\'on.

\section*{Cap\'itulo 3: C\'alculo de la Homolog\'ia de Koszul}

En este cap\'itulo introducimos una t\'ecnica nueva para el c\'alculo de la homolog\'is de Koszul, los n\'umeros de Betti y resoluciones de ideales de monomios. Esta t\'ecnica se basa en las sucesiones y \'arboles de Mayer-Vietoris de un ideal de monomios. 

En la primera secci\'on del cap\'itulo se introducen los \'erboles de Mayer.Vietoris, a partir de las sucesiones del mismo nombre introducidas en el primer cap\'itulo. Los \'arboles de Mayer-Vietoris pueden considerarse por un lado una descripci\'on del ideal junto con la parte relevante de su ret\'iculo de m\'inimos com\'un m\'ultiplos ($lcm$-lattice en el texto), y por otro como un algoritmo que nos permite el c\'alculo sencillo de los invariantes homol\'ogicos del ideal. En concreto permite obtener cotas para los n\'umeros de Betti multigraduados sin calcular la resoluci\'on m\'inima del ideal. En esta secci\'on se introducen tambi\'en als familias de ideales de Mayer-Vietoris, que son aquellos apra los que estas cotas son alcanzadas. Cada \'arbol de Mayer-Vietoris de un ideal nos proporciona una resoluci\'on del mismo, y por tanto nos da una expresi\'ond e su serie de Hilbert multigraduada. Adem\'as, se dan procesos para calcular explicitamente resoluciones y la homolog\'ia de Koszul del ideal a partir de estos \'arboles.

La segunda secci\'on analiza distintos tipos de \'arboles de Mayer-Vietoris que ser\'an aplicados en el cap\'itulo cuarto para ciertas familias de ideales. Esta secci\'on incide en los \'arboles de Mayer-Vietoris como descripciones del ideal, y muestran c\'omo el an\'alisis de estos \'arboles puede usarse para obtener los principales invariantes homol\'ogicos del ideal y nos permiten un estudio m\'as detallado del ideal considerado.

Finalmente, la tercera secci\'on se ocupa de los \'arboles de Mayer-Vietoris desde un punto de vista algor\'itmico. Se expone en ella el algoritmo b\'asico de construcci\'on de \'arboles de Mayer-Vietoris as\'i como distintas cuestiones relativas a su implementaci\'on y diferentes versiones. Se exponen los resultados de una implementaci\'on usando la librer\'ia $\cocoalib$, y las razones del uso de esta lbrer\'ia, junto con una breve descripci\'on de algunas cuestiones t\'ecnicas. En la \'ultima parte de la secci\'on se muestran algunos experimentos realizados sobre ditintos tipos de ideales y su comparaci\'on con otros algoritmos usados para calcular invariantes homl\'ogicos de ideales, implementados en los principales sistemas de \'algebra computacional dirigidos a este tipo de c\'alculos. En concreto se compara el algoritmo basado en \'arboles de Mayer-Vietoris con los algoritmos usados para calcular resoluciones minimales y series de Hilbert multigraduadas en \cocoa, \singular  y \macaulay. 

\section*{Cap\'itulo 4: Aplicaciones}
Este cap\'itulo est\'a dedicado a mostrar una colecci\'on de aplicaciones de las t\'ecnicas y c\'alculos expuestos en los cap\'itulos anteriores, en particular de los c\'alculos y analisis de ideales realizados con ayude de \'arboles de Myer-Vietoris, o de la homolog\'ia de Koszul de ideales de monomios.

La primera secci\'on aplica el an\'alisis de ideales monomiales mediante \'arboles de Mayer-Vietoris a distintas familias de ideales. Los ideales estudiados son: En primer lugar, ideales de tipo \emph{Borel fixed}, ideales estables o segmentos. Estos ideales tienen importancia desde el punto de vista te\'orico, en particular por su relaci\'on con la funci\'on de Hilbert. En segundo lugar, estudiamos ideales gen\'ericos, en el sentido de \cite{MSY00,MS04}, con ayuda de los \'arboles de Mayer-Vietoris obtenemos resoluciones m\'inimals de este tipo de ideales. En tercer y cuarto lugar estudiamos ideales monomiales que tienen importancia por su aplicai\'on en geometr\'ia algebr\'aica y teor\'ia de grafos, se trata de los ideales de Valla \cite{V04} y Ferrers \cite{CN06,CN07}. el an\'alisis de este tipo de ideales usando las t\'ecnicas introducidas en los cap\'itulos anteriores permite dar resoluciones m\'iniams de estos ideales, e incluso f\'ormulas para sus n\'umeros de Betti, as\'i como descomposiciones expl\'icitas en ideales irreducibles y/o primarios. Finalmente, otro tipo de ideales estudiados son los ideales quasi-estables, introducidos en \cite{S02b}, que tienen gran importancia en el contexto de las bases involutivas y sus aplicaciones. En relaci\'on a estos ideales estudiamos la conexi\'on entre su homolog\'ia de Koszul y el procedimiento para la obtenci\'on de sus bases de Pommaret; tambi\'en se atiende el problema de la $\delta$-regularidad.

La segunda secci\'on se ocupa de aplicaciones de la homolog\'ia de Koszul en la teor\'ia formal de ecuaciones diferenciales. Tras una breve introducci\'on al tema, se explica el papel de la homolog\'ia de Koszul en el an\'alisis de dos problemas importantes en esta teor\'ia: por un lado la caracteriaci\'on de la involutividad e integrabilidad formal de un sistema diferencial, y por otro los problemas de valor inicial. En relaci\'on al primer problema, partiendo de la caracterizaci\'on de involutividad basada en un resultado de Serre que la relaciona con la anulaci\'on de ciertos $Tor$, se ve el papel de la homolog\'ia de Koszul para esta caracterizaci\'on y se dan algunos resultados y ejemplos. En relaci\'on con el segundo problema, la formulaci\'on correcta de un problema de valor inicial para un sistema diferencial se puede estudiar en t\'erminos de descomposiciones de Stanley, que como se han visto en el cap\'itulo segundo, pueden obtenerse a partir de la homolog\'ia de Koszul; en esta seci\'on damos algunos ejemplos de la aplicaci\'on de estas t\'ecnicas.

Finalmente, la tercera secci\'on estudia la aplicaci\'on de los ideales monomiales a la teor\'ia de fiabilidad. En \cite{GW04} se establece una correspondencia entre ideales de monomios y sistemas coherentes, que permite calcular la fiabilidad de dichos sistemas mediante el c\'alculo de n\'umeros de Betti multigraduados y series de Hilbert multigraduadas de los correspondientes ideales monomiales asociados. En esta secci\'on se aplican las t\'ecnicas vistas, en particular los complejos de Koszul simpliciales y en amyor medida los \'arboles de Mayer-Vietoris a distintos sistemas de particualr importancia, tanto redes (en particular redes en serie-paralelo) como a no-redes (en particular sistemas de tipo \emph{$k$-out-of-$n$} ordinarios y consecutivos). El an\'alisis de los ideales correspondientes a estos sistemas mediante las t\'ecnicas expuestas permite no s\'olo obtener algoritmos eficaces para el c\'alculo de su fiabilidad sino incluso f\'ormulas expl\'icitas.

\fancyhead[ER]{\itshape Conclusiones y trabajo futuro} 
\fancyhead[OL]{\itshape Conclusiones y trabajo futuro}
\chapter*{Conclusiones}
\spanishchapter{Conclusiones}

Esta tesis ha estado centrada en c\'alculos expl\'icitos y aplicaciones de la homolog\'ia de Koszul y los n\'umeros de Betti de ideales monomiales. Con este inter\'es presente, los objetivos principales han sido:
\begin{itemize}
\item Analizar la homolog\'ia de Koszul de ideales monomiales y aplicarla a la descripci\'on de la estructura de dichos ideales.
\item Describir algoritmos para realizar c\'alculos eficaces de los invariantes homol\'ogicos de ideales de monomios, en particular sus n\'umeros de Betti, resoluciones libres, homolog\'ia de Koszul y serie de Hilbert.
\item Aplicar la teor\'ia de ideales monomiales a problemas dentro y fuera de las matem\'aticas, haciendo uso, en particular, de los invariantes homol\'ogicos de estos ideales.
\end{itemize}
La tesis ha comenzado con una descripci\'on del concepto de homolog\'ia de Koszul y las propiedades de los ideales monomiales, junto con algunas t\'ecnicas combinatorias \'utiles para el estudio de estos ideales. En el primer cap\'itulo se ha mostrado que la naturaleza combinatoria de los ideales monomiales permite el estudio de su estructura algebraica y homol\'ogica con herramientas sencillas. En particular se han mostrado alguans t\'ecnicas dedicadas al estudio de resoluciones y homolog\'ia de Koszul, introduciendo, como una herramienta nueva las sucesiones de Mayer-Vietoris asociadas a un ideal monomial.

%

En el segundo cap\'itulo hemos analizado la estructura algebraica y homol\'ogica de los ideales monomiale sutilizando la homolog\'ia de Koszul. Hemos empezado con el problema de la obtenci\'on de una resoluci\'on libre minimal de un ideal de monomios a partir del conocimiento de la homolog\'ia de Koszul del ideal. Este problema hab\'ia sido resuleto por Aramova y Herzog \cite{AH95,AH96}, pero las construcciones expl\'icitas eran complicadas incluso para casos sencillos. En \cite{RS06,RSS06} las construcciones de Aramova y Herzog fueron transformadas en algoritmos mediante la \emph{homolog\'ia efectiva}, de modo que este problema queda resuelto tanto desde un punto de vista te\'orico como algor\'itmico. Existe todav\'ia una carencia de m\'etodos expl\'icitos simples que permitan obtener los diferenciales de una resoluci\'on m\'inima de un ideal de monomios a partir de un conjunto de generadores de su homolog\'ia de Koszul. 


Otro problema tratado en este cap\'itulo ha sido el de encontrar descomposiciones combinatorias del anillo cociente $R/I$, donde $I\subset R=\kb[x_1,\dots,x_n]$ es un ideal monomial. En particular hemos obtenido un procedimiento para obtener una descomposici\'on de Stanley a partir del conocimiento de la homolog\'ia de Koszul de $I$. Hemos visto que los generadores de la $(n-1)$-\'esima homolog\'ia juegan un papel muy importante en este problema. Esto es debido a la relaci\'on de estos generadores con el borde de un ideal monomial, en particular con las llamadas \emph{esquinas m\'aximas}, que definen la descomposici\'on de Stanley en el caso artiniano. El caso no artiniano se ha estudiado en t\'erminos del artiniano.


Tambi\'en se han tratado en este cap\'itulo las descomposiciones en ideales irreducibles y primarios. Se dan procedimientos para calcular una descomposici\'on en ideales irreducibles a partir de la homolog\'ia de Koszul de un ideal monomial o de la de su clausura artiniana. primero tratamos el caso artiniano usando la $(n-1)$-\'esima homolog\'ia de Koszul y el concepto de esquina m\'axima. Este acercamiento es distinto a otros \cite{R07,MS04} que tambi\'en se basan en conceptos equivalentes al de esquina m\'axima. Mientras que otros acercamientos se basan en la clausura artiniana para el caso de dimensi\'on amyor que cero, hemos desarrollado otro procedimiento que usa exclusivamente la homolog\'ia de Koszul del propio ideal. Esto es un m\'etodo alternativo a los existentes. Con respecto a descomposiciones primarias irredundantes se presenta un procedimiento sencillo para obtener una descomposici\'on primaria irredundante en concreto, a partir de la $(n-1)$-\'esima homolog\'ia de Koszul de la clausura artiniana del ideal.


El tercer cap\'itulo se centra en los c\'alculos, y presenta una de las contribuciones centrales de esta tesis: los \'arboles de Mayer-Vietoris. Dado un ideal monomial hemos presentado sus \'arboles de Mayer-Vietoris asociados y las propiedades de los mismos. Los \'arboles de Mayer-Vietoris constituyen una herramienta computacional basada en las sucesiones de Mayer-Vietoris introducidas en el cap\'itulo primero. Nos permiten llevar a cabo c\'alculos homol\'ogicos sobre los ideales monomiales a los que estan asociados. Todo \'arbol de Mayer-Vietoris produce una resoluci\'on multigraduada del ideal corespondiente, y por tanto tambi\'en una expresi\'on de su serie de Hilbert multigraduada. Adem\'as, proporcionan un subconjunto y un superconjunto de los multigrados de los generadores de la homolog\'ia de Koszul en cada grado homol\'ogico. Se han dado tambi\'en condiciones bajo las cuales estos subconjuntos y\o superconjuntos proporcionan los n\'umeros de Betti multigraduados reales del ideal, y por tanto el \'arbol de Mayer-Vietoris correspondiente nos da la resoluci\'on libre m\'inima. las familias de ideales para las que podemos obtener los n\'umeros de Betti multigraduados directamente de los \'arboles de Mayer-Vietoris se han llamado \emph{ideales de Mayer-Vietoris}, y se han dividido en varios tipos. Se han dado tambi\'en algunas familias importantes de ideales de Mayer-Vietoris.


Los \'arboles de Mayer-Vietoris no son solamente una herramienta computacional sino tambi\'en un medio para analizar ideales monomiales y su estructura homol\'ogica, inclu\'idas resoluciones libres, n\'umeros de Betti y serie de Hilbert. En su art\'iculo \cite{PS07}, I. Peeva y M. Stillmann proponen varias conjeturas y problemas abiertos en torno a sizigias y funciones de Hilbert. En relaci\'on con los ideales monomiales, el principal objetivo es enunciado como

\noindent{\bf Problema 3.9.1 \cite{PS07}: }Introducir nuevas construcciones e ideas sobre resoluciones minimales.

En este contexto, los \'arboles de Mayer-Vietoris son una herramienta sencilla pero computacionalmente eficiente para analizar y obtener resoluciones monomiales. En aprticular, proporcionan una mirada cercana a las relaciones de divisibilidad en el ret\'iculo de $mcm$s del ideal, con la atenci\'on puesta en la homolog\'ia de Koszul del mismo, y por tanto en las resoluciones libres del ideal y en sus n\'umeros de Betti. En este sentido, se han presentado varios ejemplos de este tipo de an\'alisis en esta secci\'on y en el cap\'itulo cuarto.

%
%

Tambi\'en en el cap\'itulo tercero se han tratado algunas cuestiones algor\'itmicas. Se ha dado una descripci\'on del algoritmo b\'asico para el c\'alculo de \'arboles de Mayer-Vietoris junto con algunos detalles de su implementaci\'on. se ha hecho una implementaci\'on de este algoritmo usando la librer\'ia de \cpp $\cocoalib$. Esta librer\'ia, en desarrollo por el equipo de $\cocoa$ \cite{cocoa}, une la eficiencia de \cpp y la posibilidad de realizar c\'alculos en \'algebra conmutativa y comunicarse con otros sistemas de \'algebra computacional. se han realizado varios experimentos con esta implementaci\'on para demostrar el funcionamiento y capacidad de distintas estrategias usadas para construir los \'arboles de Mayer-Vietoris, los tiempos de c\'alculo comparados con otros algoritmos que calculan resoluciones, series de Hilbert multigraduadas y n\'umeros de Betti en otros sistemas, as\'i como algunas mejoras del algoritmo b\'asico. Estos experimentos muestran que el c\'alculo de los \'arboles de Mayer-Vietoris es un acercamiento eficaz para el c\'alculo de resoluciones y series de Hilbert. En particular, cuando el n\'umero de variables crece o cuando la resoluci\'on m\'inima es grande, los \'arboles de Mayer-Vietoris pueden proporcionarnos la informaci\'on que se necesita en la mayor parte de la aplicaciones, en aprticular los n\'umeros de Betti (multigraduados) o al menos cotas para ellos. Tambi\'en son una alternativa eficiente para el c\'alculo de series de Hilbert.


El cuarto cap\'itulo est\'a dedicado a las aplicaciones. Se han dado varias aplicaciones de las t\'ecnicas, procedimientos y algoritmos presentados en los cap\'itulos previos. Los hemos aplicado a diferentes clases de ideales monomiales, a otras \'areas de las matem\'aticas, a saber, la teor\'ia formal de ecuaciones diferenciales, y finalmente, fuera de las matem\'aticas, a la teor\'ia de fiabilidad.


En primer lugar se han estudiado varios tipos de ideales. Hemos analizado los ideales de tipo Borel-fijo, ideales estables y segmento con ayuda de los \'arboles de Mayer-Vietoris, y se ha probado que son una buena herramienta para realizar afirmaciones sobre estos tipos de ideales. Hemos dado demostraciones alternativas de varios resultados importantes sobre la estructura homol\'ogica de estos ideales. Lo mismo se ha aplicado a ideales de tipo Scarf, entre los que destacan los ideale sgen\'ericos. Todos estos tipos de ideales tienen importantes caracter\'isticas te\'oricas. En segundo lugar se han estudiado tres tipos de ideales con aplicaciones en otras \'areas de las matem\'aticas. Los ideales de Valla \cite{V04} tienen aplicaiones en geometr\'ia algebraica, se ha demostrado que son Mayer-Vietoris de tipo $B2$ y usando \'arboles de Mayer-Vietoris hemos dado f\'ormulas expl\'icitas para sus n\'umeros de Betti graduados y descomposiciones irreducibles y primarias irredundantes de ellos. Los ideales de Ferrers \cite{CN06,CN07} tienen aplicaciones en la teor\'ia de grafos, entre otras; son Mayer-Vietoris de  tipo $A$ y a partir de sus \'arboles de Mayer-Vietoris hemos obtenido f\'ormulas expl\'icitas para sus n\'umeros de Betti y descomposiciones irreducibles y primarias de ellos, as\'i como otros invariantes. Finalmente, hemos estudiado ideales cuasi-estables, que guardan una estrecha relaci\'on can bases involutivas (en particular las de Pommaret) y aparecen tambi\'en en la teor\'ia formal de ecuaciones diferenciales. Aqu\'i hemos usado su homolog\'ia de Koszul para desarrollar un procedimiento que completa un ideal de este tipo hasta una base de Pommaret y decide tambi\'en la $\delta$-regularidad de las coordenadas en las que el ideal est\'a dado, proporcionando as\'i un acercamiento alternativo a los existentes en la literatura \cite{S02b,HS02}.


Las aplicaciones a la \emph{teor\'ia formal de sistemas diferenciales} se basa en la dualidad entre la cohomolog\'iad e Spencer y la homolog\'ia de Koszul. El marco para esta aplicaci\'on es una cercamiento geom\'etrico a las EDPs que hace posible una an\'alisis algebraico de las mismas. El papel de la homolog\'ia de Koszul en este contexto tiene importancia en relaci\'on con la involuci\'on y la integrabilidad formal. Conociendo la homolog\'ia de Koszul y\o los n\'umeros de Betti del s\'imbolo de un sistema diferencial nos da un modo de detectar la involutividad y cuasi-regularidad de un sistema diferencial dado. Esto se basa en la correspondencia en tre el grado de involuci\'on y el concepto algebraico de regularidad de Castelnuovo-Mumford \cite{M03,S07}. Otra aplicaci\'on en este \'area es la formulaci\'on de problemas de valor inicial correctamente enunciados. Esta formulaci\'on esta estrechamente relacionada con descomposiciones de Stanley del correspondiente anillo cociente, y aplicamos el procedimeto presentado en el cap\'itulo segundo para obtener dichas descomposiciones a partir de la homolog\'iad e Koszul del ideal correspondiente.


Finalmente usamos ideales monomiales en teor\'ia de fiabilidad. En \cite{GW04} se estableci\'o una conexi\'on entre ideales monomiales y sistemas coherentes. En particular, ciertas expresiones de la serie de Hilbert de un ideal monomial proporcionan cotas para la fiabilidad del sistema correspondiente. Hemos usado diferentes t\'ecnicas para obtener cotas ajustadas o incluso f\'ormulas para calcular la fiabilidad de varios tipos importantes de sistemas. Un tipo muy natural de redes, las llamadas redes \emph{en serie-paralelo} se ha tratado usando \'arboles de Mayer-Vietoris. Se ha demostrado que los \'arboles asociados a estas redes son Mayer-Vietoris de tipo $A$. Por tanto, las resoluciones de Mayer-Vietoris de este tipo de ideales son m\'inimas, y podemos encontrar cotas ajustadas para la fiabilidad de las redes correspondientes. Hemos dado adem\'as procedimientos expl\'icitos para calcular la fiabilidad de estas redes sin calcular resoluciones, y para calcular sus n\'umeros de Betti. Hemos dado tambi\'en algunos subtipos de estas redes para los que proporcionamos f\'ormulas expl\'icitas para los n\'umeros de Betti de los ideales asociados. Tambi\'en se han estudiado otros tipos de sistemas, en particular varios sistemas de tipo $k$-entre-$n$, a saber, el usual, el consecutivo y el multiestado. En estos casos, bien la homolog\'ia simplicial de Koszul o los \'arboles de Mayer-Vietoris han sido usados para producir f\'ormulas para los n\'umeros de Betti de los ideales asociados, y por tanto para la fiabilidad del sistema. En algunos casos se dan demostraciones sencillas de algunos resultados conocidos, y en otros casos se han obtenido nuevos resultados utilizando nuestras t\'ecnicas.


Las direcciones exploradas en esta tesis abren varios caminos en los que realizar trabajos futuros, tanto desde el punto de vista te\'orico como desde el aplicado, as\'i como en un marco computacional.

La homolog\'ia de Koszul tiene un car\'acter combinatorio que ha sido estudiado mediante algunas tecnicas sencilals provinientes de la topolog\'ia algebraica y mediante \'algebra homol\'ogica b\'asica. Esta direcci\'on debe ser estudiada con mayor profundidad usando diferentes t\'ecnicas que nos permitan obtener nuevos resultados. Dos direcciones a seguir en primer lugar en este contexto son la \emph{homolog\'ia relativa} y las \emph{sucesiones espectrales}, de las que podemos esperar nuevos resultados sobre la homolog\'ia de Koszul combinatoria. Tambi\'en puede proporcionar nueva informaci\'on el estudio de los objetos estudiados en esta tesis desde un punto de vista cohomol\'ogico y mediante el uso de los funtores $Ext$ y $Hom$.

%

En cuanto al estudio de los ideales monomiales pueden seguirse varias direcciones. Un \'area amplia en la etor\'ia de ideales monomiales se refiere a los ideales asociados a grafos (\emph{edge ideals}) y a los \emph{facet ideals} asociados a complejos simpliciales, v\'ease por ejemplo \cite{V01}, \cite{FV07,HV07} y \cite{F02}. Algunos de los resultados de esta tesis pueden interpretarse en este contexto, pero debe hacerse un estudio riguros de la homolog\'ia de Koszul de ideales \emph{edge} y \emph{facet} y de la aplicaci\'on de los \'arboles de Mayer-Vietoris a estos ideales. Hay tambi\'en varios tipos de ideales que pueden ser estudiados usando nuestras t\'ecnicas y que tienen importancia desde un punto de vista te\'orico, como los ideales $LPP$ \cite{FR07}, los ideales \emph{pretty clean} \cite{JZ07} y otros. Tambi\'en hay algunos conceptos como polarizaci\'on \cite{F05} o splitting \cite{EK90} de ideales monomiales que podr\'ian ser tratados con estas t\'ecnicas, de lo que cabe esperar nuevos an\'alisis de estos ideales.


Con respecto a las aplicaciones, los ideales monomiales son un concepto que aparece en muchos lugares. Las aplicaciones presentadas en esta tesispueden verse completadas mediante un estudio m\'as profundo. Las primeras direcciones a seguir a partir de aqu\'i son las relaciones entre la homolog\'ia de Koszul y las bases de Pommaret, los \'arboles de Mayer-Vietoris de ideales cuasi-estables y diferentes sistemas coherentes de la teor\'ia de fiabilidad, tales como los isistemas $k$-entre-$n$ consecutivos multidimensionales. Tambi\'en se pueden estudiar mediante nuestras t\'ecnicas otras aplicaciones distintas de los ideales monomiales, en particular las relacionadas con las ciencias de la vida, como la relaci\'on entre ideales monomiales y \'arboles filogen\'eticos, ver por ejemplo \cite{GLRRV06}.


Finalmente, desde el punto de vista algor\'itmico, el primer paso a dar es una implementaci\'on completa de los procedimientos para obtener resoluciones m\'inimas usando ideales de Mayer-Vietoris. Tambi\'en la mejora en la implementaci\'on de los algoritmos ya implementados y la adici\'on de nuevos procedimientos para detectar la homolog\'ia en los multigrados relevantes repetidos es un punto a completar. Esto derivar\'a en un nuevo algoritmo alternativo para el c\'alculo de invariantes homol\'ogicos de ideales monomiales. Unido a las t\'ecnicas de homolog\'ia efectiva estudiadas en el cap\'itulo segundo, esto puede ser aplicado a ideales polinomiales m\'as generales.


\fancyhead[ER]{\itshape Publicaciones} 
\fancyhead[OL]{\itshape Publicaciones}
\chapter*{Publicaciones}

Incluimos a continuaci\'on los res\'umenes de las publicaciones a las que ha dado lugar el trabajo presentado en esta memoria.

\begin{itemize}

\item {\bf E. S\'aenz de Cabez\'on. Computing Koszul homology for monomial ideals. En Proceedings of the Xth Rhine Workshop of Computer Algebra, Basel, 2006.
}

En este trabajo se presentan t\'ecnicas y algoritmos para el c\'alculo de la homolog\'ia de Koszul y n\'umeros de Betti espec\'ificamente dirigidos a ideals de monomios. Se tiene en cuenta el car\'acter combinatorio de estos ideales para el dise\~no de t\'ecnicas que exploten esta caracter\'istica. Se usan complejos simpliciales o el c\'alculo de los m\'odulos de sizigias dimensi\'on a dimensi\'on. Finalmente se aplican estos resultados al c\'alculo en ideales polinomiales.

\item {\bf E. S\'aenz de Cabez\'on. Mayer-Vietoris trees of monomial ideals. In Global Integrability of Field Theories, Proceedings of GIFT 2006,
pp. 311--333. Universit\"atsverlag Karlsruhe, 2006.}

En este art\'iculo se presentan la definici\'on y principales propiedades de los \'arboles de Mayer-Vietoris y su aplicaci\'on al c\'alculo de la homolog\'ia de Koszul de ideales monomiales. Se presenta un algoritmo basado en estos \'arboles para el c\'alculo de n\'umeros de Betti de ideales monomiales y se muestran los resultados de la aplicaci\'on de dicho algoritmo a distintos ejemplos y familias de ideales. 

\item {\bf E. S\'aenz de Cabez\'on. Homological Description of Monomial Ideals using Mayer-Vietoris Trees. In Proceedings of MEGA 2007, Strobl am Wolfgangsee (Austria), 2007.}

En este art\'iculo se emplean los \'arboles de Mayer-Vietoris para obtener resultados acerca de distintos tipos de ideales. Se definen los ideales de Mayer-Vietoris y se muestra que algunos tipos de ideales que aparecen en la literatura (Valla, Ferrers, genericos, etc...) son ideales de Mayer-Vietoris.

\item {\bf E. S\'aenz de Cabez\'on, H.P. Wynn. Betti numbers and minimal free resolutions for multi-state system reliability bounds. In Proceedings of MEGA 2007, Strobl am Wolfgangsee (Austria), 2007.}

En este art\'iculo se decribe la aplicaci\'on de las t\'ecnicas de c\'alculo de funciones de Hilbert multigraduadas y homolog\'ia de Koszul de ideales monomiales a la evaluaci\'on de la fiabilidad de sistemas coherentes. Se hace uso de distintas t\'ecnicas entre las que est\'an los \'arboles de Mayer-Vietoris y el complejo de Koszul simplicial. Se dan f\'ormulas para el c\'alculo de la fiabilidad de sistemas muy relevantes, como los \emph{k-out-of-n}, \emph{consecutive k-out-of-n} o redes \emph{series-paralelo}.

\item {\bf J. Rubio, F. Sergeraert, E. S\'aenz de Cabez\'on. Minimal free resolutions, enviado para su publicaci\'on, 2006.}

El art\'iculo est\'a enfocado a proporcionar un algoritmo sencillo para el c\'alculo de resoluciones m\'inimas de $A_0$-m\'odulos de tipo finito, en el caso en el que $A_0$ es un anillo de polinomios ordinario $A_0=\kb[x_1,\dots,x_n]_0$ localizado en $0\in\kb^n$. Los m\'etodos usados se basan en la homolog\'ia efectiva, desarrollada por el primer y segundo autor. Se dan los resultados de algunos algoritmos implementados por el primer autor.

\item {\bf E. S\'aenz de Cabez\'on. Multigraded Betti numbers without computing minimal free resolutions, enviado para su publicaci\'on, 2007.}

En este art\'iculo se usan los \'arboles Mayer-Vietoris para calcular n\'umeros de Betti multigraduados de ideales monomiales sin calcular sus resoluciones minimales. El m\'etodo proporciona no s\'olo un algoritmo competitivo para realizar \'estos c\'alculos sino tambi\'en una nueva herramienta para el n\'alisis de la estructura homol\'ogica de los ideales monomiales.

\end{itemize}

\end{otherlanguage}


\end{document}